\documentclass[11pt,a4paper,fleqn]{article}
\usepackage{amsfonts,amsmath}
\usepackage{latexsym}
\usepackage{amssymb}
\usepackage{euscript}
\usepackage{graphicx}
\usepackage{bm}
\usepackage{a4wide}
\usepackage{lscape}
\usepackage{float}

\numberwithin{equation}{section}
\newcommand{\noi}{\noindent}
\newcommand{\dis}{\displaystyle}

\newcommand{\qed}{\nopagebreak\hspace*{\fill}
{\vrule width6pt height6ptdepth0pt}\par}

\begin{document}

\begin{center}{\bf \LARGE A global view of Brownian penalisations}\end{center}

\begin{center} {\Large J. Najnudel$^{(1)}$, B. Roynette$^{(2)}$, M. Yor$^{(3)}{}^{(4)}$} \end{center}

\begin{center} {\large 30/01/2009} \end{center}

\bigskip

\hbox{\noi$^{(1)}$ Institut f\"{u}r Mathematik, Universit\"{a}t Z\"{u}rich \\}
\hbox{\hspace*{0.55cm} Winterthurerstrasse 190, CH 8057 Z\"{u}rich \\}
\hbox{\noi$^{(2)}$Institut Elie Cartan, Universit\'{e} Henri Poincar\'e,    \\}
\hbox{\hspace*{0.55cm} B.P. 239, 54506 Vandoeuvre les Nancy Cedex\\}

\hbox{\noi$^{(3)}$ Laboratoire de Probabilit\'es et Mod\`{e}les al\'eatoires, \\}
\hbox{\hspace*{0.55cm}Universit\'es de Paris VI et VII, 4 Place Jussieu, Case 188 \\}
\hbox{\hspace*{0.55cm}F - 75252 Paris Cedex 05\\}

\hbox{\noi$^{(4)}$ Institut Universitaire de France \\}

\vskip 10pt 

\begin{center}{\large Preface}\end{center}

\vskip 10pt 

\noi {\bf 1)} Let $\dis \big(\Omega = \mathcal{C} (\mathbb{R}_{+} \to \mathbb{R}), (X_{t}, \mathcal{F}_{t})_{t \ge 0}, \mathcal{F}_{\infty} = \mathop{\vee}_{t \ge 0}^{} \mathcal{F}_{t}, W_{x} (x \in \mathbb{R})\big)$ denote the canonical realisation of one-dimensional Brownian motion. With the help of Feynman-Kac type penalisation results for Wiener measure, we have, in  [RY, M], constructed on $(\Omega, \mathcal{F}_{\infty})$ a positive and $\sigma$-finite measure {\bf W}. The aim of this second monograph, in particular \underline{Chapter 1}, is to deepen our understanding of {\bf W}, as we discuss there other remarkable properties of this measure. 

\noi For pedagogical reasons, we have chosen to take up here again the construction of {\bf W} found in [RY, M], so that the present monograph may be read, essentially, independently from our previous papers, including [RY, M]. 

\noi Among the main properties of {\bf W} presented here, let us cite :

\noi $\bullet$ the close links between {\bf W} and probabilities obtained by penalising Wiener measure by certain functionals : see Theorems 1.1.2, 1.1.11 and 1.1.11' ;

\noi $\bullet$ the existence of integral representation formulae for the measure {\bf W} : see Theorems 1.1.6 and 1.1.8. These formulae allow to express {\bf W} in terms of the laws of Brownian bridges and of the law of the 3-dimensional Bessel process \big(see formula (1.1.43)\big). They also allow to express {\bf W} in terms of the law of Brownian motion stopped at the first time when its local time at 0 reaches level $l$, $l$ varying, and of the law of the 3-dimensional Bessel process \big(see formula (1.1.40)\big). One may observe that these representation formulae are close to those obtained by Biane and Yor in [BY] for some different $\sigma$-finite measures on Wiener path space.

\noi $\bullet$ the existence, for every $F \in L_{+}^{1}(\mathcal{F}_{\infty}, {\bf W} )$, of a $\big( (\mathcal{F}_{t}, \; t \ge 0), W\big)$ martingale $\big(M_{t} (F), \; t \ge 0\big)$ which converges to 0, as $t \to \infty$ (see Theorem 1.2.1). Many examples of such martingales are given (see Chap. 1, Examples 1 to 7). The Brownian martingales of the form $\big(M_{t} (F), \; t \ge 0)$ are characterized among the set of all Brownian martingales (see Corollary 1.2.6) and a decomposition theorem of every positive Brownian  supermartingale involving the martingales $\big(M_{t} (F), \; t \ge 0\big)$ is established in Theorem 1.2.5. In the same spirit, we show (see Theorem 1.2.11) that every martingale $\big(M_{t} (F), \; t \ge 0\big)$ with $F \in L^{1} (\mathcal{F}_{\infty}, {\bf W}), \; F$ not necessarily $\ge 0$, may be decomposed in a canonical manner into the sum of two quasi-martingales which enjoy some remarkable properties. In particular, this result allows to obtain a characterization of the martingales $\big(M_{t} (F), \; t \ge 0\big)$, with $F \in L^{1} (\mathcal{F}_{\infty}, {\bf W})$ which vanish on the zero set of the process $(X_{t}, \; t \ge 0)$. This is Theorem 1.2.12.

\smallskip

\noi $\bullet$ a general penalisation Theorem, for Wiener measure, which is valid for a large class $\mathcal{C}$ of penalisation functionals $(F_{t}, \;t \ge 0)$ and whose proof hinges essentially upon some remarkable properties of {\bf W} : this is the content of Subsection 1.2.5 and particularly Theorem 1.2.14 and Theorem 1.2.15.

\noi $\bullet$ the existence of invariant measures, which are intimately related with {\bf W}, for several Markov processes taking values in function spaces (see Section 1.3). Chapter 1 of this monograph is devoted to the results we have just described.

\smallskip

\noi {\bf 2)} The results relative to the 1-dimensional Brownian motion are extended, in \underline{Chapter 2} of this monograph to 2-dimensional Brownian motion (we identify $\mathbb{R}^{2}$ to $\mathbb{C}$, and use complex notation). In this framework, the role of the measure {\bf W} is played by a positive and $\sigma$-finite measure, which we denote ${\bf W}^{(2)}$ on $\big( \Omega = \mathcal{C} (\mathbb{R}_{+} \rightarrow \mathbb{C}), \; \mathcal{F}_{\infty}\big)$. The properties of ${\bf W}^{(2)}$ are, mutatis mutandis, analogous to those of {\bf W}. However, in the set-up of the $\mathbb{C}$-valued Brownian motion $(X_{t}, \; t\ge 0)$, it is of interest to consider the winding process $(\theta_{t}, \; t \ge 0 )$ :
$$(\theta_{t}, \; t \ge 0) = \left(\theta_{0} + {\rm Im} \int_{0}^{t} \frac{d X_{s}}{X_{s}}, \; t \ge 0\right)$$
\noi We study this process under ${\bf W}^{(2)}$. We then obtain a Spitzer type limit theorem about the asymptotic behavior in distribution for $\theta_{t}$, adequately normalized, as $t \to \infty$. This is Theorem 2.3.1. (see also Remark 2.3.2).

\noi {\bf 3)} \underline{Chapter 3} of this Monograph is devoted to the transcription of several of the preceding results to a more general framework, that of a certain class of linear diffusions (taking values in $\mathbb{R}_{+}$). This class is described in Section 3.1. It is in fact the class of the linear diffusions studied by Salminen, Vallois and Yor \big(see [SVY]\big). These are diffusions taking values in $\mathbb{R}_{+}$, and associated with a speed measure $m$ and a scale function $S$, both of which have adequate properties. Fundamental examples of such diffusions are the Bessel processes with dimension $d=2(1-\alpha)$ for $0<d<2$. \big(We also refer to $\alpha \in ]0,1[$, or to the index $-\alpha \in ]-1,0[\big)$. The case $d=1 \left({\rm or} \; \alpha = \dis \frac{1}{2}\right)$ is that of reflected Brownian motion.

\smallskip

\noi We particularize, in Section 3.3, for these examples, the general results obtained for this class of linear diffusions (see Theorem 3.3.1). The analogue, for the Bessel process with index $(-\alpha)$, of the measure {\bf W}, is denoted ${\bf W}^{(-\alpha)}$. Then, still in this framework of the Bessel process of index $(-\alpha)$, we establish some link between, on one hand, the measure ${\bf W}^{(- \alpha)}$ and, on the other hand, a Feynman-Kac type penalisation of a Bessel process with index $(- \alpha)$ (see Remark 3.3.2 and 3.3.3). Finally, in Section 3.4, we give a new description of the measure ${\bf W}^{(-\alpha)}$ restricted to $\mathcal{F}_{g}$, with $g := \sup \{t \ge 0 \;;\; X_{t}=0\}$. This is Theorem 3.4.1. This description is the transcription in our situation of results of Pitman-Yor \big(see [PY2]\big). In some sense, this description of ${\bf W}^{(-\alpha)}$ restricted to $\mathcal{F}_{g}$ resembles the description due to D. Williams (see [Wi]) of the It\^{o} measure of Brownian excursions.

\smallskip

\noi {\bf 4)} \underline{Chapter 4} of this monograph consists in obtaining, this time in the framework of Markov chains taking values in a countable set, the analogue of the preceding results. Section 4.1 is devoted to the definition of the measures $(\mathbb{Q}_{x}, \; x \in E)$ which play here the role of the measures ${\bf W}_{x}$ in the precedings chapters. Also as in the preceding chapters, certain martingales are associated to these measures $(\mathbb{Q}_{x}, \; x \in E)$; see the description in Corollary 4.2.2. In this new framework of Markov chains, the $\sigma$-finite measures $(\mathbb{Q}_{x}, \; x \in E)$ depend, from our construction, on a point $x_{0} \in E$ and on a function $\phi$. This dependence with respect to $x_{0}$ and $\phi$ is studied in subsections 4.2.3 and 4.2.4. Section 4.3 is devoted to the study of many examples ; in particular, for random walks on trees, it appears that there may exist a whole family of different measures $(\mathbb{Q}_{x}, \; x \in E)$. All results found in this Chapter 4 are due solely to J. Najnudel. \\ Finally, a very concise summary of some of the results found in this Monograph is presented, without proofs, in our Comptes Rendus de l'Acad\'emie des Sciences Note [NRY]. 

 
\newpage

{\bf{\LARGE Table of contents}}

\bigskip

\noi {\bf Chap. 1} \fbox{Existence and main properties of $\bf W$}

\smallskip

\begin{itemize}
\item[1.0] \underline{Introduction}
\smallskip
\item[1.1] \underline{Existence of {$\bf W$} and first properties}
\item[1.1.1] A few more notations
\item[1.1.2] A Feynman-Kac penalisation result
\item[1.1.3] Definition of $\bf W$
\item[1.1.4] Study of the canonical process under $W_{\infty}^{(\lambda \delta_{0})}$
\item[1.1.5] Some remarkable properties of $\bf W$
\item[1.1.6] Another approach to Theorem 1.1.6
\item[1.1.7] Relation between $\bf W$ and other penalisations (than the Feynman-Kac ones)
\smallskip
\item[1.2] \underline{$W$-Brownian martingales associated to $\bf W$}
\item[1.2.1] Definition of the martingales $\big(M_{t} (F), \; t \ge 0\big)$
\item[1.2.2] Examples of martingales $\big(M_{t} (F), \; t \ge 0\big)$
\item[1.2.3] A decomposition Theorem for positive Brownian supermartingales
\item[1.2.4] A decomposition result for the martingale $\big(M_{t} (F), \; t \ge 0\big)$
\item[1.2.5] A penalisation Theorem, for functionals in class $\mathcal{C}$
\item[1.2.6] Some other results about the martingales $(M_t(F), \; t \geq 0)$
\item[1.3] \underline{Invariant measures related to ${\bf W}_{x}$ and \boldmath${\Lambda}$\unboldmath$_{x}$}
\item[1.3.1] The process $(\mathcal{X}_{t}, \; t \ge 0)$
\item[1.3.2] The measure \boldmath$\Lambda$\unboldmath$_{x}$
\item[1.3.3] Invariant measures for the process $\big((X_{t}, L_{t}^{\bullet}), \; t \ge 0\big)$
\item[1.3.4] Invariant measures for the process $(L_{t}^{X_t - \bullet}, \; t \ge 0)$
\end{itemize}

\bigskip

\noi {\bf Chap. 2} \fbox{Existence and properties of $\bf W^{(2)}$}

\smallskip

\begin{itemize}
\item[2.1] \underline{Existence of $\bf W^{(2)}$}
\item[2.1.1] Notation and Feynman-Kac penalisations in two dimensions
\item[2.1.2] Existence of the measure $\bf W^{(2)}$
\smallskip
\item[2.2] \underline{Properties of $\bf W^{(2)}$}
\item[2.2.1] Some notation
\item[2.2.2] Description of the canonical process under $W_{\infty}^{(2,q_{0})}$
\item[2.2.3] Another description of the measure $\bf W^{(2)}$
\item[2.3] \underline{Study of the winding process under $\bf W^{(2)}$}
\item[2.3.1] Spitzer's Theorem
\item[2.3.2] An analogue of Spitzer's Theorem
\item[2.4] \underline{$W^{(2)}$-martingales associated to $\bf W^{(2)}$}
\item[2.4.1] Definition of $\big(M_{t}^{(2)} (F), \; t \ge 0\big)$
\item[2.4.2] A decomposition Theorem for positive $W^{(2)}$-supermartingales
\item[2.4.3] A decomposition Theorem of the martingales $\big(M_{t}^{(2)} (F), \; t \ge 0\big)$
\end{itemize}
\bigskip

\noi {\bf Chap. 3} \fbox{The analogue of the measure $\bf W$ for a class of linear diffusions}

\smallskip

\begin{itemize}
\item[3.1] \underline{Main hypotheses and notations}
\item[3.1.1] Our framework is that of Salminen-Vallois-Yor
\item[3.1.2] The semi-group of $(X_{t}, \; t \ge 0)$
\item[3.1.3] The local time process 
\item[3.1.4] The process $X$ conditioned not to vanish
\item[3.1.5] A useful Proposition
\item[3.2] \underline{The $\sigma$-finite measure ${\bf W}^{*}$}
\item[3.2.1] Definition of ${\bf W}^{*}$
\item[3.2.2] Some properties of ${\bf W}^{*}$
\item[3.2.3] Relation between the measure ${\bf W}^{*}$ and penalisations
\smallskip
\item[3.3] \underline{The example of Bessel processes with dimension $d$ $(0<d<2)$}
\item[3.3.1] Transcription of our notation in the context of Bessel processes
\item[3.3.2] The measure $\bf W^{(- \alpha)}$
\item[3.3.3] Relations between $\bf W^{(-\alpha)}$ $\big(d = 2(1-\alpha)\big)$ and Feynman-Kac penalisations
\item[3.4] \underline{Another description of ${\bf W}^{(-\alpha)}$ (and of ${\bf W}^{*}_g$)}
\item[3.5] \underline{Penalisations of $\alpha$-stable symmetric L\'evy process $(1<\alpha \le 2)$}
\item[3.5.1] Notation and classical results
\item[3.5.2] Definition of the $\sigma$-finite measure ${\bf P}$
\item[3.5.3] The martingales $\big(M_{t} (F), \; t \ge 0\big)$ associated with ${\bf P}$
\item[3.5.4] Relations between ${\bf P}$ and penalisations
\end{itemize}

\bigskip

\noi {\bf Chap. 4} \fbox{An analogue of the measure $\bf W$ for discrete Markov chains}

\smallskip

\begin{itemize}
\item[4.0] \underline{Introduction}
\item[4.1] \underline{Construction of the $\sigma$-finite measures $(\mathbb{Q}_{x}, \; x \in E)$ }
\item[4.1.1] Notation and hypothesis
\item[4.1.2] A family of new measures
\item[4.1.3] Definition of the measures $(\mathbb{Q}_{x}, \; x \in E)$
\item[4.2] \underline{Some properties of $(\mathbb{Q}_{x}, \; x \in E)$}
\item[4.2.1] Martingales associated with $(\mathbb{Q}_{x}, \; x \in E)$
\item[4.2.2] Properties of the canonical process under $(\mathbb{Q}_{x}, \; x \in E)$
\item[4.2.3] Dependence of $\mathbb{Q}_{x}$ on $x_0$
\item[4.2.4] Dependence of $\mathbb{Q}_{x}$ on $\phi$
\item[4.3] \underline{Some examples}
\item[4.3.1] The standard random walk
\item[4.3.2] The "bang-bang random walk"
\item[4.3.3] The random walk on a tree
\item[4.3.4] Some more general conditions for the existence of $\phi$
\item[4.3.5] The standard random walk on $\mathbb{Z}^{2}$
\end{itemize}


\newpage

\noi {\bf{\Large Chapter 1. On a remarkable $\sigma$-finite measure ${\bf W}$ on path space, which rules penalisations for linear Brownian motion}}

\medskip

{\bf \large 1.0 Introduction.}

\smallskip

\noi {\bf 1.0.1} $\big(\Omega, (X_{t}, \mathcal{F}_{t}), \; t \ge 0, \mathcal{F}_{\infty}, W_{x} (x \in \mathbb{R})\big)$ denotes the canonical realisation of 1-dimensional Brownian motion. $\Omega = \mathcal{C} (\mathbb{R}_{+} \to \mathbb{R}), \; (X_{t}, t \ge 0)$ is the coordinate process on this space and $(\mathcal{F}_{t}, \; t \ge 0)$ denotes its natural filtration ; $\mathcal{F}_{\infty} = \dis \mathop{\vee}_{t \ge 0}^{} \mathcal{F}_{t}$. For every $x \in \mathbb{R}, \; W_{x}$ denotes Wiener measure on $(\Omega, \mathcal{F}_{\infty})$ such that $W_{x} (X_{0}=x)=1$. We write $W$ for $W_{0}$ and if $Z$ is a r.v. defined on $(\Omega, \mathcal{F}_{\infty})$, we write $W_{x}(Z)$ for the expectation of $Z$ under the probability $W_{x}$.

\smallskip

\noi {\bf 1.0.2} In a series of papers \big([RVY, $i], \; i=I, II, \cdots, X \big)$ we have studied various penalisations of Wiener measure with certain positive functionals $(F_{t}, \; t \ge 0)$ ; that is for each functional $(F_{t}, \; t \ge 0)$ in a certain class, we have been able to show the existence of a probability $W_{\infty}^{F}$ on $(\Omega, \mathcal{F}_{\infty})$ such that : for every $s \ge 0$ and every $\Gamma_{s} \in b(\mathcal{F}_{s})$, the space of bounded $\mathcal{F}_{s}$ measurable variables :
$$
\mathop{\rm lim}_{t \to \infty}^{} \; \frac{W(\Gamma_{s} F_{t})}{W(F_{t})} = W_{\infty}^{F} (\Gamma_{s}) \eqno(1.0.1)
$$

\noi In this paper, we shall construct a positive and $\sigma$-finite measure $\bf W$ on $(\Omega, \mathcal{F}_{\infty})$ which, in some sense, "rules all these penalisations jointly".

\smallskip

\noi {\bf 1.0.3} In Section 1.1 of this chapter, we show the existence of $\bf W$ and we describe some of its properties.

\noi In Section 1.2, we show how to associate to $\bf W$ a family of $\big((\mathcal{F}_{t}, \; t \ge 0), W\big)$ martingales $\big(M_{t}(F), \; t \ge 0\big) \; \big(F \in L_{+}^{1} (\mathcal{F}_{\infty}, {\bf W})\big)$. We study the properties of these martingales and give many examples.

\noi In Section 1.3, we describe links between $\bf W$ and a $\sigma$-finite measure \boldmath$\Lambda$\unboldmath $\,$  which is defined as the "law" of the total local time of the canonical process under ${\bf W}$ in Chapter 3 of [RY, M]. 

\noi In particular, we construct an invariant measure \boldmath$\widetilde{\Lambda}$\unboldmath $\,$ for the Markov process $\big((X_{t}, L_{t}^{\bullet}), \; t \ge 0\big)$ (and \boldmath$\widetilde{\Lambda}$\unboldmath$\,$ is intimately related to \boldmath$\Lambda$\unboldmath). Here, $L_{t}^{\bullet}$ denotes the local times process $(L_{t}^{x}, x \in \mathbb{R}_{+})$, so that this Markov process $(X,L^{\bullet})$ takes values in $\mathbb{R} \times \mathcal{C} (\mathbb{R} \longrightarrow \mathbb{R}_{+})$.

\smallskip

\noi {\bf 1.0.4} \underline{Notations} : As certain $\sigma$-finite measures play a prominent role in our paper, we write them, as a rule, in bold characters. Thus, no confusion should arise between the $\sigma$-finite measure ${\bf W}_{x}$ and the Wiener measure $W_{x}$.

\bigskip

{\bf \large 1.1 Existence of W and first properties.}

\smallskip

\noi Our aim in this section is to define, via Feynman-Kac type penalisations, a positive and $\sigma$-finite measure ${\bf W}$ on $(\Omega, \mathcal{F}_{\infty})$. Moreover, independently from this penalisation procedure, we give several remarkable descriptions of ${\bf W}$. 

\smallskip

\noi {\bf 1.1.1} \underline {A few more notations.}

\smallskip

\noi $\big(\Omega, (X_{t}, \mathcal{F}_{t})_{t \ge 0}, \mathcal{F}_{\infty}, W_{x} (x \in \mathbb{R})\big)$ denotes the canonical realisation of 1-dimensional Brownian motion.

\noi We denote by $\mathcal{I}$ the set of positive Radon measures $q(dx)$ on $\mathbb{R}$, such that :
$$
0 < \int_{0}^{\infty} \big(1+|x|\big) \; q(dx) < \infty \eqno(1.1.1)
$$

\noi For every $q \in \mathcal{I}, \; (A_{t}^{(q)}, \; t \ge 0)$ denotes the additive functional defined by :
$$
A_{t}^{(q)} := \int_{\mathbb{R}} L_{t}^{y} \; q (dy) \eqno(1.1.2)
$$

\noi where $(L_{t}^{y}, \; t \ge 0, y \in \mathbb{R})$ denotes the jointly continuous family of local times of Brownian motion $(X_{t}, \; t \ge 0)$. When the Radon measure $q$ admits a density with respect to the Lebesgue measure on $\mathbb{R}$ (and then we denote again this density by $q$) the density of occupation formula yields :
$$
A_{t}^{(q)} = \int_{\mathbb{R}} L_{t}^{y} \; q(dy) = \int_{0}^{t} q(X_{s}) ds \eqno{(1.1.3)}
$$

\noi We denote by $b(\mathcal{F}_{s})$ \big(resp. $b_{+}(\mathcal{F}_{s})\big)$ the vector space of bounded and real valued (resp. the cone of bounded and positive) $\mathcal{F}_{s}$ measurable r.v.'s.

\noi As our means to construct $\bf W$, we use a penalisation result obtained in [RVY, I] \big(see also [RY, M]\big). In the next subsection, we recall this result.

\bigskip

\noi {\bf 1.1.2} \underline{A Feynman-Kac penalisation result.}

\smallskip

\noi {\bf Theorem 1.1.1.} {\it Let $q \in \mathcal{I}$ and :  
$$
\begin{array} {rrlrl}
\hspace*{1cm}&D_{x,t}^{(q)}& := W_{x} \left({\rm exp} \left(-\frac{1}{2} \; A_{t}^{(q)} \right)\right) &\hspace*{7cm}&{(1.1.4)} \\
&W_{x,t}^{(q)}& :=  \frac{\dis{\rm exp} \left(-\frac{1}{2} \; A_{t}^{(q)}\right)}{\dis D_{x,t}^{(q)}} \cdot W_{x} &&{(1.1.5)}
\end{array} $$

\noi {\bf 1)} For every $s�\ge 0$ and $\Gamma_{s} \in b (\mathcal{F}_{s}), \; W_{x,t}^{(q)} (\Gamma_{s})$ admits a limit as $t \longrightarrow \infty$, denoted by $W_{x, \infty}^{(q)} (\Gamma_{s})$, i.e. :
$$
W_{x,t}^{(q)} (\Gamma_{s}) \mathop{\longrightarrow}_{t \to \infty}^{} W_{x, \infty}^{(q)} (\Gamma_{s}) \eqno(1.1.6)
$$

\noi We express this property by writing that $W_{x,t}^{(q)}$ converges, as $t \longrightarrow \infty$, to $W_{x, \infty}^{(q)}$ \underline{along the} \underline{filtration $(\mathcal{F}_{s}, \; s \ge 0)$}.

\smallskip

\noi {\bf 2)} $W_{x, \infty}^{(q)}$ induces a probability on $(\Omega, \mathcal{F}_{\infty})$ such that : 
$$
{W_{x, \infty}^{(q)}}_{ |_{\mathcal{F}_{s}}} = M_{x,s}^{(q)} {\cdot W_{x}}_{|_{\mathcal{F}_{s}}} \eqno(1.1.7)
$$

\noi where $(M_{x,s}^{(q)}, \; s \ge 0)$ is the $\big((\mathcal{F}_{s}, \; s \ge 0), \; W_{x}\big)$ martingale defined by :
$$
M_{x,s}^{(q)} := \frac{\varphi_{q} (X_{s})}{\varphi_{q} (x)} \;\; {\rm exp} \left(- \frac{1}{2} \; A_{s}^{(q)}\right) \eqno(1.1.8)
$$

\noi In particular, $M_{x,0}^{(q)} = 1 \quad W_{x}$ a.s. 

\smallskip

\noi The function $\varphi_{q} : \mathbb{R} \longrightarrow \mathbb{R}_{+}$ which is featured in (1.1.8) is strictly positive, continuous, convex and satisfies : 
$$
\varphi_{q} (x) \mathop{\sim}_{|x| \to \infty}^{} |x| \eqno(1.1.9)
$$

\noi {\bf 3)} $\varphi_{q}$ may be defined via one or the other of the two following properties : 
\begin{itemize}
\item[\it i)] $\varphi_{q}$ is the unique solution of the Sturm-Liouville equation :
$$
\varphi'' = \varphi \cdot q \quad \hbox{(in the sense of distributions)} \eqno(1.1.10)
$$

which satisfies the boundary conditions :
$$
\varphi' (+ \infty) =- \varphi' (-\infty) =1 \eqno(1.1.11)
$$

\item[{\it ii)}] $\dis \; \sqrt{\frac{\pi t}{2}} \; W_{x} \left({\rm exp} \left(- \frac{1}{2} \; A_{t}^{(q)}\right) \right) \mathop{\longrightarrow}_{t \to \infty}^{} \varphi_{q} (x) \hfill (1.1.12) $
\end{itemize}

\noi {\bf 4)} Under the family of probabilities $(W_{x, \infty}^{(q)}, \; x \in \mathbb{R})$, the canonical process $(X_{t}, \; t \ge 0)$ is a transient time homogeneous diffusion. More precisely, there exists a $\big(\Omega, (\mathcal{F}_{t}, \; t \ge 0), W_{\infty}^{(q)}\big)$ Brownian motion $(B_{t}, \; t \ge 0)$ such that :
$$
X_{t} = x + B_{t} + \int_{0}^{t} \frac{\varphi'_{q} (X_{s})}{\varphi_{q} (X_{s})} \; ds \eqno(1.1.13)
$$

\noi In particular, this diffusion process $(X_{t}, \; t \ge 0)$ admits the following function $\gamma_{q}$ as its scale function :
$$
\gamma_{q} (x) := \int_{0}^{x} \frac{dy}{\varphi_{q}^{2} (y)} \eqno(1.1.14)
$$

\noi \big(and : $\big|\gamma_{q} (\pm \infty)\big| < \infty\big)$.  }

\smallskip

\noi We note that the function $\varphi_{q}$ featured in Theorem 1.1 is not exactly the one found in [RY, M]. It differs from it by the factor $\dis \sqrt{\frac{\pi}{2}}$ ; we have made this slight change in order to simplify some further formulae.

\noi {\bf 1.1.3}  \underline{Definition of {\bf W}.}

\smallskip

\noi We now use Theorem 1.1.1 to construct the $\sigma$-finite measure {\bf W}. In fact, we define, for every $x \in \mathbb{R}$, a positive and $\sigma$-finite ${\bf W}_{x}$ which is deduced from {\bf W} via the following simple translation by $x$ :
$$
{\bf W}_{x} \big(F (X_{s}, \; s \ge 0)\big) = {\bf W} \big(F (x + X_{s}, \; s \ge 0)\big) \eqno(1.1.15)
$$

\noi for every positive functional $F$. This formula (1.1.15) explains why, most of the time, we may limit ourselves to consider ${\bf W}_{0}$, which we denote simply by {\bf W}.

\bigskip

\noi {\bf Theorem 1.1.2.} {\it (Existence of {\bf W})

\noi There exists, on $(\Omega, \mathcal{F}_{\infty})$ a positive and $\sigma$-finite measure {\bf W}, with infinite total mass, such that, for every $q \in \mathcal{I}$ :
$$
{\bf W} = \varphi_{q} (0) \; {\rm exp} \left(\frac{1}{2} \; A_{\infty}^{(q)}\right) \cdot W_{\infty}^{(q)} \eqno(1.1.16)
$$

\noi or

$$
W_{\infty}^{(q)}  = \frac{1}{\varphi_{q} (0)} \; {\rm exp} \left(- \frac{1}{2} \; A_{\infty}^{(q)}\right) \cdot {\bf W} \eqno(1.1.16')
$$

\noi In other terms, the RHS of (1.1.16) does not depend on $q \in \mathcal{I}$. In particular :
$$
{\bf W} \left({\rm exp} \left(- \frac{1}{2} \; A_{\infty}^{(q)} \right) \right) = \varphi_{q} (0) \eqno(1.1.17)
$$   }

\noi or more generally, from (1.1.15) : 
$$
{\bf W}_x \left({\rm exp} \left(- \frac{1}{2} \; A_{\infty}^{(q)} \right) \right) = \varphi_{q} (x). \eqno(1.1.17')
$$   
\smallskip

\noi As we shall soon see, the measure $\bf W$ is such that, for every $t > 0$ and for every r.v. $\Gamma_{t} \in b_{+} (\mathcal{F}_{t}), \; {\bf W} (\Gamma_{t})$ equals $0$ or $+ \infty$ depending whether $W(\Gamma_{t}) = 0$ or is strictly positive. Thus, the measure $\bf W$, although, as we show later, it is $\sigma$-finite on $(\Omega, \mathcal{F}_{\infty})$, is not $\sigma$-finite on either of the measurable spaces $(\Omega, \mathcal{F}_t), \; t > 0$. 

\smallskip

\noi {\bf Proof of Theorem 1.1.2.}

\noi {\it i)} \underline{We shall establish that}, for every $q \in \mathcal{I}$, the measure on $(\Omega, \mathcal{F}_{\infty})$ :
\begin{equation*}
\varphi_{q} (0) \; {\rm exp} \left(\frac{1}{2} \; A_{\infty}^{(q)}\right) \cdot W_{\infty}^{(q)}
\end{equation*}

\noi does not depend on $q$, which allows to define $\bf W$ from formula (1.1.16). Then, we shall prove that $\bf W$, thus defined, is $(\Omega, \mathcal{F}_{\infty}) \; \sigma$-finite.

\smallskip

\noi {\it ii)} {\bf Lemma 1.1.3.} {\it For every $q \in \mathcal{I}$ and every $x \in \mathbb{R}$ :
\begin{eqnarray*}
\hspace*{1cm}1) & {\rm if} \; \lambda < 1 & W_{x, \infty}^{(q)} \left({\rm exp} \, \frac{\lambda}{2} \; A_{\infty}^{(q)} \right) < \infty \hspace*{7cm}{(1.1.18)} \\
2) & {\rm if} \; \lambda \ge1 & W_{x, \infty}^{(q)} \left({\rm exp} \, \frac{\lambda}{2} \; A_{\infty}^{(q)} \right) = + \infty \hspace*{6.7cm}(1.1.19)
\end{eqnarray*}  }

\noi {\bf Proof of Lemma 1.1.3.}

\noi From (1.1.7), for every $\lambda \in ]0, 1[$ :
\begin{gather*}
W_{x,\infty}^{(q)} \left({\rm exp} \, \frac{\lambda}{2} \; A_{t}^{(q)}\right) = W_{x} \left( \frac{\varphi_{q} (X_{t})} {\varphi_{q} (x)} \; {\rm exp} \left(-\left(\frac{1- \lambda}{2 }\right) A_{t}^{(q)}\right) \right) \nonumber \\
\quad = \frac{\varphi_{(1-\lambda)q} (x)}{\varphi_{q} (x)} \; W_{x} \left( \frac{\varphi_{q} (X_{t})}{\varphi_{(1 - \lambda)q} (X_{t})} \; \frac{\varphi_{(1-\lambda)q} (X_{t})} {\varphi_{(1-\lambda)q}(x)} \; {\rm exp} \left( - \left( \frac{1- \lambda}{2 } \right) A_{t}^{(q)}\right)\right) \hspace*{1.3cm}(1.1.20)
\end{gather*}

\noi We have been able to write (1.1.20) because the functions $\varphi_{q}$ and $\varphi_{(1-\lambda) q}$ are strictly positive. On the other hand, since for every $q \in \mathcal{I}$, $\varphi_{q} (x) \dis \mathop{\sim}_{|x| \to \infty}^{} |x|$, there exist two constants :
$$
0 < C_{1} (\lambda, q) \le C_{2} (\lambda, q) < \infty
$$

\noi such that :
$$
C_{1} (\lambda, q) \le \mathop{\rm inf}_{y \in \mathbb{R}}^{} \; \frac{\varphi_{q}(y)}{\varphi_{(1-\lambda)q} (y)} \le \mathop{\rm sup}_{y \in \mathbb{R}}^{} \; \frac{\varphi_{q} (y)}{\varphi_{(1 - \lambda)q} (y)} \le C_{2} (\lambda, q) \eqno(1.1.21)
$$

\noi Thus, from (1.1.20) :
$$
W_{x, \infty}^{(q)} \left({\rm exp} \, \frac{\lambda}{2} \; A_{t}^{(q)} \right) \le \frac{\varphi_{(1-\lambda) q}(x)}{\varphi_{q} (x)} \; \mathop{\rm sup}_{y \in \mathbb{R}}^{} \; \frac{\varphi_{q}(y)}{\varphi_{(1-\lambda) q} (y)} \; W_{x, \infty}^{((1- \lambda)q)} (1) \le \frac{C_{2}(\lambda, q)}{C_{1} (\lambda, q)} \eqno(1.1.22)
$$

\noi We now let $t \longrightarrow \infty$ and we use the monotone convergence Theorem to obtain point 1) of Lemma 1.1.3.

\noi We now write relation (1.1.20) with $\lambda =1$ :
$$
W_{x, \infty}^{(q)} \left({\rm exp} \, \frac{1}{2} \; A_{t}^{(q)}�\right) = W_{x} \left( \frac{\varphi_{q} (X_{t})}{\varphi_{q} (x)} \right) \mathop{\sim}_{t \to \infty}^{} k(x) \sqrt{t} \eqno(1.1.23)
$$

\noi with $\dis k(x)= \frac{1}{\varphi_{q}(x)} \cdot \sqrt{\frac{2}{\pi}} > 0$, since $\varphi_{q} (x) \dis \mathop{\rm \sim}_{|x| \to \infty}^{} |x|$. It then remains to let $t \longrightarrow \infty$ in (1.1.23), then to apply once again the monotone convergence Theorem to obtain point 2) of Lemma 1.1.3.

\smallskip

\noi {\it iii)} \underline{Formula (1.1.16) is then a consequence of} : 

\noi {\bf Lemma 1.1.4.} {\it The measure $\dis \varphi_{q} (x) {\rm exp} \left(\frac{1}{2} \; A_{\infty}^{(q)}\right) \cdot W_{x, \infty}^{(q)}$ does not depend on $q \in \mathcal{I}$. }

\smallskip

\noi We note that the measure $\dis \varphi_{q} (x) {\rm exp} \left(\frac{1}{2} \; A_{\infty}^{(q)}\right) \cdot W_{x, \infty}^{(q)}$ is well defined since, from point 1) of Lemma 1.1.3, the r.v. $A_{\infty}^{(q)}$ is $W_{x, \infty}^{(q)}$ a.s. finite. On the other hand, the measure

\noi $\dis \varphi_{q} (x) {\rm exp} \left(\frac{1}{2} \; A_{\infty}^{(q)}\right) \cdot W_{x, \infty}^{(q)}$ has infinite total mass from point 2) of Lemma 1.1.3.

\smallskip

\noi {\bf Proof of Lemma 1.1.4.}

\noi Let $q_{1}, q_{2} \in \mathcal{I}$. Then, from (1.1.7), we have for every $\Gamma_{u} \in b_{+} (\mathcal{F}_{u})$, with $u \le t$ :
\begin{eqnarray*} 
W_{x,\infty}^{(q_{1})} \left(\Gamma_{u} \varphi_{q_{1}} (x) {\rm exp} \left(\frac{1}{2} \;A_{t}^{(q_{1})}\right) \right) 
&=& W_{x} \big(\Gamma_{u} \varphi_{q_{1}} (X_{t})\big) \nonumber \\
&=& W_{x} \left( \Gamma_{u} \frac{\varphi_{q_{1}} (X_{t})}{\varphi_{q_{2}} (X_{t})} \; \varphi_{q_{2}} (X_{t}) \right) \nonumber \\�
&=& W_{x, \infty}^{(q_{2})} \left( \Gamma_{u} \varphi_{q_{2}} (x) \frac{\varphi_{q_{1}} (X_{t})} {\varphi_{q_{2}} (X_{t})} {\rm exp} \left( \frac{1}{2} \; A_{t}^{(q_{2})} \right) \right) \hspace*{1cm} (1.1.24)
\end{eqnarray*}

\noi Since the relation (1.1.24) takes place for every $\Gamma_{u} \in b_{+} (\mathcal{F}_{u})$ for any $u \le t$, we may replace $\Gamma_{u}$ by $\Gamma_{u} {\rm exp} \,(- \varepsilon A_{t}^{(q_{1}+q_{2})}) \quad (\varepsilon >0)$. We obtain :
\begin{gather*} 
W_{x, \infty}^{(q_{1})} \left[ \Gamma_{u} \varphi_{q_{1}} (x) \; {\rm exp} \left(\left(\frac{1}{2} - \varepsilon\right) A_{t}^{(q_{1})}\right) \cdot {\rm exp} \big(-\varepsilon A_{t}^{(q_{2})}\big) \right]  \nonumber \\
\quad = W_{x, \infty}^{(q_{2})} \left[ \Gamma_{u} \varphi_{q_{2}} (x) \frac{\varphi_{q_{1}} (X_{t})}{\varphi_{q_{2}} (X_{t})} \, {\rm exp} \left(\left(\frac{1}{2} - \varepsilon\right) A_{t}^{(q_{2})} \right) \cdot {\rm exp}  \big(- \varepsilon A_{t}^{(q_{1})}\big) \right] \hspace*{2cm}(1.1.25)
\end{gather*}

\noi However --- this is point 4) of Theorem 1.1.1 --- $\;\; |X_{t}| \dis \mathop{\longrightarrow}_{t \to \infty}^{} \infty$, $W_{x, \infty}^{(q_{2})}$ a.s. and the function $\dis x \longrightarrow \frac{\varphi_{q_{1}} (x)} {\varphi_{q_{2}} (x)}$ is bounded and tends to 1 when $|x| \longrightarrow \infty$. The dominated convergence Theorem - which we may apply thanks to Lemma 1.1.3 - implies then, by letting $t \longrightarrow \infty$ in (1.1.25) :
\begin{eqnarray*}
\lefteqn{  \varphi_{q_{1}} (x) W_{x, \infty}^{(q_{1})} \left[ \Gamma_{u} \, {\rm exp} \left(\left( \frac{1}{2} -\varepsilon \right)  A_{\infty}^{(q_{1})}\right) {\rm exp} \big(- \varepsilon A_{\infty}^{(q_{2})} \big) \right]   }\\
&= &\varphi_{q_{2}} (x) W_{x, \infty}^{(q_{2})} \left[ \Gamma_{u} \left({\rm exp} \left( \left( \frac{1}{2} - \varepsilon \right)\right. A_{\infty}^{(q_{2})} \right) \cdot {\rm exp} \big(-\varepsilon A_{\infty}^{(q_{1})} \big)\right] \hspace*{3cm} (1.1.26) 
\end{eqnarray*}

\noi Since (1.1.26) holds for every $\Gamma_{u} \in b_{+}( \mathcal{F}_{u})$ the monotone class Theorem implies that (1.1.26) is still true when we replace $\Gamma_{u} \in b_{+} (\mathcal{F}_{u})$ by $\Gamma \in b_{+} (\mathcal{F}_{\infty})$. It then remains to let $\varepsilon \longrightarrow 0$ and to use the monotone convergence Theorem to obtain : for every $\Gamma \in b_{+} (\mathcal{F}_{\infty})$ :
$$
\varphi_{q_{1}} (x) W_{x, \infty}^{(q_{1})} \left(\Gamma \; {\rm exp} \left(\frac{1}{2} \; A_{\infty}^{(q_{1})}\right) \right) = \varphi_{q_{2}} (x) W_{x, \infty}^{(q_{2})} \left(\Gamma  \, {\rm exp} \left(\frac{1}{2} \; A_{\infty}^{(q_{2})} \right) \right)
$$

\noi This is Lemma 1.1.4 and point 1) of Theorem 1.1.2.

\smallskip

\noi {\it iv)} \underline{We now show that $\bf W$ has infinite mass, but is $\sigma$-finite on $\mathcal{F}_{\infty}$ .}

\noi Firstly, it is clear, from point 2) of Lemma 1.1.3, that : 
$$
{\bf W} (1) = \varphi_{q} (0) W_{\infty}^{(q)} \left({\rm exp} \left(\frac{1}{2} \; A_{\infty}^{(q)} \right) \right) = + \infty \eqno(1.1.27)
$$

\noi On the other hand, from point 1) of Lemma 1.1.3, $A_{\infty}^{(q)} < \infty \quad W_{\infty}^{(q)}$ a.s. Hence :
$$ 1_{ A_{\infty}^{(q)} \le n } \uparrow 1 \qquad W_{\infty}^{(q)} \quad {\rm a.s. } $$

\noi Thus :
$$
{\bf W}(A_{\infty}^{(q)} \le n) = \varphi_{q} (0) \; W_{\infty}^{(q)} \left(\left({\rm exp} \left(\frac{1}{2} \; A_{\infty}^{(q)} \right) \right) \cdot 1_{ A_{\infty}^{(q)}  \le n} \right) \le \varphi_{q} (0) e^{\frac{n}{2}} \eqno(1.1.28) $$

\noi which proves that $\bf W$ is $(\Omega, \mathcal{F}_{\infty})$ $\sigma$-finite.

\smallskip

\noi {\it v)} \underline{We now show that, for every $\Gamma_{t} \in b_{+} (\mathcal{F}_{t}), \; {\bf W} (\Gamma_{t}) = 0$ or $+ \infty$.}

\noi By definition of $\bf W$, we have :
\begin{eqnarray*} 
{\bf W} (\Gamma_{t}) 
&=& \varphi_{q} (0) \; W_{\infty}^{(q)} \left(\Gamma_{t} \; {\rm exp} \left(\frac{1}{2} \; A_{\infty}^{(q)}\right) \right) \nonumber \\
&=& \varphi_{q} (0) \; W_{\infty}^{(q)} \left(\Gamma_{t} \; {\rm exp} \left(\frac{1}{2} \; A_{t}^{(q)}\right)  W_{X_{t}, \infty}^{(q)} \left({\rm exp} \left(\frac{1}{2} \; A_{\infty}^{q}\right) \right)\right) \hspace*{2cm}(1.1.29)
\end{eqnarray*} 

\noi from the Markov property. But, from Lemma 1.1.3, $\dis W_{x, \infty}^{(q)} \left({\rm exp} \left(\frac{1}{2} \; A_{\infty}^{(q)}\right) \right) = + \infty$ for every $x \in \mathbb{R}$. Thus, ${\bf W} (\Gamma_{t})$ equals $0$ or $+ \infty$ according to whether $W_{\infty}^{(q)} (\Gamma_{t})$ is $0$ or is strictly positive, i.e. according to whether $W(\Gamma_{t})$ equals $0$ or is strictly positive since, from (1.1.7) and (1.1.8), the probabilities $W$ and $W_{\infty}^{(q)}$ are equivalent on $\mathcal{F}_{t}$. 

\noi The careful reader may have been surprised about our use in the proof of Lemma 1.1.4 of the r.v. ${\rm exp} \big( - \varepsilon \; A_{t}^{(q_{1} + q_{2})}\big)$. This is purely technical and "counteracts" the fact that $\bf W$ takes only the values $0$ and $+ \infty$ on $\mathcal{F}_{t}$.

\smallskip

\noi We shall now give several other descriptions of the measure $\bf W$. In order to obtain these descriptions we use a particular case of Theorem 1.1.1, which shall play a key role in our study. This particular case is that of $q = \delta_{0}$ (or more generally $q = \lambda \delta_{0})$, the Dirac measure in $0$. We begin by recalling a result in this case. 

\smallskip

\noi {\bf 1.1.4}  \underline{Study of the canonical process under $W_{\infty}^{(\lambda \delta_{0})}$.}

\smallskip

\noi Theorem 1.1.5 below has been obtained in [RVY, II], Theorem 8, p. 339, with $\dis h^{+} (x) = h^{-} (x) = {\rm exp} \left(- \frac{\lambda x}{2}\right) \quad (\lambda, x \ge 0)$.

\smallskip

\noi {\bf Theorem 1.1.5.} {\it  \big(A particular case of Theorem 1.1.1, with $q = \lambda \delta_{0}$, hence $A_{t}^{(q)} = \lambda L_{t}, \; t \ge~0$, where $(L_{t}, \; t \ge 0)$ is the Brownian local time at $0$.\big)

\smallskip

\noi {\bf 1)} The function $\varphi_{\lambda \delta_{0}}$ defined by (1.1.10), (1.1.11) equals :
$$
\varphi_{\lambda \delta_{0}} (x) = |x| + \frac{2}{\lambda}\;;\; \quad {\it hence,} \quad \varphi_{\lambda \delta_{0}} (0) = \frac{2}{\lambda} \eqno(1.1.30) $$

\noi while the martingale $(M_{s}^{(\lambda \delta_{0})}, \; s \ge 0)$ \big(see (1.1.8)\big) equals :
$$
M_{s}^{(\lambda \delta_{0})} = \left(1 + \frac{\lambda}{2} |X_{s}|\right) \; {\rm exp} \left(- \frac{\lambda} {2} \; L_{s}\right) \eqno(1.1.31) $$

\noi {\bf 2)} Under $W_{\infty}^{(\lambda \delta_{0})}$ :

\smallskip

\noi {\it i)} The r.v. $g := {\rm sup}\{u \ge 0 \;;\; X_{u}=0\}$ is $W_{\infty}^{(\lambda \delta_{0})}$ a.s. finite and $L_{\infty} (=L_{g})$ has density :
$$
f_{L_{\infty}}^{W_{\infty}^{(\lambda \delta_{0})}} (l) = \frac{\lambda}{2} \; e^{- \frac{\lambda}{2}l} \;\;1_{[0, \infty[} \; (l) \eqno(1.1.32) $$

\noi {\it ii)} The processes $(X_{u}, \; u \le g)$ and $(X_{g +u}, u \ge 0)$ are independent.

\smallskip

\noi {\it iii)} The process $(X_{g +u}, u \ge 0)$ is distributed with $P_{0}^{(3, \; {\rm sym})}$ where :
$$
P_{0}^{(3, \; {\rm sym})} = \frac{1}{2} \; (P_{0}^{(3)} + \widetilde{P}_{0}^{(3)}) \eqno(1.1.33) $$

\noi with $P_{0}^{(3)}$ (resp. $\widetilde{P}_{0}^{(3)}$) denoting the law of 3-dimensional Bessel process (resp. its opposite) starting from 0.

\smallskip

\noi {\it iv)} Conditionally on $L_{\infty} (=L_{g}) =l, \; (X_{u}, \; u \le g)$ is a Brownian motion starting from $0$, considered until its local time at $0$ reaches level $l$, that is up to the stopping time : 
$$
\tau_{l} := {\rm inf} \{t \ge 0 \;;\; L_{t} > l\} \eqno(1.1.34) $$

\noi We write $W_{0}^{\tau_{l}}$ for the law of this process. }

\smallskip

$$
 {\bf 3)} \hspace*{4cm}W_{\infty}^{(\lambda \delta_{0})} = \frac{\lambda}{2} \int_{0}^{\infty} e^{- \frac{\lambda}{2} l} \big(W_{0}^{\tau_{l}} \circ P_{0}^{(3, \; {\rm sym})}\big) dl \hspace*{3cm}\eqno(1.1.35) $$

\noi In (1.1.35), we write $W_{0}^{\tau_{l}} \circ P_{0}^{(3, \; {\rm sym})}$ for the image of the probability $W_{0}^{\tau_{l}} \otimes P_{0}^{(3, \; {\rm sym})}$ by the concatenation operation $\circ$ : 
$$
\circ : \Omega \times \Omega \longrightarrow \Omega $$

\noi defined by (note that $X_{\tau_{l} (\omega)} =0$) :
$$ 
\hspace*{4cm} X_{t} (\omega \circ \widetilde{\omega}) =
\left\{ \begin{array}{rl}
 X_{t} (\omega) \hspace*{1cm} {\rm if} \qquad t \le \tau_{l} ( \omega)  \nonumber \\
 X_{t - \tau_{l} (\omega)} (\widetilde{\omega})   \hspace*{1cm}  {\rm if} \qquad t \ge \tau_{l} (\omega)\\ 
  \end{array}\right.  \hspace*{2cm}(1.1.36) $$

\noi Such a notation $\circ$ has been used by Biane-Yor [BY] to whom we refer the reader. Let us note that formula (1.1.35) is nothing else but the translation of the results of point 2) of Theorem 1.1.5.

\smallskip

\noi{\bf 1.1.5}  \underline{Some remarkable properties of $\bf W$.}

\smallskip

\noi We may now describe the measure {\bf W} independently from any penalisation. We introduce : 
$$
\begin{array} {rrl}
\hspace*{1cm}&g := {\rm sup} \{t \;;\; X_{t} =0\}, \qquad g_{a} := {\rm sup} \{t \;;\; X_{t} =a\} \hspace*{4.5cm}& (1.1.37) \\
&\sigma_{a,b} := {\rm sup} \big\{t, \; X_{t} \in [a,b] \big\} \qquad \, (a < b)\hfill & (1.1.38) \\
&\sigma_{a} := {\rm sup} \big\{t, \; X_{t} \in [-a,a]\big\} \qquad (a \ge 0) \hfill & (1.1.39)
\end{array}$$ 

\newpage

\noi {\bf Theorem 1.1.6.} {\it The following identities hold : 

\noi$ {\bf 1)}  \qquad {\bf W} = \int_{0}^{\infty} dl \; (W_{0}^{\tau_{l}} \circ P_{0}^{(3, {\rm sym})} )\hspace*{8.6cm} (1.1.40) $

\noi {\bf 2)} {\it i)} For every $(\mathcal{F}_{t}, \; t \ge 0)$ stopping time $T$ and for any r.v. $\Gamma_{T}$ which is positive and $\mathcal{F}_{T}$ measurable :
$$
{\bf W} (\Gamma_{T} \, 1_{g < T} \, 1_{T < \infty}) = W \big(\Gamma_{T} |X_{T}| 1_{T < \infty}\big) \eqno(1.1.41) $$

\noi {\it ii)} The law of $g$ under {\bf W} is given by : 
$$
{\bf W} (g \in dt) = \frac{dt}{\sqrt{ 2 \pi t}} \qquad (t \ge 0) \eqno(1.1.42) $$

\noi {\it iii)} Conditionally on $g=t$, the process $(X_{u}, \; u \le g)$ under $\bf W$ is a Brownian bridge with length $t$. We denote by $\Pi_{0,0}^{(t)}$ the law of this bridge.

\smallskip

\noi $ {\it iv)} \qquad {\bf W} = \int_{0}^{\infty} \frac{dt}{\sqrt{2 \pi t}} \; \big(\Pi_{0,0}^{(t)} \circ P_{0}^{(3, {\rm sym})}\big) \hspace*{8cm}(1.1.43)$

\noi {\it v)} For every previsible and positive process $(\phi_{s}, \; s \ge 0)$ we have :
$$
{\bf W} (\phi_{g}) = W \left(\int_{0}^{\infty} \phi_{s} d L_{s} \right) \eqno(1.1.44) $$

\noi {\bf 3)} {\it i)} For every $(\mathcal{F}_{t}, \; t \ge 0)$ stopping time $T$, the law under $\bf W$ of $L_{\infty} - L_{T}$, on $ T < \infty$ is given by :
\begin{eqnarray*}
{\bf W} (L_{\infty} - L_{T} \in dl, \; {T < \infty})
& = & W \big(T < \infty \big) \,  1_{[0, \infty]}(l) dl + W \big(|X_{T}| 1_{T < \infty} \big) \delta_{0} (dl) \\
& =  &W \big(T < \infty \big) \, 1_{[0, \infty]} (l) dl + {\bf W} (g \le T < \infty) \delta_{0} (dl) \end{eqnarray*}

\noi In particular, for $T=t$ : 
$$
{\bf W} (L_{\infty} - L_{t} \in dl) = 1_{[0, \infty[} (l)dl + \sqrt{\frac{2t}{\pi}} \; \delta_{0} (dl) \eqno(1.1.45) $$
\begin{gather*}
 \hspace*{-1cm} {\it ii)} \; \hbox{For every} \;\; l >0, \; \hbox{conditionally on} \; L_{\infty} - L_{T} =l, \; T<\infty, \; (X_{u}, \; u \le T) \; \hbox{is a Brownian motion} \nonumber \\
\hbox{\hspace*{-1,1cm}  indexed by} \; [0,T] \hspace*{11,6cm}(1.1.46)
\end{gather*}
\noi {\it iii)} The density of $(g, L_{\infty})$ under $\bf W$ equals :
$$
f_{g, L_{\infty}}^{\bf W} (u,l) = \frac{l \; {\rm exp} \dis \; \left(- \frac{l^{2}}{2u}\right)}{\sqrt{2 \pi u^{3}}} \; 1_{[0, \infty[} (u) 1_{[0, \infty[} (l) \eqno(1.1.47)$$       }

\noi {\bf Remark 1.1.7.} 

\noi {\bf 1)} We deduce from formulae (1.1.43) and (1.1.17) that :
\begin{eqnarray*} 
\varphi_{q} (0)
&=& {\bf W} \left({\rm exp} \left( - \frac{1}{2} \; A_{\infty}^{(q)}\right) \right) \nonumber \\
&=& \int_{0}^{\infty} \frac{dt}{\sqrt{2 \pi t}} \; \Pi_{0,0}^{(t)} \left({\rm exp} \left(- \frac{1}{2} \; A_{t}^{(q)} \right) \right) \cdot P_{0}^{(3, {\rm sym)}} \left({\rm exp} \left(- \frac{1}{2} \; A_{\infty}^{(q)} \right) \right) \hspace*{1cm} (1.1.48)
\end{eqnarray*}

\noi {\bf 2)} It is proven in Biane-Yor \big([BY], see also [Bi]\big) that : 
$$
\int_{0}^{\infty} dl \; W_{0}^{\tau_{l}} = \int_{0}^{\infty} \frac{dt}{\sqrt{2 \pi t}} \; \Pi_{0,0}^{(t)} $$

\noi Thus, from this identity, we deduce easily that (1.1.40) implies (1.1.43).

\smallskip

\noi {\bf 3)} Formula (1.1.41) \big(see also formulae (1.1.52), (1.1.54), (1.1.55), (1.1.56), (1.1.73)\big) yields a "representation" of the Brownian sub-martingale $\big(|X_{t}|, \; t \ge 0\big)$ in terms of the increasing process $(1_{g \le t}, \; t \ge 0)$. \big(By a "representation" of a $\big(P, \; (\mathcal{F}_{t}, \; t \ge 0)\big)$ submartingale $(Z_{t}, \; t \ge 0)$, we mean a couple $\big(Q, (C_{t}, \; t \ge 0)\big)$ where $Q$ is a $\sigma$-finite measure and $(C_{t}, \; t \ge 0)$ is a increasing process such that, for every $\Gamma_{t} \in b (\mathcal{F}_{t}) \;:\; Q (\Gamma_{t} \cdot C_{t}) = E_{P} [\Gamma_{t} \cdot Z_{t}]$.\big) Here, $({\bf W}, 1_{g \le t}, t \geq 0)$ is a representation of the submartingale $\big(|X_{t}|, \; t \ge 0)$.

\smallskip 

\noi Before we prove Theorem 1.1.6, we present a slightly different version of it. We shall not prove this version, whose proof relies on close arguments to those we needed to obtain Theorem 1.1.6.

\smallskip

\noi {\bf Theorem 1.1.8.} {\it Let $a \ge 0$ ; the following formulae hold : 

\noi {\bf 1)} For every $(\mathcal{F}_{t}, \; t \ge 0)$ stopping time $T$ and for every r.v. $\Gamma_{T}$ positive and $\mathcal{F}_{T}$ measurable :
$$
{\bf W} (\Gamma_{T} \, 1_{(\sigma_{a} < T < \infty)} \big) = W \big(\Gamma_{T} \big(|X_{T}| -a \big)_{+} 1_{T < \infty}\big) \eqno(1.1.49)$$

\noi {\bf 2)} $\; {\it i)} \qquad {\bf W} (\sigma_{a} \in dt) = \frac{ e^{- \frac{a^{2}}{2t} } }{\sqrt{2 \pi t}} \; dt \quad (t \ge 0) \hspace*{7,4cm}(1.1.49')$

\noi \quad\; ${\it ii)} \qquad \, \;{\bf W} = \int_{0}^{\infty} \frac{ dt }{\sqrt{2 \pi t}} \; e^{- \frac{a^{2}}{2t}} \frac{1}{2} \big(\Pi_{0,a}^{(t)} \circ P^{(a,3)} + \Pi_{0,-a}^{(t)} \circ P^{(-a,3)}\big) \hspace*{4,1cm}(1.1.50)$

\noi where $\Pi_{\alpha, \beta}^{(t)}$ denotes the law of the Brownian bridge of length $t$ starting from $\alpha$ and ending in $\beta$ and where $P^{(a,3)}$ (resp. $P^{(-a,3)}$) is the law of the process $(a+ R_{t}, \; t \ge 0)$ \big(resp. $(-a-R_{t},$ $t \ge 0)$\big) where $(R_{t},\; t \ge 0)$ is a 3-dimensional Bessel process starting from 0. In particular, the law of $(X_{u}, \; u \le \sigma_{a})$, conditionally on $\sigma_{a} = t$ is $\dis \frac{1}{2} (\Pi_{0,a}^{(t)} + \Pi_{0,-a}^{(t)})$

\smallskip

{\it iii)} For every positive and previsible process $(\phi_{u}, u \ge 0)$, we have :
$$
{\bf W} (\phi_{\sigma_{a}}) = W \left(\int_{0}^{\infty} \phi_{u} \; d_{u} (\widetilde{L}_{u}^{a}) \right) \eqno(1.1.51) $$

\noi with $\widetilde{L}_{u}^{a} := \dis \frac{1}{2} (L_{u}^{a} + L_{u}^{-a})$. }

\smallskip

\noi We note that points 1) and 2) of Theorem 1.1.6 are particular cases of the corresponding ones in Theorem 1.1.8 when $a=0$. On the other hand, in the same spirit as for (1.1.49) we have, with the same kind of notation :
$$
W \big(\Gamma_{T} (X_{T}-a)_{+} 1_{T<\infty}\big) = {\bf W}^{+} (\Gamma_{T} 1_{g_{a} < T<\infty}) \eqno{(1.1.52)_{+}} $$

$$
\hspace*{-0,5cm} W \big(\Gamma_{T} (X_{T}-a)_{-} 1_{T<\infty}\big) = {\bf W}^{-} (\Gamma_{T} 1_{g_{a} < T<\infty})\eqno{(1.1.52)_{-}} $$

\noi where : 
$$
{\bf W}^{+} = \frac{1}{2} \int_{0}^{\infty} \frac{dt}{\sqrt{2 \pi t}} \; \Pi_{0,0}^{(t)} \circ P_{0}^{(3)} \eqno{(1.1.53)_{+}} $$

$$
\hspace*{-0,3cm}{\bf W}^{-} = \frac{1}{2} \int_{0}^{\infty} \frac{dt}{\sqrt{2 \pi t}} \; \Pi_{0,0}^{(t)} \circ \widetilde{P}_{0}^{(3)} \eqno{(1.1.53)_{-}} $$

\noi Adding (1.1.52)$_{+}$ and (1.1.52)$_{-}$ yields :
$$
W(\Gamma_{T} |X_{T}-a| 1_{T<\infty}) = {\bf W} (\Gamma_{T} 1_{g_{a}<T<\infty}) \eqno(1.1.54) $$

\noi and also, with $a<b$ : 
$$
W\big(\Gamma_{T} \big((X_{T}-b)_{+} + (a - X_{T})_{+}\big) 1_{T<\infty}\big) = {\bf W} (\Gamma_{T} 1_{\sigma_{a,b} < T < \infty}) \eqno(1.1.55) $$

\noi and
$$
W�\big(\Gamma_{T} \big(|X_{T}|-a\big)_{+} 1_{T < \infty}\big) = {\bf W} (\Gamma_{T} 1_{\sigma_{a} < T < \infty}) \qquad (a \ge 0) \eqno(1.1.56) $$

\noi {\bf Proof of Theorem 1.1.6.}

\noi Here is the plan of our proof. We shall use formula (1.1.16) with $q = \delta_{0}$ :
$$ {\bf W} = \varphi_{\delta_{0}} (0) \; e^{\frac{1}{2} \; L_{\infty}} \; . \; W_{\infty}^{(\delta_{0})} = 2 \; e^{\frac{1}{2} \; L_{\infty}} \; . \; W_{\infty}^{(\delta_{0})} \eqno(1.1.57) $$

\noi (from (1.1.30)), as well as the properties of $W_{\infty}^{(\delta_{0})}$ recalled in Theorem 1.1.5.

\bigskip

\noi {\it i)} \underline{ We prove (1.1.40) }.

\noi Let $F$ and $G$ be two positive functionals. We have, from (1.1.57) :
\begin{eqnarray*}
\lefteqn{{\bf W} \big(F (X_{s}, \; s \le g) \cdot G (X_{g+s}, \; s \ge 0)\big) } \nonumber \\
&=&2 \, W_{\infty}^{(\delta_{0})} \big(e^{\frac{1}{2} \, L_{\infty}} F (X_{s}, \; s \le g) \, G (X_{g+s}, \; s \ge 0)\big)\\
&= & 2 \, W_{\infty}^{(\delta_{0})} \big(e^{\frac{1}{2} \, L_{g}} F (X_{s}, \; s \le g) \, G (X_{g+s}, \; s \ge 0)\big)\\ 
&&(\text{since} \;\; L_{\infty} = L_{g}) \\
&=& 2 \, W_{\infty}^{(\delta_{0})} \big(e^{\frac{1}{2} \, L_{g}} F (X_{s}, \; s \le g)\big) \cdot P_{0}^{(3, {\rm sym})} \big(G (X_{s}, \; s \ge 0)\big) \\
&& \text{\big(from Point 2)ii) of Theorem 1.1.5 and from (1.1.33)\big)} \\
&= & \left(2 \int_{0}^{\infty} W_{\infty}^{(\delta_{0})} \big(e^{\frac{1}{2} \, L_{g}} F(X_{s}, \; s \le g) \big|{L_{g}=l}\big) \frac{1}{2} \; e^{- \frac{l}{2}} dl\right) \cdot P_{0}^{(3, {\rm sym})} \big( G (X_{s}, \; s \ge 0)\big) \\
&& \text{\big(from (1.1.32)\big)} \\
&=& \left(2 \int_{0}^{\infty} e^{ \frac{l}{2}} W \big(F (X_{s}, \; s \le \tau_{l}) \big) \frac{1}{2} \; e^{-\frac{l}{2}} dl \right)  \cdot P_{0}^{(3, {\rm sym})} \big( G (X_{s}, \; s \ge 0)\big) \\
&= & \int_{0}^{\infty} dl (W_{0}^{\tau_{l}} \circ  P_{0}^{(3, {\rm sym})}) \; \big(F (X_{s}, \; s \le g) \cdot G (X_{g+s}, \; s \ge 0) \big)
\end{eqnarray*}

\noi from point 2, {\it iv)} of Theorem 1.1.5.

\smallskip

\noi {\it ii)} \underline{We now prove (1.1.41)}.

\noi For this purpose, we apply formula 
 with $q= \lambda \delta_{0}$. Thus :
$$
A_{t}^{(q)} = \lambda L_{t} \quad \hbox{and, from (1.1.30)}, \quad \varphi_{\lambda \delta_{0}} (x)= \frac{2}{\lambda} + |x|.$$

\noi Thus, from (1.1.7), (1.1.30), (1.1.31), (1.1.16) and Doob's optional stopping Theorem :
\begin{eqnarray*}
W \left(\Gamma_{T} \left( \frac{2}{\lambda} + |X_{T}|\right) 1_{T < \infty}\right)
&=& \frac{2}{\lambda} \, W_{\infty}^{(\lambda \delta_{0})} \big(e^{\frac{\lambda}{2} \, L_{T}}\Gamma_{T}  \, 1_{T < \infty}\big) \nonumber \\
&=& {\bf W} (\Gamma_{T} \, 1_{g \le T < \infty}\big) + {\bf W} \big(\Gamma_{T} \, 1_{g > T} e^{- \frac{\lambda}{2} (L_{\infty} -L_{T})} \big) \quad (1.1.58)
\end{eqnarray*}

\noi We then let $\lambda \longrightarrow \infty$ in (1.1.58) and note that $L_{\infty} - L_{T} > 0$ on $g>T$. The monotone convergence Theorem implies : 
$$
W\big(\Gamma_{T} |X_{T}| 1_{T < \infty}\big) = {\bf W} (\Gamma_{T} \, 1_{g \le T < \infty}) $$

\noi This is precisely relation (1.1.41). Relation (1.1.42) is an easy consequence of (1.1.41).

\smallskip

\noi {\it iii)} \underline{We prove (1.1.45) and (1.1.46)}.

\noi We note that (1.1.41) and (1.1.58) imply :
\begin{eqnarray*}
\frac{2}{\lambda} \, W (\Gamma_{T} \, 1_{T < \infty}) 
&=& {\bf W} \left(\Gamma_{T} \, 1_{g > T} \; {\rm exp} \left(- \frac{\lambda}{2} (L_{\infty} - L_{T}) \right) \right) \hspace*{4cm}(1.1.59)\\
&=& W(\Gamma_{T} \, 1_{T < \infty}) \left(\int_{0}^{\infty} e^{- \frac{\lambda}{2}l} dl \right) \nonumber
\end{eqnarray*}

\noi Thus, by injectivity of the Laplace transform, for every function $\psi$ : $\mathbb{R}_{+} \longrightarrow \mathbb{R}_{+}$ Borel and integrable :
$$
W(\Gamma_{T} \, 1_{T < \infty}) \; \left(\int_{0}^{\infty} \psi (l)dl\right) = {\bf W} \big(\Gamma_{T} \, \psi (L_{\infty}-L_{T}) 1_{g >T}\big) \eqno(1.1.60) $$

\noi and
$$
W\big(|X_{T}| 1_{T < \infty}\big) = {\bf W} (g \le T < \infty)= {\bf W} (L_{\infty} -L_{T}=0, \; T < \infty) \eqno(1.1.61) $$

\noi In other terms, we have : 
$$
{\bf W} (L_{\infty} - L_{T} \in dl, \; {T < \infty}) = W \big( T < \infty \big) 1_{[0, \infty[} (l)dl + W \big(|X_{T}| 1_{T < \infty}\big) \delta_{0} (dl) $$

\noi and, under ${\bf W}$, conditionally on $L_{\infty} - L_{T} = l \;\; (l >0), \;\; (X_{u}, \; u \le T)$ is a Brownian motion indexed by $[0,T]$. This is (1.1.45) and (1.1.46).

\noi {\it iv)} \underline{We now prove point 2, {\it iii)} of Theorem 1.1.6.}

\noi For this purpose, we write (1.1.41), choosing for $\Gamma_{t}$ a r.v. of the form $\Phi_{g^{(t)}}$, where $(\Phi_{u}, \; u \ge 0)$ is a previsible positive process, and where $g^{(t)} := {\rm sup} \{ s \le t, \; X_{s} =0\}$.

\noi The RHS of (1.1.41) becomes, with $T=t$ :
\begin{eqnarray*}
\hspace*{1cm}W\big(|X_{t}| \Phi_{g^{(t)}} \big) 
& =& W \left(\int_{0}^{t} \Phi_{s} \, dL_{s}\right) \\
&&\text{\Big(from the balayage formula \big(cf [ReY], Chap. VI, p. 260\big)\Big)} \\
&=& \int_{0}^{t} W \big(\Phi_{s} |X_{s}=0\big) W (dL_{s}) \\
&=& \int_{0}^{t} W \big(\Phi_{s} | X_{s}=0\big) \frac{ds}{\sqrt{2 \pi s}} \hspace*{5,5cm}(1.1.62)\\
&&\Biggl({\rm since} \;\; W(L_{s})
= W \big(|X_{s}|\big) = \sqrt{\frac{2 s}{\pi}} \Biggl)
\end{eqnarray*}

\noi The LHS of (1.1.41) writes : 
\begin{eqnarray*}
\hspace*{1cm}{\bf W} (\Phi_{g^{(t)}} 1_{g \le t}) 
&=& {\bf W} (\Phi_{g} \, 1_{g \le t}) \\
&& \text{ \big(since} \; g=g^{(t)} \; \text{on the set} \; \{g \le t\}\big) \\
&=& \int_{0}^{t} {\bf W} \big(\Phi_{g} | g=s\big) \frac{ds}{\sqrt{2 \pi s}} \hspace*{4cm} (1.1.63) \\
&& \text{from (1.1.42). Thus :} \\
\int_{0}^{t} W \big(\Phi_{s} | X_{s} = 0\big) \frac{ds}{\sqrt{2 \pi s}} 
&=& \int_{0}^{t} {\bf W} \big(\Phi_{g}| g=s\big) \frac{ds}{\sqrt{2 \pi s}}
\end{eqnarray*}

 \noi This relation implies ${\bf W} \big(\Phi_{s} | g=s\big) =W\big(\Phi_{s}|X_{s}=0\big)$, i.e. point 2, {\it iii)} of Theorem 1.1.6.

\noi We also note that we deduce from the equality between (1.1.62) and (1.1.63) : 
$$
W \left( \int_{0}^{t} \Phi_{s} \, dL_{s}\right) = \int_{0}^{t} {\bf W} \big(\Phi_{g} | g=s\big) \frac{ds}{\sqrt{2 \pi s}} $$

\noi that :
\begin{eqnarray*} 
W \left(\int_{0}^{\infty} \Phi_{s} \, dL_{s}\right)
&=& \int_{0}^{\infty} {\bf W} \big(\Phi_{g} | g=s\big) {\bf W} (g \in ds) \\
&=& {\bf W} (\Phi_{g}) \hspace*{9cm}(1.1.64)
\end{eqnarray*}

\noi i.e. point 2, {\it v)} of Theorem 1.1.6.

\smallskip

\noi {\it v)} \underline{We now prove point 2, {\it iv)} of Theorem 1.1.6.}

\noi We obtain, with the help of (1.1.57), for two positive functionals $F$ and $G$ :
\begin{eqnarray*}
\lefteqn{{\bf W} \big(F (X_{s}, \; s \le g) \; G (X_{g+s}, \; s \ge 0)\big) } \\
&=& 2 \, W_{\infty}^{(\delta_{0})} \big(F (X_{s}, \; s \le g) \, e^{\frac{1}{2} \, L_{g}} G (X_{g+s, \, s \ge 0}) \big)\\
&=& 2 \, W_{\infty}^{(\delta_{0})} \big(F (X_{s}, \; s \le g) \, e^{\frac{1}{2} \, L_{g}}\big) \; P_{0}^{(3, {\rm sym})} \big(G (X_{s}, \; s \ge 0)\big) \\
\lefteqn{\hbox{(from point 2 {\it ii)} and 2 {\it iii)} of Theorem 1.1.5)}} \\
&=& {\bf W} \big(F (X_{s}, \; s \le g)\big) \cdot P_{0}^{(3, {\rm sym})} \big(G (X_{s}, \; s \ge 0)\big) \\
\lefteqn{\hbox{\big(using once again (1.1.57)\big)}}  \\
&=& \left(\int_{0}^{\infty} {\bf W} \big(F (X_{s}, \; s \le g) | g=t\big) \frac{dt}{\sqrt{2 \pi t}} \right)  \cdot P_{0}^{(3, {\rm sym})} \big(G (X_{s}, \; s \ge 0)\big) \\
\lefteqn{\hbox{\big(from (1.1.42)\big)}} \\
&=& \int_{0}^{\infty} \frac{dt}{\sqrt{2 \pi t}} \; \Pi_{0,0}^{(t)} \big(F (X_{s}, \; s \le t)\big) \cdot P_{0}^{(3, {\rm sym})} \big(G (X_{s}, \; s \ge 0)\big) \\
\lefteqn{\hbox{(from point 2 {\it iii)} of Theorem 1.1.6)}} \\
&=& \int_{0}^{\infty} \frac{dt}{\sqrt{2 \pi t}}  \big(\Pi_{0,0}^{(t)} \circ P_{0}^{(3, {\rm sym})}\big)\;\; \big(F (X_{s}, \; s \le g) \; G (X_{g+s}, \; s \ge 0)\big)\cdot
\end{eqnarray*}

\noi {\it vi)} Formula (1.1.47) is a consequence of (1.1.42), (1.1.43) and the fact that : 

\noi Under $\Pi_{0,0}^{(t)}, \; L_{t}$ is distributed as $\sqrt{2 t \mathfrak{e}}$, where $\mathfrak{e}$ is a standard exponential r.v.

\noi {\bf Remark 1.1.9.}

\noi {\bf 1)} We have, from (1.1.16) and Theorem 1.1.5 : 
$$
\frac{\lambda}{2} \; e^{- \frac{\lambda}{2} \, L_{\infty}} \cdot {\bf W} = W_{\infty}^{(\lambda \delta_{0})} \eqno(1.1.65) $$

\noi But, from Theorem 1.1.1, under $W_{\infty}^{(\lambda \delta_{0})}$ :
$$
X_{t} = B_{t} + \int_{0}^{t} \frac{{\rm sgn} \, X_{s}}{ \frac{2}{\lambda} + |X_{s}|} \; ds $$

\noi Hence \big(see [RY, M], Chap. 4\big) : $ \dis W_{\infty}^{(\lambda \delta_{0})} \mathop{\longrightarrow}_{\lambda \to \infty}^{} P_{0}^{(3,\,{\rm sym})}$.

\noi Thus  $\quad \dis \frac{\lambda}{2} \; (e^{- \frac{\lambda}{2} \, L_{\infty}}) \, {\bf W}  \mathop{\longrightarrow}_{\lambda \to \infty}^{}  P_{0}^{(3,\,{\rm sym})} \hspace*{8,2cm}(1.1.66)$

\noi This convergence holds in the sense of weak convergence with respect to the topology of uniform convergence on compacts in $\mathcal{C} \big([0, \infty[ \to \mathbb{R})$.

\smallskip

\noi {\bf 2)} Formula (1.1.41) may be proven in a different manner than by the way we have  indicated. Indeed, from (1.1.57) (where, to simplify, we choose $T=t$)
\begin{eqnarray*}
\hspace*{0,7cm}{\bf W} (\Gamma_{t} \, 1_{g \le t})
&=& 2 \, W_{\infty}^{(\delta_{0})} (\Gamma_{t} \, 1_{g \le t} \, e^{\frac{1}{2} \, L_{\infty}} ) \\
&=& 2 \, W_{\infty}^{(\delta_{0})} (\Gamma_{t} \, 1_{g \le t} \, e^{\frac{1}{2} \, L_{t}} ) \\
&& \big(\text{since} \;  L_{\infty} = L_{t} \; \text{ on the set} \; (g \le t)\big) \\
&=&2 \, W_{\infty}^{(\delta_{0})} (\Gamma_{t} \, e^{\frac{1}{2} \, L_{t}} W_{\infty}^{(\delta_{0})} \big(1_{g \le t} | \mathcal{F}_{t}\big) \big) \hspace*{5cm}(1.1.67)
\end{eqnarray*}

\noi But
\begin{eqnarray*}
\hspace*{2cm}W_{\infty}^{(\delta_{0})} \big(1_{g \le t} | \mathcal{F}_{t} \big)
&=& W_{\infty}^{(\delta_{0})} \; \big(T_{0} \circ \theta_{t} = \infty | \mathcal{F}_{t} \big) \\
&=& W_{X_{t}, \infty}^{(\delta_{0})} \; (T_{0} = \infty) \hspace*{6cm}(1.1.68)
\end{eqnarray*}

\noi with $T_{0} = {\rm inf} \{ t \ge 0 \;;\; X_{t} =0\}$, by Markov property. But, from (1.1.14), the scale function of the process $(X_{t}, \; t \ge 0)$ under $(W_{x, \infty}^{(\delta_{0})}, \; x \in \mathbb{R})$ equals :
$$
\gamma_{\delta_{0}} (x) = \frac{x}{2 \big(2 + |x|\big)} \eqno(1.1.69) $$

\noi We deduce from (1.1.69) :
$$
W_{x, \infty}^{(\delta_{0})} (T_{0} = \infty) = \frac{|x|}{2+|x|} \eqno(1.1.70) $$

\noi Plugging (1.1.70) and (1.1.68) in (1.1.67), we obtain :
\begin{eqnarray*}
{\bf W} (\Gamma_{t} \, 1_{g \le t}) 
&=& 2 \, W_{\infty}^{(\delta_{0})} \left( \Gamma_{t} \, e^{\frac{1}{2} \, L_{t}} \, \frac{|X_{t}|}{2 + |X_{t}|} \right) \\
&=& 2 \, W \left(\Gamma_{t} \; e^{\frac{1}{2} \, L_{t}} \;  \frac{|X_{t}|}{2 + |X_{t}|} \; \frac{2+|X_{t}|}{2} \; e^{- \frac{1}{2} \, L_{t}}\right) \\
&& \textrm{ (from (1.1.31), with $\lambda = 1$, and (1.1.7))} \\
&=& W \big(\Gamma_{t} |X_{t}|\big)
\end{eqnarray*}

\noi Formulae (1.1.54), (1.1.56), (1.1.57) may be proven following the same arguments.

\smallskip

\noi {\bf 3)} Let $q \in \mathcal{I}$ such that the convex hull of its support equals the interval $[a,b] \; (a \le b)$. From (1.1.7) and (1.1.6) we have : 
\begin{eqnarray*} 
\hspace*{1cm}W\big(\varphi_{\lambda q} (X_{t}) \cdot \Gamma_{t}\big)
&=& \varphi_{\lambda q} (0) \, W_{\infty}^{(\lambda q)} \big(\Gamma_{t} \, e^{- \frac{\lambda}{2} \, A_{t}^{(q)}}\big) \\
&=& {\bf W} \big(\Gamma_{t} \, e^{- \frac{\lambda}{2} \, (A_{\infty}^{(q)}- A_{t}^{(q)})}\big)  \\
&=& {\bf  W} (\Gamma_{t} \, 1_{\sigma_{a,b} \le t}) + {\bf W} \big(\Gamma_{t} e^{- \frac{\lambda}{2} \, (A_{\infty}^{(q)}- A_{t}^{(q)})} 1_{\sigma_{a,b} > t}\big) \hspace*{2cm}(1.1.71)
\end{eqnarray*}

\noi On the other hand, we have proven in [RY, IX] \big(see also [RY, M], Chap. 2\big) that there exists, for every $x \in \mathbb{R}$, a positive and $\sigma$-finite measure ${\bm \nu}_{x}^{(q)}$ such that :
$$
\int_{0}^{\infty} e^{- \frac{\lambda y}{2}} {\bm \nu}_{x}^{(q)} (dy) = \varphi_{\lambda q} (x) \eqno(1.1.72) $$

\noi It remains to let $\lambda \to \infty$ in (1.1.72) to obtain, since $A_{\infty}^{(q)} - A_{t}^{(q)} > 0$ on the set $(\sigma_{a\; b} > t)$ :
$$
W \big(\Gamma_{t} \, {\bm \nu}_{X_{t}}^{(q)} (\{0\})= {\bf W} (\Gamma_{t} \, 1_{\sigma_{a,b} \le t}) \eqno(1.1.73)$$

\noi Hence, ${\bm \nu}_{x}^{(q)} (\{0\})$ depends only on supp$(q)$ and $\big({\bm \nu}_{X_{t}}^{(q)} \, (\{0\}), \; t \ge 0\big)$ is a sub-martingale. Formula (1.1.55) (with $T=t)$ is a particular case of (1.1.73), since :
$$
{\bm \nu}_{x}^{(\delta_{a} + \delta_{b})} (\{0\}) = (x-b)_{+} + (a-x)_{+} \eqno(1.1.74) $$

\noi \big(see [RY, IX]\big).

\smallskip

\noi {\bf 1.1.6.} \underline{Another approach to Theorem 1.1.6.}

\smallskip

\noi Let, for $q \in \mathcal{I}$, the probability $W_{\infty}^{(q)}$ be defined by (1.1.7). Then :
$$
{W_{\infty}^{(q)}}_{|_{\mathcal{F}_{t}}} = \frac{\varphi_{q} (X_{t})}{\varphi_{q} (0)} \; e^{- \frac{1}{2} \, A_{t}^{(q)}} {\cdot W}_{|_{\mathcal{F}_{t}}} \eqno(1.1.75) $$

\noi In Theorem 1.1.2, we have defined the measure $\bf W$ from the formula : 
$$
{\bf W} = \varphi_{q} (0) \; e^{\frac{1}{2} \, A_{\infty}^{(q)}} W_{\infty}^{(q)} \eqno(1.1.76) $$

\noi then, we have shown that : 
$$
\hspace*{1,5cm}{\bf W} = \int_{0}^{\infty} \frac{dt}{\sqrt{2 \pi t}} \; (\Pi_{0,0}^{(t)} \circ P_{0}^{(3,\, {\rm sym})})  \eqno(1.1.77)$$

\noi \big(cf Theorem 1.1.6, relation (1.1.43)\big). We now "forget" our previous results and proceed in a reverse way. For this purpose, we define, for the time being, the measure :
$$
\hspace*{1,4cm}\widetilde{\bf W} = \int_{0}^{\infty} \frac{dt}{\sqrt{2 \pi t}} (\Pi_{0,0}^{(t)} \circ P_{0}^{(3,\, {\rm sym})}) \eqno(1.1.78) $$

\noi We shall show that, for every $q \in \mathcal{I}$ : 
\begin{equation*}
\frac{1}{\varphi_{q} (0)} \; e^{- \frac{1}{2} \, A_{\infty}^{(q)}} \cdot \widetilde{\bf W} = W_{\infty}^{(q)}
\end{equation*}

\noi {\bf Theorem 1.1.10.}

\noi {\it Let $\widetilde{\bf W}$ be defined by (1.1.78) and $W_{\infty}^{(q)}$ be defined by (1.1.75). Then, for every $q \in \mathcal{I}$ : 
$$
\frac{1}{\varphi_{q} (0)} \; e^{- \frac{1}{2} \, A_{\infty}^{(q)}} \cdot \widetilde{\bf W} = W_{\infty}^{(q)} \eqno(1.1.79) $$ }

\newpage

\noi {\bf Proof of Theorem 1.1.10.}

\noi We compute the value of $W_{\infty}^{(q)}$ when integrating the following general class of functionals which are $\mathcal{F}_{t}$-measurable and positive : 
$$
F(X_{u}, \; u \le g^{(t)}) \cdot G (X_{g^{(t)}+ u}, \; u \le t-g^{(t)}) \eqno(1.1.80) $$

\noi We have : 
\begin{eqnarray*} 
\lefteqn{W_{\infty}^{(q)} \big(F (X_{u}, \; u \le g^{(t)}) \; G (X_{g^{(t)} +u} \;;\; u \le t-g^{(t)})\big)} \nonumber \\
&=& \frac{1}{\varphi_{q} (0)} \; W \left[ F (X_{u}, \; u \le g^{(t)}) \; G (X_{g^{(t)} +v} \;;\; v \le t-g^{(t)}) \, {\rm exp} \left(- \frac{1}{2} \, A_{t}^{(q)} \right) \varphi_{q} (X_{t}) \right]  \\
\lefteqn{\rm \big(from (1.1.75)\big)}  \\
&=& \frac{1}{\varphi_{q} (0)} \; W \left[ F(X_{u}, \; u \le g^{(t)}) \, {\rm exp} \left(- \frac{1}{2} \,A_{g^{(t)}}^{(q)} \right) \cdot G (X_{g^{(t)}+u}, \; u \le t-g^{(t)}) \right. \\
&& \hspace*{2cm}\left. \cdot \varphi_{q} (X_{t})\; {\rm exp} \left( - \frac{1}{2} \big(A_{t}^{(q)} - A_{g^{(t)}}^{(q)} \big) \right) \right] \hspace*{5cm}(1.1.81)
\end{eqnarray*}

\noi We now consider the probability $W$ restricted to $\mathcal{F}_{t}$, denoted as $W^{(t)}$, which we disintegrate with respect to the law of $g^{(t)}$ : 
\begin{eqnarray*} 
W^{(t)} 
&= & \int_{0}^{t} \frac{du}{\pi \sqrt{u (t-u)}} \; \big(\Pi_{0,0}^{(u)} \circ M^{(t-u, {\rm sym})}\big) \hspace*{5cm}(1.1.82) \\
&& \text{with :} \\
W(g^{(t)} \in du) 
&=& \frac{du}{\pi \sqrt{u(t-u)}} \qquad u \le t
\end{eqnarray*}

\noi and where $\Pi_{0,0}^{(u)}$ denotes the law of the Brownian bridge with length $u$ and $M^{(t, \, {\rm sym})}$ is the law of a symmetric Brownian meander of length $t$. Thus, (1.1.81) becomes :
\begin{gather*}
W_{\infty}^{(q)} \big[ F (X_{u}, \; u \le g^{(t)}) \; G (X_{g^{(t)}+v}, \; v \le t -g^{(t)})\big] \\
\qquad = \frac{1}{\varphi_{q}(0)} \; \int_{0}^{t} \frac{du}{\pi \sqrt{u (t-u)}} \; \Pi_{0,0}^{(u)} \big(F (X_{v}, \; v \le u) e^{- \frac{1}{2} \, A_{u}^{(q)}} \big) \\
\hspace*{2,5cm} \cdot M^{(t-u, \, {\rm sym})} \big(\varphi_{q} (X_{t-u}) e^{- \frac{1}{2} \, A_{t-u}^{(q)}} \cdot G(X_{l}, \; l \le t-u) \big)
\end{gather*}

\noi Using now Imhof's relation \big(see [RY, M], Chap. 1, Item $G$\big) :
$$
M^{(t, \, {\rm sym})} = \sqrt{\frac{\pi t}{2}} \; \frac{1}{|X_{t}|} \; P_{0}^{(3, \, {\rm sym})} |_{\mathcal{F}_{t}}\eqno(1.1.83) $$

\noi we obtain :
\begin{eqnarray*} 
\lefteqn{W_{\infty}^{(q)} \big[F (X_{u}, \; u \le g^{(t)}) \, G (X_{g^{(t)} + v}, \; v \le t-g^{(t)})\big] }\nonumber \\
&=& \frac{1}{\varphi_{q}(0)} \; \int_{0}^{t} \frac{du}{\pi \sqrt{u(t-u)}} \; \Pi_{0,0}^{(u)} \big(F (X_{v}  \; v \le u) \, e^{- \frac{1}{2} \, A_{u}^{(q)}} \big) \nonumber \\
&& \qquad \cdot P_{0}^{(3, {\rm sym})} \left(\varphi_{q} (X_{t-u})  \frac{G(X_{l}, \; l \le t-u)}{|X_{t-u}|}  \; \sqrt{\frac{\pi}{2} (t-u)} \; e^{- \frac{1}{2} \, A_{t-u}^{(q)}} \right) \hspace*{1,5cm}(1.1.84)
\end{eqnarray*}

\noi We observe that the factor $\sqrt{t-u}$ simplifies on the RHS of (1.1.84). We then let $t \longrightarrow \infty$ in (1.1.84), by using the fact that $\varphi_{q} (x) \dis \mathop{\sim}_{|x| \to \infty}^{} |x|$. We obtain, since $g^{(t)} \dis \mathop{\longrightarrow}_{t \to \infty}^{} g$ a.s. under $W_{\infty}^{(q)}$ (cf Theorem 1.1.1) : 
\begin{eqnarray*}
\lefteqn{W_{\infty}^{(q)} \big(F (X_{u}, \; u \le g) \cdot G(X_{g+v}, \; v \ge 0)\big) } \nonumber \\
&=& \frac{1}{\varphi_{q} (0)} \; \left( \int_{0}^{\infty} \frac{du}{\sqrt{2 \pi u}} \; \Pi_{0,0}^{(u)} \big(F (X_{v},\; v \le u) e^{- \frac{1}{2} \, A_{u}^{(q)}} \right) \cdot P_{0}^{(3, {\rm sym})} \big( G (X_{l}, \; l \ge 0) \; e^{- \frac{1}{2} \, A_{\infty}^{(q)}} \big) \\
&=& \frac{1}{\varphi_{q} (0)} \; \widetilde{\bf W} \big( e^{- \frac{1}{2} \, A_{\infty}^{(q)}} F(X_{u}, \; u \le g) \, G(X_{g+l}, \; l \ge 0) \big)
\end{eqnarray*}

\noi This is the statement of Theorem 1.1.10.

\smallskip

\noi{\bf 1.1.7} \underline{Relations between {\bf W} and other penalisations (than the Feynman-Kac ones).}

\smallskip

\noi We have shown - this is Theorem 1.1.2 - that for every $q \in \mathcal{I}$ :
\begin{eqnarray*} 
{\bf W}
&=& \varphi_{q} (0) \, {\rm exp} \left(\frac{1}{2} \, A_{\infty}^{(q)}\right) \cdot W_{\infty}^{(q)} \nonumber \\
&=& {\bf W} \left( {\rm exp} \left(- \frac{1}{2} \, A_{\infty}^{(q)}\right) \right) \cdot {\rm exp} \left( \frac{1}{2} \, A_{\infty}^{(q)}\right) \cdot W_{\infty}^{(q)} \hspace*{5cm}(1.1.85)
\end{eqnarray*}

\noi Of course, this formula is very much linked with the penalisation of the Wiener measure by the multiplicative functional $\dis \left(F_{t} = {\rm exp} \left(- \frac{1}{2} \, A_{t}^{(q)}\right), \; t \ge 0\right)$. Here, we shall prove that formulae analogous to (1.1.85) are true for other penalisations than these Feynman-Kac ones. 

\noi We now fix some notations : 
$$
S_{t} := \mathop{\rm sup}_{s \le t}^{} X_{s}, \quad I_{t} := \mathop{\rm inf}_{s \le t}^{} X_{s} \eqno(1.1.86) $$

\vspace*{-0.50cm}

$$
\Gamma_{+} := \{ \omega \in \Omega \;;\; \mathop{\rm lim}_{t \to \infty}^{} X_{t} (\omega) = + \infty\}, \; \Gamma_{-} := \{ \omega \in \Omega, \mathop{\rm lim}_{t \to \infty}^{} X_{t} (\omega)  = -�\infty\} \eqno(1.1.87)$$

\vspace*{-0.50cm}

$$
{\bf W}^{+} := 1_{\Gamma_{+}} \cdot {\bf W} , \quad {\bf W}^{-} := 1_{\Gamma_{-}} \cdot {\bf W} \eqno(1.1.88) $$

\vspace*{-0.50cm}

$$
\theta_{+} := {\rm sup} \{ t \;;\; S_{t} < S_{\infty}\}, \qquad \theta_{-} := {\rm sup} \{t \;;\; I_{t} > I_{\infty}\} \eqno(1.1.89) $$

\smallskip

\noi Let $\psi_{+}$ (resp. $\psi_{-}$) a Borel and integrable function from $\mathbb{R}_{+}$ to $\mathbb{R}_{+}$ (resp. from $\mathbb{R}_{-}$ to $\mathbb{R}_{+}$). We denote by $(M_{s}^{\psi_{+} (S)}, \; s \ge 0)$ \big(resp. $(M_{s}^{\psi_{-}(I)}, \; s \ge 0)\big)$ the 
Az\'ema-Yor martingale defined by :
$$
M_{s}^{\psi_{+} (S)} := \frac{1}{\dis \left(\int_{0}^{\infty} \psi_{+} (y) dy\right)} \; \left(\psi_{+} (S_{s}) \, (S_{s} - X_{s}) + \int_{S_{s}}^{\infty} \psi_{+} (y) dy\right) \eqno(1.1.90)$$

$$
M_{s}^{\psi_{-} (I)} := \frac{1}{\left(\dis \int_{- \infty}^{0} \psi_{-} (y)dy\right)} \left(\psi_{-} (I_{s}) \, (X_{s} - I_{s}) + \int_{-\infty}^{I_{s}} \psi_{-} (y) dy\right) \eqno(1.1.91)$$

\noi Let $W_{\infty}^{\psi_{+} (S)}$ (resp $W_{\infty}^{\psi_{-}(I)})$ denote the probability on $(\Omega, \mathcal{F}_{\infty})$ characterized by : 
$$
{W_{\infty}^{\psi_{+} (S)}}_{|_{\mathcal{F}_{t}}} = {M_{t}^{\psi_{+}(S)} \cdot W}_{|_{\mathcal{F}_{t}}} , \qquad {W_{\infty}^{\psi_{-}(I)} \cdot }{|_{\mathcal{F}_{t}}} = M_{t}^{\psi_{-}(I)} \cdot W_{|_{\mathcal{F}_{t}}}
\eqno(1.1.92) $$

\noi \big(see [RVY, II] for more informations about these probabilities\big).

\bigskip

\noi The analogue of formulae (1.1.85) and (1.1.41) is here :

\smallskip

\noi {\bf Theorem 1.1.11.} {\it Let $\psi_{+}, \psi_{-}$ as above, with $\psi_{+} (\infty) =\psi_{-} (- \infty) =0$.}

\smallskip

\noi $ {\bf 1)} \qquad W_{\infty}^{\psi_{+} (S)} = \frac{1}{{\bf W}
 \big(\psi_{+} (S_{\infty})\big)} \cdot 
\psi_{+} (S_{\infty}) \cdot {\bf W}^-  \hfill(1.1.93) $

\smallskip

\noi $ \qquad \;\;\; \, W_{\infty}^{\psi_{-} (I)} = \frac{1}{{\bf W}
 \big(\psi_{-} (I_{\infty})\big)} \cdot 
\psi_{-} (I_{\infty}) \cdot {\bf W}^+  \hfill (1.1.94)$
 
 \noi {\bf 2)}{\it For every $t \ge 0$ and $\Gamma_{t} \in b_{+} (\mathcal{F}_{t})$ :
 
 \smallskip

\noi $ \qquad \;\;\; \; W \big(\Gamma_{t} (S_{t} - X_{t})\big) = {\bf W}^{-} (\Gamma_{t} \; 1_{\theta_{+} \le t}) \hfill(1.1.95) $ 

\smallskip

\noi $ \qquad \;\;\; \; W\big(\Gamma_{t} (X_{t} - I_{t})\big) = {\bf W}^{+} (\Gamma_{t} \; 1_{\theta_{-} \le t})
\hfill(1.1.96) $ 

\smallskip}

\noi {\bf Proof of Theorem 1.1.11.}

\noi {\it i)} We have, from (1.1.85), for $q \in \mathcal{I}$, and $\Gamma_{t} \in b_{+} (\mathcal{F}_{t})$ :
\begin{eqnarray*} 
\lefteqn{{\bf W} (e^{- \frac{1}{2} \, A_{\infty}^{(q)}} \cdot \Gamma_{t}) } \\
&=& \varphi_{q} (0) \; W_{\infty}^{(q)} (\Gamma_{t}) \\
&=& W \big(\Gamma_{t} \varphi_{q} (X_{t}) e^{- \frac{1}{2} \, A_{t}^{(q)}}\big) \\
&& \big(\text{from (1.1.7) and (1.1.8)\big)} \\
&=& \left( \int_{0}^{\infty} \psi (y) dy \right) \cdot W_{\infty}^{\psi(S)} \left( \Gamma_{t} \; \frac{\varphi_{q} (X_{t}) e^{- \frac{1}{2} \, A_{t}^{(q)}}} {\psi (S_{t}) (S_{t}-X_{t}) + \dis \int_{S_{t}}^{\infty} \psi (y) dy} \right) \hspace*{2cm}(1.1.97) 
\end{eqnarray*}

\noi from (1.1.92) and (1.1.90), and we have written, to simplify, $\psi$ for $\psi_{+}$. Formula (1.1.97) being true for every $\Gamma_{t} \in b_{+} (\mathcal{F}_{t})$, we may take $\Gamma_{t}=\Gamma_{u} \, 1_{\psi(S_{t}) >a} \cdot 1_{S_{t}-X_{t} > b|X_{t}|} \cdot 1_{\int_{S_{t}}^{\infty} \psi (y) dy > c}$ with $0<b<1$,   $a,c > 0$ for any $\Gamma_{u} \in \mathcal{F}_{u}, \; u \le t$. We obtain thus :
\begin{gather*} 
{\bf W} \big[\Gamma_{u} e^{- \frac{1}{2} \, A_{\infty}^{(q)}} \; 1_{\psi (S_{t}) >a} 1_{S_{t}-X_{t} >b|X_{t}|} \; 1_{\int_{S_{t}}^{\infty} \psi (y) dy > c}\big]  \\
= \left(\int_{0}^{\infty} \psi(y) dy\right) \cdot W_{\infty}^{\psi(S)} \left[ \Gamma_{u} \frac{\varphi_{q} (X_{t}) e^{-\frac{1}{2} \, A_{t}^{(q)}}} {\psi (S_{t}) (S_{t}-X_{t}) + \int_{S_{t}}^{\infty} \psi (y) dy} \; 1_{\psi(S_{t}) >a}  \right.\\
\qquad \qquad \qquad \qquad \qquad \qquad  \left. 1_{S_{t}-X_{t} >b |X_{t}|} \; 1_{\int_{S_{t}}^{\infty} \psi (y) dy > c} \right] \hspace*{4,2cm}(1.1.98) \end{gather*}

\noi We shall now let $t \to \infty$ in (1.1.98) with $u$ being fixed. On the LHS, we have : 
\begin{gather*} 
{\bf W}^{+} \; {\rm a.s.} \qquad 1_{\psi(S_{t}) >a} \mathop{\longrightarrow}_{t \to \infty}^{} 0 \qquad \hbox{\big(since} \; S_{t} \longrightarrow + \infty \;\; {\rm and} \;\; \psi (S_{t}) \mathop{\longrightarrow}_{ t \to \infty}^{} 0, \; \textrm{since} \; \psi(+ \infty) = 0 \big)\\
{\bf W}^{-} \; {\rm a.s.} \qquad 1_{\psi(S_{t}) >a} \mathop{\longrightarrow}_{t \to \infty}^{} 1_{\psi(S_{\infty}) > a}  \\
\hspace*{2,2cm} 1_{S_{t}-X_{t} >b |X_{t}]} \mathop{\longrightarrow}_{t \to\infty}^{} 1 \hspace*{8,1cm}{(1.1.99)} \\  
\hspace*{2,2cm} 1_{\int_{S_{t}}^{\infty} \psi (y) dy > c} \longrightarrow 1_{\int_{S_{\infty}}^{\infty} \psi (y) dy > c}
\end{gather*}

\noi Thus, from Lebesgue's dominated convergence Theorem, the LHS of (1.1.98) converges, as $t \to \infty$, towards $L$, with :
$$
L = {\bf W}\big(\Gamma_{u} 1_{\Gamma_{-}} e^{- \frac{1}{2} \, A_{\infty}^{q}} 1_{\psi (S_{\infty})>a} \;
1_{\int_{S_{\infty}}^{\infty} \psi (y) dy > c} \big) \eqno(1.1.100)$$

\noi We now consider the RHS of (1.1.98). On the set :
$$
\big(\psi (S_{t}) > a\big) \cap \big(S_{t} - X_{t} > b |X_{t}|\big) \cap \left(\int_{S_{t}}^{\infty} \psi (y) dy > c�\right), $$

\noi we have :
$$
\frac{\varphi_{q} (X_{t})}{\psi (S_{t}) (S_{t}-X_{t}) + \int_{S_{t}}^{\infty} \psi (y) dy} \le \frac{d + |X_{t}|}{ab |X_{t}| + c} \le k $$

\noi since $|\varphi_{q} (x) | \le d+|x|$ ; thus, we may apply the dominated convergence Theorem to obtain, since under $W_{\infty}^{\psi (S)}$ \big(see 
[RVY, II]\big) : $X_{t} \dis\mathop{\longrightarrow}_{t \to \infty}^{} - \infty$ and $S_{t} \dis\mathop{\longrightarrow}_{t \to \infty}^{} S_{\infty}$ a.s., the convergence of the RHS of (1.1.98) to $R$, with :
$$
R = \left(\int_{0}^{\infty} \psi (y)dy\right) \cdot W_{\infty}^{\psi (S)} \left(\frac{\Gamma_u}{\psi (S_{\infty})} \; e^{- \frac{1}{2} \, A_{\infty}^{(q)}} 1_{\psi(S_{\infty}) >a} \; 1_{\int_{S_{\infty}}^{\infty} \psi (y) dy > c}\right) \eqno(1.1.101) $$

\noi \big(since $\varphi_{q} (x) \dis \mathop{\sim}_{|x| \to \infty}^{}|x|\big)$. Hence, letting $a,c \to 0$ and applying the monotone class Theorem, the equality between (1.1.100) and (1.1.101) implies, for every $\Gamma \in b_{+} (\mathcal{F}_{\infty})$ : 
$$
{\bf W}^{-} (\Gamma e^{-\frac{1}{2} \, A_{\infty}^{(q)}}) = \left(\int_{0}^{\infty} \psi (y) dy \right) W_{\infty}^{\psi (S)} \left(\frac{\Gamma}{\psi (S_{\infty})} \;e^{- \frac{1}{2} \, A_{\infty}^{(q)}} \right) $$

\noi then, replacing $\Gamma e^{- \frac{1}{2} \, A_{\infty}^{(q)}}$ by $\Gamma$ :
\begin{eqnarray*}
{\bf W}^{-} (\Gamma)
&=& \left( \int_{0}^{\infty} \psi (y) dy\right) \cdot W_{\infty}^{\psi (S)} \left(\frac{\Gamma}{\psi (S_{\infty})}\right) \\
&=& {\bf W}^{-} \big(\psi (S_{\infty})\big) \, W_{\infty}^{\psi(S)} \left(\frac{\Gamma}{\psi (S_{\infty}) }\right) \\
&=& {\bf W} \big(\psi (S_{\infty})\big) \, W_{\infty}^{\psi(S)} \left(\frac{\Gamma}{\psi (S_{\infty}) }\right)
\end{eqnarray*}

\noi since $\psi (\infty) = 0$ and $S_{\infty} = + \infty \quad {\bf W}^{+}$ a.s. 

\noi We note that there is no problem to divide by $\psi (S_{\infty})$ since $\psi (S_{\infty}) > 0 \quad W_{\infty}^{\psi (S)}$ a.s. \big(under $W_{\infty}^{\psi (S)}, \; S_{\infty}$ admits $\psi$ as density \big(see [RVY, II]\big)\big).

\noi We have proven (1.1.93), and the proof of (1.1.94) is similar.

\smallskip

\noi {\it ii)} \underline{We now prove (1.1.95)}.

\noi For this purpose, we use the penalisation by $(e^{- \frac{\lambda}{2} \, S_{t}}, \; t \ge 0)$ i.e. (1.1.90) and (1.1.92), with
$\psi_+(x) = e^{- \lambda x}$. We obtain :
$$
M_{t}^{\psi_+ (S)} = \left(1 + \frac{\lambda}{2}\; (S_{t} - X_{t})\right) e^{- \frac{\lambda}{2} \, S_{t}} \eqno(1.1.102) $$

\noi Hence, for every $t \ge 0$ and $\Gamma_{t} \in b_{+} (\mathcal{F}_{t})$ :
\begin{eqnarray*} 
W \left( \Gamma_{t} \left( \frac{2}{\lambda} + (S_{t} -X_{t})\right) \right)
&= & \frac{2}{\lambda} \; W_{\infty}^{\psi_+ (S)} \big[ e^{\frac{\lambda \, S_{t}}{2}} \Gamma_{t}\big] \\
&=&{\bf W}^{-} \big(e^{-\frac{\lambda}{2} \, (S_{\infty}-S_{t})} \Gamma_{t} \big) \qquad {\rm \big(from} \; (1.1.93)\big)  \\
&=& {\bf W}^{-} (\Gamma_{t} \, 1_{\theta_{+} \le t}) + {\bf W}^{-} (\Gamma_{t} \, e^{- \frac{\lambda}{2} \, (S_{\infty} -S_{t})} 1_{\theta_{+} >t}) \hspace*{1cm}(1.1.103)
\end{eqnarray*}

\noi We then let $\lambda \to + \infty$ in (1.1.103) by noting that $S_{\infty} - S_{t} > 0$ on $(\theta_{+} > t)$. We obtain : 
$$
W \big(\Gamma_{t} (S_{t} - X_{t})\big) = {\bf W}^{-} (\Gamma_{t} \, 1_{\theta_{+} \le t}) $$

\noi This is (1.1.95). By symmetry, (1.1.96) now follows.

\smallskip

\noi {\bf Remark 1.1.12} We deduce from(1.1.103) and (1.1.95) that :
$$
W ( \Gamma_{t}) \; \left(\int_{0}^{\infty} e^{- \frac{\lambda}{2} \, y} dy \right) = {\bf W}^{-} \big(\Gamma_{t} \, e^{-\frac{\lambda}{2} (S_{\infty}-S_{t})} 1_{\theta_{+} > t}\big) $$

\noi and operating as in the proof of point 3), {\it i)} of Theorem 1.1.6, we obtain : 
\begin{eqnarray*} 
{\bf W}^{-} (S_{\infty}-S_{t} \in dl) 
&=& 1_{[0, \infty[} (l) dl + W (S_{t} - X_{t}) \, \delta_{0} (dl)  \\
&=& 1_{[0,\infty[} (l)dl + \sqrt{\frac{2t}{\pi}} \, \delta_{0} (dl) \hspace*{6cm}(1.1.104)
\end{eqnarray*}

\noi and, conditionally on $S_{\infty} - S_{t}=l, \; l>0, \; (X_{u}, \; u \le t)$ is, under ${\bf W}^-$, a Brownian motion indexed by $[0,t]$. Theorem 1.1.11 is the prototype of similar results which we may obtain for other penalisations. Here are, without proof, some examples.

\smallskip

\noi {\bf Theorem 1.1.11'.}

\noi {\bf 1)} {\it Let $h^{+},  h^{-} \;:\; \mathbb{R}_{+} \longrightarrow \mathbb{R}_{+}$ such that $\dis \int_{0}^{\infty} (h^{+} +  h^{-}) (y) dy < \infty$. Let $W_{\infty}^{h^{+},  h^{-}}$ denote the probability defined by \big(see [RVY, II]\big) :
$$
{ W_{\infty}^{h^{+},  h^{-} } }_{| \mathcal{F}_{t}} = M_{t}^{h^{+},  h^{-}} \cdot W_{|\mathcal{F}_{t}} \eqno(1.1.105) $$

\noi with
\begin{gather*}
M_{t}^{h^{+},  h^{-}} = \dis \frac{1} { \dis \frac{1}{2} \, \int_{0}^{\infty} (h^{+}+ h^{-}) (y) dy} \; \Bigl\{ \Bigl( X_{t}^{+} \, h^{+} (L_{t}) + X_{t}^{-} \, h^{-} (L_{t}) \\
\hspace*{5cm} \left. \left. + \int_{L_{t}}^{\infty} \frac{1}{2} ( h^{+} +  h^{-}) (y) dy \right) \right \} \hspace*{3,4cm}{(1.1.106)}
\end{gather*}

\noi Then :
$$
{\bf W} = \Big\{{\bf W}^{+} \big(\big(h^{+} (L_{\infty})\big) + {\bf W}^{-} (h^{-}(L_{\infty})\big)\Big\} \left(1_{\Gamma_{+}} \; \frac{1}{h^{+} (L_{\infty})} + 1_{\Gamma_{-}} \; \frac{1}{h^{-} (L_{\infty})}\right) W_{\infty}^{h^{+}, h^{-}} \eqno(1.1.107) $$

\noi In other words :
$$
1_{\Gamma_{+}} . W_{\infty}^{h^{+}, h^{-}}  = \frac{1} {{\bf W}^{+} \big((h^{+} 
(L_{\infty}) \big)+ {\bf W}^{-} (h^{-} (L_{\infty})\big) } 
 \; h^{+} (L_{\infty}) \; . {\bf W}^+ \eqno(1.1.108) $$

$$
1_{\Gamma_{-}} . W_{\infty}^{h^{+}, h^{-}}  = \frac{1} {{\bf W}^{+} \big((h^{+} 
(L_{\infty}) \big)+ {\bf W}^{-} (h^{-} (L_{\infty})\big) } 
 \; h^{-} (L_{\infty}) \; . {\bf W}^- \eqno(1.1.109) $$
 
 \noi {\bf 2)} Let $\psi : \mathbb{R}_{+} \longrightarrow \mathbb{R}_{+}$ be Borel and integrable, and let us define :
$$
M_{t}^{\psi (S_{g})} := \left(\frac{1}{2} \; \psi (S_{g^{(t)}}) |X_{t}| + \psi (S_{t}) (S_{t}-X_{t}^{+}) + \int_{S_{t}}^{\infty} \psi (y) dy\right) \cdot \frac{1}{\int_{0}^{\infty} \psi (y) dy} \eqno(1.1.110) $$

\noi with $g^{(t)} := \sup \{ s \le t, \; X_{s} =0\}$. If $W_{\infty}^{\psi (S_{g})}$ is given by :
$$
{W_{\infty}^{\psi (S_{g})}}_{|�\mathcal{F}_{t}} = M_{t}^{\psi (S_{g})} \cdot W_{|\mathcal{F}_{t}} \eqno(1.1.111) $$   

\noi \big(see [RY, VIII]\big), then : 

 \noi {\it i)} $\qquad  W_{\infty}^{\psi (S_{g})}  = \frac{ \psi (S_{g}) }{ {\bf W} \big(\psi (S_{g})\big)} 
\; \cdot {\bf W} \hfill(1.1.112) $
 
 \noi {\it ii)} If $\rho := \sup \{u \le g, \; S_{g^{(u)}} < S_{g}\}$, then, for all $t$ and for all $\Gamma_{t} \in b_{+} (\mathcal{F}_{t})$ :
$$
 {\bf W} (\Gamma_{t} \, 1_{\rho \le t}) = W \left(\Gamma_{t} \left(\frac{1}{2} \; |X_{t}| + (S_{t} - X_{t}^{+}) 1_{S_{t} = S_{g}^{(t)}}\right) \right) \eqno(1.1.113) $$   }

\noi We could also present analogous results for penalisations associated to the numbers of downcrossings \big(see [RVY, II]\big) or the length of the longest excursion before $g^{(t)}$ \big(see [RVY, VII]\big), etc...

\noi We use, in Section 2, Theorem 1.1.11 and 1.1.11' to give explicit examples of martingales $\big(M_{t} (F), \; t \ge 0\big), \; F \in L_{+}^{1} ({\bf W})$. These martingales are defined in Theorem 1.2.1

\bigskip

{\bf \large 1.2 $W$-Brownian martingales associated to $\bf W$.}

\smallskip

\noi The notation in this Section 1.2 is the same as in Section 1.1. Our aim here is to associate to every r.v. in $L_{+}^{1} (\Omega, \mathcal{F}_{\infty}, {\bf W})$ a $W$-martingale and to study a few of its properties. Thus, {\bf W} appears as "a machine to construct $W$-martingales". We shall also prove (see Theorem 1.2.5) a decomposition Theorem which is valid for every positive Brownian supermartingale.

\smallskip

\noi{\bf 1.2.1}  \underline{Definition of the martingales $\big(M_{t} (F), \; t \ge 0\big)$.}
\\ \\
\noi {\bf Theorem 1.2.1.} {\it Let $F \in L_{+}^{1} ( \Omega, \mathcal{F}_{\infty}, {\bf W})$. There exists a  positive (necessarily continuous) $\big( (\mathcal{F}_{t}, \; t \ge 0), W\big)$ martingale $\big(M_{t} (F), \; t \ge 0\big)$ such that :

\smallskip

\noi {\bf 1)} For every $t \ge 0$ and $\Gamma_{t} \in b(\mathcal{F}_t)$ : 
$$
{\bf W} (F \cdot \Gamma_{t}) = W \big(M_{t} (F) \cdot \Gamma_{t}\big) \eqno(1.2.1) $$

\noi In particular, for every $t \ge 0$ :
$$
 {\bf W} (F) = W \big(M_{t} (F)\big) \eqno(1.2.2)$$

\noi {\bf 2)} $(M_t(F), t \geq 0)$ may be computed via the "characteristic formula" :
 $$ \qquad M_{t} (F) =  \widehat{\bf W}_{X_{t}} \big(F (\omega_{t}, \widehat{\omega}^{t}) \big) \eqno(1.2.3)$$

\noi (cf Point 1 of Remark 1.2.2 for this notation)

\smallskip

\noi {\bf 3)} $\qquad \dis M_{t} (F) \mathop{\longrightarrow}_{t \to \infty}^{} 0 \quad W \; {\rm a.s.} \hfill(1.2.4) $

\noi In particular, the martingale $\big(M_{t}(F), \; t \ge 0\big)$ is not uniformly integrable.

\smallskip

\noi {\bf 4)} For every $q \in \mathcal{I}$ :
$$
M_{t} (F) = \varphi_{q} (0) \; M_{t}^{(q)} \, W_{\infty}^{(q)} (F \, e^{\frac{1}{2} \, A_{\infty}^{(q)}} |\mathcal{F}_{t}) \eqno(1.2.5)$$

\noi where $M_{t}^{(q)}, \; \varphi_{q}$ and $W_{\infty}^{(q)}$ are defined in Theorem 1.1.1. }

\smallskip

\noi {\bf Remark 1.2.2.} 

\noi 1. We now give some explanation about the notation in (1.2.3). If $\omega \in \mathcal{C} (\mathbb{R}_{+} \to \mathbb{R})$, then $\omega_{t}$ (resp. $\omega^{t}$) denotes the part of $\omega$ before $t$ (resp. after $t$) : 
\begin{equation*}
\omega =(\omega_{t}, \omega^{t})
\end{equation*}

\noi that is, precisely : 
\begin{gather*}  
 X_{u} (\omega) =
\left\{ \begin{array}{rl}
X_{u}(\omega_{t}) \quad\quad {\rm if} \quad u \le t \hfill \nonumber \\
 X_{u-t} (\omega^{t}) \quad {\rm if} \quad u \ge t \hfill  \\ 
  \end{array}\right.  
\end{gather*}

\noi and our notation $\widehat{\bf W}_{X_{t} } \big(F (\omega_{t}, \widehat{\omega}^{t})\big)$ stands for the expectation of $F(\omega_{t}, \bullet)$ with respect to ${\bf W}_{X_{t} (\omega)}$.

\smallskip

\noi 2. To every r.v. $G$ in $L_{+}^{1}(\Omega, \mathcal{F}_{\infty}, W)$ we may of course associate the positive martingale $\big(\widetilde{M}_{t} (G) := W(G|\mathcal{F}_{t}), \; t \ge 0\big)$. But, contrarily to the description for $M_{t}(F)$ given in Theorem 1.2.1, this is a uniformly integrable martingale.

\smallskip

\noi 3. Formula (1.2.5) may seem ambiguous, since the r.v. $W_{\infty}^{(q)} \big(F \, e^{\frac{1}{2} \, A_{\infty}^{(q)}} |�\mathcal{F}_{t}\big)$ is only defined $W_{\infty}^{(q)}$ a.s. But since from (1.1.7), the probability $W_{\infty}^{(q)}$ is absolutely continuous on $\mathcal{F}_{t}$ with respect to $W$, there is in fact no ambiguity. On the other hand, from (1.1.16) : 
$$
{\bf W} (F) = \varphi_{q} (0) \, W_{\infty}^{(q)} \left(F \exp \left(\frac{1}{2} \, A_{\infty}^{(q)}\right)\right) < \infty
\eqno(1.2.6) $$

\noi as soon as $F \in L^{1} ({\bf W})$. Thus, the $\big((\mathcal{F}_{t}, \; t \ge 0), W_{\infty}^{(q)}\big)$ martingale $\left(W_{\infty}^{(q)} \left(F \exp (\frac{1}{2} \, A_{\infty}^{(q)}) | \mathcal{F}_{t}\right), t \ge 0 \right)$ is $W_{\infty}^{(q)}$-uniformly integrable. 

\smallskip

\noi 4. Of course, $\big(M_{t} (F), \; t \ge 0\big)$ is continuous, as it is a $\big((\mathcal{F}_{t}, \; t \ge 0), W \big)$ martingale.

\smallskip

\noi 5. On the injectivity of $F \longrightarrow \big(M_{t} (F), \; t \ge 0\big)$ : assume that, for $F_{1}$ and $F_{2}$ belonging to $L^{1}(\Omega, \mathcal{F}_{\infty}, {\bf W})$ we have : $M_{t} (F_{1})=M_{t}(F_{2})$ a.s., for every $t \ge 0$. Then $F_{1}=F_{2} \; \; {\bf W}$ a.s. Indeed, from (1.2.1) :
$$
W\big(\Gamma_{t}(M_{t}(F_{1})-M_{t}(F_{2})\big)=0 
= {\bf W} \big(\Gamma_{t} (F_{1}-F_{2})\big)  $$

\noi As this relation is true for every $t \ge 0$ and $\Gamma_{t} \in b(\mathcal{F}_{t})$, the monotone class Theorem implies that, for every $\Gamma \in b (\mathcal{F}_{\infty})$ :
\begin{equation*}
{\bf W} \big(\Gamma (F_{1}-F_{2})\big) = 0, \;\; {\rm i.e.} \;\; F_{1}=F_{2} \;\; {\bf W} \; a.s.
\end{equation*}

\noi Later in this Section (see Lemma 1.2.8), we shall obtain a more direct "construction" of $F$ from $\big(M_{t} (F), \; t \ge 0\big)$.

\smallskip

\noi {\bf Proof of Theorem 1.2.1.}

\noi {\it i)} \underline{We show point 1}.

\noi We denote by $W^{F}$ the finite positive measure on $(\Omega, \mathcal{F}_{\infty})$ defined by :
$$
W^{F} (G) := {\bf W} (F \cdot G) \eqno(1.2.7) $$

\noi Let $\Gamma_{t} \in b_{+}(\mathcal{F}_{t})$ such that $W (\Gamma_{t})=0$. From (1.1.7), for every $q \in \mathcal{I}, \; W_{\infty}^{(q)} (\Gamma_{t})=0$ hence, from (1.1.16) :
$$
W^{F} (\Gamma_{t}) = {\bf W} (F \cdot \Gamma_{t})= \varphi_{q} (0) \; W_{\infty}^{(q)} (F\, e^{\frac{1}{2} \, A_{\infty}^{(q)}} \Gamma_{t}) = 0 $$

\noi from (1.2.6). Thus :
$$
W^{F}_{| \mathcal{F}_{t}} \ll W_{|\mathcal{F}_{t}} $$

\noi Consequently, from the Radon-Nikodym Theorem, there exists a $W$ integrable, positive r.v. $M_{t} (F)$, such that 
$$
W^{F}_{| \mathcal{F}_{t}} = M_{t} (F) \cdot W_{| \mathcal{F}_{t}} \eqno(1.2.8) $$

\noi This is a rewriting of formula (1.2.1). Formula (1.2.2) is obtained from (1.2.1) by taking $\Gamma_{t} \equiv 1$. The fact that $\big(M_{t} (F), \; t \ge 0\big)$ is a $\big((\mathcal{F}_{t}, \; t \ge 0), \; W\big)$ martingale follows from (1.2.8). We also note that, as every Brownian martingale, the process $\big(M_{t} (F), \; t \ge 0\big)$ admits a continuous version (which we shall always consider).

\bigskip

\noi {\it ii)} \underline{We show point 4}.

\noi From (1.2.1), (1.1.16) and (1.1.7), we have for every $\Gamma_{t} \in b_{+} (\mathcal{F}_{t})$
\begin{eqnarray*}
{\bf W}(\Gamma_{t} F)
&=& W \big(\Gamma_{t} M_{t} (F)\big) \\
&=& \varphi_{q} (0) W_{\infty}^{(q)} (\Gamma_{t} F \, e^{\frac{1}{2} \, A_{\infty}^{(q)}} ) \qquad {\rm \big(from } \;(1.1.16)\big) \\
&=& \varphi_{q} (0) W_{\infty}^{(q)} \big(\Gamma_{t} W_{\infty}^{(q)} (F \, e^{\frac{1}{2} \, A_{\infty}^{(q)}} | \mathcal{F}_{t})\big) \\
&=& \varphi_{q} (0) W \big(\Gamma_{t} M_{t}^{(q)} W_{\infty}^{(q)} \big(F \, e^{\frac{1}{2} \, A_{\infty}^{(q)}} | \mathcal{F}_{t} \big) \big) \qquad  \qquad {\rm \big(from } \;(1.1.7)\big)
\end{eqnarray*}

\noi (1.2.5) follows.

\smallskip

\noi {\it iii)} \underline{We show point 3}.

\noi $\bullet$ For every $s \ge 0$ and $\Gamma_{s} \in b (\mathcal{F}_{s})$, we have for $s \le t$ from (1.2.1) :
$$
{\bf W} (\Gamma_{s} \cdot F) = W \big(\Gamma_{s} \cdot M_{t} (F)\big) \eqno(1.2.9) $$

\noi Since the $\big((\mathcal{F}_{t}, \; t > 0), W\big)$ martingale $\big(M_{t} (F), \; t \ge 0\big)$ is positive, it converges $W$ a.s. towards $M_{\infty} (F)$. Letting $t \to \infty$ in (1.2.9) and using Fatou's Lemma, we have :
\begin{equation*}
W \big(\Gamma_{s} M_{\infty} (F)\big) \le {\bf W} (\Gamma_{s} \cdot F)
\end{equation*}

\noi Choosing $\Gamma_{s} = 1_{g^{(s)} \ge a} $, with $g^{(s)} := \sup \{ u \le s, X_{u}=0\}$ we obtain :
$$
W \big(1_{g^{(s)} \ge a} \cdot M_{\infty} (F)\big) \le {\bf W} (1_{g^{(s)} \ge a} \cdot F) \eqno(1.2.10) $$

\noi Letting $s \to \infty$ in (1.2.10) and noting that : 
\begin{gather*}
1_{g^{(s)} \ge a} \longrightarrow 1 \qquad W \;\; {\rm a.s.} \\
1_{g^{(s)} \ge a} \longrightarrow 1_{g \ge a} \;\; {\bf W} \;\; {\rm a.s.}.
\end{gather*}

\noi we obtain :
$$
W\big(M_{\infty} (F)\big) \le {\bf W} (1_{g \ge a} \cdot F) $$

\noi Now, from Theorem 1.1.6 we know that $g < \infty$ {\bf W} a.s., hence we get : ${\bf W} (1_{g \ge a} \cdot F) \dis \mathop{\longrightarrow}_{a \to \infty}^{} 0$. Thus :
$$
W \big(M_{\infty} (F)\big) =0 \qquad {\rm and} \qquad M_{\infty} (F) = 0 \qquad W \; {\rm a.s.} $$

\noi $\bullet$ Another way to prove point 3. consists in writing, for $s \le t$ : 
\begin{eqnarray*}
\hspace*{1cm}W \big(\Gamma_{s} M_{t} (F)\big) 
&=&  \varphi_{q} (0) \; W \big(\Gamma_{s} M_{t}^{(q)} W_{\infty}^{(q)} \big( F \, e^{\frac{1}{2} \, A_{\infty}^{(q)}} | \mathcal{F}_{t})\big) \quad \text{\big(from (1.2.5)\big)}\\
&= &\varphi_{q} (0) W_{\infty}^{(q)} \big(\Gamma_{s} W_{\infty}^{(q)} \big( F \, e^{\frac{1}{2} \, A_{\infty}^{(q)}} | \mathcal{F}_{t})\big) \quad \text{\big(from (1.1.7)\big)}\hspace*{1cm} (1.2.11) 
\end{eqnarray*}

\noi But, since the $W_{\infty}^{(q)}$ martingale $\big(W_{\infty}^{(q)} \big( F \, e^{\frac{1}{2} \, A_{\infty}^{(q)}} | \mathcal{F}_{t}), \; t \ge 0\big)$ is uniformly integrable it converges a.s. and in $L^{1} (W_{\infty}^{(q)})$ towards $F \, e^{\frac{1}{2} \, A_{\infty}^{(q)}}$ as $t \to�\infty$. Thus, letting $t \to \infty$ in (1.2.11) and using again Fatou's Lemma, we obtain :
$$
W\big(\Gamma_{s} M_{\infty} (F)\big) \le \varphi_{q} (0) \; W_{\infty}^{(q)} (\Gamma_{s} F \, e^{\frac{1}{2} \, A_{\infty}^{(q)}}) \eqno(1.2.12) $$

\noi We then choose $\Gamma_{s} = {\bf{1}}_{\{A_{s}^{(q)} \ge a\}}$ and obtain 
$$
W \big(1_{(A_{s}^{(q)} \ge a)} \, M_{\infty} (F)\big) \le \varphi_{q} (0) \, W_{\infty}^{(q)} \big(1_{A_{s}^{(q)} \ge a} \, F \, e^{- \frac{1}{2} \, A_{\infty}^{(q)}}\big) \eqno(1.2.13) $$

\noi We then let $s \to \infty$ and note that :
\begin{gather*}
1_{A_{s}^{(q)} \ge a} \longrightarrow 1 \qquad \;\;\; W \; \hbox{\rm a.s. (since Brownian motion is recurrent)} \\
1_{A_{s}^{(q)} \ge a} \longrightarrow 1_{A_{\infty}^{(q)} \ge a} \; W_{\infty}^{(q)} \; {\rm a.s.}
\end{gather*}

\noi Hence :
$$
W \big(M_{\infty} (F)\big) \le \varphi_{q} (0) \, W_{\infty}^{(q)}\big(1_{A_{\infty}^{(q)} \ge a} \, F \, e^{- \frac{1}{2} \, A_{\infty}^{(q)}}\big) \eqno(1.2.14) $$

\noi It now suffices to let $a \to \infty$, using the fact that $A_{\infty}^{(q)} < \infty \; W_{\infty}^{(q)}$ a.s., and that $F \, e^{- \frac{1}{2} \, A_{\infty}^{(q)}} \in L^{1} (W_{\infty}^{(q)})$ \big(from (1.2.6)\big) to obtain :
$$
W\big(M_{\infty} (F)\big) = 0 \quad \hbox{and hence :} \quad M_{\infty} (F) =0 \quad W \;�{\rm a.s.} $$

\noi {\it iv)} \underline{We prove point 2, i.e. the "characteristic formula" for $M_t(F)$}.

\noi We have, from (1.2.5) : 
\begin{eqnarray*}
M_{t} (F)
&=& \varphi_{q} (0) \, M_{t}^{(q)} \, W_{\infty}^{(q)} \big(F \, e^{- \frac{1}{2} \, A_{\infty}^{(q)}} | \mathcal{F}_{t}\big) \\
&=& \varphi_{q} (X_{t}) e^{- \frac{1}{2} \, A_{t}^{(q)}} W_{\infty}^{(q)} \big(F \, e^{- \frac{1}{2} \, A_{\infty}^{(q)}} | \mathcal{F}_{t}\big) \\
&& \text{\big(from the definition (1.1.8) of} \; M_{t}^{(q)}\big) \\
&=& \varphi_{q} (X_{t}) \widehat{W}_{ X_{t}, \infty}^{(q)}  \big( e^{\frac{1}{2} (A_{\infty}^{q}-A_{t}^{q}) } F( \omega_{t}, \widehat{\omega}^{t})\big) \\
&& \text{(from the Markov property)} \\
&=& \widehat{\bf W}_{X_{t}} \big(F ( \omega_{t}, \widehat{\omega}^{t})\big), \qquad {\rm from} \; (1.1.16)
\end{eqnarray*}

\noi{\bf 1.2.2} \underline{ Examples of martingales $\big(M_{t} (F), \; t \ge 0\big)$.}

\smallskip

\noi Formula (1.2.3) which provides an "explicit" expression for $M_{t} (F)$ is not always, practically, easy to compute.

\noi {\bf 1.2.2.1}  \underline{A first method to obtain examples of $\big(M_{t} (F), \; t \ge 0\big)$.}

\noi To begin with, we present a "computation principle" to obtain $M_{t} (F)$.

\noi "\underline{\bf Computation principle}"

\noi Let $(N_{t}, \; t \ge 0)$ denote a $\big((\mathcal{F}_{t}, \; t \ge 0\big), \; W\big)$ positive martingale such that $N_{0}=1$. Let $W_{\infty}^{N}$ be the probability on $(\Omega, \mathcal{F}_{\infty})$ which is characterized by : 
$$
{W_{\infty}^{N}}_{ | \mathcal{F}_{t}} = N_{t} \cdot W_{| \mathcal{F}_{t}} \eqno(1.2.15) $$

\noi Let us assume that there exists a r.v. $F \in L_{+}^{1} (\Omega, \mathcal{F}_{\infty}, {\bf W})$ such that :
$$
F \cdot {\bf W} = {\bf W} (F) \cdot W_{\infty}^{N} \eqno(1.2.16) $$

\noi Then
$$
M_{t} (F) = {\bf W} (F) \cdot N_{t} \eqno(1.2.17) $$

\noi {\bf Proof of the "Computation principle".}

\noi We have, for every $t \ge 0$ and $\Gamma_{t} \in b (\mathcal{F}_{t})$, from (1.2.1) :
$$
{\bf W} (F \cdot \Gamma_{t}) = W \big(M_{t} (F) \cdot \Gamma_{t} \big) $$

\noi On the other hand, from the hypothesis (1.2.16) : 
$$
{\bf W} (F \cdot \Gamma_{t}) = {\bf W}(F) \, W_{\infty}^{N} (\Gamma_{t}) $$

\noi Hence, this quantity also equals :
$$
 {\bf W} (F) \, W (\Gamma_{t} \cdot N_{t}) \eqno(1.2.18) $$

\noi from (1.2.15). Since $\Gamma_{t}$ denotes any $\mathcal{F}_{t}$ measurable set in (1.2.18), one obtains :
$$
M_{t} (F) = {\bf W} (F) \cdot N_{t} \quad W \;\; {\rm a.s.} $$

\noi {\bf Example 1.} Let $q \in \mathcal{I}$ and $N_{t} := \dis \frac{\varphi_{q} (X_{t})}{\varphi_{q} (0)}$ exp$\dis \left(-\frac{1}{2} \, A_{\infty}^{(q)}\right)$.

\noi From (1.1.16) and (1.1.7), the hypotheses of the "Computation principle" are satisfied with $\dis F = \exp \left(-\frac{1}{2} \, A_{\infty}^{(q)}\right)$. Thus :
\begin{eqnarray*}
M_{t} (e^{- \frac{1}{2} \, A_{\infty}^{(q)}})
& = &{\bf W} (e^{- \frac{1}{2} \, A_{\infty}^{(q)}}) \cdot \frac{\varphi_{q} (X_{t})}{\varphi_{q} (0)} \; \exp \left(-\frac{1}{2} \, A_{t}^{(q)}\right) \nonumber \\
&=& \varphi_{q} (X_{t}) \, \exp \left(-\frac{1}{2} \, A_{t}^{(q)}\right) \ \hspace*{6,5cm}(1.2.19)
\end{eqnarray*}

\noi since, from (1.1.17), ${\bf W} \left(\exp \left( - \frac{1}{2} \, A_{\infty}^{(q)}\right) \right) = \varphi_{q} (0)$.

\smallskip

\noi {\bf Example 2.} Let $h : \mathbb{R}_{+} \to \mathbb{R}_{+}$ Borel and integrable and :
$$
N_{t} := \frac{1}{\int_{0}^{\infty} h(y)dy} \cdot \left(h (L_{t}) |X_{t}| + \int_{L_{t}}^{\infty} h(y)dy\right) \eqno(1.2.20) $$

\noi From Theorem 1.1.11', the hypotheses of the "Computation principle" are satisfied with $F=h(L_{\infty}) $ \big(we note from point 3, {\it i)} of Theorem 1.1.6 : ${\bf W} \big(h(L_{\infty})\big)= \dis \int_{0}^{\infty} h(l)dl < \infty\big)$. Thus :
$$
M_{t} \big(h (L_{\infty})\big) = h(L_{t}) |X_{t}| + \int_{L_{t}}^{\infty} h(y)dy \eqno(1.2.21) $$

\noi \big(cf [RVY, II] for the use of this martingale\big).

\smallskip

\noi {\bf Example 3.} Let $\dis S_{t} := \mathop{\sup}_{s \le t}^{} X_{s}$ and $\psi : \mathbb{R}_{+} \to \mathbb{R}_{+}$ Borel and integrable, such that $\psi (+ \infty)=0$. Due to Theorem 1.1.11, the "Computation principle" applies with $F= \psi (S_{\infty})$ and 
$$
N_{t} := \frac{1}{\int_{0}^{\infty} \psi (y)dy} \; \left(\psi (S_{t}) (S_{t} - X_{t}) + \int_{S_{t}}^{\infty} \psi (y) dy \right) $$

\noi We note that, from (1.1.104) (taken with $t=0$) :
$$
{\bf W} \big(\psi (S_{\infty})\big) = \int_{0}^{\infty} \psi (l) dl < \infty. \eqno(1.2.22) $$

\noi Thus : 
$$
M_{t} \big(\psi (S_{\infty})\big) = \psi (S_{t}) (S_{t}-X_{t}) + \int_{S_{t}}^{\infty} \psi (y) dy \eqno(1.2.23) $$

\noi Another manner to obtain (1.2.23) may be to invoke L\'evy's Theorem :
\begin{equation*}
\big((S_{t}, S_{t}-X_{t}), \; t \ge 0\big) \mathop{=}_{}^{(\rm law)} \big((L_{t}, |X_{t}|), \; t \ge 0\big)
\end{equation*}

\noi then to use (1.2.21).

\noi The reader may refer to [RVY, II] for links between the Az\'ema-Yor martingale $\dis \biggl(\psi (S_{t}) (S_{t} - X_{t})$

\noi $\dis \left. + \int_{S_{t}}^{\infty} \psi (y) dy, \; t \ge 0 \right)$ and the penalisation problem with the process $\big(\psi (S_{t}), \; t \ge 0\big)$.

\smallskip

\noi {\bf Example 4.} Let $\psi : \mathbb{R}_{+} \to \mathbb{R}_{+}$ a Borel, integrable function with $\psi (\infty) = 0$. The "Computation principle", with the help of Theorem 1.1.11', yields to, with $F=\psi (S_{g})$ :
\begin{eqnarray*}
\hspace*{1cm}M_{t} \big(\psi (S_{g})\big)
& =& \dis \frac{1}{2} \; \psi (S_{g^{(t)}}) |X_{t}| + \psi (S_{t}) (S_{t}-X_{t}^{+}) + \int_{S_{t}}^{\infty} \psi (y) dy \hspace*{2cm}(1.2.24) \\
&=& \dis \frac{1}{2} \; \psi (S_{g^{(t)}}) \cdot X_{t} + M_{t} \big(\psi (S_{\infty})\big)
\end{eqnarray*}

\noi where $M_{t} \big(\psi (S_{\infty}\big) \big)$ is defined by (1.2.23). We note that, from (1.1.104), since $\psi (+\infty)=0$ :
$$
{\bf W}^{-} \big(\psi (S_{g})\big) = {\bf W} \big(\psi (S_{\infty})\big) = \int_{0}^{\infty} \psi (l) dl < \infty \eqno(1.2.25) $$

\noi {\bf Example 5.} Let $a<b$ and :
\begin{gather*}
\begin{array}{cccc}
T^{(1)} \hfill & :=  \inf \{ t \ge 0 \;;\; X_{t} > b\},  \hfill & T^{(2)} \hfill &:= \inf \{ t \ge T^{(1)} \;;\; X_{t} < a \} \hfill \\
T ^{(2n+1)}& :=  \inf \{ t \ge T^{(2n)} \;;\; X_{t} > b \}, \hfill & T^{(2n+2)} &:= \inf \{ t \ge T^{(2n+1)} \;;\; X_{t} < a\}
\end{array}
\end{gather*}

\noi Define :
\begin{equation*}
D_{t}^{[a,b]} := \sum_{n \ge 1} \; 1_{(T^{(2n)} \le t)}
\end{equation*}

\noi $D_{t}^{[a,b]}$ is the number of down-crossings on the interval $[a,b]$ before time $t$. Let $h: \mathbb{N} \to \mathbb{R}_{+}$ such that $h$ is decreasing, $h(0)=1, h(+\infty)=0$ and denote $\Delta h(n) := h(n)-h(n+1)$. The "Computation principle" and an extension to this situation of Theorem 1.1.11' lead to :
\begin{gather*} 
\hspace*{-1cm}M_{t} \big(\Delta h (D_{\infty}^{[a,b]})\big) = \sum_{n \ge 0}  \left\{ 1_{[T^{(2n)}, T^{(2n+1)}[} (t) \left[ \frac{h(n)}{2} \left( 1+ \frac{b-X_{t}}{b-a}\right) +  \frac{h(n+1)}{2} \left(\frac{X_{t}-a}{b-a}\right) \right]  \right.\\
\hspace*{2cm} \left.+ 1_{[T^{(2n+1)}, T^{(2n+2)}[} (t) \left[ \frac{h(n+1)}{2} \left( 1+ \frac{b-X_{t}}{b-a}\right) + \frac{h(n)}{2} \left(\frac{X_{t}-a}{b-a}\right) \right] \right\} \;(1.2.26)
\end{gather*}

\noi The reader may refer to [RVY, II] for more information relative to this martingale.

\smallskip

\noi {\bf Example 6.} Let $\Sigma_{g^{(t)}}$ denote the length of the longest excursion of Brownian motion $(X_{u}, \; u \ge~0)$ before $g^{(t)} := \sup \{ s \le t \;;\; X_{s}=0\}$. Let $h: \mathbb{R}_{+} \to \mathbb{R}_{+}$ such that $\dis \int_{0}^{\infty}  z|h'(z)|dz < \infty$. Then, the "Computation principle" and an extension of Theorem 1.1.11', lead to :
\begin{gather*}
M_{t}\big(\sqrt{h} (\Sigma_{g})\big) = \sqrt{h} (\Sigma_{g^{(t)}}) \cdot |X_{t}| + h_{1} (A_{t}) \Phi \left( \frac{|X_{t}|} {\sqrt{(\Sigma_{g^{(t)}}-A_{t})_{+}}}\right) \nonumber \\
\hspace*{3cm} + \int_{0}^{ \frac{|X_{t}|} {\sqrt{(\Sigma_{g^{(t)}}-A_{t})_{+}}}} h_{1} \left(A_{t} + \frac{X_{t}^{2}}{v^{2}} \right) \left(\exp \left(-\frac{v^{2}}{2}\right)\right) dv \hspace*{2cm}(1.2.27)
\end{gather*}

\noi with
\begin{gather*}
A_{t} := t-g^{(t)}, \; \Phi (x) := \int_{x}^{\infty} \exp \left(-\frac{v^{2}}{2}\right) dv \\
h_{1} (x) := - \int_{\sqrt{x}}^{\infty} z h'(z) dz
\end{gather*}

\noi \big(see [RY, VIII] or [RY, M], Chap. 3\big).

\smallskip

\noi{\bf 1.2.2.2} \underline{A $\;$second $\;$manner $\;$to compute explicitly martingales of the form $\big(M_{t} (F), \; t \ge 0\big).$}

\noi This method hinges upon the following Theorem 1.2.3. $(F_{u}, \; u \ge 0)$ denotes a positive predictable process such that : 
$$
{\bf W} (F_{g}) < \infty \eqno(1.2.28) $$

\noi We note that from Theorem 1.1.6, this condition is equivalent to :
$$
W \left(\int_{0}^{\infty} F_{s} \, dL_{s}�\right) < \infty \eqno(1.2.29) $$

\noi or equivalently after the change of variable $l=L_{s}$, to :
$$
\int_{0}^{\infty} W (F_{\tau_{l}}) dl = W \left(\int_{0}^{\infty} F_{s} \, dL_{s}�\right) < \infty \eqno(1.2.30) $$

\noi with :
$$
\tau_{l} := \inf \{ t > 0 \;;\; L_{t} > l\} \eqno(1.2.31) $$

\noi {\bf Theorem 1.2.3.} {\it Let $(F_{u}, \; u \ge 0)$ denote a positive predictable process such that : 
$$
{\bf W} (F_{g}) = \int_{0}^{\infty} W (F_{\tau_{l}}) dl < \infty \eqno(1.2.32) $$

\noi Then, the martingale $\big(M_{t} (F_{g}), \; t \ge 0\big)$ may be expressed as : 
\begin{eqnarray*}
M_{t} (F_{g})
&=& F_{g^{(t)}} \cdot |X_{t}| + \int_{t}^{\infty} p_{u-t} (X_{t}) \; \Pi_{0,0}^{(u)} (F_{u}|\mathcal{F}_{t}) \,du \hspace*{5,2cm}(1.2.33) \\
&=& F_{g^{(t)}} \cdot |X_{t}| + \int_{L_{t}}^{\infty} W (F_{\tau_{l}} | \mathcal{F}_{t}) \,dl \hspace*{7,1cm}(1.2.34) \\
&=& \int_{0}^{t} F_{g^{(s)}} \, {\rm sgn}(X_{s}) \, dX_{s} + W \left(\int_{0}^{\infty} F_{\tau_{l}} \, dl | \mathcal{F}_{t} \right) \hspace*{5cm}(1.2.35)
\end{eqnarray*}

\noi In (1.2.33), $\Pi_{0,0}^{(u)}$ denotes the law of Brownian bridge of length $u$ and :
$$
p_{s} (x) := \frac{1}{\sqrt{2 \pi s}} \; e^{- \frac{x^{2}}{2s}} \eqno(1.2.36) $$  }

\noi {\bf Proof of Theorem 1.2.3.}

\noi {\it i)} \underline{We first prove (1.2.33)}.

\noi For every $t \ge 0$ and $\Gamma_{t} \in b(\mathcal{F}_{t})$ we have by (1.2.1) :
\begin{eqnarray*}
{\bf W} (\Gamma_{t} F_{g}) 
&=& W \big(\Gamma_{t} \, M_{t}(F_{g})\big) \\
&=& {\bf W} (\Gamma_{t} \, F_{g} \, 1_{g \le t}) + {\bf W} (\Gamma_{t} \, F_{g} \, 1_{g > t}) \\
&=& {\bf W} (\Gamma_{t} \, F_{g^{(t)}} \, 1_{g \le t} ) + {\bf W} (\Gamma_{t} \, F_{g} \, 1_{g > t}) \\
&& \text{\big(since} \; g=g^{(t)} \; \text{on the set} \; (g \le t)\big)
\end{eqnarray*} 
\begin{equation*}
\hspace*{1,8cm} := (1_{t}) + (2_{t}) \hspace*{9,4cm}(1.2.37) \end{equation*}

\noi We study successively $(1_{t})$ and $(2_{t})$ :
\begin{eqnarray*}
(1_{t})
&=& {\bf W} (\Gamma_{t} \, F_{g^{(t)}} \, 1_{g \le t}) \\
&=& W \big(\Gamma_{t} \, F_{g^{(t)}} \, |X_{t}|\big)\hspace*{9,1cm}(1.2.38) \\
&&\text{\big( from point 2, {\it i)} of Theorem 1.1.6.\big)} \\
(2_{t})
&=& {\bf W} (\Gamma_{t} \, F_{g} \, 1_{g > t}) \\
&=& \int_{t}^{\infty} \frac{du}{\sqrt{2 \pi u}} \; {\bf W} \big(\Gamma_{t} \, F_{g} | g=u\big) \qquad \hbox{(from (1.1.42))} \\
&=& \int_{t}^{\infty} \frac{du}{\sqrt{2 \pi u}} \;  \Pi_{0,0}^{(u)} (\Gamma_{t} \, F_{u}) \qquad \hbox{(from point 2)iii) of Theorem 1.1.6)} \\
&=& \int_{t}^{\infty} \frac{du}{\sqrt{2 \pi u}} \;  \Pi_{0,0}^{(u)}  \big(\Gamma_{t} \, \Pi_{0,0}^{(u)} (F_{u} | \mathcal{F}_{t})\big)
\end{eqnarray*}

\noi We now use the (partial) absolute continuity formula for the law of the Brownian bridge with respect to that of Brownian motion :
$$
{\Pi_{0,0}^{(u)}}_{|\mathcal{F}_{t}} = \frac{p_{u-t}(X_{t})}{p_{u} (0)} \cdot W_{| \mathcal{F}_{t}} \qquad (u > t)  \eqno(1.2.39) $$

\noi to obtain : 
\begin{eqnarray*}
(2_{t})
&=& \int_{t}^{\infty} \frac{du}{\sqrt{2 \pi u}} \; W \left(\frac{\Gamma_{t} \, p_{u-t}(X_{t})}{p_{u}(0)} \; \Pi_{0,0}^{(u)} \big(F_{u} | \mathcal{F}_{t}\big) \right)\\
 &=& \int_{t}^{\infty} du \,  W  \big(\Gamma_{t} \, p_{u-t} (X_{t}) \, \Pi_{0,0}^{(u)} \big(F_{u} | \mathcal{F}_{t}\big) \big) \hspace*{6cm}(1.2.40)
\end{eqnarray*}

\noi since $p_{u} (0) = \dis \frac{1}{\sqrt{2 \pi u}}\cdot$ Gathering (1.2.40), (1.2.37) and (1.2.38), we obtain (1.2.33).

\bigskip

\noi {\it ii)} \underline{We now prove (1.2.34)}.

\noi Of course, (1.2.34) is equivalent to : 
$$
\int_{t}^{\infty} p_{u-t} (X_{t}) \, \Pi_{0,0}^{(u)} \big(F_{u} | \mathcal{F}_{t}\big) du = \int_{L_{t}}^{\infty} W \big(F_{\tau_{l}} | \mathcal{F}_{t}\big) dl \eqno(1.2.41) $$

\noi or to :
$$
W \left(\Gamma_{t} \int_{t}^{\infty} p_{u-t} (X_{t}) \; \Pi_{0,0}^{(u)} \big(F_{u} | \mathcal{F}_{t}\big) du\right) = W \left(\Gamma_{t} \cdot \int_{L_{t}}^{\infty} F_{\tau_{l}} \, dl \right) \eqno(1.2.42) $$
\noi for any $\Gamma_{t} \in b(\mathcal{F}_{t})$. But we have : 
\begin{eqnarray*}
W \left(\Gamma_{t} \cdot \int_{L_{t}}^{\infty} F_{\tau_{l}} dl \right) 
&= & W \left(\Gamma_{t} \int_{t}^{\infty} F_{u} \; dL_{u}\right) \\
&& \text{(after the change of variable} \;  l=L_{u}) \\
&=& \int_{t}^{\infty} \frac{du}{\sqrt{2 \pi u}} \;  \Pi_{0,0}^{(u)} \, (F_{u} \, \Gamma_{t}) \\
&=& \int_{t}^{\infty} \frac{du}{\sqrt{2 \pi u}} \;  \Pi_{0,0}^{(u)}  \, \left(\Gamma_{t} \, \Pi_{0,0}^{(u)} \big(F_{u} | \mathcal{F}_{t} \big)\right) \\
&= &\int_{t}^{\infty} \frac{du}{\sqrt{2 \pi u}} \; W \left(\Gamma_{t} \; \frac{p_{u-t} (X_{t})}{p_{u} (0)} \; \Pi_{0,0}^{(u)} \big(F_{u} | \mathcal{F}_{t} \big)\right) \\
&& \text{\big(by the absolute continuity formula (1.2.39)\big)} \\
&= &W \left(\Gamma_{t} \int_{t}^{\infty} p_{u-t} (X_{t}) \, \Pi_{0,0}^{(u)} \big(F_{u} | \mathcal{F}_{t} \big) du \right)
\end{eqnarray*}

\noi {\it ii')} \underline{We give now a direct proof} - i.e. without using (1.2.33) - of (1.2.34). We have, for every $t \ge 0$ and $\Gamma_{t} \in b_{+} (\mathcal{F}_{t})$ :
\begin{eqnarray*}
{\bf W} (F_{g} \, \Gamma_{t})
&=& {\bf W} (F_{g} \, \Gamma_{t} \, 1_{g \le t}) + {\bf W} (F_{g} \, \Gamma_{t} \, 1_{g > t} ) \\
&=& {\bf W} (F_{g^{(t)}} \, \Gamma_{t} \, 1_{g \le t}) + {\bf W} (\widetilde{\Gamma}_{g} \, F_{g})
\end{eqnarray*}

\noi (since $g=g^{(t)}$ on the set $(g \le t)$), and we have used the notation :
\begin{gather*}
(\widetilde{\Gamma}_{u}, \; u \ge 0) := \big(\Gamma_{t} \, 1_{]t, \infty[} (u), \; u \ge 0\big) \\
\quad = W (\Gamma_{t} \, F_{g^{(t)}} |X_{t}| + W \left( \int_{0}^{\infty} \widetilde{\Gamma}_{\tau_{l}} \, F_{\tau_{l}} \, dl \right)
\end{gather*}

\noi \big(from point 2 {\it i)} of Theorem 1.1.6 and from formula (1.1.44)\big). 

\smallskip

\noi  Hence :
\begin{eqnarray*}
{\bf W} (F_{g} \, \Gamma_{t})
&=& W \big(\Gamma_{t} \, F_{g^{(t)}} |X_{t}|\big) + W \left(\Gamma_{t} \int_{0}^{\infty} 1_{t < \tau_{l}} \, F_{\tau_{l}} \, dl \right) \\
&=& W \big(\Gamma_{t} \, F_{g^{(t)}} |X_{t}|\big) + W \left(\Gamma_{t} \int_{L_{t}}^{\infty} F_{\tau_{l}} \, dl \right) \\
&=& W \big(\Gamma_{t} \, F_{g^{(t)}} |X_{t}|\big) + W \left(\Gamma_{t} \int_{L_{t}}^{\infty} W \big( F_{\tau_{l}} | \mathcal{F}_{t} \big) dl \right)
\end{eqnarray*}

\noi which implies (1.2.34).

\smallskip

\noi {\it iii)} \underline{We now prove (1.2.35)}.

\noi To go from (1.2.34) to (1.2.35), we use the balayage formulae, which yields :
\begin{equation*}
F_{g^{(t)}} \cdot |X_{t}| = \int_{0}^{t} F_{g^{(s)}} \, {\rm sgn} (X_{s}) \, dX_{s} + \int_{0}^{t} F_{u} \, dL_{u}
\end{equation*}

\noi and we add this expression to $\dis \int_{L_{t}}^{\infty} W \big(F_{\tau_{l}} | \mathcal{F}_{t}\big) \,dl =W \left(\int_{t}^{\infty} F_{u} \, dL_{u} | \mathcal{F}_{t}\right)$ on the RHS. It is now clear that (1.2.34) implies (1.2.35).

\smallskip

\noi {\bf Corollary 1.2.4.} 

\noi {\bf 1)} {\it Formula (1.2.34) expresses the martingale $\big(M_{t} (F_{g}), \; t \ge 0\big)$ as the sum of a submartingale $(F_{g^{(t)}} \cdot |X_{t}|, \; t \ge 0)$ and a supermartingale $\dis \left( W \left(\int_{0}^{\infty} F_{\tau_{l}} \, 1_{\tau_{l} > t} \, dl |�\mathcal{F}_{t}\right), \; t \ge 0 \right)$ both of which converge to 0 a.s., as $t \to \infty$. 

\smallskip

\noi {\bf 2)} The variable $\dis \int_{0}^{\infty} F_{g^{(u)}}^{2} du$ is finite a.s. but it satisfies :
$$
W\left( \left( \int_{0}^{\infty} F_{g^{(u)}}^{2} \, du \right)^{\frac{1}{2}}\right) = + \infty \eqno(1.2.43) $$

\noi unless $F_{g}=0, \;\; {\bf W}$ a.s.  }

\smallskip

\noi {\bf Proof of Corollary 1.2.4.}

\noi The first statement is obvious since $F_{g^{(t)}} |X_{t}|$ is the absolute value of the martingale $F_{g^{(t)}} \cdot X_{t}$. Moreover, $|F_{g^{(t)}} \cdot X_{t}| \le M_{t} (F_{g})$, hence since $M_{t} (F_{g}) \dis \mathop{\longrightarrow}_{t \to \infty}^{} 0$ a.s. (see Theorem 1.2.1) the same is true for $F_{g^{(t)}} \cdot X_{t}$. To prove the second item, assume that :
$$
W \left(\left( \int_{0}^{\infty} F_{g^{(u)}}^{2} \, du\right)^{\frac{1}{2}}\right) < \infty $$

\noi Then, the martingale $\dis \left(\int_{0}^{t} F_{g^{(s)}} \, {\rm sgn}(X_{s}) \, dX_{s}, \; t \ge 0\right)$ would be in $H^{1}$ ; a fortiori it would be uniformly integrable. From (1.2.35), since $W \dis \left( \int_{0}^{\infty} F_{\tau_{l}} \, dl\right) < \infty, \; \big(M_{t} (F_{g}), \; t \ge 0\big)$ would also be uniformly integrable ; but this is only possible, since this martingale converges a.s. to 0 (see Theorem 1.2.1) if it is identically equal to 0, that is $F_{g} =0 \;\; {\bf W}$ a.s. (see point 5 of Remark 1.2.2).

\bigskip

\noi Of course, if we want to compute $\big(M_{t} (F), \; t \ge 0\big)$ in a completely explicit manner, we need to compute $\Pi_{0,0}^{(u)} \big(F_{u}�|�\mathcal{F}_{t}\big)$, for $t \le u$ $\biggl($or $W \dis \left(\int_{0}^{\infty} F_{\tau_{l}}\, dl | \mathcal{F}_{t}\right)\biggl)$. This is what has been done in the Examples 4 and 6 above. Here is an example where this computation is immediate.

\smallskip

\noi {\bf Example 7.} Let $\psi : \mathbb{R}_{+} \to \mathbb{R}_{+}$ Borel such that : 
$$
\int_{0}^{\infty} \psi (t) \; \frac{dt}{\sqrt{2 \pi t}} < \infty \eqno(1.2.44) $$

\noi Then : 
$$
M_{t} \big(\psi (g)\big) = \psi (g^{(t)}) |X_{t}| + \int_{0}^{\infty} \; \frac{du}{\sqrt{2 \pi u}} \; e^{- \frac{X_{t}^{2}} {2u}} \psi (t+u) \eqno(1.2.45) $$

\noi To obtain (1.2.45), we apply Theorem 1.2.3 with the (deterministic) process $(F_{u}, \; u \ge 0) := \big(\psi (u), \; u \ge 0\big)$ and we use : 
$$
\Pi_{0,0}^{(u)} \big(F_{u} | \mathcal{F}_{t}\big) = \Pi_{0,0}^{(u)} \big( \psi (u) | \mathcal{F}_{t} \big) = \psi (u)
$$

\noi We then make the change of variable $u-t=v$ in (1.2.33).

\noi More generally (see Theorem 1.1.8), with $g_{a} := \sup \{t \;;\; X_{t}=a\}$, we have :
$$
M_{t} \big[\psi (g_{a})\big] = \psi (g_{a}^{(t)}) |X_{t}-a| + \int_{0}^{\infty} \frac{du}{\sqrt{2 \pi u}} \; e^{- \frac{( X_{t}-a)^{2}}{2u}} \psi (t+u) \eqno(1.2.46) $$

\noi with :
$$
g_{a}^{(t)} := \sup \{ s \le t \;;\; X_{s}=a\} \eqno(1.2.47) $$

\noi {\bf Back to Example 2.} Formula (1.2.21) is a particular case of (1.2.34). Indeed, if we apply (1.2.34) with $(F_{u}, \; u \ge 0) := \big(h (L_{u}), \; u \ge 0\big)$, we obtain :
\begin{eqnarray*}
M_{t} \big(h (L_{\infty})\big) 
&=& M_{t} \big(h (L_{g})\big) \\
&=& h (L_{g^{(t)}}) |X_{t}| + W \left( \int_{L_{t}}^{\infty} h (L_{\tau_{l}})dl | \mathcal{F}_{t}\right) \\
&=& h(L_{t}) |X_{t}| + \int_{L_{t}}^{\infty} h(l) \, dl
\end{eqnarray*}

\noi since $L_{g^{(t)}}  = L_{t}$ and $L_{\tau_{l}} = l$.

\noi In the same spirit, for $h : \mathbb{R}_{+} \times \mathbb{R}_{+} \to \mathbb{R}_{+}$ Borel such that :
$$
\int_{0}^{\infty} \int_{0}^{\infty} h (l,u) \; \frac{ l\, e^{  - \frac{ l^{2} }{2u}}} { \sqrt{ 2 \pi u^{3} } }  \; dl \, du < \infty
\eqno(1.2.48) $$

\noi then \big(see (1.1.47)\big)  ${\bf W} \big(h(L_{\infty}, g)\big) < \infty$  and
\begin{gather*}
M_{t} \big(h (L_{\infty},g)\big) = h (L_{g^{(t)}}, g^{(t)}) \cdot |X_{t}| + W \left(\int_{L_{t}}^{\infty} h (L_{\tau_{l}} , \tau_{l}) \, dl | \mathcal{F}_{t} \right)\nonumber \\
\quad = h (L_{t}, g^{(t)}) \cdot |X_{t}| + \widehat{W}_{X_{t}} \left(\int_{0}^{\infty} h (L_{t} + \widehat{L}_{v}, t+v) d \widehat{L}_{v}\right) \hspace*{3,8cm}(1.2.49)
\end{gather*}

\noi{\bf 1.2.2.3}  \underline{A third manner to obtain explicit examples of martingales $\big(M_{t} (F), \; t \ge 0\big)$.}

\smallskip

\noi $\bullet$ We begin with a definition. We shall say that a family of r.v.'s $(F_{t}, \; t \ge 0)$ converges, as $t \to \infty$, towards $F \;\; {\bf W}$ a.s. if for some $G > 0$, $G \in L_{+}^{1} (\mathcal{F}_{\infty}, {\bf W})$ $F_{t} \dis \mathop{\longrightarrow}_{t \to \infty}^{} F\;\; W^{G}$ a.s. We recall : $W^{G} (\Gamma):= {\bf W} (G \; \Gamma), \; \Gamma \in b (\mathcal{F}_{\infty})$. Clearly, this definition does not depend on the r.v. $G$ chosen in the above class. In particular, it may be convenient to take for $G$ the r.v. exp$\dis \left(-\frac{1}{2} \, A_{\infty}^{(q)}\right)$ for some $q \in \mathcal{I}$ ; hence, the a.s. ${\bf W}$-convergence is precisely the $W_{\infty}^{(q)}$ a.s. convergence.

\noi This definition may seem complicated. However, its aim is to take care of the difficulty arising from the fact that for every $\Gamma_{t} \in b_{+} (\mathcal{F}_{t}), \; {\bf W} (\Gamma_{t})$ equals either 0 or $+ \infty$ \big(see point {\it v)} of the proof of Theorem 1.1.2\big). 

\noi Equivalently, $\dis F_{t} \mathop{\longrightarrow}_{t \to \infty}^{} F$ $\quad {\bf W}$ a.s. if and only if ${\bf W} (\Delta) =0$ with $\Delta = \big\{ \omega \;;\; F_{t} (\omega) \dis \mathop{\not\longrightarrow}_{t \to \infty}^{} F(\omega)\big\}$

\smallskip

\noi $\bullet$ In Section 1.2.3 below we shall obtain the following result : (it is a Corollary of Theorem 1.2.5, in the same Section 1.2.3)

\noi {\bf Corollary 1.2.6.} {\it A positive $\big((\mathcal{F}_{t}, \; t \ge 0), W\big)$ martingale $M_{t}, \; t \ge 0$ is of the form $\big(M_{t} (F),$ 

\noi $t \ge 0\big)$ for some $F \in L_{+}^{1}(\mathcal{F}_{\infty}, {\bf W})$ if and only if :
\begin{equation*}
\mathop{\rm lim}_{t \to \infty}^{} \; \frac{M_{t}}{1+|X_{t}|} \; \textrm{exists} \; {\bf W}\textrm{-a.s.} 
\end{equation*}
and 
\begin{equation*}
M_{0} = {\bf W} \left(\mathop{\rm lim}_{t \to \infty}^{} \; \frac{M_{t}}{1+|X_{t}|} \right)
\end{equation*}
\noi and, in this case : }
\begin{equation*}
F = \mathop{\rm lim}_{t \to \infty}^{} \; \frac{M_{t}}{1+|X_{t}|} \quad {\bf W} \; {\rm a.s.}
\end{equation*}

\noi $\bullet$ We now illustrate with 3 examples how due to this Corollary, we may compute explicitly $\big(M_{t} (F), \; t \ge 0\big)$ for some $F \in L_{+}^{1}(\mathcal{F}_{\infty}, {\bf W})$.   

\smallskip

\noi {\bf Back to Example 1.} Let $q \in \mathcal{I}$ and $M_{t} := \varphi_{q} (X_{t}) \; \exp \dis \left(-\frac{1}{2} \, A_{t}^{(q)}\right)$. Since
\begin{equation*}
\varphi_{q} (x) \mathop{\rm \sim}_{|x| \to \infty}^{} |x| \;\; {\rm and} \;\; |X_{t}| \mathop{\longrightarrow}_{t \to \infty}^{} \infty \quad {\bf W} \; {\rm a.s.}
\end{equation*}

\noi we have :
\begin{equation*}
\frac{M_{t}}{1+|X_{t}|} \; \mathop{\longrightarrow}_{t \to \infty}^{} \exp \left(-\frac{1}{2} \, A_{\infty}^{(q)}\right)
 := F \quad {\bf W} \; \textrm{a.s.}
\end{equation*}

\noi On the other hand, 
\begin{equation*}
M_{0} = \varphi_{q} (0) = {\bf W} \left(\exp \left(-\frac{1}{2} \, A_{\infty}^{(q)}\right)\right) \qquad \hbox{\big(from (1.1.17)\big)}
\end{equation*}

\noi Thus, from Corollary 1.2.6. : 
\begin{equation*}
M_{t}  \left(\exp \left(-\frac{1}{2} \, A_{\infty}^{(q)}\right)\right) = \varphi_{q} (X_{t}) \, \exp \left( - \frac{1}{2} \, A_{t}^{(q)}\right)
\end{equation*}

\noi {\bf Back to Example 2.} Let $h:\mathbb{R}_{+} \to \mathbb{R}_{+}$ Borel and integrable and :
\begin{equation*}
M_{t} := h(L_{t})|X_{t}| + \int_{L_{t}}^{\infty} h(y)dy
\end{equation*}

\noi It is clear that :
\begin{equation*}
\frac{M_{t}}{1+|X_{t}|} \; \mathop{\longrightarrow}_{t \to \infty}^{} h(L_{\infty}) \quad {\bf W} \; {\rm a.s.}
\end{equation*}

\noi and that from point 3)i) of Theorem 1.1.6. : 
\begin{equation*}
M_{0} = \int_{0}^{\infty} h(y)dy = {\bf W} \big(h (L_{\infty})\big)
\end{equation*}

\noi Thus, from Corollary 1.2.6 :
\begin{equation*}
M_{t} \big(h (L_{\infty})\big) = h (L_{t}) |X_{t}| + \int_{L_{t}}^{\infty} h(y)dy
\end{equation*}

\noi {\bf Back to Example 3.} Let $\psi : \mathbb{R}_{+} \to \mathbb{R}_{+}$ Borel and integrable, with $\psi (\infty) =0$. Let
\begin{equation*}
M_{t} := \psi (S_{t}) (S_{t}-X_{t}) + \int_{S_{t}}^{\infty} \psi (y)dy
\end{equation*}

\noi Then :
\begin{equation*}
\frac{M_{t}}{1 + |X_{t}|} \; \mathop{\longrightarrow}_{t \to \infty}^{} \psi (S_{\infty}) \qquad {\bf W} \; {\rm a.s.} \;\; \big({\rm see} \; (1.1.99)\big)
\end{equation*}

\noi From (1.2.22) :
\begin{equation*}
{\bf W} \big(\psi (S_{\infty})\big) = \int_{0}^{\infty} \psi (l)dl = M_{0}
\end{equation*}

\noi Hence :
\begin{equation*}
M_{t} \big(\psi (S_{\infty})\big) = \psi (S_{t}) (S_{t}-X_{t}) + \int_{S_{t}}^{\infty} \psi (y) dy
\end{equation*}

\noi{\bf 1.2.3} \underline{A decomposition Theorem for positive Brownian supermartingales.}

\noi Here is the most inportant result of this Section 1.2.

\noi {\bf Theorem 1.2.5.} {\it Let $(Z_{t}, \; t \ge 0)$ denote a positive $\big((\mathcal{F}_{t}, \; t \ge 0), W\big)$ supermartingale. We denote $Z_{\infty} := \dis \mathop{\rm lim}_{t \to \infty}^{} Z_{t} \quad (W$ a.s.). Then : 

\noi {\bf 1)}  $$z_{\infty} := \mathop{\rm lim}_{t \to \infty}^{} \; \frac{Z_{t}}{1+|X_{t}|}
\; \operatorname{\it exists} \;  {\bf W} \; {\rm a.s.} \eqno(1.2.50) $$

$$
 and \qquad {\bf W} (z_{\infty}) < \infty \eqno(1.2.51) $$

\noi {\bf 2)} $(Z_{t}, \; t \ge 0)$ decomposes in a unique manner in the form :
$$
Z_{t} = M_{t}(z_{\infty}) + W\big(Z_{\infty}|\mathcal{F}_{t}\big) + \xi_{t}, \qquad t \ge 0 \eqno(1.2.52) $$

\noi where $\big(M_{t} (z_{\infty}), \; t \ge 0\big)$ and $\big(W \big(Z_{\infty}|\mathcal{F}_{t}\big), \; t \ge 0\big)$ denote two $\big((\mathcal{F}_{t}, \; t \ge 0), W\big)$ martingales and :
\begin{equation*}
(\xi_{t}, \; t \ge 0) \;\; \hbox{is a} \;\; \big((\mathcal{F}_{t}, \; t \ge 0), W\big) \;\; \hbox{positive supermartingale}
\end{equation*}

\noi such that :

\noi {\it i)} $Z_{\infty} \in L_{+}^{1} (\mathcal{F}_{\infty}, W)$, hence $W \big(Z_{\infty}|\mathcal{F}_{t}\big)$ converges $W$ a.s. and in $L^{1} (\mathcal{F}_{\infty},W)$ towards $Z_{\infty}$.

\smallskip

\noi  {\it ii)} $ \qquad \frac{W \big(Z_{\infty}|\mathcal{F}_{t}\big) + \xi_{t}}{1+|X_{t}|} \; \underset{t \rightarrow \infty}{\longrightarrow} 0 \qquad {\bf W} \; {\rm a.s.} \hfill(1.2.53) $

\smallskip

\noi ${\it iii)} \qquad M_{t} (z_{\infty}) + \xi_{t}  \; \underset{t \rightarrow \infty}{\longrightarrow} 0 \qquad  W \; {\rm a.s.} \hfill(1.2.54) $

\bigskip   }

\noi After proving Theorem 1.2.5, we shall give a number of examples of $\big((\mathcal{F}_{t}, \; t \ge 0), W\big)$ supermartingales for which we can compute explicitly the decomposition (1.2.52).

\noi We refer the reader to subsection 1.2.2.3 for the definition of the a.s. ${\bf W}$ convergence.

\smallskip

\noi {\bf Corollary 1.2.6.} {\it \big(Characterisation of martingales of the form $\big(M_{t} (F), \;\; t \ge 0)\big)$. 

\noi A $\big((\mathcal{F}_{t}, \; t�\ge 0), \; W\big)$ positive martingale $(Z_{t}, \; t \ge 0)$ is equal to $\big(M_{t} (F), \; t \ge 0\big)$ for an $F \in L_{+}^{1} (\mathcal{F}_{\infty}, {\bf W})$ if and only if :
$$
Z_{0} = {\bf W} \left(\mathop{\rm lim}_{t \to \infty}^{} \; \frac{Z_{t}}{1+|X_{t}|}\right) \eqno(1.2.55) $$

\noi Note that $\dis \mathop{\rm lim}_{t \to \infty}^{} \; \frac{Z_{t}}{1+|X_{t}|}$ exists $\bf W$ a.s. from (1.2.50). }

\bigskip

\noi {\bf Proof of Corollary 1.2.6.} 

\noi We write, from (1.2.52) :
\begin{equation*}
Z_{t} = M_{t} (z_{\infty}) + W(Z_{\infty} | \mathcal{F}_{t}) + \xi_{t}
\end{equation*}

\noi \big(where, in this situation, $(\xi_{t}, \; t \ge 0)$ is a positive martingale\big). Hence :
\begin{equation*}
Z_{0} = W \big(M_{0} (z_{\infty})\big) + W \big(W (Z_{\infty}| \mathcal{F}_{0})\big) + W(\xi_{0})
\end{equation*}

\noi i.e., from (1.2.55) and (1.2.2) :
\begin{equation*}
Z_{0} = {\bf W}(z_{\infty}) = {\bf W}(z_{\infty}) + W(Z_{\infty}) + W (\xi_{0})
\end{equation*}

\noi hence :
\begin{equation*}
W(Z_{\infty}) = W(\xi_{0}) = 0 \quad {\rm and} \quad W(Z_{\infty} | \mathcal{F}_{t}) = \xi_{t} = 0, \quad {\rm i.e.} \;\; Z_{t} = M_{t}(z_{\infty})
\end{equation*}

\noi {\bf Proof of Theorem 1.2.5.} 

\noi This proof hinges on the three following Lemmas. 

\smallskip

\noi {\bf Lemma 1.2.7.} {\it Let $F, G \in L_{+}^{1} (\mathcal{F}_{\infty}, {\bf W})$ and $G > 0 \quad {\bf W}$ a.s. Then :
$$
\frac{ M_{t}(F)}{M_{t} (G)} = W^{G} \left( \frac{F}{G} \Big| \mathcal{F}_{t}\right) \quad W^{G} \; {\rm a.s.}
\eqno(1.2.56) $$

\noi Consequently : }
$$
\frac{ M_{t}(F)}{M_{t} (G)} \; \mathop{\longrightarrow}_{t \to \infty}^{} \frac{F}{G} \quad W^{G} \; {\rm a.s.}
\quad ({\rm hence} \; {\bf W} \; {\rm a.s.}) \eqno(1.2.57) $$

\noi {\bf Lemma 1.2.8.} {\it Let $F \in L_{+}^{1} (\mathcal{F}_{\infty}, {\bf W})$. Then :
$$  
\frac{M_{t} (F)}{1 + |X_{t}|} \; \mathop{\longrightarrow}_{t \to \infty}^{} F \qquad {\bf W} \; {\rm a.s.} \eqno(1.2.58) $$  }

\noi {\bf Lemma 1.2.9.} {\it Let $(Z_{t}, \; t \ge 0)$ denote a positive $\big((\mathcal{F}_{t}, \; t \ge 0), W\big)$ supermartingale. Then :

\smallskip

\noi ${\bf 1)} \qquad z_{\infty} := \mathop{\rm lim}_{t \to \infty}^{} \; \frac{Z_{t}}{1+|X_{t}|} \qquad
 {\rm exists} \;\; {\bf W} \; {\rm a.s.} \hfill(1.2.59) $

\smallskip 

\noi Furthermore :
$$
{\bf W} (z_{\infty}) < \infty \eqno(1.2.60) $$

\noi {\bf 2)} For every  $ t \ge 0 : \quad M_{t}(z_{\infty}) \le Z_{t} \qquad W $  a.s. \hfill{\rm (1.2.61)}

\smallskip }

\noi {\bf Proof of Lemma 1.2.7.}

\noi We have, for every $t \ge 0$ and $\Gamma_{t} \in b (\mathcal{F}_{t})$ :
\begin{eqnarray*}
W^{G} \left( \Gamma_{t} \; \frac{M_{t} (F)}{M_{t} (G)} \right)
&=& {\bf W} \left( \Gamma_{t} \, G \; \frac{M_{t} (F)}{M_{t} (G)} \right)  \hspace*{1,2cm}\big(\hbox{by definition of} \; W^{G}\big) \\
&=& W \left(\Gamma_{t} \, M_{t} (G) \;  \frac{M_{t} (F)}{M_{t} (G)} \right) \ \hspace*{0,3cm} \big(\hbox{by definition of} \;  M_{t}(G)\big) \\
&=& W \big(\Gamma_{t} \, M_{t} (F)\big) \\
&=& {\bf W} (\Gamma_{t} \, F) \ \hspace*{2,8cm} \hbox{\big(by definition of} \; M_{t} (F)\big) \\
&=& W^{G} \left(\Gamma_{t} \; \frac{F}{G}\right) \ \hspace*{2,1cm} \hbox{(by definition of} \; W^{G}) \\
&=& W^{G} \left( \Gamma_{t} \, W^{G} \left(\frac{F}{G} \Big| \mathcal{F}_{t}\right) \right)
\end{eqnarray*}

\noi This is (1.2.56). Now, (1.2.57) is an immediate consequence of (1.2.56) since $\dis \frac{F}{G} \in L^{1} (W^{G})$. Indeed :  $\dis W^{G} \left(\frac{F}{G}\right) = {\bf W} \left(G \cdot \frac{F}{G} \right) = {\bf W} (F) < \infty$.

\smallskip

\noi {\bf Proof of Lemma 1.2.8.}

\noi {\it i)} We first apply Lemma 1.2.7 with $\dis G := {\rm exp} \left(- \frac{1}{2} \, A_{\infty}^{(q)}\right)$, for any $q \in \mathcal{I}$. Then, recall that (Example 1) $M_{t} (G) = \varphi_{q} (X_{t}) \, {\rm exp} \dis \left( - \frac{1}{2}\, A_{t}^{(q)}\right)$ and, since $\varphi_{q} (x) \sim |x|$ as $|x| \to \infty$, we get :
\begin{equation*}
\frac{M_{t} (G)}{1+ |X_{t}|} \mathop{\longrightarrow}_{t \to \infty}^{} {\rm exp} \left( - \frac{1}{2}\, A_{\infty}^{(q)}\right) = G \qquad {\bf W} \; {\rm a.s.}
\end{equation*}

\noi which is the statement of Lemma 1.2.8 with $F = {\rm exp} \dis \left( - \frac{1}{2}\, A_{\infty}^{(q)}\right)$.

\smallskip

\noi {\it ii)} For a general $F \in L_{+}^{1} (\mathcal{F}_{\infty}, {\bf W})$, we write : 
\begin{equation*}
\frac{M_{t} (F)}{1+ |X_{t}|} = \frac{M_{t} (F)}{M_{t} (G)} \cdot \frac{M_{t} (G)}{1+|X_{t}|} \mathop{\longrightarrow}_{t \to \infty}^{} \frac{F}{G} \cdot G \qquad {\bf W} \; {\rm a.s.}
\end{equation*}

\noi by applying Lemma 1.2.7, and the result of point {\it i)} above.

\smallskip

\noi {\bf Proof of Lemma 1.2.9.}

\noi {\it i)} We begin with an argument similar to the one we  used to prove Lemma 1.2.8, that is : 

\noi we write : 
\begin{equation*}
\frac{Z_{t}}{1+|X_{t}|} = \frac{Z_{t}}{M_{t}(G)} \; \frac{M_{t} (G)}{1+|X_{t}|}
\end{equation*}

\noi We now use the fact that $\dis \left(\frac{Z_{t}}{M_{t}(G)}, \; t \ge 0\right)$ is a $\big((\mathcal{F}_{t}, \; t \ge 0), \; W^{G}\big)$ positive supermartingale; hence it converges $W^{G}$ a.s. to a r.v. $\zeta$ ; consequently : 
\begin{equation*}
z_{\infty} := \mathop{\rm lim}_{t \to \infty}^{} \; \frac{Z_{t}}{1+|X_{t}|} \qquad \hbox{exists} \; W^{G} \; {\rm a.s.}
\end{equation*}

\noi and we have : 
\begin{equation*}
z_{\infty} = \zeta \cdot G
\end{equation*}

\noi {\it ii)} Since $\dis \zeta := \mathop{\rm lim}_{t \to \infty}^{} \frac{Z_{t}}{M_{t}(G)}, \; W^{G}$ a.s., is the limit as $t \to \infty$ of a $W^{G}$ supermartingale, we have~: 
\begin{gather*}
W^{G} (\zeta) \le \frac{Z_{0}}{M_{0} (G)} \qquad {\rm hence :} \\
{\bf W}(z_{\infty}) = W^{G} (\zeta) \le \frac{Z_{0}}{M_{0}(G)} < \infty
\end{gather*}

\noi {\it iii)} For any $t \ge 0$ and $\Gamma_{t} \in b_{+} (\mathcal{F}_{t})$, we have :
\begin{eqnarray*}
{\bf W} (\Gamma_{t}\, z_{\infty})
&=& {\bf W} \left( \Gamma_{t} \mathop{\rm lim}_{u \to \infty}^{} \frac{Z_{u}}{1 + |X_{u}|}\right) \\
&=& {\bf W} \left( \Gamma_{t} \mathop{\rm lim}_{u \to \infty}^{} \frac{Z_{u}}{1 + |X_{u}|} \cdot 1_{g \le u} \right) \\
&\le& \mathop{\rm \underline{lim}}_{u \to \infty}^{} {\bf W} \left(\Gamma_{t} \, \frac{Z_{u}}{1+|X_{u}|} \; 1_{g \le u} \right) \quad \text{(from Fatou's Lemma)} \\
&=& \mathop{\rm \underline{lim}}_{u \to \infty}^{} \quad W \left(\Gamma_{t} \, \frac{Z_{u}}{1+ |X_{u}|} |X_{u}|\right) \quad \text{\big(from point 2 {\it i)} of Theorem 1.1.6\big)} \\
&\le& \mathop{\rm \underline{lim}}_{u \to \infty}^{} \quad W (\Gamma_{t} \, Z_{u}) \quad \left(\text{since} \; \dis \frac{|X_{u}|}{1+|X_{u}|} \le 1\right)\\
& \le& W (\Gamma_{t} \, Z_{t})
\end{eqnarray*}

\noi since $(Z_{t}, \; t \ge 0)$ is a supermartingale. Finally : 
\begin{equation*}
{\bf W} (\Gamma_{t} \, z_{\infty}) = W\big(\Gamma_{t} \, M_{t} (z_{\infty})\big) \le W (\Gamma_{t} \cdot Z_{t})
\end{equation*}

\noi which is equivalent to point 2 of Lemma 1.2.9.

\smallskip

\noi \underline{We may now end the proof of Theorem 1.2.5.} 
\begin{equation*}
 \hspace*{-1cm} {\rm Let } \quad \widetilde{Z}_{t} := Z_{t} - M_{t} (z_{\infty}) \qquad (t \ge 0)
\end{equation*}

\noi Since $\big(M_{t} (z_{\infty}), \; t \ge 0\big)$ is a $\big((\mathcal{F}_{t}, \; t \ge 0), \; W\big)$ martingale, the process $(\widetilde{Z}_{t}, \; t \ge 0)$ is still a $\big((\mathcal{F}_{t}, \; t \ge 0), \; W\big)$ positive \big(from (1.2.61)\big) supermartingale, and since $\dis M_{t} (z_{\infty}) \mathop{\longrightarrow}_{t \to \infty}^{} 0 \quad W$ a.s. from Theorem 1.2.1, we obtain : 
\begin{equation*}
\widetilde{Z}_{t} \mathop{\longrightarrow}_{t \to \infty}^{} Z_{\infty} \qquad W \; {\rm a.s.}
\end{equation*}

\noi Since $(\widetilde{Z}_{t}, \; t \ge 0)$ is a positive supermartingale, we obtain :
\begin{equation*}
W(Z_{\infty} | \mathcal{F}_{t}) \le \widetilde{Z}_{t}
\end{equation*}

\noi We now write : 
\begin{equation*}
\xi_{t} := \widetilde{Z}_{t} - W (Z_{\infty} | \mathcal{F}_{t}) \qquad t \ge 0
\end{equation*}

\noi This is a positive supermartingale such that $\dis \mathop{\rm lim}_{t \to \infty}^{} \xi_{t} =0 \quad W$ a.s. On the other hand, $\bf W$ a.s.~: 
\begin{equation*}
 \mathop{\rm lim}_{t \to \infty}^{} \; \frac{\xi_{t}}{1+|X_{t}|} =  \mathop{\rm lim}_{t \to \infty}^{} \;\frac{\widetilde{Z}_{t}}{1+|X_{t}|} = z_{\infty} - z_{\infty} = 0
\end{equation*}

\noi The uniqueness of decomposition (1.2.52) being immediate, Theorem 1.2.5 is proven.

\bigskip

\noi{\bf 1.2.4} \underline{ A decomposition result for the martingale $\big(M_{t} (F), \; t \ge 0\big)$.}

\noi A difference with the preceding subsection is that the r.v.'s $F$ which we now consider belong to $L^{1} (\mathcal{F}_{\infty}, W)$, but are not necessarily positive.

\noi We shall now prove a decomposition result of the $\big((\mathcal{F}_{t}, \; t \ge 0), \; W\big)$ martingale $\big(M_{t} (F), \; t \ge 0\big)$. For this purpose, we shall use the following lemma.

\smallskip

\noi {\bf Lemma 1.2.10.} {\it Let  $F \in L^{1} (\mathcal{F}_{\infty}, {\bf W})$ 

\noi {\bf 1)} There exists a predictable process $\big(k_{s} (F), \; s \ge 0\big)$ which is defined $dL_{s}(\omega) W(d \omega)$ a.s., and is positive if $F$ is positive, such that :
$$
W \left( \int_{0}^{\infty} |k_{s} (F)|dL_{s}\right) = {\bf W} \big(|k_{g} (F)|\big) \le {\bf W} \big(|F|\big) < \infty
\eqno(1.2.62) $$

\noi and such that for every bounded predictable process $(\Phi_{s}, \; s \ge 0)$
\begin{eqnarray*}
\hspace*{2,5cm}{\bf W} (\Phi_{g} F)
&=& W \left(\int_{0}^{\infty} \Phi_{s} \, k_{s} (F) dL_{s}\right)\hspace*{5,5cm}(1.2.63) \\
&=& {\bf W} \big(\Phi_{g} \, k_{g} (F)\big) \hspace*{7,2cm}(1.2.64) \\
 \hspace*{-1cm}Thus :   \hspace*{1cm}{\bf W} (F|\mathcal{F}_{g}) &=& k_{g} (F) \hspace*{8,5cm}(1.2.65) \\
\end{eqnarray*} 

\noi {\bf 2)} We have ${\bf W} \big(|k_{g} (F)|\big) < \infty$ \big(from (1.2.62)\big)
$$
{\bf W} \big(|k_{g} (F)|\big) \le {\bf W} (|F|) < \infty \eqno(1.2.66) $$

\noi and
$$
\big(k_{s} (k_{g} (F), \; s \ge 0) \big) = (k_{s} (F), \; s \ge 0\big) \quad dL_{s} (\omega) \; W (d \omega) \qquad {\rm a.s.} \eqno(1.2.67) $$

\noi {\bf 3)} If $(h_{s}, \; s \ge 0)$ is a predictable process such that  ${\bf W} (|h_{g}|) < \infty$,  then : 
$$
\big(k_{s} (h_{g}), \; s \ge 0\big) = (h_{s}, \; s \ge 0) \quad dL_{s} (\omega) W (d�\omega) \qquad {\rm a.s.} 
\eqno(1.2.68) $$   }

\noi {\bf Proof of Lemma 1.2.10.}

\noi It suffices, by linearity, to prove this Lemma when $F \ge 0$.

\smallskip

\noi {\it i)} Formula (1.2.64), written for $F \equiv 1$ and $k_{s} (F) \equiv 1$ :
$$
{\bf W} (\Phi_{g}) = W \left(\int_{0}^{\infty} \Phi_{s} \, dL_{s}\right) \eqno(1.2.69) $$

\noi is formula (1.1.44). Let us define the measure $\mu_{F}$, on the predictable $\sigma$-field, and more generally on the set of positive predictable processes by : 
$$
\mu_{F} (\Phi) = {\bf W} (\Phi_{g} \cdot F) \eqno(1.2.70) $$

\noi Clearly, $\mu_{F}$ is absolutely continuous, on the predictable $\sigma$-field, with respect to $\mu_{1}$, which is the measure $\mu_{F}$ for $F \equiv 1$. Thus, from (1.2.69), $\mu_{F}$ is absolutely continuous on the predictable $\sigma$-field with respect to the measure $dL_{s} (\omega) W (d \omega)$. Thus, there exists, from the Radon-Nikodym Theorem, a process $\big(k_{s} (F), \; s \ge 0\big)$ which is predictable such that, for every $\Phi \ge 0$ predictable : 
\begin{equation*}
\mu_{F} (\Phi) = {\bf W} (\Phi_{g} \cdot F) = W \left(\int_{0}^{\infty} \Phi_{s} \, k_{s} (F) dL_{s}\right)
\end{equation*}

\noi This is relation (1.2.64). The further relations (1.2.65) and (1.2.66) follow immediately.

\smallskip

\noi {\it ii)} The other points of Lemma 1.2.10 are elementary. We show, for example, (1.2.68). We have, from (1.2.63) and (1.2.69), for every predictable and bounded process $\Phi$ : 
\begin{eqnarray*}
{\bf W} (\Phi_{g} \, h_{g}) 
&=& W \left(\int_{0}^{\infty} \Phi_{s} \, k_{s} (h_{g}) dL_{s}\right) \\
&=& W \left(\int_{0}^{\infty} \Phi_{s} \, h_{s} \, dL_{s}\right)
\end{eqnarray*}

\noi Hence, $\Phi$ being arbitrary, (1.2.68). Relation (1.2.67) is obtained by application of (1.2.68) with $(h_{s}, \; s \ge 0) = \big(k_{s} (F), \; s \ge 0\big)$. 

\smallskip

\noi Here is now the announced decomposition Theorem.

\smallskip

\noi {\bf Theorem 1.2.11.} {\it Let $F \in L^{1} ( \mathcal{F}_{\infty}, {\bf W})$. There exist two continuous positive processes $\dis \big(\Sigma_{t} (F)$, $t \ge 0\big)$ and $\big(\Delta_{t} (F), \; t \ge 0\big)$ such that, for every $t \ge 0$  :
$$
M_{t} (F) = \Sigma_{t} (F) + \Delta_{t} (F) \qquad (t \ge 0) \eqno(1.2.71) $$

\noi with : 

\smallskip

\noi {\bf 1)}{\it i)} For every $t \ge 0$ and $\Gamma_{t} \in b (\mathcal{F}_{t})$ :
$$
{\bf W} (\Gamma_{t} \, 1_{g \le t} \, F) = W\big(\Gamma_{t} \Sigma_{t} (F)\big) \eqno(1.2.72) $$

$\hspace*{-0,5cm}$ {\it ii)} $\big(\Sigma_{t} (F), \; t \ge 0\big)$ is a quasimartingale (a positive submartingale if $F \ge 0$) which vanishes on the zero set of $(X_{u}, \; u \ge 0)$. Its Doob-Meyer decomposition writes :
$$
\Sigma_{t} (F) = - M^{\Sigma (F)}_{t} + \int_{0}^{t} k_{s} (F) \, dL_{s} \eqno(1.2.73) $$

\noi In particular, the bounded variation part of this decomposition is absolutely continuous with respect to $dL_{s}$. In (1.2.73), $\big(M^{\Sigma (F)}_{t} , \; t \ge 0\big)$ is a $\big( (\mathcal{F}_{t}, \; t \ge 0), \quad W\big)$ martingale satisfying, if $F \ge 0$ : 
\begin{gather*}
\hspace*{1,5cm}\mathop{\rm sup}_{s \le t}^{} \, M^{\Sigma (F)}_{s} = \int_{0}^{t} k_{s} (F) dL_{s} \hspace*{7,4cm}(1.2.74) \\
\hspace*{1,5cm}\mathop{\rm lim}_{t \to \infty}^{} M^{\Sigma (F)}_{t} := M^{\Sigma (F)}_{\infty} = \int_{0}^{\infty} k_{s} (F) \, dL_{s} = \mathop{\rm sup}_{t \ge 0}^{} M^{\Sigma (F)}_{t} \hspace*{3cm}(1.2.75)
\end{gather*}

\noi In particular, this martingale is not uniformly integrable. 

$\hspace*{-0,5cm}$ {\it iii)} We have the "explicit formula" :
$$
\Sigma_{t} (F) = |X_{t}| \cdot \widehat{E}_{X_{t}}^{(3)} \big(F (\omega_{t}, \widehat{\omega}_{t})\big)
\eqno(1.2.76)$$ 

\noi (see point 1 of Remark 1.2.2 for such a notation).

\noi In (1.2.76), the expectation is taken with respect to $\widehat{\omega}_{t}$, the letter $\omega_{t}$, and $X_{t}$, being frozen ; $\widehat{E}_{X_{t}}^{(3)}$ denotes the expectation relatively to a 3-dimensional Bessel process starting from $X_{t}$, if $X_{t} > 0$, and the expectation with respect to the opposite of a 3-dimensional Bessel process, if $X_{t} < 0$.

\smallskip

$\hspace*{-0,5cm}$ {\it iv)} The application $F \to \big(\Sigma_{t} (F), \; t \ge 0\big)$ is injective since :
$$
\frac{\Sigma_{t} (F)}{1 + |X_{t}|} \; \mathop{\longrightarrow}_{t \to \infty}^{} F \qquad {\bf W} \; {\rm a.s.}
\eqno(1.2.77) $$

$\hspace*{-0,5cm}$ {\it v)} We have, for every $t \ge 0$ :
$$
W  \big\{ \Sigma_{t} (F) - \Sigma_{t} \big(k_{g} (F)\big) \big| \mathcal{F}_{g^{(t)}} \big\} = 0  \eqno(1.2.78) $$

\noi {\bf 2)}{\it i)} For every $t \ge 0$ and $\Gamma_{t} \in b (\mathcal{F}_{t})$ : 
$$
{\bf W} (\Gamma_{t} \, 1_{g > t} \, F) = W \big(\Gamma_{t} \, \Delta_{t} (F)\big) \eqno(1.2.79) $$

$\hspace*{-0,5cm}$ {\it ii)} $\big(\Delta_{t} (F), \; t \ge 0\big)$ is a quasimartingale (a positive supermartingale if $F \ge 0$). Its Doob-Meyer decomposition writes :
$$
\Delta_{t} (F) = M^{\Delta (F)}_{t} -  \int_{0}^{t} k_{s} (F) \, dL_{s} \eqno(1.2.80) $$

\noi where $\big(M^{\Delta (F)}_{t}, \; t \ge 0\big)$ is the $\big((\mathcal{F}_{t}, \; t \ge 0), \quad W\big)$ martingale given by :
$$
M^{\Delta (F)}_{t} = W \left(\int_{0}^{\infty} k_{s} (F) \, dL_{s} | \mathcal{F}_{t} \right) \eqno(1.2.81) $$

\noi In particular, since from (1.2.62), $\dis \int_{0}^{\infty} k_{s} (F) \, dL_{s} \in L^{1}  (\mathcal{F}_{\infty}, W)$, this martingale is uniformly integrable. 

\smallskip

$\hspace*{-0,5cm}$ {\it iii)}  The application $F \to \big(\Delta_{t} (F), \; t \ge 0\big)$ is not injective since :
$$
 \big(\Delta_{t} (F), \; t \ge 0\big) = \big(\Delta_{t} (k_{g} (F)\big), \; t \ge 0\big) \eqno(1.2.82) $$ 
 
 \noi \big(and $k_{g} (F) \neq F$ when $F$ is not $\mathcal{F}_{g}$ measurable\big).     
 
 \bigskip
 
 \noi {\bf 3)} The martingale $\big(M_{t} (F), \; t \ge 0\big)$ satisfies : 
$$
 \big(W \big(M_{t} (F) | \mathcal{F}_{g^{(t)}}\big), \; t \ge 0\big) = \big(W \big(M_{t} (k_{g} (F)\big) | \mathcal{F}_{g^{(t)}}\big), \; t \ge 0\big) \eqno(1.2.83) $$     }
 
\noi The following Theorem is an important consequence of Theorem 1.2.11.   

\newpage

\noi {\bf Theorem 1.2.12.} {\it Let $F \in L^{1} (\mathcal{F}_{\infty}, {\bf W})$. 

\noi Then, the $\big((\mathcal{F}_{t}, \; t \ge 0), \, W\big)$ martingale $\big(M_{t} (F), \; t \ge~0\big)$ vanishes on the zeros of $(X_{u}, \; u \ge 0)$ if and only if $k_{g} (F) =0$.    }

\smallskip

\noi {\bf Remark 1.2.13}

\noi {\bf 1)} If $F=F_{g}$, with $(F_{u}, \; u \ge 0)$ a positive previsible process Theorem 1.2.3 implies, in this particular case : 
\begin{equation*}
\Sigma_{t} (F_{g}) = F_{g^{(t)}} \cdot |X_{t}|, \quad \Delta_{t} (F_{g}) = \int_{L_{t}}^{\infty} W (F_{\tau_{l}} |
\mathcal{F}_{t}) \, dl.
\end{equation*}

\noi {\bf 2)} If $F \ge 0$, the supermartingale $\big(\Delta_{t} (F), \; t \ge 0\big)$ satisfies :
\begin{equation*}
\Delta_{t} (F) \mathop{\longrightarrow}_{t \to \infty}^{} 0 \qquad W \; {\rm a.s.}, \quad {\rm since} \;\; 0 \le \Delta_{t} (F) \le M_{t} (F)
\end{equation*}

\noi and 
\begin{equation*}
\frac{\Delta_{t} (F)}{1+|X_{t}|} = \frac{M_{t} (F)}{1+|X_{t}|} -  \frac{\Sigma_{t} (F)}{1+|X_{t}|} \mathop{\longrightarrow}_{t \to \infty}^{} F-F=0 \qquad {\bf W} \; {\rm a.s.}
\end{equation*}

\noi from Lemma 1.2.8 and (1.2.77). Hence, in the decomposition (1.2.52) of the supermartingale $\Delta_{t} (F)$, there remains uniquely the term $(\xi_{t}, \; t \ge 0)$.

\smallskip

\noi {\bf 3)} When $F \ge 0$, gathering the terms (1.2.71), (1.2.73), (1.2.80) and (1.2.81), we have : 
\begin{equation*}
M_{t} (F) = - M^{\Sigma (F)}_{t} + W \left(\int_{0}^{\infty} k_{s} (F) dL_{s} | \mathcal{F}_{t}\right)
\end{equation*}

\noi This formula implies \big(from (1.2.75)\big) that $(M_{t}^{\Sigma (F)}, \; t \ge 0)$ is not uniformly integrable since if it were, then $\big(M_{t} (F), \; t \ge 0\big)$ would be null.

\smallskip

\noi {\bf 4)} From relation (1.2.83) there exists an application 
\begin{eqnarray*}
m : \qquad L^{1} (\mathcal{F}, {\bf W})
&\longrightarrow & \mathcal{M} \big((\mathcal{F}_{g^{(t)}}, \; t \ge 0 ), \quad W \big) \\
F
&\longrightarrow & \big(m_{t} (F), \; t \ge 0\big)
\end{eqnarray*}

\noi where $\mathcal{M} \big((\mathcal{F}_{g^{(t)}}, \; t \ge 0 ), \; W\big)$ denotes the set of $\big((\mathcal{F}_{g^{(t)}}, \; t \ge 0 ), \; W\big)$ martingales ; this application $m$ is such that : 
$$
m_{t} (F) = W \big(M_{t} (k_{g} (F)) | \mathcal{F}_{g^{(t)}} \big) \eqno(1.2.84) $$

\noi with
\begin{equation*}
m_{t} (F) := \sigma_{t} (F) + \delta_{t} (F)
\end{equation*}

\noi and
\begin{eqnarray*}
\sigma_{t} (F) &=& \sqrt{\frac{\pi}{2}} \; k_{g^{(t)}} \, (F) \, \sqrt{t-g^{(t)}} \\
\delta_{t} (F) &=& W \left( \int_{0}^{\infty} k_{s} (F) \, dL_{s} | \mathcal{F}_{g^{(t)}}\right)
\end{eqnarray*}

\noi If $F \ge 0$, $\big(\sigma_{t} (F), \; t \ge 0\big)$ resp. $\big(\delta_{t} (F), \; t \ge 0\big)$ is a $\big( (\mathcal{F}_{g^{(t)}}, \; t \ge 0), \quad W\big)$ submartingale \big(resp. $\big((\mathcal{F}_{g^{(t)}}, \; t \ge 0), \quad W\big)$ supermartingale \big).

\smallskip

\noi {\bf 5)} We recall that by definition, a process $(Z_{t}, \; t \ge 0)$ is a quasimartingale if, for every $t \ge 0$ : 
\begin{equation*}
{\rm sup} \;\; W \left(\sum_{i=1}^{n-1}\big| W \big(Z_{t_{i+1}} - Z_{t_{i}} \big) \big| \, \, \big| {\mathcal{F}_{t_{i}}}  \right) < \infty
\end{equation*}

\noi the sup being taken over the set of subdivisions $0 \le t_{1} <\cdots < t_{n}<t$. In fact, such a process is the difference of two supermartingales \big(see [R]\big). On the other hand, the F\"{o}llmer measure \big(see [F]\big) $\mu_{Z}$ - with finite mass - of a supermartingale $(Z_{t}, \; t \ge 0)$ (or of a quasimartingale) is the measure defined on the predictable $\sigma$-field and characterised by : 
\begin{equation*}
\mu_{Z} (\Gamma_{t} \, 1_{]t, \infty]}) = W (\Gamma_{t} \cdot Z_{t}) \qquad \big(\Gamma_{t} \in b(\mathcal{F}_{t})\big)
\end{equation*}

\noi Hence formulae (1.2.65), (1.2.70) and (1.2.79) imply that the measure $\mu_{F}$ defined by (1.2.70) is the F\"{o}llmer measure of the quasimartingale $\big(\Delta_{t} (F), \; t \ge 0\big)$.

\smallskip

\noi {\bf Proof of Theorem 1.2.11.}

\noi {\it i)} \underline{We define $\Sigma_{t} (F)$ via} :
$$
\Sigma_{t} (F) = M_{t} (F \, 1_{g \le a})_{\big| a=t} \eqno(1.2.85) $$

\noi Hence, for every $\Gamma_{t} \in b (\mathcal{F}_{t})$ :
$$
{\bf W} (\Gamma_{t} \, 1_{g \le t} \cdot F) = W \big(\Gamma_{t} \, \Sigma_{t} (F)\big) \eqno(1.2.86)$$

\noi It is easy to deduce from (1.2.86) that $\dis \big(\Sigma_{t} (F) = \Sigma_{t} (F^{+}) - \Sigma_{t} (F^{-}), \; t \ge 0\big)$ is a semimartingale, as the difference of two submartingales and we shall show below (see  point {\it vi)} of this proof) that it is in fact a quasimartingale which admits a continuous version.

\smallskip

\noi {\it ii)} \underline{We show (1.2.73).}

\noi By linearity, it suffices to prove (1.2.73) for $F \ge 0$. From (1.2.86), we have  for $s \le t$ and  $\Gamma_{s} \in b (\mathcal{F}_{s})$ :
\hspace*{1cm}\begin{eqnarray*}
{\bf W} ({\Gamma_{s}} \, 1_{s\le g \le t} \, F)
&=& W \Big(\Gamma_{s} \big(\Sigma_{t} (F) - \Sigma_{s} (F)\big)\Big)  \\
&=& W \left(\Gamma_{s} \cdot \int_{s}^{t} k_{u} (F) \, dL_{u} \right) \hspace*{5,5cm}(1.2.87)
\end{eqnarray*}

\noi by using Lemma 1.2.10 with $(\Phi_{u} := \Gamma_{s} \, 1_{]s,t]} \, (u), \; u \ge 0)$. (1.2.73) follows immediately from (1.2.87).

\smallskip

\noi {\it iii)} \underline{We show (1.2.74) and (1.2.75).}

\noi Since, if  $F \ge 0$, then $\Sigma_{s} (F) \ge 0$, we have :
\begin{gather*}
\mathop{\rm sup}_{s \le t}^{} \, M^{\Sigma (F)}_{s} \le \int_{0}^{t} k_{u} (F) dL_{u} \qquad {\rm and} \\
\mathop{\rm sup}_{s \le t}^{} \, M^{\Sigma (F)}_{s} \ge \mathop{\rm sup}_{s \le g^{(t)}}^{} \; M^{\Sigma (F)}_{s} = \int_{0}^{g^{(t)}} k_{u} (F) dL_{u} = \int_{0}^{t} k_{u} (F) dL_{u}
\end{gather*}

\noi since $\Sigma_{g^{(t)}} (F) =0$ from (1.2.76) (which is proven below). 

\noi On the other hand, since $0 \le \Sigma_{t} (F) \le M_{t} (F)$ and since $M_{t} (F)$ $\dis \mathop{\longrightarrow}_{t \to \infty}^{} 0 \quad W$ a.s. from Theorem 1.2.1, we have $\dis \Sigma_{t} (F) \mathop{\longrightarrow}_{t \to \infty}^{} 0 \quad W$ a.s.,  and thus, from (1.2.73) :
\begin{equation*}
\mathop{\rm lim}_{t \to \infty}^{} \, M^{\Sigma (F)}_{t} := M_{\infty}^{\Sigma (F)} = \int_{0}^{\infty} k_{s} (F) \, dL_{s} = \mathop{\rm sup}_{t \ge 0}^{} \, M^{\Sigma (F)}_{t}
\end{equation*}

\noi which, in particular, proves,  that $(M^{\Sigma (F)}_{t}, \; t \ge 0)$ is not uniformly integrable.

\smallskip

\noi {\it iv)} \underline{We show (1.2.76).}

\noi For this purpose, we shall use the notation and results of subsection 1.1.4. We have, for every $t \ge 0$ and $\Gamma_{t} \in b (\mathcal{F}_{t})$ :
\begin{eqnarray*}
\lefteqn{ W\big(\Gamma_{t} \, \Sigma_{t} (F)\big) } \\
&=& {\bf W} (\Gamma_{t} \, 1_{g \le t} \, F)  \qquad \hbox{\big(from (1.2.86)\big)}  \\
&=& 2 \, W_{\infty}^{(\delta_{0})} (\Gamma_{t} \, 1_{g \le t} \, F \, e^{\frac{1}{2} \, L_{\infty}}) \\
&=& 2 \, W_{\infty}^{(\delta_{0})} (\Gamma_{t} \, e^{\frac{1}{2} \, L_{t}} \, 1_{g \le t} \, F) \quad
\big({\rm since} \;  L_{\infty} = L_{t} \; \hbox{on the set} \; (g \le t)\big)  \\
&=& 2 \, W_{\infty}^{(\delta_{0})} (\Gamma_{t} \, e^{\frac{1}{2} \, L_{t}} W_{\infty}^{(\delta_{0})} \big(1_{g \le t} \cdot F | \mathcal{F}_{t})\big)  \\
&=& 2 \, W_{\infty}^{(\delta_{0})} \big(\Gamma_{t} \, e^{\frac{1}{2} \, L_{t}} W_{\infty}^{(\delta_{0})} (1_{T_{0} \circ \theta_{t} = \infty} \cdot F | \mathcal{F}_{t})\big)  \quad
 \big({\rm since} \; (g \le t) = (T_{0} \circ \theta_{t} = \infty)\big)  \\
&=& 2 \, W_{\infty}^{(\delta_{0})} \big(\Gamma_{t} \, e^{\frac{1}{2} \, L_{t}} \, \widehat{W}_{X_{t}, \infty}^{(\delta_{0})} \big(1_{T_{0} = \infty} F (\omega_{t}, \widehat{\omega}^{t}) \big)  \\
&& \hbox{(by the Markov property} )  \\
&=& 2 \, W_{\infty}^{(\delta_{0})} \big(\Gamma_{t} \, e^{\frac{1}{2} \, L_{t}} \, \widehat{W}_{X_{t}, \infty}^{(\delta_{0})} \big(F(\omega_{t}, \widehat{\omega}^{t}) | T_{0} = \infty\big) \cdot W_{X_{t}, \infty}^{(\delta_{0})} (T_{0} = \infty)\big) \hspace*{2cm}(1.2.88)
\end{eqnarray*}

\noi But, from (1.1.70) :
\begin{equation*}
W_{X_{t}, \infty}^{(\delta_{0})} (T_{0} = \infty) = \frac{|X_{t}|}{2+|X_{t}|}
\end{equation*}

\noi and, from Theorem 1.1.5, conditionally on $(T_{0} = \infty)$, $W_{\infty, x}^{(\delta_{0})}$ is the law of a Bessel (3) process (resp. of the opposite of a Bessel (3) process) started at $x$ if $x > 0$ (resp. if $x < 0$). Then :
\begin{eqnarray*}
W \big(\Gamma_{t} \, \Sigma_{t} (F)\big)
&=& 2 \, W_{\infty}^{(\delta_{0})} \left( \Gamma_{t} \, e^{\frac{1}{2} \, L_{t}} \frac{|X_{t}|}{2+|X_{t}|} \; \widehat{E}_{X_{t}}^{(3)} \, \big(F (\omega_{t}, \widehat{\omega}^{t})\big) \right) \\
&=& W \left( \Gamma_{t} \, e^{\frac{1}{2} \, L_{t}} \, \frac{|X_{t}|}{2+|X_{t}|} \; \widehat{E}_{X_{t}}^{(3)} \big( F (\omega_{t}, \widehat{\omega}^{t})\big) \, e^{-\frac{1}{2} \, L_{t}} \big(2+|X_{t}|\big) \right)
\end{eqnarray*}

\noi \big(from (1.1.31) and (1.1.7)\big). 

\noi Finally  $ \dis W \big(\Gamma_{t} \, \Sigma_{t} (F)\big) = W \Big(\Gamma_{t} |X_{t}| \, \widehat{E}_{X_{t}}^{(3)} \big(F (\omega_{t}, \widehat{\omega}^{t})\big) \Big)$

\noi It is relation (1.2.76). Observe that this relation implies $\big(\Sigma_{t} (F), \; t \ge 0\big)$ vanishes on the zeros of $(X_{t}, \; t \ge 0)$. On the other hand, (1.2.76) implies (1.2.77), since, under ${\bf W}, \; \dis |X_{t}| \; \mathop{\longrightarrow}_{t \to \infty}^{} \infty$ a.s.

\smallskip

\noi {\it v)} \underline{We show (1.2.83) and (1.2.78).}

\noi For every positive, bounded and predictable process $(\Phi_{u}, \; u \ge 0)$, we have : 
$$
W \big(\Phi_{g^{(t)}} \, M_{t} (F)\big) = {\bf W} (\Phi_{g^{(t)}} \cdot F) \eqno(1.2.89) $$

\noi by definition of $M_{t} (F)$. But, the $\sigma$-algebra $\mathcal{F}_{g^{(t)}}$ is contained in $\mathcal{F}_{g}$. Hence the RHS of (1.2.89) equals from (1.2.64) : 
$$
{\bf W} \big(\Phi_{g^{(t)}} \, k_{g} (F)\big) = W\Big(\Phi_{g^{(t)}} \, M_{t} \, \big(k_{g} (F)\big)\Big) $$

\noi Finally :
$$
W\big(\Phi_{g^{(t)}} \, k_{g} (F)\big) = W\Big(\Phi_{g^{(t)}} \, M_{t} \, \big(k_{g} (F)\big)\Big) $$

\noi Thus $\;\; W\big(M_{t}(F) - M_{t}(k_{g} (F)) \big| \mathcal{F}_{g^{(t)}}\big) = 0$ i.e. (1.2.83) is satisfied. (1.2.78) is proven by using the same arguments.

\smallskip

\noi {\it vi)} \underline{We show (1.2.79).}

\noi We define $\; \Delta_{t} (F) \;$ by : 
$$ \Delta_{t} (F) := M_{t} (F \, 1_{g > a}) |_{a=t} $$

\noi It is clear that :
\begin{equation*}
M_{t} (F) = \Sigma_{t} (F) + \Delta_{t} (F)
\end{equation*}

\noi and that, for every $t \ge 0$ and $\Gamma_{t} \in b(\mathcal{F}_{t})$ :
$$ {\bf W} (\Gamma_{t} \, 1_{g > t} \, F) = W \big(\Gamma_{t} \, \Delta_{t} (F)\big) $$

\noi Then writing $\Delta_{t}(F)=\Delta_{t} (F^{+}) - \Delta_{t} (F^{-})$ we deduce easily from this formula that $\big(\Delta_{t} (F^{\pm}),$

\noi $t \ge 0\big)$ are two positive supermartingales and then $\big(\Delta_{t} (F), \; t \ge 0\big)$ is a quasimartingale. Since $\Sigma_{t} (F) = M_{t} (F) - \Delta_{t} (F) = M_{t} (F^{+}) - M_{t} (F^{-}) - \Delta_{t} (F^{+}) + \Delta_{t} (F^{-})$, it is clear that $\big(\Sigma_{t} (F), \; t \ge~0\big)$ is still a quasimartingale. Formula (1.2.80) then results from (1.2.73) and (1.2.71). Finally, thanks to (1.2.80) and (1.2.73), $\big(\Delta_{t} (F), \; t \ge 0\big)$ and $\big(\Sigma_{t} (F), \; t \ge 0\big)$ admit continuous versions.

\smallskip

\noi {\it vii)} \underline{We show (1.2.81).}

\noi We have, from (1.2.79), for every $\Gamma_{t} \in b(\mathcal{F}_{t})$

\begin{eqnarray*}
{\bf W} (\Gamma_{t} \, 1_{g>t} \, F) &=& W\big(\Gamma_{t} \, \Delta_{t} (F)\big) \\
&=& {\bf W} (\widetilde{\Gamma}_{g} \cdot F) \\
\Big({\rm with} \quad (\widetilde{\Gamma}_{u}, \; u \ge 0 ) 
&:=& \big(\Gamma_{t} \, 1_{]t, \infty[} (u), \; u \ge 0\big)\Big) \\
&=& W \left(\int_{0}^{\infty} \widetilde{\Gamma}_{u} \, k_{u} (F) dL_{u}\right) \;\; \hbox{(from Lemma 1.2.10)} \\
&=& W \left( \Gamma_{t} \cdot \int_{t}^{\infty} k_{u} (F) dL_{u} \right) \\
&=& W \left(\Gamma_{t} \, W \left(\int_{t}^{\infty} k_{u} (F) dL_{u} | \mathcal{F}_{t} \right) \right) \\
{\rm Hence} \; : \; \Delta_{t} (F)
&=& W \left(\int_{t}^{\infty} k_{u} (F) dL_{u} | \mathcal{F}_{t}\right) \\
&=& W \left( \int_{0}^{\infty} k_{u} (F) dL_{u} \Big| \mathcal{F}_{t} \right) - \int_{0}^{t} k_{u} (F) dL_{u}
\end{eqnarray*}

\noi This equality implies (1.2.80) and (1.2.81).

\smallskip

\noi {\it viii)} \underline{We show (1.2.82).}

\noi It suffices, to prove (1.2.82), to show that for every $t \ge 0$ and $\Gamma_{t} \in b(\mathcal{F}_{t})$, we have :
$$ W \big(\Gamma_{t} \, \Delta_{t} (F - k_{g} (F))\big) = 0 $$

\noi But :
\begin{gather*}
W\big(\Gamma_{t} \, \Delta_{t} (F- k_{g} (F))\big) = {\bf W} \big(\Gamma_{t} \, 1_{g>t} (F-k_{g} (F)\big) \\
\qquad = {\bf W} \big(\widetilde{\Gamma}_{g} (F- k_{g} (F)\big) 
\end{gather*}

\noi \Big(with $\quad \big(\widetilde{\Gamma}_{u} := \Gamma_{t} \, 1_{]t, \infty[} \, (u), \; u \ge 0\big)$\Big)
\begin{gather*}
\qquad = W \left(\Gamma_{t} \int_{t}^{\infty} \big(k_{u} (F) - k_{u} (k_{g} (F))\big) \, dL_{u}\right) \quad \hbox{\big(from (1.2.63)\big)} \\
\qquad =0
\end{gather*} 

\noi since $\quad k_{u} (F) = k_{u} (k_{g} (F))$, from (1.2.67).

\smallskip

\noi {\it ix)} \underline{Observe that}, by using (1.2.82), (1.2.83) is a consequence of (1.2.78). Indeed :
\begin{eqnarray*}
M_{t} \big(F-k_{g} (F)\big)
&=& \Sigma_{t} \big(F-k_{g} (F)\big) + \Delta_{t} \big(F-k_{g} (F)\big) \\
&=& \Sigma_{t} \big(F-k_{g} (F)\big) \qquad \hbox{from (1.2.82)}
\end{eqnarray*}

\noi Thus :
\begin{equation*}
W \big(M_{t} \big(F-k_{g} (F)\big) \big| \mathcal{F}_{g^{(t)}}\big) = W \big(\Sigma_{t} \big(F-k_{g} (F)\big) \big| \mathcal{F}_{g^{(t)}} \big) = 0 \quad�\hbox{from (1.2.78)}
\end{equation*}

\noi This ends the proof of Theorem 1.2.11.

\bigskip

\noi {\bf Proof of Theorem 1.2.12.}

\noi For this purpose, we need the following result, due to Az\'ema and Yor \big(see [AY2]\big) : a $\big( (\mathcal{F}_{t}, \; t \ge~0), \; W\big)$ martingale $(M_{t}, \; t \ge 0)$ vanishes on the zeros of $(X_{u}, \; u \ge 0)$ if and only if for every $t \ge 0$ :
$$ W(M_{t} | \mathcal{F}_{g^{(t)}}) = 0. \eqno(1.2.90) $$

\noi \underline{Suppose} $k_{g} (F) =0$

\noi From (1.2.83), we have : $\;\; W(M_{t} (F)| \mathcal{F}_{g^{(t)}}) = W \big(M_t(k_{g} (F)) | \mathcal{F}_{g^{(t)}}\big)=0$. Thus, from (1.2.90), \hbox{$\big(M_{t} (F), \; t \ge 0\big)$} vanishes on the zeros of $(X_{u}, \; u \ge 0)$.

\smallskip

\noi \underline{Conversely}, suppose that $\big(M_{t} (F), \; t \ge~0\big)$ vanishes on the zeros of $(X_{u} \ge 0)$. Then we have from (1.2.90) and (1.2.83), for every $s$ and $t$, $s \le t$ and $\Gamma_{s} \in b (\mathcal{F}_{s})$, since $\Gamma_{s} \, 1_{s \le g^{(t)}}$ is a $\mathcal{F}_{g^{(t)}}$ measurable r.v. : 
\begin{eqnarray*}
0
&=& W \big(\Gamma_{s} \, 1_{s <g^{(t)}} \, M_{t} (k_{g} (F))\big) \\
&=& {\bf W} \big(\Gamma_{s} \, 1_{s <g^{(t)}} \, k_{g} (F) \big) \mathop{\longrightarrow}_{t \to \infty}^{} {\bf W} \big( \Gamma_{s} \; 1_{s \le g} \, k_{g} (F)\big)
\end{eqnarray*}

\noi since $\dis g^{(t)} \mathop{\longrightarrow}_{t \to \infty}^{} g \quad {\bf W}$ a.s. 

\smallskip

\noi Thus : 
$$ {\bf W} \big(\Gamma_{s} \, 1_{s \le g} \, k_{g} (F) \big) = 0 $$

\noi We deduce from the monotone class Theorem that, for every bounded $\mathcal{F}_{g}$ measurable r.v. $\Phi$~:
$$ {\bf W} \big(\Phi \, k_{g} (F)\big) = 0. \eqno(1.2.91) $$

\noi i.e. $\;k_{g} (F) =0\;\;$ since $\;\;k_{g} (F)$ is $\mathcal{F}_{g}$-measurable. $\hfill \blacksquare$

\bigskip

\noi \underline{We end this subsection with some examples of decomposition (1.2.71).}

\smallskip

\noi {\bf Example 8.} Let $F := {\rm exp} \dis \left(- \frac{\lambda}{2} \, L_{\infty}\right)$. We have shown (Example 2) that :
$$ M_{t} (F) = \left( \frac{2}{\lambda} + |X_{t}|\right) \, e^{- \frac{\lambda}{2} \, L_{t}} \eqno(1.2.92) $$

\noi We then have : 
\begin{gather*}
M_{t} (F) = \Sigma_{t} (F) + \Delta_{t} (F) \qquad {\rm with} \nonumber \\
\Sigma_{t} (F) = |X_{t}| \, e^{- \frac{\lambda}{2} \, L_{t}},\; \Delta_{t} (F) = \frac{2}{\lambda} \, e^{- \frac{\lambda}{2} \, L_{t}} \hspace*{6,8cm} (1.2.93)
\end{gather*}

\noi Indeed, from (1.2.72) :
\begin{eqnarray*}
{\bf W} (\Gamma_{t} \, 1_{g \le t} \, e^{- \frac{\lambda}{2} \, L_{t}})
&=& W \big( \Gamma_{t} \, \Sigma_{t} (e^{-\frac{\lambda}{2} \, L_{\infty}})\big) \\
&=& {\bf W} (\Gamma_{t} \, 1_{g \le t} \, e^{-\frac{\lambda}{2} \, L_{t}}) \;\; \big({\rm since}  \; L_{\infty} = L_{t} \; \hbox{on the set} \; (g \le t)\big) \\
&=& W (\Gamma_{t} |X_{t}| e^{- \frac{\lambda}{2} \, L_{t}})
\end{eqnarray*}

\noi from point 2 {\it i)} of Theorem 1.1.6. 

\noi Thus :  
$$ \Sigma_{t} (e^{- \frac{\lambda}{2} \, L_{\infty}}) = |X_{t}| \, e^{- \frac{\lambda}{2} \, L_{t}} $$

\smallskip

\noi {\bf Example 9.} This example generalises Example 8. Let $q \in \mathcal{I}$ and $\dis F := \exp \left(- \frac{1}{2} \, A_{\infty}^{(q)}\right)$. We know (see Example 1) that : 
$$
M_{t} (e^{- \frac{1}{2} \, A_{\infty}^{(q)}}) = \varphi_{q} (X_{t}) \, \exp \left(- \frac{1}{2} \, A_{t}^{(q)}\right)
\eqno(1.2.94) $$

\noi Then :
$$
\Sigma_{t} (e^{- \frac{1}{2} \, A_{\infty}^{(q)}}) = \psi_{q} (X_{t}) \, e^{-\frac{1}{2} \, A_{t}^{(q)}}, \; \Delta_{t} (e^{ - \frac{1}{2} \, A_{\infty}^{(q)}}) = (\varphi_{q} - \psi_{q}) (X_{t}) e^{- \frac{1}{2} \, A_{t}^{(q)}}  \eqno(1.2.95)$$

\noi with $\psi_{q}$ solution of : 
\begin{gather*} 
\psi'' = q \, \psi \quad {\rm on} \; \mathbb{R} \backslash \{0\} \nonumber \\
\psi  (x) \mathop{\sim}_{|x| \to \infty}^{} |x|, \quad \psi (0) = 0 \hspace*{9cm}(1.2.96) 
\end{gather*}

\noi {\bf Proof of (1.2.95).} We have :
\begin{eqnarray*}
{\bf W} (\Gamma_{t} \, 1_{g \le t} \, e^{- \frac{1}{2} \, A_{\infty}^{(q)}}) 
&=& \varphi_{q} (0) \;W_{\infty}^{(q)} (\Gamma_{t} \, 1_{g \le t}) \qquad \big({\rm from} \; (1.1.16)\big), \\
&& \hbox{(with the notation of Theorems 1.1.1 and 1.1.2)}  \\
&=& \varphi_{q} (0) \, W_{\infty}^{(q)} \big(\Gamma_{t} \, W_{\infty}^{(q)} (1_{g \le t}| \mathcal{F}_{t})\big)  \\
&=& \varphi_{q} (0) \, W_{\infty}^{(q)} \big(\Gamma_{t} \, W_{\infty, X_{t}}^{(q)} (T_{0} = \infty)\big) \hspace*{4,4cm} (1.2.97)
\end{eqnarray*}

\noi But, by using the scale function $\gamma_{q}$ of the Markov process $(X_{t}, \; t \ge 0)$ under $W_{\infty}^{(q)}$, we have, with $\gamma_{q}$ given by (1.1.14) :
\begin{eqnarray*}
W_{\infty, x}^{(q)} (T_{0}= \infty)
&=& \frac{\gamma_{q} (x) - \gamma_{q} (0)}{\gamma_{q} (\infty) - \gamma_{q} (0)} \qquad {\rm if} \; x > 0 \\
&=& \frac{\gamma_{q} (0) - \gamma_{q} (x)}{\gamma_{q} (0) - \gamma_{q} (- \infty)} \quad {\rm if} \; x < 0 \hspace*{5cm} (1.2.98) \\
&:=& \lambda_{q} (x) 
\end{eqnarray*}

\noi Hence, by definition of $\Sigma_{t} (e^{-\frac{1}{2} \, A_{\infty}^{(q)}})$ :
\begin{eqnarray*}
\hspace*{1.2cm}W \big(\Gamma_{t} \, \Sigma_{t} \, (e^{-\frac{1}{2} \, A_{\infty}^{(q)}})\big)
&=& \varphi_{q} (0) \, W_{\infty}^{(q)} \big(\Gamma_{t} \, \lambda_{q} (X_{t})\big) \\
&=& W \big(\Gamma_{t} \, \varphi_{q} (X_{t}) \lambda_{q} (X_{t}) e^{-\frac{1}{2} \, A_{t}^{(q)}})\big) \hspace*{3cm}(1.2.99)
\end{eqnarray*}

\noi Thus :
$$ \Sigma_{t} (e^{- \frac{1}{2} \, A_{\infty}^{(q)}}) = \psi_{q} (X_{t}) \, e^{\frac{1}{2} \, A_{t}^{(q)}} $$

\noi with
$$ \psi_{q} (x)  := \lambda_{q} (x) \, \varphi_{q} (x) \eqno(1.2.100) $$

\noi It is clear, from (1.2.100), (1.2.98) and since $\dis \varphi_{q} (x) \mathop{\sim}_{|x| \to \infty}^{} |x|$ that : 
$$ \psi_{q} (x) \mathop{\sim}_{|x| \to \infty}^{} |x| \quad {\rm and} \quad \psi_{q} (0) =0. $$

\noi On the other hand, the relation $\psi''_{q} = q \, \psi_{q}$ on $\mathbb{R}$ is the consequence of direct calculation using the explicit form of $\gamma_{q}$ given by (1.1.14) (see Lemma 1.3.3 below for such a computation). We deduce from (1.2.95) and from It\^{o}-Tanaka :
\begin{eqnarray*}
\Sigma_{t} (e^{- \frac{1}{2} \, A_{\infty}^{(q)}}) 
&=& \int_{0}^{t} \psi'_{q} (X_{s}) \, e^{- \frac{1}{2} \, A_{s}^{(q)}} dX_{s}  + \frac{1}{2} \, \int_{0}^{t} \big(\psi'_{q} (0_{+}) - \psi'_{q} (0_{-})\big) \, e^{- \frac{1}{2} \, A_{s}^{(q)}} dL_{s} \\
{\rm i.e.}
&& M_{t}^{\Sigma (F)} = - \int_{0}^{t} \psi'_{q} (X_{s}) e^{- \frac{1}{2} \, A_{s}^{(q)}} \, dX_{s} \\
&& k_{s} (e^{- \frac{1}{2} \, A_{\infty}^{(q)}}) = \frac{1}{2} \, \big(\psi'_{q} (0_{+}) - \psi'_{q}  (0_{-})\big) \, e^{- \frac{1}{2} \, A_{s}^{(q)}}.
\end{eqnarray*}

\noi {\bf Example 10.} Let $\psi : \mathbb{R}_{+} \to \mathbb{R}_{+}$ Borel and integrable with $\psi (\infty) =0$ and $F:= \psi (S_{\infty})$. We know (see Example 3) that :
$$ M_{t} \big(\psi (S_{\infty})\big) = \psi (S_{t}) (S_{t} - X_{t}) + \int_{S_{t}}^{\infty} \psi (y) dy \eqno(1.2.101) $$

\noi We have : 
$$
\Sigma_{t} \big(\psi (S_{\infty})\big) = \psi (S_{t}) X_{t}^{-}, \quad \Delta_{t} \big(\psi (S_{\infty})\big) = \psi (S_{t}) (S_{t} - X_{t}^{+}) + \int_{S_{t}}^{\infty} \psi (y) dy \eqno(1.2.102) $$

\noi Indeed :
\begin{equation*}
{\bf W} \big(\Gamma_{t} \, 1_{g \le t} \, \psi(S_{\infty}\big) \big) = {\bf W}^{-} \big(\Gamma_{t} \, 1_{g \le t} \, \psi(S_{g}\big) \big)
\end{equation*}
\begin{eqnarray*}
\big({\rm since} \;\; \psi (\infty) = 0,  S_{\infty} 
=&& \hspace*{-0,5cm}\infty \;\; {\rm on} \;\; \Gamma_{+}, \; S_{\infty} = S_{g} \;\; {\rm on} \;\; \Gamma_{-}\big). 
\end{eqnarray*}
\begin{eqnarray*}
&=& {\bf W}^{-} \big(\Gamma_{t} \, 1_{g \le t} \, \psi (S_{g^{(t)}})\big) \hspace*{1cm} \big({\rm since} \;\; g^{(t)} =g \; {\rm on} \; (g \le t)\big) \\
&=& W \big(\Gamma_{t} \, \psi (S_{g^{(t)}}) X_{t}^{-}\big) \hspace*{1,5cm} (\hbox{from (1.1.52)}) \\
&=& W \big(\Gamma_{t} \, \psi (S_{t}) \, X_{t}^{-}\big) \hspace*{1,8cm} ({\rm since} \; S_{g^{(t)}} = S_{t} \quad {\rm if} \;\; X_{t} < 0) \\
{\rm Thus} && \; \Sigma_{t}  \big(\psi (S_{\infty})\big) = \psi (S_{t}) \, X_{t}^{-}
\end{eqnarray*}

\noi {\bf Example 11.} In some sense, the present example stands midway between Examples 9 and 10. Let $q : \mathbb{R} \to \mathbb{R}_{+}$ such that $q (x) =0$ if $x < 0, \; q (x)>0$ if $x >0$, $\dis \mathop{\underline{\rm lim}}_{x \to \infty}^{} \; q(x) >0$. We have shown, in [RY, IX] \big(see also [RY, M]\big) the existence for every $x \in \mathbb{R}$ of a $\sigma$-finite measure ${\bm \nu}_{x}^{(q)}$, on $\mathbb{R}_{+}$ such that :
$$
M_{t} \big(h (A_{\infty}^{(q)})\big) = \int_{\mathbb{R}_{+}} h (A_{t}^{(q)} + y) \,  {\bm \nu}_{X_{t}}^{(q)} (dy)
\eqno(1.2.103) $$

\noi for $h : \mathbb{R}_{+} \to \mathbb{R}_{+}$ sub-exponential at infinity.

\noi We then have :
\begin{eqnarray*}
\hspace*{1.2cm}\Sigma_{t} \big(h (A_{\infty}^{(q)})\big) 
&=& h (A_{t}^{(q)}) \cdot X_{t}^{-} \hspace*{7cm} (1.2.104) \\
\Delta_{t} \big(h (A_{\infty}^{(q)})\big) 
&=& \int_{\mathbb{R}_{+}} h(A_{t}^{(q)} + y) \big({\bm \nu}_{X_{t}}^{(q)} (dy) - X_{t}^{-} \delta_{0} (dy)\big) \\
&=&  \int_{\mathbb{R}_{+}} h(A_{t}^{(q)} + y) \, {\bm \nu}_{X_{t}}^{(q), a} (dy) \hspace*{5cm} (1.2.105)
\end{eqnarray*}

\noi where ${\bm \nu}_{X_{t}}^{(q), a}$ denotes the absolute continuous part of ${\bm \nu}_{X_{t}}^{(q)} $. Relation (1.2.104) is obtained from the same arguments as those used for relation (1.2.102) by noting that $\; 1_{X_{t} \le 0} \;\; d A_{t}^{(q)} =0$ and (1.2.105) results from : 
\begin{eqnarray*}
{\rm if} \; x < 0, \;\; {\bm \nu}_{x}^{(q)} (dy) &=& {\bm \nu}_{x}^{(q),a} (dy) + x^{-} \delta_{0} (dy) \\
{\rm if} \; x > 0, \;\; {\bm \nu}_{x}^{(q)} (dy) &=& {\bm \nu}_{x}^{(q),a} (dy) \hspace*{3cm} \hbox{\big(see [RY, IX]\big)}
\end{eqnarray*}

\noi {\bf Example 12.} Let $q : \mathbb{R} \to \mathbb{R}_{+}$ such that : 
\begin{equation*}
\int_{- \infty}^{0} \big(1 +|x|\big) q(x) dx < \infty \;;\; \mathop{\underline{\rm lim}}_{x \to \infty}^{} \, x^{2 \alpha} q(x) > 0 \quad {\rm for \; some} \; \alpha < 1 
\end{equation*}

\noi and $\dis A_{t}^{(q)} := \int_{0}^{t} q (X_{s}) ds$.  Let $\varphi_{q}$ the solution of $\varphi'' =q \, \varphi, \quad \varphi' (-\infty) =-1, \quad \varphi (+ \infty) =0$. Then, we have : 
\begin{eqnarray*}
\hspace*{2cm}&&M_{t} \left({\rm exp} - \frac{1}{2} \, A_{\infty}^{(q)}\right) = \varphi_{q} (X_{t}) \; {\rm exp} \left( - \frac{1}{2} \, A_{t}^{(q)}\right) \hspace*{2,5cm}(1.2.106) \\
{\rm and}
&& e^{- \frac{1}{2} \, A_{\infty}^{(q)}} \cdot {\bf W}^{-} = {\bf W} (e^{- \frac{1}{2} \, A_{\infty}^{(q)}}) \cdot W_{\infty}^{(q)} \hspace*{5cm}(1.2.107)
\end{eqnarray*}

\noi where the probability $W_{\infty}^{(q)}$ is characterised by 
$$
W_{\infty}^{(q)} |_{\mathcal{F}_{t}} = \frac{\varphi_{q} (X_{t})}{\varphi_{q} (0)} \; {\rm exp} \left(-\frac{1}{2} \, A_{t}^{(q)}\right) \cdot W|_{\mathcal{F}_{t}} \eqno(1.2.108) $$

\noi \big(see [RVY, I], the one-sided case, p. 209\big). We then have :
$$ \Sigma_{t} (e^{- \frac{1}{2} \, A_{\infty}^{(q)}}) = \psi_{q} (X_{t}) e^{- \frac{1}{2} \, A_{t}^{(q)}} \eqno(1.2.109) $$

\noi with
\begin{gather*}
\psi_{q} (x) = 0 \quad {\rm if} \;\; x \ge 0 \\
\psi_{q} (x) \mathop{\sim}_{x \to - \infty}^{} |x| \quad {\rm and} \quad \psi''_{q}=q \, \psi_{q} \;\; {\rm on} \; \mathbb{R}_{-}
\end{gather*}

\noi Hence
\begin{gather*} 
\Delta_{t} (e^{- \frac{1}{2} \, A_{\infty}^{(q)}}) =
\left\{ \begin{array}{rl}
M_{t} (e^{- \frac{1}{2} \, A_{\infty}^{(q)}}) \hspace*{2cm} {\rm if} \; X_{t} \ge 0 \hfill \\
(\varphi_{q} - \psi_{q}) (X_{t}) \, e^{- \frac{1}{2} \, A_{t}^{(q)}} \quad {\rm if} \; X_{t} \le 0 
\end{array} \right. 
\end{gather*}

\noi (1.2.109) is obtained by following the same arguments as those in Example 9. What changes is that, under the probability $W_{x, \infty}^{(q)}$, we have $X_{t} \to - \infty$ a.s., for every $x$ \big(see Theorem 5.1 in [RVY, I]\big).

\smallskip

\noi {\bf Example 13.} Let $\psi : \mathbb{R}_{+} \to \mathbb{R}_{+}$ Borel and integrable, such that $\dis \int_{0}^{\infty} \psi (y) dy=1$. Then we have, from (1.2.24)
$$
M_{t} \big(\psi (S_{g})\big) = \frac{1}{2} \, \psi (S_{g^{(t)}}) |X_{t}| + \psi (S_{t}) (S_{t}-X_{t}^{+}) + \int_{S_{t}}^{\infty} \psi(y) dy \eqno(1.2.110) $$

\noi \big(see Example 4, (1.2.24) and (1.2.25)\big) ; 
$$ {\bf W}^{-} \big(\psi (S_{g})\big) = {\bf W} \big(\psi (S_{\infty})\big) = \int_{0}^{\infty} \psi (l) dl $$

\noi On the other hand, we have : 
\begin{eqnarray*}
{\bf W} \big(\Gamma_{t} \, 1_{g \le t} \, \psi (S_{g})\big)
&=& W \big(\Gamma_{t} \, \Sigma_{t} \big(\psi (S_{g})\big)\big) \\
&=& {\bf W} \big(\Gamma_{t} \, 1_{g \le t} \, \psi (S_{g^{(t)}})\big) \quad \big({\rm since} \;\; g = g^{(t)}\;\; {\rm  on} \;\; (g \le t)\big) \\
&=&  W \big(\Gamma_{t} \, \psi (S_{g^{(t)}}) |X_{t}|\big) \quad \; \hbox{(from point 2 {\it i)} of Theorem 1.1.6.) }
\end{eqnarray*}

\noi Hence : 
$$ \Sigma_{t} \big(\psi (S_{g})\big) = \psi (S_{g^{(t)}}) |X_{t}| \eqno(1.2.111) $$

\noi and, from (1.2.110) : 
$$
\Delta_{t} \big(\psi (S_{g})\big) = \psi (S_{t})(S_{t} - X_{t}^{+}) + \int_{S_{t}}^{\infty} \psi (y) dy \eqno(1.2.112) $$

\bigskip

\noi {\bf 1.2.5} \underline{A penalisation Theorem, for functionals in class $\mathcal{C}$}

\noi In Section 1 of this Chapter, we constructed the measure {\bf W} from the penalisation results, and more particularly from Feynman-Kac type penalisations. We shall now operate in a reverse order : starting from the existence and the properties of the measure {\bf W} which we just established, we shall obtain penalisation results.

\noi Here is the class of functionals $(F_t, \; t \geq 0)$ for which we shall obtain such a penalisation result.

\noi {\bf Definition 1.2.13. } Let $(F_t, \; t \geq 0)$ denote an adapted, positive process. We shall say 
that this process belongs to the class $\mathcal{C}$ if 

\noi {\it i)} $(F_t, \; t \geq 0)$ is a decreasing process, i.e. if $s \leq t$ : 
$$ 0 \leq F_t \leq F_s \quad W \; \textrm{a.s.} \eqno(1.2.113)$$
\noi In particular, since $0 \leq F_t \leq F_0$ and since $F_0$ is a.s. constant, this process
is bounded by a constant $C=F_0$.







\noi {\it ii)} There exists $a \ge 0$ such that for every $t \ge \sigma_{a}$, with :
$$ \sigma_{a} := \sup \{ t \ge 0 \;;\; X_{t} \in [-a,a]\}  $$

\noi we have :
$$ F_{t} = F_{\sigma_{a}} = F_{\infty} \eqno(1.2.114) $$

\noi {\it iii)}
$$ {\bf W} (F_{\infty}) = {\bf W} ( F_{\sigma_{a}}) < \infty \eqno(1.2.115) $$

\smallskip

\noi One of the advantages of this class $\mathcal{C}$ is that it contains a large number of processes $(F_{t}, \; t \ge 0)$ for which we have already obtained a penalisation result. More precisely, let $\varphi : \mathbb{R}^{n} \to \mathbb{R}_{+}$ Borel. Then :
\begin{equation*}
F_{t} := \varphi (L_{t}^{a_{1}}, \cdots L_{t}^{a_{r}}, A_{t}^{(q_{1})}, \cdots A_{t}^{(q_{s})}, D_{t}^{[\alpha_{1}, \beta_{1}]}, \cdots D_{t}^{[\alpha_{u}, \beta_{u}]}, S_{g^{(t)}}, -I_{g^{(t)}})
\end{equation*}

\noi (see Examples 1 to 9 for these notations) belongs to the class $\mathcal{C}$ \big(if (1.2.115) is satisfied\big) as soon as $q_{1}, \cdots, q_{s}$ are elements of $\mathcal{I}$ with compact support (if we choose $a$ large enough) and $\varphi$ is a function which is decreasing with respect to each of its arguments. We may add $S_{t}$ and $(-I_{t})$ to the list of the arguments of $\varphi$, if $\varphi$ has compact support in these arguments.

\noi One can give some examples of functionals $(F_t, \; t \geq 0)$ which are not in the class $\mathcal{C}$
and for which the statement of Theorem 1.2.14 below does not apply. One of these examples is the functional : 
$$ \left( F_t = \exp \left(- \int_{- \infty}^{\infty} (L_t^y)^2 \, dy \right), \; t \geq 0 \right)$$
(see [N3] for a study of this functional). 

\noi Here is the first step towards a penalisation result.

\smallskip

\noi {\bf Theorem 1.2.14.} {\it Let $(F_{t}, \; t \ge 0)$ be a process which belongs to $\mathcal{C}$. Then :  }

\noi {\bf 1)}
$$
\sqrt{\frac{\pi t}{2}} \; W(F_{t}) \mathop{\longrightarrow}_{t \to \infty} {\bf W} (F_{\infty}) \eqno(1.2.116) $$
\noi {\bf 2)}

$$ W (F_t. |X_t|) \underset{t \rightarrow \infty}{\longrightarrow} {\bf W} (F_{\infty}) \eqno(1.2.117) $$

\noi {\bf Proof of Theorem 1.2.14.}

\noi {\bf 1)} \underline{We start with the proof of point 1)}

\noi We write $F_{t}$ in the form :
\hspace*{1.2cm}\begin{gather*} 
F_{t} = F_{t} \, \frac{(|X_{t}|-a)_{+}}{1+|X_{t}|} + F_{t} \, \frac{1+ |X_{t}|-(|X_{t}|-a)_{+}}{(1+ |X_{t}|)^{2}} \,
(|X_{t}|-a)_{+}  \\
\qquad + F_{t} \, \frac{\big[(1+|X_{t}|)-(|X_{t}|-a)_{+}\big]^{2}}{(1+|X_{t}|)^{2}} := F_{t}^{(1)} +F_{t}^{(2)} + F_{t}^{(3)} \hspace*{4cm}(1.2.118)
\end{gather*}

\noi and we study each term of this decomposition of $F_{t}$. 

\smallskip

\noi {\it i)} \underline{Study of $W(F_{t}^{(1)})$.}

\noi For $\lambda > 0$, we have : 
\begin{gather*}
\int_{0}^{\infty} e^{- \lambda t} W(F_{t}^{(1)})dt = \int_{0}^{\infty} e^{- \lambda t} W\left( F_{t} \, \frac{(|X_{t}|-a)_{+})}{1+ |X_{t}|}\right) dt  \\
\qquad = \int_{0}^{\infty} e^{-\lambda t} {\bf W} \left(\frac{F_{t}}{1+|X_{t}|} \; 1_{\sigma_{a} \le t}\right) dt \\
\hbox{\big(by Theorem 1.1.8, relation (1.1.49) \big)}  \\
\qquad = \int_{0}^{\infty} e^{-\lambda t} {\bf W} \left(F_{\sigma_{a}} \; \frac{1_{\sigma_{a}\le t}}{1+|X_{t}|} \right) dt  \qquad \hbox{\big(from (1.2.114)\big)}  \\
\qquad = {\bf W} \left( F_{\sigma_{a}} \, e^{-\lambda \sigma_{a}} \int_{0}^{\infty} e^{-\lambda u} \frac{du}{1 + |X_{\sigma_{a} +u}|} \right)  \\
\hbox{(after the change of variable} \;\; t = \sigma_{a}+u) \\
\qquad = {\bf W} (F_{\infty} e^{- \lambda \sigma_{a}}) \; E_{0}^{(3)} \left(\int_{0}^{\infty} e^{- \lambda u} \frac{du}{1+a+R_{u}} \right) \hspace*{5cm}(1.2.119)
\end{gather*}

\noi from point 2 of Theorem 1.1.8, where in (1.2.119) $\; (R_{u}, \; u \ge 0)$ denotes a Bessel process of dimension 3 started at 0. But
$$E_{0}^{(3)} \left[\frac{1}{1+a+R_{u}}\right] \mathop{\sim}_{u \to \infty}^{} \sqrt{\frac{2}{\pi u}}  $$

\noi and is a decreasing function of $u$. By the (easy part of the) Tauberian Theorem \big(see [Fe]) :
$$
\int_{0}^{\infty} e^{-\lambda t} \, W (F_{t}^{(1)}) dt \mathop{\sim}_{\lambda \to 0}^{} {\bf W} (F_{\infty}) \sqrt{\frac{2}{\lambda}} \eqno(1.2.120) $$

\noi {\it ii)} \underline{Study of $W(F_{t}^{(2)})$.}

\noi For $\lambda > 0$, we have :
\begin{equation*}
\int_{0}^{\infty} e^{-\lambda t} \, W(F_{t}^{(2)}) dt \le (1+a) \int_{0}^{\infty} e^{- \lambda t} W \left(F_{t} \; \frac{(|X_{t}|-a)_{+}}{(1+|X_{t}|)^{2}}\right) dt
\end{equation*}

\noi \big(from (1.2.118) and since : $\dis 0 \le 1 + |X_{t} | - (|X_{t}| -a)_{+} \le 1+a \big)$
\begin{gather*}
\qquad = (1+a) {\bf W} (F_{\infty} \, e^{- \lambda \sigma_{a}} ) \, E_{0}^{(3)} \left( \int_{0}^{\infty} e^{-\lambda u} \frac{du}{ (1+a+R_{u})^{2} }\right)  \\
\hbox{\big(by using the same argument as in point {\it i)}\big)}  \\
\qquad \le (1+a) \, {\bf W} (F_{\infty}) \int_{0}^{\infty} e^{-\lambda u} E_{0}^{(3)} \left( \frac{1}{ (1+a+R_{u})^{\frac{3}{2}} }\right) du  \\
\qquad \le (1+a) \, {\bf W}(F_{\infty}) O \left(\frac{1}{\lambda^{\frac{1}{4}}}\right) = o \left( \frac{1}{ \sqrt{\lambda}} \right) \qquad (\lambda \to 0) \hspace*{4cm}(1.2.121)
\end{gather*}

\noi {\it iii)} \underline{Study of $W(F_{t}^{(3)})$.}
$$ W(F_{t}^{(3)}) \le (1+a)^{2} C \, W \left( \frac{1}{1+|X_{t}|^{2}}\right) $$

\noi from (1.2.118). Hypothesis {\it i)} : $0 \le F_{t} \le C$ imply :
$$ \sqrt{\frac{\pi t}{2}} \; W(F_{t}^{(3)}) \mathop{\longrightarrow}_{t \to \infty}^{} 0  \eqno(1.2.122) $$

\noi Thus :
$$ \int_{0}^{\infty} e^{- \lambda t} W(F_{t}^{(3)})dt = o \left(\frac{1}{\sqrt{\lambda}} \right)\qquad (\lambda \to 0) \eqno(1.2.123) $$

\noi Gathering (1.2.120), (1.2.121) and (1.2.123) we obtain :
$$ 
\int_{0}^{\infty} e^{-\lambda t} W(F_{t}) dt \mathop{\sim}_{\lambda \to 0}^{} \sqrt{\frac{2}{\lambda}} \; {\bf W} (F_{\infty}) \eqno(1.2.124) $$

\noi $W(F_{t})$ being by hypothesis a decreasing function in $t$, the Tauberian Theorem implies :
$$ \sqrt{\frac{\pi t}{2}} \; W (F_{t}) \mathop{\longrightarrow}_{t \to\infty}^{} {\bf W} (F_{\infty}) $$

\noi This is precisely the statement of point 1) of Theorem 1.2.14.

\smallskip

\noi {\bf 2)} \underline{We now prove point 2 of Theorem 1.2.14}

\noi We write
\begin{eqnarray*}
W \big( F_{t} \cdot |X_{t}|\big)
&=& W\big(F_{t} \big(|X_{t}|-a\big)_{+}\big) + W \big(F_{t} \big(|X_{t}|-\big(|X_{t}|-a\big)_{+}\big) \\
&:=& (1_{t}) + (2_{t})
\end{eqnarray*}

\noi and we study successively $(1_{t})$ and $(2_{t})$.
\begin{eqnarray*}
\cdot (1_{t})
&=& W\big(F_{t}(|X_{t}|-a)_{+}\big)  = {\bf W} (F_{t} \, 1_{\sigma_{a} \le t}) \qquad \hbox{(from Theorem 1.1.8)} \\
&=& {\bf W} (F_{\infty} \, 1_{\sigma_{a} \le t}) \qquad \hbox{(from (1.2.114))} \\
&\dis \mathop{\longrightarrow}_{t \to \infty}^{}& {\bf W} (F_{\infty}) \qquad \hbox{(by the monotone convergence Theorem)} \\
&& \big({\rm since} \; F_{\infty} \in L^{1}_{+} (\mathcal{F}_{\infty}, {\bf W})\big) \\
\cdot (2_{t}) 
&=& W \big(F_{t} \big(|X_{t}|-\big(|X_{t}|-a\big)_{+}\big) \le a \; W(F_{t})
\end{eqnarray*}

\noi We now write :
\begin{eqnarray*}
\hspace*{1cm}W(F_{t})
&=& W \left( F_{t} \; \frac{ 1 + |X_{t}| - \big(|X_{t}|-a \big)_{+} }{1+ |X_{t}|}  \right) + W  \left( F_{t} \; \frac{\big(|X_{t}|-a\big)_{+}}{1+|X_{t}|} \right) \hspace*{1cm} (1.2.125)\\
&=& (3_{t}) + (4_{t}) \qquad \hbox{and we have} \\
(3_{t})
&=& W \left( F_{t} \; \frac{1+|X_{t}|-\big(|X_{t}|-a\big)_{+}}{1+ |X_{t}|} \right) \le (1+a) W \left(\frac{F_{t}}{1+ |X_{t}|}\right) \\
&& \le (1+a) \, C \; W \;\left(\frac{1}{1+ |X_{t}|}\right) \mathop{\longrightarrow}_{t \to \infty}^{} 0 \\
\hbox{since}
&& (F_{t}, \; t \ge 0) \; \hbox{is bounded} \\
(4_{t})
&=& W \left( F_{t} \; \frac{\big(|X_{t}|-a\big)_{+}}{1+|X_{t}|} \right) = {\bf W} \left( \frac{F_{t}}{1+|X_{t}|} \; 1_{\sigma_{a} \le t} \right) \qquad \hbox{(from Th. 1.1.8.)} \\
&=& {\bf W} \left( \frac{F_{\infty}} {1+ |X_{t}|} \; 1_{\sigma_{a} \le t}\right) \qquad \hbox{(from (1.2.114))} \\
&\dis \mathop{\longrightarrow}_{t \to \infty}^{}& 0 \qquad {\rm since} \; |X_{t}| \mathop{\longrightarrow}_{t \to \infty}^{}  + \infty \quad  {\bf W} \;\; \hbox{a.s. and we apply the dominated convergence} \\
&&\hbox{ Theorem.}
\end{eqnarray*}

\noi This ends the proof of Theorem 1.2.14. We are now able to state the announced penalisation Theorem. 

\noi {\bf Theorem 1.2.15.} {\it (General penalisation Theorem)

\noi Let $(F_{t}, \; t \ge 0)$ be a process belonging to $\mathcal{C}$. Then, for every $s \ge 0$ and $\Gamma_{s} \in b(\mathcal{F}_{s})$ :

\noi {\bf 1)} The limit, as $\dis t \to \infty$, $\dis \frac{W(\Gamma_{s} F_{t})}{W (F_{t})}$ exists \hfill{\rm (1.2.126)}

\noi {\bf 2)} This limit equals :
$$
\mathop{\lim}_{t \to \infty}^{} \; \frac{W(\Gamma_{s} F_{t})}{W(F_{t})} = \frac{W\big(\Gamma_{s} M_{s}(F_{\infty})\big)}{{\bf W} (F_{\infty})} := W_{\infty}^{F} (\Gamma_{s}) \eqno(1.2.127) $$

\noi The probability $W_{\infty}^{F}$, which is characterised by (1.2.127) satisfies :}
$$  W_{\infty}^{F} = \frac{F_{\infty}}{{\bf W} (F_{\infty})} \cdot {\bf W} \eqno(1.2.128) $$

\noi By comparing (1.2.128) with (1.1.16'), (1.1.93), (1.1.94), (1.1.108), (1.1.109) and (1.1.112), 
one can see that Theorem 1.2.15 is a general Theorem which implies many results given in Section 1.1 of 
this monograph, for example Theorems 1.1.1, 1.1.2, 1.1.11 and 1.1.11'. 

\noi {\bf Proof of Theorem 1.2.15.}

\noi {\it i)} We shall use the following notations : let $\omega_{s} \in \mathcal{C} \big([0,s] \to \mathbb{R}\big)$ and $(F_{t}^{(\omega_{s})}, \; t \ge 0)$ the functional defined by :
$$
F_{t}^{(\omega_{s})} (X_{u}, \; u \ge 0) := F_{t +s} \big(\omega_{s} \circ \big(\omega_{s} (s) + X_{u}, \; u \ge 0\big)\big) \eqno(1.2.129) $$

\noi With this notation, we have the following Lemma.

\smallskip

\noi {\bf Lemma 1.2.16.} {\it If $(F_{t}, \; t \ge 0) \in \mathcal{C}$, then, for $W$-almost every $\omega_{s} \in \mathcal{C} \big([0,s] \to \mathbb{R}\big)$ $(F_{t}^{(\omega_{s})}, \; t \ge 0) \in \mathcal{C}$.}

\smallskip

\noi {\bf Proof of Lemma 1.2.16.}

\noi {\it i)} It is clear that $(F_{t}^{(\omega_{s})}, \; t \ge 0)$ is a monotone function of $t$ and that, from (1.2.129) and (1.2.114) we have, for $t \ge \sigma_{|\omega_{s} (s)|+a}$ :
\begin{equation*}
F_{t}^{(\omega_{s})} (X_{u}, \; u \ge 0) = F_{\sigma_{|\omega_{s} (s)|+a}}^{(\omega_{s})} (X_{u}, \; u \ge 0) = F_{\infty}^{(\omega_{s})} (X_{u}, \; u \ge 0)
\end{equation*}

\noi {\it ii)} We need to prove that ${\bf W} (F_{\infty}^{(\omega_{s})}) < \infty$. We note that : 
\begin{eqnarray*}
{\bf W}(F_{\infty}^{(\omega_{s})})
&=& {\bf W} \big(F_{\infty} \big(\omega_{s} \circ \big(\omega_{s} (s)+ X_{u}, \; u \ge 0\big)\big) \\
&=& M_{s} (F_{\infty}) (\omega_{s}) \qquad \big({\rm from} (1.2.3)\big)
\end{eqnarray*}

\noi Hence : 
\begin{equation*}
W\big({\bf W} (F_{\infty}^{(\omega_{s})})\big) = W \big(M_{s} (F_{\infty})\big) = {\bf W}(F_{\infty}) < \infty \qquad \text{\big(from (1.2.2)\big)}
\end{equation*}

\noi In particular :
\begin{equation*}
{\bf W} (F_{\infty}^{(\omega_{s})}) < \infty \qquad W \;\; {\rm a.s.}
\end{equation*}

\noi This is Lemma 1.2.16.

\smallskip

\noi {\it ii)} We may now end the proof of Theorem 1.2.15. We have, for $t \ge s$ :
\begin{eqnarray*}
\hspace*{0,2cm}\frac{W \big(F_{t} | \mathcal{F}_{s}\big)}{W (F_{t})}
&=& \frac{\widehat{W} (F_{t-s}^{(\omega_{s})})}{W(F_{t})} \quad \hbox{(from the Markov property)} \nonumber \\
&=& \frac{ \sqrt{\frac{\pi t}{2}} \; \widehat{W}(F_{t-s}^{(\omega_{s})})}{ \sqrt{\frac{\pi t}{2}} \; W (F_{t})} \mathop{\longrightarrow}_{t \to \infty}^{} \frac{{\bf W}(F_{\infty}^{(\omega_{s})})}{{\bf W}(F_{\infty})} \qquad {\rm a.s.} \hspace*{3,3cm} (1.2.130)
\end{eqnarray*}

\noi \big(from Theorem 1.2.14 applied to $(F_{t}, \; t \ge 0)$ and to $(F_{t}^{(\omega_{s})}, \; t \ge 0)$ due to Lemma 1.2.16.\big)
\begin{equation*}
\hspace*{2.2cm} = \frac{M_{s} (F_{\infty})}{{\bf W} (F_{\infty})}
\end{equation*}

\noi (from point 2 of Theorem 1.2.1.)

\smallskip

\noi To show Theorem 1.2.15, it now suffices to see that the convergence is (1.2.130) also holds in $L^{1} (\mathcal{F}_{\infty}, W)$. However, from Scheff\'e's Lemma \big(see [M], T. 21\big) this is implied by the equality : $W \left(\frac{M_{s}(F_{\infty})}{{\bf W}(F_{\infty})} \right) =1$ for every $s \ge 0$, which follows immediately from Theorem 1.2.1 (equality (1.2.2)).

\smallskip

\noi {\bf Remark 1.2.17.}

\noi Let $\varphi : \mathbb{R}_{+} \to \mathbb{R}_{+}$ Borel such that : $\dis \int_{0}^{\infty} \varphi(x) (1+ x^{2}) dx < \infty$ and let :
\begin{eqnarray*}
&&F_{t}^{(1)} := \varphi (S_{t}) \;1_{(X_{t} >0)} \qquad (t \ge 0) \\
&& F_{t}^{(2)} := \varphi (S_{d_{t}}) 1_{(X_{t} > 0)} \quad  \;\; \, (t \ge 0)
\end{eqnarray*}

\noi It is shown in [RY, VIII] that :

\noi {\it i)} $\quad \dis E \big(F_{t}^{(1)}\big) \mathop{\sim}_{t \to \infty} \frac{3}{2} \sqrt{\frac{2}{\pi t^{3}}} \int_{0}^{\infty} \varphi (x) x^{2} dx \hfill(1.2.131)$

$\; \dis E \big(F_{t}^{(2)}\big) \mathop{\sim}_{t \to \infty} \sqrt{\frac{2}{\pi t^{3}}} \int_{0}^{\infty} \varphi (x) x^{2} dx \hfill(1.2.132)$

\noi {\it ii)} for every $s \ge 0$ and $\Gamma_{s} \in b (\mathcal{F}_{s})$ 

$\dis \frac{E [\Gamma_{s} \, F_{t}^{(i)}]}{E (F^{(i)}_{t})} \mathop{\longrightarrow}_{t \to \infty} E(\Gamma_{s} \, M_{s}^{\psi}) \qquad (i=1,2) \hfill(1.2.133)$

\noi where the martingale $\big(M_{s}^{\psi}, \; s \ge 0\big)$ is defined by :
\begin{eqnarray*}
&& M_{s}^{\psi} = \psi (S_{s}) (S_{s} - X_{s}) + \int_{S_{s}}^{\infty} \psi (y) dy \\
{\rm and}
&& \psi (x) := \varphi (x) x^{2} + 2 \int_{x}^{\infty} \varphi (y) y \, dy \qquad (x \ge 0)
\end{eqnarray*}

\noi We now inspect Theorem 1.2.15 in the light of this result. If we assume that $\dis \mathop{\rm lim}_{y \to + \infty} \varphi (y) =0$, we obtain :
$$   \mathop{\rm lim}_{t \to \infty}  \; F_{t}^{(i)} = 0 \qquad {\bf W} \qquad {\rm a.s.} $$

\noi and, from (1.2.131) and (1.2.132).
$$\dis \mathop{\rm lim}_{t \to \infty}  \sqrt{t} \, E [F_{t}^{(i)}]= 0 \quad (i=1,2)$$

\noi Thus, we are working here in a degenerate case of Theorem 1.2.15 and of Theorem 1.2.14, i.e. : in a case where $F_{\infty} \equiv 0$. However, from (1.2.133), this situation is not so "degenerate", since it allows to obtain a non-trivial penalisation.

\bigskip

\noi {\bf 1.2.6} \underline{Some other results about the martingales $(M_t(F), \; t \geq 0)$.}

\smallskip

\noi Let us first state the following definition : 

\noi {\bf Definition 1.2.18. } Let $(F_t, \; t \geq 0)$ denote an adapted, positive process. We shall say that 
this process belongs to the class $\widetilde{\mathcal{C}}$ if :

\noi {\it i)} $(F_t, \; t \geq 0)$ is a decreasing process, i.e., if $s \leq t$ : 
$$ 0 \leq F_t \leq F_s \quad W \; \textrm{a.s.} \eqno(1.2.134) $$
\noi In particular, since $0 \leq F_t \leq F_0$ and since $F_0$ is a.s. constant, this process is bounded by
a constant $C= F_0$. 

\noi {\it ii)} There exists $a \geq 0$ such that, for every $t \geq \sigma_a$, with 
$$ \sigma_a := \sup \{ t \geq 0; \; X_t \in [-a,a] \}$$
$$ F_t = F_{\sigma_a} = F_{\infty} \eqno(1.2.135)$$
and there exists $k > 0$ such that 
$$\underset{x \in [-a,a]}{\sup} \, {\bf W}_x (F_{\infty})  \leq k \eqno(1.2.136)$$
\noi {\it iii)} For every random time $T < \infty$ a.s. and every $u \geq 0$ :
$$ F_{T+u} (\omega) \leq F_u (\theta_T  \omega) \eqno(1.2.136').$$
\noi Of course, there is the inclusion $\widetilde{\mathcal{C}} \subset \mathcal{C}$. As the class $\mathcal{C}$, the 
class $\widetilde{\mathcal{C}}$ contains many interesting functionals $(F_t, \; t \geq 0)$. The following result holds :

\noi {\bf Theorem 1.2.19. } {\it Let $(F_t, \; t \geq 0)$ be a process in the class $\widetilde{\mathcal{C}}$
and $F_{\infty} := \underset{t \rightarrow \infty}{\lim} \, F_t$. Then, there exists a bounded process
 $(Y_t, \; t \geq 0)$ : 
$$0 \leq |Y_t| \leq c \eqno(1.2.137)$$
\noi such that : 
$$M_t (F_{\infty} ) = F_t |X_t| + Y_t \quad W \; \textrm{a.s.} \eqno(1.2.138)$$ }
\noi {\bf Examples}

\noi {\bf 1) } Let $(F_t:=h(L_t), \; t \geq 0)$ with $h : \mathbb{R}_+ \rightarrow \mathbb{R}_+$ Borel, such 
that $$ \int_0^{\infty} h(y) dy = 1.$$
\noi Then, (see (1.2.21)) :
$$M_t (h(L_{\infty})) = h(L_t) |X_t| + \int_{L_t}^{\infty} h(y) dy$$
\noi i.e. this is (1.2.138) with $$Y_t = \int_{L_t}^{\infty} h(y) dy.$$
\noi {\bf 2) } Let $\left(F_t := \exp \left( - \frac{1}{2} A_t^{(q)} \right), \; t \geq 0 \right)$ 
with $q \in \mathcal{I}$ and $q$ with compact support. Then (see (1.2.19)) : 
\begin{align*}
 M_t \left( e^{-\frac{1}{2} A_{\infty}^{(q)} } \right) & = \varphi_q (X_t) \, e^{-\frac{1}{2} A_{t}^{(q)}}
\\ & = e^{-\frac{1}{2} A_t^{(q)} } \, |X_t| + e^{-\frac{1}{2} A_t^{(q)} } \left( \varphi_q(X_t) - |X_t| \right)
\end{align*}
\noi i.e. (1.2.138) with 
$$ Y_t =  e^{-\frac{1}{2} A_t^{(q)} } \left( \varphi_q(X_t) - |X_t| \right)$$
\noi and we note that 
$$0 \leq |Y_t| \leq |\varphi_q(X_t) - |X_t| | \leq k$$
\noi since $\varphi_q$ is convex and $\varphi_q (x)$ is equivalent to $|x|$ as $|x|$ goes to infinity. 

\noi {\bf Proof of Theorem 1.2.19.}

\noi {\it i) } \underline{It is sufficient to prove} 
$$M_t (F_{\infty}) = F_t . (|X_t|-a)_+ + \widetilde{Y}_t \eqno(1.2.139)$$
\noi with $|\widetilde{Y}_t| \leq c'$. Indeed, if (1.2.139) is satisfied, then :
\begin{align*}
|M_t (F_{\infty} ) - F_t . |X_t|| & = | F_t (|X_t| - a)_+ + \widetilde{Y}_t - F_t. |X_t|| \\
& = |\widetilde{Y}_t + F_t ((|X_t| - a)_+ - |X_t| ) | \\ 
& \leq |\widetilde{Y}_t| + ak \leq c' + ak = c''.
\end{align*}
\noi {\it ii) } \underline{We now prove (1.2.139)}

\noi From point 2) of Theorem 1.2.1, we know that :
\begin{align*}
M_t (F_{\infty}) & = \widehat{\bf W}_{X_t} \big( F_{\infty} (\omega_t, \hat{\omega}^t) \big) \hspace*{7cm}(1.2.140)  \\
& = \widehat{\bf W}_{X_t} \big( F_{\infty} (\omega_t, \hat{\omega}^t) \, 1_{\sigma_a (\omega_t, 
\hat{\omega}^t ) < t } \big) \\
& + \widehat{\bf W}_{X_t} \big( F_{\infty} (\omega_t, \hat{\omega}^t) \, 1_{\sigma_a (\omega_t, 
\hat{\omega}^t ) > t } \big) \hspace*{5cm}(1.2.141) \\
 & := (1)_t + (2)_t 
\end{align*}

\noi \underline{Study of $(1)_t$}
\begin{align*}
(1)_t & = \widehat{\bf W}_{X_t} \big( F_{\infty} (\omega_t, \hat{\omega}^t) \, 1_{\sigma_a (\omega_t, 
\hat{\omega}^t ) < t } \big) \\
 & = \widehat{\bf W}_{X_t} \big( F_{\infty} (\omega_t) \, 1_{\sigma_a (\omega_t, 
\hat{\omega}^t ) < t } \big) 
\end{align*}
\noi since, on $\sigma_a < t$, $F_{\infty} = F_t$ (from (1.2.135)). Hence : 
$$ (1)_t = F_t (\omega_t) \widehat{\bf W}_{X_t} \big( 1_{\sigma_a (\omega_t, \hat{\omega}^t) < t} \big).$$
\noi But one can easily check that : 
$$\widehat{\bf W}_{X_t} \big( 1_{\sigma_a (\omega_t, \hat{\omega}^t) < t} \big) = (|X_t|-a)_+. \eqno(1.2.142)$$
\noi Indeed, we have :
$${\bf W}_x (\sigma_a = 0) = (|x|-a)_+. \eqno(1.2.142')$$
\noi Since, from (1.1.17') and relation (1.1.30) of Theorem 1.1.5, we have :
$$ {\bf W}_x \left( \exp \left( -\frac{1}{2} \, \lambda \, L_{\infty} \right) \right) = \frac{2}
{\lambda} + |x|.$$
\noi Letting $\lambda \rightarrow \infty$, we have $${\bf W}_x( T_0 = \infty) = |x|.$$
\noi Hence, we obtain (1.2.142') by translation. By (1.2.142), we deduce : 
$$(1)_t = \widehat{\bf W}_{X_t} \big( F_{\infty} (\omega_t, \hat{\omega}^t) \, 1_{\sigma_a (\omega_t, 
\hat{\omega}^t ) < t } \big) = F_t. ( |X_t| - a)_+). \eqno(1.2.143)$$
\noi \underline{Study of $(2)_t$}
$$ (2)_t = \widehat{\bf W}_{X_t} \big( F_{\infty} (\omega_t, \hat{\omega}^t) \, 1_{\sigma_a (\omega_t, 
\hat{\omega}^t ) > t } \big) $$
\begin{itemize}
\item If $|X_t| < a$, then $\sigma_a (\omega_t, 
\hat{\omega}^t ) > t$ and 
\begin{align*}
\widehat{\bf W}_{X_t} \left( F_{\infty} (\omega_t, \hat{\omega}^t) 1_{\sigma_a (\omega_t, \hat{\omega}^t) > t} \right)
& = \widehat{\bf W}_{X_t} \left( F_{\infty} (\omega_t, \hat{\omega}^t) \right) \\
& \leq \underset{x \in [-a,a]}{\sup} {\bf W}_x (F_{\infty}) \leq k \hspace*{4cm}(1.2.144) 
\end{align*}
\noi (from (1.2.136)). 
\item If $X_t \, (=x) \, \notin [-a,a]$ : 
$$ \widehat{\bf W}_{x} \left( F_{\infty} (\omega_t, \hat{\omega}^t) 1_{\sigma_a (\omega_t, \hat{\omega}^t) > t}
\right) 
= \widehat{\bf W}_{x} \left( F_{\infty} (\omega_t, \hat{\omega}^t) 1_{\widehat{T}_a < \infty} \right)$$
\noi where $\widehat{T}_a$ is the hitting time for $\hat{\omega}^t$ of $a$ or $-a$ (it does not depend on $\omega_t$).
\end{itemize}
\noi Hence 
$$ \widehat{\bf W}_{x} \left( F_{\infty} (\omega_t, \hat{\omega}^t) 1_{\sigma_a (\omega_t, \hat{\omega}^t) > t}
\right) \leq \widehat{\bf W}_{x} \left( F_{\infty} \left( \theta_{\widehat{T}_a} (\hat{\omega}^t) \right) 
1_{\widehat{T}_a < \infty} \right)$$
\noi since, from (1.2.136') : $$F_{\infty} (\omega_t, \hat{\omega}^t)  \leq F_{\infty} 
\big(\theta_{\widehat{T}_a} (\hat{\omega}^t) \big)$$
\noi on the event $\{(\sigma_a (\omega_t, \hat{\omega}^t)  > t ) \cap \widehat{T}_a (\hat{\omega}^t) < \infty\}$.
Hence,
\begin{align*}
& \widehat{\bf W}_{x} \left( F_{\infty} (\omega_t, \hat{\omega}^t) 1_{\sigma_a (\omega_t, \hat{\omega}^t) > t}
\right) \\
\leq & \, \varphi_{\delta_0} (x) \, W_{x, \infty}^{(\delta_0)} \left( e^{\frac{1}{2} L_{\infty} } 
F( \theta_{T_a} \omega) \, 1_{T_a < \infty} \right) \\
& \big(\textrm{from (1.1.57)}\big) \\
= & \, \varphi_{\delta_0} (x) \, W_{x, \infty}^{(\delta_0)} \left( e^{\frac{1}{2} L_{T_a} } 
 1_{T_a <\infty} W_{X_{T_a}, \infty}^{(\delta_0)} (e^{\frac{1}{2} L_{\infty}} F_{\infty}) \right) \\
& \textrm{(from the Markov property)} \\
= & \, \frac{\varphi_{\delta_0} (x) }{\varphi_{\delta_0} (a)} \, W_{x, \infty}^{(\delta_0)} \, 
\left( 1_{T_a < \infty}  \, {\bf W}_a (F_{\infty}) \right) \\
& \, \big(\textrm{from (1.1.57) and since} \; L_{T_a} = 0 \; W_{x, \infty}^{(\delta_0)} \; \textrm{a.s. for} 
\; |x| > a) \\
= & \, {\bf W}_a (F_{\infty}) \, \frac{\varphi_{\delta_0} (x)}{\varphi_{\delta_0} (a)}  \;
W_{x, \infty}^{(\delta_0)} (T_a < \infty).
\end{align*}
\noi But, $\varphi_{\delta_0} (x) = 2 + |x|$ and 
$$ W_{x, \infty}^{(\delta_0)} (T_a < \infty) \underset{|x| \rightarrow \infty}{\sim} \, \frac{2}{2 + |x|}$$
\noi (see (1.1.70)). 

\noi Hence : 
$$ \underset{x \in [-a,a]}{\sup} \, \varphi_{\delta_0}(x) \, W_{x, \infty}^{(\delta_0)} 
\left( F(\omega_t, \hat{\omega}^t) 1_{\sigma_a (\omega_t, \hat{\omega}^t) > t} \right) \leq c'' \eqno(1.2.145)$$
\noi Gathering (1.2.145), (1.2.144) and (1.2.143), we obtain Theorem 1.2.19.

\noi {\bf Corollary 1.2.20.}

\noi {\it Let $(F_t, \; t \geq 0)$ and $(G_t, \; t \geq 0)$ be two processes in $\widetilde{\mathcal{C}}$. Then : }

\noi {\bf 1) } 
$$ W \left( \frac{M_t (F_{\infty}) M_t (G_{\infty}) } {1 + |X_t|} \right) \underset{t \rightarrow \infty}
{\longrightarrow} \, {\bf W} (F_{\infty}. G_{\infty})  \hspace*{5cm}(1.2.146) $$
\noi {\bf 2) }
$$ \frac{1}{2} \sqrt{\frac{\pi}{2t}} \, W \big( M_t (F_{\infty}) M_t (G_{\infty}) \big)  \underset
{t \rightarrow \infty}{\longrightarrow} \, {\bf W} (F_{\infty}. G_{\infty}) \eqno(1.2.147)$$
\noi Note that, since $(F_t, \; t \geq 0)$ and $(G_t, \; t \geq 0)$ are in $\widetilde{\mathcal{C}}$, 
one has : 
$$ {\bf W} (F_{\infty}. G_{\infty} ) \leq k \, {\bf W} (G_{\infty}) < \infty.$$
\noi {\bf Proof of Corollary 1.2.20.} 

\noi {\bf 1) } \underline{We start with point 1) }

\noi We have : 
\begin{align*}
& W \left( \frac{M_t(F_{\infty}) M_t (G_{\infty})}{1+ |X_t|} \right) \\
= & {\bf W} \left( F_{\infty} \,  \frac{M_t (G_{\infty}) }{1+ |X_t|} \right) \\
& \textrm{(from Theorem 1.2.1)} \\
= & {\bf W} \left( F_{\infty} \, \frac{G_t |X_t| + Y_t^G}{1 + |X_t|} \right) \\
& \textrm{(from Theorem 1.2.19)} \\
 & \underset{t \rightarrow \infty}{\longrightarrow} \, {\bf W} (F_{\infty} G_{\infty} )
\end{align*}
\noi since $ \frac{G_t |X_t|}{1+ |X_t|} \, \leq \, G_t \, \leq \, k$, $F_{\infty} \in L^1 ({\bf W})$, 
$G_t$ decreases to $G_{\infty}$ when $t \rightarrow \infty$, and $|Y_t^G| < c$. 

\noi {\bf 2) } \underline{We now prove point 2) (briefly)}

\noi By polarization, it is sufficient to prove (1.2.147) when $F_{\infty} = G_{\infty}$. In this case, 
$t \rightarrow W( M_t^2 (F_{\infty}) )$ is an increasing function of $t$ and one can apply the Tauberian
Theorem. Let us compute :
$$ \int_0^{\infty} e^{- \lambda t} W (M_t^2 (F_{\infty})) dt 
= \int_0^{\infty} e^{- \lambda t} W \left[ (F_t. |X_t| + Y_t)^2 \right] \, dt.$$
\noi It is not difficult to see that in this expression, the terms $Y_t^2$ and $F_t|X_t|. Y_t$ are negligible, so 
we only need to deal with the term $F_t^2 |X_t|^2$. By doing as in the proof of the point 1) of 
Theorem 1.2.14, one has : 
\begin{align*}
F_t^2 |X_t|^2 & = F_t^2 \, \frac{|X_t|^2 (|X_t|-a)_+ }{1+ |X_t|} \\
 & + F_t^2 |X_t|^2 \, \frac{ 1+ |X_t| - (|X_t|-a)_+ }{( 1+ |X_t|)^2} \, (|X_t|-a)_+ \\
& + F_t^2 \, \frac{ |X_t|^2 (1 + |X_t| - (|X_t| - a)_+)}{(1 + |X_t|)^2} \\
& = (\tilde{1})_t + (\tilde{2})_t + (\tilde{3})_t \hspace*{5cm}(1.2.148) 
\end{align*}
\noi Now : 
\begin{align*}
(1)_t & :=  \int_{0}^{\infty} e^{- \lambda t} W ((\tilde{1})_t) dt \\
 & = \int_{0}^{\infty} e^{- \lambda t} W \left( F_t^2 \, \frac{|X_t|^2 (|X_t|-a)_+}{1 + |X_t|} \right) dt \\
 & = \int_{0}^{\infty} e^{- \lambda t} {\bf W} \left( \frac{F_t^2 |X_t|^2}{1 + |X_t|} \, 1_{\sigma_a < t} \right)
dt \\
& \textrm{(from Theorem 1.1.16)}  \\
& = \int_{0}^{\infty} e^{- \lambda t} {\bf W} \left( F_{\sigma_a}^2 \, \frac{|X_t|^2}{1 + |X_t|} 
\, 1_{\sigma_a < t} \right) dt \\
& \big(\textrm{from (1.2.135)}\big) \\
& = {\bf W} \left( F_{\sigma_a}^2 e^{- \lambda \sigma_a} \int_0^{\infty} e^{- \lambda u} du \, 
\frac{|X_{\sigma_a + u}|^2 } {1 + |X_{\sigma_a + u}|} \right) \\
& \textrm{(after the change of variables $t = \sigma_a + u$)} \\
 & = {\bf W} \left( F_{\sigma_a}^2 e^{- \lambda \sigma_a} \right) 
E_0^{(3)} \left( \int_0^{\infty} e^{- \lambda u} \, \frac{(a + R_u)^2}{1+ a + R_u} \, du \right) \\
& \textrm{(from Theorem 1.2.1)}
\end{align*} 
\noi Since, by scaling : 
$$E_0^{(3)} \left( \int_0^{\infty} e^{- \lambda u} \, \frac{(a + R_u)^2}{1+ a + R_u} \, du \right) 
\underset{\lambda \rightarrow 0}{\sim} \, \frac{\sqrt{2}}{\lambda^{3/2}},$$
\noi one has : 
$$(1)_t \underset{\lambda \rightarrow 0}{\sim} {\bf W} (F_{\infty}) \, \frac{\sqrt{2}}{\lambda^{3/2}}.$$
\noi It is now easy, by using the same arguments as in the proof of point 1) of Theorem 1.2.14, to see
that $(2)_t$ and $(3)_t$ are, when $\lambda$ tends to zero, negligible with respect to $(1)_t$. Finally, 
from Tauberian Theorem : 
$$\frac{1}{2} \sqrt{\frac{\pi}{2 t} } \, W \left( M_t (F_{\infty}) M_t (G_{\infty}) \right) 
\underset{t \rightarrow \infty}{\longrightarrow} \, {\bf W} (F_{\infty}. G_{\infty}). $$

\noi {\bf Remark 1.2.21.}

\noi {\bf 1) } By using the same arguments as in Corollary 1.2.20, one can see that if $F^{(1)},...,F^{(k)}$ are
$k$ processes in the class $\widetilde{\mathcal{C}}$, then :
$$ W \left( \frac{\prod_{i=1}^k \, F_{\infty}^{(i)} }{(1+ |X_t|)^{k-1}} \right) 
\underset{t \rightarrow \infty}{\longrightarrow} {\bf W} \left( \prod_{i=1}^k 
 F_{\infty}^{(i)} \right) \eqno(1.2.149)$$
and
$$ t^{(k-1)/2} \, W \left( \prod_{i=1}^k F_{\infty}^{(i)} \right) 
\underset{t \rightarrow \infty}{\longrightarrow} c_k \,{\bf W} \left( \prod_{i=1}^k 
 F_{\infty}^{(i)} \right) \eqno(1.2.150)$$
\noi where $c_k$ is a universal constant. 

\noi Note that, at first sight, (1.2.149) and (1.2.150) seem quite strange since one knows (from Theorem 1.2.1) that 
$M_t (F_{\infty}^{(i)}) \underset{t \rightarrow \infty}{\longrightarrow} 0$, $W$ a.s. for all 
$i = 1,...,k$.

\noi {\bf 2) } Let $(F_t, \; t \geq 0)$ and $(G_t, \; t \geq 0)$ be two processes in $\widetilde{\mathcal{C}}$. 
We penalise Wiener measure by the process $(F_t, \; t \geq 0)$ (see Theorem 1.2.15) and we denote by $W_{\infty}
^F$ the probability obtained with this penalisation. Now, let us penalise the probability $W_{\infty}^F$
by $(G_t, \; t \geq 0)$ : we obtain the probability $W_{\infty}^{F,G}$. On the other hand, 
if we penalise Wiener measure by the functional $(F_t. G_t, \; t \geq 0)$, we obtain the probability
$W_{\infty}^{F.G}$. It is not difficult to see, by using Theorem 1.2.19, that $W_{\infty}^{F,G} = 
W_{\infty}^{F.G}$. 

\noi {\bf 3) } Let $(F_t, \; t \geq 0)$ be an adapted, positive and increasing process, such that, for
some $\lambda_0 > 0$, $(e^{- \lambda_0 F_t}, \; t \geq 0 )$ is in $\mathcal{C}$, and such that for all $x$, 
${\bf W}_x (e^{- \lambda_0 F_{\infty}}) < \infty$. Then, for all $x \in \mathbb{R}$, there exists a 
positive and $\sigma$-finite measure ${\bm \nu}_x^{(F_{\infty})}$, carried on $\mathbb{R}_+$, and such that 
for all continuous functions $h$ with compact support : 
$$\sqrt{t} \, W_x [h(F_t)] \underset{t \rightarrow \infty}{\longrightarrow} \int_{\mathbb{R}_+}
h(y) {\bm \nu}_x^{(F_{\infty}) } (dy) \eqno(1.2.151) $$
\noi This Theorem is a generalization of a result in [RY, IX]. In [RY, IX], it is obtained when 
$(F_t, \; t \geq 0)$ is an additive functional. In fact, the measure ${\bm \nu}_x^{(F_{\infty})}$ 
is the image of ${\bf W}_x$ by $F_{\infty} : \Omega \rightarrow \mathbb{R}_+$. The proof of (1.2.151)
is essentially a consequence of Theorem 1.2.14. 

\bigskip

{\large \bf 1.3 Invariant measures related to ${\bf W}_{x}$ and \boldmath$\Lambda$\unboldmath$_{x}$.}

\noi We shall now show that the measure {\bf W}, and the measure \boldmath$\Lambda$\unboldmath $\,$ which we shall define very soon, are closely related to invariant measures of some Markov process taking values in certain functional spaces. 

\noi{\bf 1.3.1} \underline{ The process $(\mathcal{X}_{t}, \; t \ge 0)$.}

\smallskip

\noi As before, $\big(\Omega, (X_{t}, \mathcal{F}_{t})_{t \ge 0}, \mathcal{F}_{\infty}, W_{x} (x \in \mathbb{R})\big)$ denotes the canonical realisation of Brownian motion, starting at zero. Let $\mathcal{X}_{0} \in \Omega = \mathcal{C} (\mathbb{R}_{+} \to \mathbb{R})$. We define the process $(\mathcal{X}_{t}, \; t \ge 0)$ taking values on $\mathcal{C} (\mathbb{R}_{+} \to \mathbb{R})$, and issued from $\mathcal{X}_{0}$, by :
\begin{gather*}  
\mathcal{X}_{t}(u) :=
\left\{ \begin{array}{rl}
\mathcal{X}_{0} (u-t) \qquad \qquad {\rm if} \hfill \quad u \ge t  \\
\mathcal{X}_{0} (0) + X_{t-u} \qquad {\rm if} \quad u \le t  \\
\end{array}\right.  \hspace*{6,6cm}\hfill(1.3.1) \\
\end{gather*} 

\noi It is easy enough to see that this process is Markov \big(we denote by $(P_{t}, \; t \ge 0)$ the semigroup associated with this Markov Process $(\mathcal{X}_{t}, \; t \ge 0)$\big) and that the measure :
$$  {\widetilde{\bf W}} := \int_{\mathbb{R}} dx \; W_{x} \eqno(1.3.2) $$

\noi is an invariant measure for this process. However, this process admits other invariant measures. More precisely : 

\smallskip

\noi {\bf Theorem 1.3.1.} {\it Let $a,b \ge 0$, with $a+b > 0$, and :
$$ {\bf W}_{x}^{a,b} := a \, {\bf W}_{x}^{+} + b \, {\bf W}_{x}^{-} \eqno(1.3.3) $$

\noi Then :
$$ \widetilde{\bf W}^{a,b} := \int_{\mathbb{R}} dx \; {\bf W}_{x}^{a,b} \eqno(1.3.4) $$

\noi is an invariant measure for the process $(\mathcal{X}_{t}, \; t \ge 0)$. Recall that ${\bf W}_{x}^{+}$ and ${\bf W}_{x}^{-}$ are defined in (1.1.88) by : }
$$
{\bf W}_{x}^{+} = 1_{\Gamma_{+}} \cdot {\bf W}_{x}, \; {\bf W}_{x}^{-} = 1_{\Gamma_{-}} \cdot {\bf W}_{x}$$

\noi {\bf Proof of Theorem 1.3.1.}

\noi By symmetry, it suffices to prove that the measure $\widetilde{\bf W}^{+}$ defined by $\widetilde{\bf W}^{+} := \dis \int_{\mathbb{R}} dx \, {\bf W}_{x}^{+}$ is invariant. For every measurable and positive functional $F : \Omega \to \mathbb{R}_{+}$, we have :
\begin{gather*}
\int_{\mathbb{R}} dx \int_{\Omega} {\bf W}_{x}^{+} (d \mathcal{X}) P_{t} F (\mathcal{X}) \\
\quad = \int_{\mathbb{R}} dx  \int_{\Omega} {\bf W}_{x}^{+} (d \mathcal{X}) W \big(F (x+X_{t-u}, \; u \le t \;;\; \mathcal{X} (v-t), \; v \ge t)\big) \qquad \hbox{\big(from (1.3.1)\big)} \\
\quad =\int_{\mathbb{R}} dx \int_{\Omega} {\bf W}_{x}^{+} (d \mathcal{X}) W \big(F (x+X_{t} -X_{u}, \; u \le t \;;\; \mathcal{X} (v-t), \; v \ge t)\big)
\end{gather*}

\noi \big(since $(X_{t-u}, \; u \le t)$ has the same law under $W$ as $(X_{t} -X_{u}, \; u \le t)\big)$ 
\begin{equation*}
\quad = \int_{\mathbb{R}} dy \; W \left( \int_{\Omega} {\bf W}_{y-X_{t}}^{+} (d \mathcal{X}) F (y-X_{u}, \; u \le t \;;\; \mathcal{X} (v-t), v \ge t)\right)
\end{equation*}

\noi \big(from Fubini and after making the change of variable $x + X_{t} =y$\big)
\begin{equation*}
\quad = \int_{\mathbb{R}} dy \; W \left( \int_{\Omega} {\bf W}_{y-X_{t}} (d \mathcal{X}) F (y-X_{u}, \; u \le t \;;\; \mathcal{X} (v-t), \; v \ge t) \, 1_{\Gamma_{+}} (\mathcal{X}) \right)
\end{equation*}

\noi (from the definition of ${\bf W}^{+}$ and since $\mathcal{X} \in \Gamma_+$ if and only if :
$\underset{v \rightarrow \infty}{\lim} \mathcal{X}(v-t) = + \infty$)

\begin{equation*}
\quad = \int_{\mathbb{R}} dy \; W_y \left( \int_{\Omega} {\bf W}_{X_t} (d \mathcal{X}) F(X_u, \, u \leq t;
 \, \mathcal{X} (v-t), \, v \geq t ) \, 1_{\Gamma_+} (\mathcal{X}) \right)
\end{equation*}
\noi since $(X_{u}, \; u \geq 0)$ and $(-X_{u}, \; u \geq 0)$ have the same law under $W_0$. We now write :
\begin{equation*}
\int_{\Omega} {\bf W}_{X_t} (d \mathcal{X}) F(X_u, \, u \leq t; \, \mathcal{X}(v-t), \, v \geq t)
 \, 1_{\Gamma_+} (\mathcal{X}) 
\end{equation*}
\begin{equation*}
 \quad = \widehat{{\bf W}}_{X_t} \left( F 1_{\Gamma_+} (\omega_t, \hat{\omega}^t) \right)
\end{equation*}

\noi where $\omega_t \in \mathcal{C} ([0,t] \, \rightarrow \mathbb{R})$, $\hat{\omega}^t \in \mathcal{C} (
\mathbb{R}_+ \, \rightarrow \mathbb{R})$, and :
\begin{equation*}
\omega_t (u) = X_u \quad \textrm{for} \; u \leq t, 
\end{equation*}
\begin{equation*}
\hat{\omega}^t (v) = X_t + \mathcal{X} (v) \quad \textrm{for} \; v \geq 0
\end{equation*}

\noi (see point 1 of Remark 1.2.2 for such a notation). In the preceding relation, $\omega_t$ is frozen and
expectation is taken with respect to $\hat{\omega}^t$. Hence, from the "characteristic formula" (1.2.3) for the
 martingale $(M_t (F 1_{\Gamma_+}), \, t \geq 0)$, we have:
\begin{equation*}
\int_{\Omega} {\bf W}_{X_t} (d \mathcal{X}) F(X_u, \, u \leq t; \, \mathcal{X} (v-t), \, v \geq t)
1_{\Gamma_+} (\mathcal{X}) 
\end{equation*}
\begin{equation*}
 \quad = M_t (F 1_{\Gamma_+}) (\omega_t).
\end{equation*}
\noi Hence:
\begin{equation*}
\int_{\mathbb{R}} dx \int_{\Omega} {\bf W}_x^+ (d \mathcal{X}) P_t F(\mathcal{X})
\end{equation*}
\begin{equation*}
 \quad = \int_{\mathbb{R}} dy \, W_y \big(M_t (F 1_{\Gamma_+}) \big)
\end{equation*}
\begin{equation*}
\quad = \int_{\mathbb{R}} dy \, W_y \big(M_0 (F 1_{\Gamma_+}) \big)
\end{equation*}
\begin{equation*}
\quad = \int_{\mathbb{R}} dy \, {\bf W}_y (F 1_{\Gamma_+})
\end{equation*}
\noi (from (1.2.2) where we replace ${\bf W}$ ($= {\bf W}_0$) by ${\bf W}_y$). 
\begin{equation*}
\quad =  \int_{\mathbb{R}} dy \, {\bf W}^+_y (F)
\end{equation*}
\noi (from the definition of ${\bf W}^+_y$).
\begin{equation*}
 \quad = \widetilde{{\bf W}}^+ (F).
\end{equation*}





\smallskip

\noi This is Theorem 1.3.1.

\noi{\bf 1.3.2} \underline{The measure  \boldmath$\Lambda$\unboldmath$_{x}$.}

\noi Let $\widetilde{\Omega} = \mathcal{C} (\mathbb{R} \to \mathbb{R}_{+})$ and $\mathcal{L} : \Omega \to \widetilde{\Omega}$, the application "total local time" defined by :
$$  \mathcal{L} (X_{t}, \; t \ge 0) = (L_{\infty}^{y}, \; y \in \mathbb{R}). \eqno(1.3.5)   $$

\noi We denote by \boldmath$\Lambda$\unboldmath$_{x}$ the image of ${\bf W}_{x}$ by $\mathcal{L}$. It is possible to give a very simple description of \boldmath$\Lambda$\unboldmath$_{x}$ \big(see [RY, M]\big). Here is this description :

\noi $\cdot$ Let $u, \alpha, \beta \in\mathbb{R}_{+}$ and $x \in \mathbb{R}$. We denote by $Q_{x,u}^{\alpha, \beta}$ the law of the process $(Y_{v}, \; v \in\mathbb{R})$ defined as follows : 
\begin{gather*}
Y_{x} = u \\
(Y_{x+t}, \; t \ge 0) \; \hbox{is the square of an} \; \alpha\hbox{-dimensional Bessel process} \\
(Y_{x-t}, \; t \ge 0) \; \hbox{is the square of a} \; \beta\hbox{-dimensional Bessel process, independent from} \\
(Y_{x+t}, \; t \ge 0). 
\end{gather*}

\noi Then :
$$  \text{\boldmath$\Lambda$\unboldmath}_{x} = \frac{1}{2} \int_{0}^{\infty}�du \; (Q_{x,u}^{0,2} + Q_{x,u}^{2,0}) \eqno(1.3.6) $$

\noi {\bf Sketch of the proof of (1.3.6).}

\noi By translation, it suffices to prove (1.3.6) for $x=0$. Then, we use (1.1.40) :
\begin{equation*}
{\bf W} = \int_{0}^{\infty} dv \, (W_{0}^{\tau_{v}} \circ P_{0}^{(3, {\rm sym})})
\end{equation*}

\noi and the following facts :

\noi $\bullet$ From the second Ray-Knight Theorem \big(see [ReY], Chap. IX\big) for Brownian motion, the process $(L_{\tau_{l}}^{y}, \; y \ge 0)$ is a 0-dimensional squared Bessel process, starting from $l$.

\noi $\bullet$ For a 3-dimensional Bessel process, starting from 0, $(L_{\infty}^{y}, \; y \ge 0)$ is a 2-dimensional squared Bessel process, starting from 0. This constitutes the "third" Ray-Knight Theorem.

\noi $\bullet$ If $(Z_{t}^{i}, \; t \ge 0), \; i =1,2$, are two squared Bessel processes with respective dimensions $d_{1}$ and $d_{2}$, starting respectively from $u_{1}$ and $u_{2}$, then $(Z_{t}^{(1)} + Z_{t}^{(2)}, \; t \ge 0)$ is a squared Bessel process with dimension $d_{1}+d_{2}$ starting from $u_{1} + u_{2}$. 

\noi Other properties about the measure \boldmath$\Lambda$\unboldmath$_{x}$ may be found in \big([RY, M], Chap. 2\big). It is easily deduced from (1.3.6) that the r.v. $L_{\infty}^{y}$, under ${\bf W}_{x}$, admits the "law" :
$$  {\bf W}_{x} (L_{\infty}^{y} \in du) = |y-x| \, \delta_{0} (du) + du \qquad (u \ge 0) \eqno(1.3.7) $$

\noi \big(see also (1.1.45)\big).

\bigskip

\noi {\bf 1.3.3}  \underline{Invariant measures for the process $\big((X_{t}, L_{t}^{\bullet}), \; t \ge 0\big)$.}

\noi The process $\big((X_{t}, L_{t}^{\bullet}), \; t \ge 0\big)$, where $L_{t}^{\bullet} = (L_{t}^{y}, \; y \in \mathbb{R})$ denotes the local times process (in the space variable) at time $t$, for Brownian motion $(X_{t}, \; t \ge 0)$ is a Markov process taking values in $\mathbb{R} \times \widetilde{\Omega} = \mathbb{R} \times \mathcal{C} (\mathbb{R} \to \mathbb{R}_{+})$. In fact, if $\mathcal{X}_0$ is a function which has a finite total local time at each level, $\big((X_{t}, L_{t}^{\bullet} + L_{\infty}^{\bullet} (\mathcal{X}_0), \; t \ge 0\big)$ is the image of the process $(\mathcal{X}_{t}, \; t \ge 0)$ (see (1.3.1)) by the application :
$$  H : \Omega \to \mathbb{R} \times \widetilde{\Omega}  $$

\noi defined by :
$$  H(Y_{t}, \; t \ge 0) = (Y_{0}, L_{\infty}^{\bullet}) \eqno(1.3.8)  $$

\noi Of course, $H$ is only defined a.s. (with respect to the law of the process $(X_t, \; t \geq 0)$), 
i.e. it is only defined for the trajectories $\omega \in \Omega$ for which local time exists. 
 As a Corollary of Theorem 1.3.1, the image of $\widetilde{\bf W}^{a,b}$ by $H$ is an invariant measure for the process $\big((X_{t}, \; L_{t}^{\bullet}), \; t \ge 0\big)$. This image, which we denote by \boldmath$\widetilde{\Lambda}$\unboldmath$^{a,b}$ is equal, from (1.3.6), to :
$$
\text{ \boldmath$\widetilde{\Lambda}$\unboldmath}^{a,b} = \frac{1}{2} \int_{\mathbb{R}} dx \int_{0}^{\infty} du \; (a \, Q_{x,u}^{2,0} + b \, Q_{x,u}^{0,2}) \eqno(1.3.9) $$

\noi Thus, we have obtained :

\smallskip

\noi {\bf Theorem 1.3.2.} {\it The measure \boldmath$\widetilde{\Lambda}$\unboldmath$^{a,b}$ is an invariant measure for the process $\big((X_{t}, \; L_{t}^{\bullet}), \; t \ge 0\big)$.}

\smallskip

\noi We shall now give a different proof of Theorem 1.3.2 than the one we have just indicated. This proof has the further advantage that it hinges on arguments which shall be useful in the sequel. We begin with the :

\smallskip

\noi {\bf Lemma 1.3.3} {\it Let $q \in \mathcal{I}$

\noi {\bf 1)} Define $\varphi_{q}^{+} (x) : = {\bf W}_{x}^{+} \big(e^{-\frac{1}{2} \, A_{\infty}^{(q)}}\big) = {\bf W}_{x} \big(e^{- \frac{1}{2} \, A_{\infty}^{(q)}} 1_{\Gamma_{+}}\big)$. Then, $\varphi_{q}^{+}$ is the unique solution of Sturm-Liouville equation} :
\begin{gather*}
\hspace*{1.2cm}\varphi'' = q \, \varphi \qquad \hbox{\it with boundary conditions :} \nonumber \\
\hspace*{1.2cm}\varphi (x) \mathop{\sim}_{x \to + \infty}^{} x \qquad \varphi (x) \mathop{\longrightarrow}_{x \to - \infty}^{} C := \frac{1}{ \dis \int_{- \infty}^{\infty} \frac{dy}{\varphi_{q}^{2} (y)}} \hspace*{4,4cm}(1.3.10)
\end{gather*}

\noi {\it {\bf 2)} Define $\varphi_{q}^{-} (x) := {\bf W}_{x}^{-} \big(e^{-\frac{1}{2} \, A_{\infty}^{(q)}}\big) = {\bf W}_{x} \big(e^{- \frac{1}{2} \, A_{\infty}^{(q)}} 1_{\Gamma_{-}}\big)$. Then, $\varphi_{q}^{-}$ is the unique solution of the Sturm-Liouville equation : }
\begin{gather*}
\hspace*{1.2cm}\varphi'' = q \, \varphi \qquad \hbox{\it with boundary conditions :}  \\
\hspace*{1.2cm}\varphi (x) \mathop{\sim}_{x \to - \infty}^{} |x| \qquad \varphi (x) \mathop{\longrightarrow}_{x \to + \infty}^{} C := \frac{1}{ \dis \int_{- \infty}^{\infty} \frac{dy}{\varphi_{q}^{2} (y)}} \hspace*{4,4cm} (1.3.11)
\end{gather*}

\noi {\bf Proof of Lemma 1.3.3.} 

\noi It suffices, by symmetry, to prove point 1. We have 
\begin{eqnarray*}
{\bf W}_{x} \big(e^{- \frac{1}{2} \, A_{\infty}^{(q)}} 1_{\Gamma_{+}}\big)
&=& \varphi_{q} (x) \, W_{x, \infty}^{(q)} (\Gamma_{+}) \qquad \hbox{\big(from (1.1.16)\big)} \\
&=& \mathop{ \mathop{\rm lim}_{a \to - \infty}^{}}_{b \to + \infty}^{} \varphi_{q} (x) \, W_{x, \infty}^{(q)} (T_{b} < T_{a})
\end{eqnarray*}

\noi But, from (1.1.14), this limit equals :
$$
\varphi_{q} (x) \; \frac{ \gamma_{q} (x) - \gamma_{q} (- \infty)}{\varphi_{q} (\infty) - \gamma_{q} (- \infty)}
:= \varphi_{q} (x) \, \frac{\gamma_{q}(x) - \alpha}{\beta - \alpha} \eqno(1.3.12) $$

\noi where $\gamma_{q}$ is given by (1.1.14). Hence :
$$
\varphi_{q}^{+} (x) = \varphi_{q} (x) \; \frac{\gamma_{q} (x)- \alpha}{\beta - \alpha} \cdot \eqno(1.3.13) $$

\noi It remains to prove that $\varphi_{q}^{+}$ satisfies the announced conditions. But (with $\gamma$ for $\gamma_{q}$) : 
\hspace*{1.2cm}\begin{eqnarray*}
(\varphi_{q}^{+})'' (x)
&=& \varphi''_{q} (x) \left(\frac{\gamma (x)-\alpha}{\beta - \alpha}\right) + 2 \varphi'_{q} (x) \frac{\gamma'(x)} {\beta - \alpha} + \varphi_{q} (x) \frac{\gamma'' (x)}{\beta- \alpha}  \\
&=& \varphi''_{q} (x) \left(\frac{\gamma (x) - \alpha}{\beta - \alpha}\right) + \frac{2 \varphi'_{q} (x)}{\beta - \alpha} \; \frac{1}{\varphi^{2}_{q} (x)} + \frac{\varphi_{q} (x)}{\beta - \alpha} \left(-2 \frac{\varphi'_{q}(x)}{\varphi_{q}^{3} (x)}\right) \qquad \hbox{\big(from (1.1.14)\big)}  \\
&=& \varphi''_{q} (x) \left(\frac{\gamma(x)-\alpha}{\beta - \alpha}\right) = q(x) \varphi_{q} (x) \frac{\gamma (x)-\alpha}{\beta - \alpha} = q(x) \varphi_{q}^{+} (x) \hspace*{2,6cm} (1.3.14)
\end{eqnarray*}

\noi On the other hand : 
\begin{gather*}
\hspace*{1.2cm}\varphi_{q}^{+} (x) = \varphi_{q} (x) \frac{ \gamma (x) -\gamma (-\infty)}{\gamma (\infty) - \gamma (- \infty) } \mathop{\sim}_{x \to \infty}^{} \varphi_{q} (x) \mathop{\sim}_{x \to \infty}^{} x \hspace*{4cm} (1.3.15) \\
\hspace*{1.2cm}\varphi_{q}^{+} (x) = \varphi_{q} (x) \frac{ \dis \int_{- \infty}^{x} \frac{dy}{\varphi_{q}^{2} (y)} }{\dis \int_{- \infty}^{\infty}  \frac{dy}{\varphi_{q}^{2} (y)}} \mathop{\sim}_{x \to -\infty}^{} C \; \frac{\varphi_{q}(x)}{|x|} \mathop{\longrightarrow}_{x \to -\infty}^{} C = \frac{1}{ \dis \int_{- \infty}^{\infty}  \frac{dy}{\varphi_{q}^{2} (y)}}
\end{gather*}
\noi (since $\varphi_q(y)$ is equivalent to $|y|$ when $y$ goes to $- \infty$). 

\noi This proves Lemma 1.3.3.  $\hfill \blacksquare$

\bigskip

\noi {\bf We now prove Theorem 1.3.2.}

\noi Of course, by symmetry, it suffices to show that the measure : \boldmath$\widetilde{\Lambda}$\unboldmath$^{+} := \dis \int_{\mathbb{R}} dx$ \boldmath$ \Lambda$\unboldmath$_{x}^{+}$, where \boldmath$\Lambda$\unboldmath$_{x}^{+}$ is the image of ${\bf W}^{+}_{x}$ by $\mathcal{L}$, is invariant for the process $\big((X_{t}, L_{t}^{\bullet}), \; t \ge 0\big)$. We note that from (1.3.6), we have :
$$
\text{\boldmath$\widetilde{\Lambda}$\unboldmath}_{x}^{+} = \frac{1}{2} \int_{0}^{\infty} du \, Q_{x,u}^{2,0}  \eqno(1.3.16)$$

\noi We denote by $(Q_{t}, \; t \ge 0)$ the semi-group which is associated to the Markov process $\big((X_{t}, L_{t}^{\bullet}), \; t \ge 0\big)$, and we consider $F : \mathbb{R} \times \widetilde{\Omega} \to \mathbb{R}_{+}$ of the form :
\begin{eqnarray*} 
\hspace*{1.2cm}F(x,l) 
&=& f(x) \; \exp \left(- \frac{1}{2} < q, l>\right) \hspace*{6cm} (1.3.17) \\
&=& f(x) \, \exp \left(- \frac{1}{2} \int_{\mathbb{R}} l(y) q(y) dy \right) 
\end{eqnarray*}

\noi for $q \in \mathcal{I}$ and $f$ Borel, bounded. Then, for such an $F$, we obtain, by definition of the process $\big((X_{t}, L_{t}^{\bullet}), \; t \ge 0\big)$ :
$$
Q_{t} F (x,l) = W \left(f (x+ X_{t}) \, \exp \left\{ - \frac{1}{2} <q,l> - \frac{1}{2} \int_{0}^{t} q (x+ X_{s}) ds \right\} \right) \eqno(1.3.18) $$

\noi Now, from the monotone class theorem, Theorem 1.3.2 shall be obtained once we show that : 
$$
\int_{\mathbb{R}} dx \int_{\widetilde{\Omega}} \text{\boldmath$\Lambda$\unboldmath}_{x}^{+} (dl) \, Q_{t} F (x,l) = \int_{\mathbb{R}} dx \int_{\widetilde{\Omega}} \text{\boldmath$\Lambda$\unboldmath}_{x}^{+} (dl) F (x,l) \eqno(1.3.19) $$

\noi for every $t \ge 0$. But, from Lemma 1.3.3, we have :
\begin{eqnarray*}
\hspace*{1.2cm}{\bf W}_{x}^{+} \left( \exp - \frac{1}{2} \, A_{\infty}^{(q)}\right)
&=& {\bf W}_{x} \left( \exp \left( - \frac{1}{2} \, A_{\infty}^{(q)}\right) \cdot 1_{\Gamma_{+}} \right) \\
&=& \int_{\widetilde{\Omega}} \text{\boldmath$\Lambda$\unboldmath}_{x}^{+} (dl) \exp \left(-\frac{1}{2} <q,l>\right) = \varphi_{q}^{+} (x) \hspace*{2,3cm} (1.3.20)
\end{eqnarray*}

\noi since \boldmath$\Lambda$\unboldmath$_{x}^{+}$ is the image of ${\bf{W}}_{x}^{+}$ by $\mathcal{L}$.

\noi Thus, the left-hand side of (1.3.19) writes : 
\begin{eqnarray*}
LHS
&=& <Q_{t} F, 1>_{\widetilde{\bf \Lambda}^{+}} \\
&=& \int_{\mathbb{R}} dx \int_{\widetilde{\Omega}} {\bf \Lambda}_{x}^{+} (dl) W \big(f (x+X_{t}) \, e^{-\frac{1}{2} <q,l> - \frac{1}{2} \int_{0}^{t} q (x +X_{s})ds}\big) \\
&& \quad \big(\textrm{from (1.3.18)}\big) \\
&=& W \left(\int_{\mathbb{R}} dx \varphi_{q}^{+} (x)  f (x + X_{t}) \, \exp  \left(-\frac{1}{2} \int_{0}^{t} q (x+ X_{s}) ds \right)\right)
\end{eqnarray*}

\noi (from Fubini and (1.3.20))
$$
\hspace*{-2cm}= \int_{\mathbb{R}} f(y) dy \, W \left(\varphi_{q}^{+} (y-X_{t}) \exp \left(- \frac{1}{2} \int_{0}^{t} q (y-X_{t} + X_{s})ds\right) \right) \eqno(1.3.21) $$

\noi after making the change of variables $x+X_{t}=y$. On the other hand, the right-hand side of (1.3.19) equals :
\begin{eqnarray*}
\hspace*{1.2cm}RHS
&=& \int_{\mathbb{R}} dy \int_{\widetilde{\Omega}} \text{\boldmath$\Lambda$\unboldmath}_{y}^{+} (dl) f(y) \, \exp \left(- \frac{1}{2} <l,q>\right) \\
&=& \int_{\mathbb{R}} f(y) \;\; \varphi_{q}^{+} (y) dy \hspace*{8cm}(1.3.22)
\end{eqnarray*}

\noi from (1.3.20). Thus, Theorem 1.3.2 is an immediate consequence of the following :

\smallskip

\noi {\bf Lemma 1.3.4.} {\it For every $q \in \mathcal{I}$, $x$ real and $t \ge 0$: 
$$
W\big(\varphi_{q}^{+} (y-X_{t}) \, \exp \left(- \frac{1}{2} \int_{0}^{t} q(y-X_{t} + X_{s})ds \right) \big) = \varphi_{q}^{+} (y) \eqno(1.3.23)  $$

\noi Furthermore, (1.3.23) is also true when $\varphi_{q}^{+}$ is replaced by $\varphi_{q}^{-}$ or $\varphi_{q}$.  }

\smallskip

\noi {\bf Proof of Lemma 1.3.4.}
\begin{gather*}
W \left(\varphi_{q}^{+} (y-X_{t}) \, \exp \left(- \frac{1}{2} \int_{0}^{t} q (y-X_{t} + X_{s}) ds \right) \right) \\
\qquad = W \left(\varphi_{q}^{+} (y-X_{t}) \, \exp \left(- \frac{1}{2}  \int_{0}^{t} q (y-X_{t} + X_{t-r}) dr \right) \right) \\
\end{gather*}

\noi (after making the change of variables $s=t-r$).
\begin{equation*}
\qquad = W \left( \varphi_{q}^{+} (y-X_{t}) \,\exp  \left(- \frac{1}{2}  \int_{0}^{t} q (y-X_{r}) dr\right)\right)
\end{equation*}

\noi \big(since the process $(X_{t} - X_{t-r}, \; 0 \le r \le t)$ has the same law as $(X_{r}, \; 0 \le r \le t)$\big)
\begin{equation*}
\qquad = W_{y} \left( \varphi_{q}^{+} (X_{t}) \, \exp \left(- \frac{1}{2}  \int_{0}^{t} q (X_{r}) dr \right) \right)
\end{equation*}

\noi \big(since $(-X_{r}, \; r \ge 0)$ has the same law as $(X_{r}, \; r \ge 0)$\big)
\begin{equation*}
\qquad = \varphi_{q}^{+} (y)
\end{equation*}

\noi because, from (1.3.10) and It\^{o}'s formula, $\dis\left(\varphi_{q}^{+} (X_{t}) \, \exp \left(- \frac{1}{2}  \int_{0}^{t} q (X_{r}) dr\right), \; t \ge 0\right)$ is a 

\noi $\big( (\mathcal{F}_{t}, \; t \ge 0)$, $W_{y}\big)$ martingale.

\smallskip

\noi {\bf Remark 1.3.5.}

\noi {\bf 1)} We denote by $\mathcal{G}$ the infinitesimal generator of the process $\big((X_{t}, , L_{t}^{\bullet}), \; t \ge 0\big)$. For a function $F$ of the form given by (1.3.17), we obtain :
\begin{eqnarray*}
\hspace*{1,5cm}\lefteqn{\mathcal{G}F (x,l)} \\
&=& \frac{d}{ds} \Big|_{s = 0} \; Q_{s} F (x,l)  \\
&=& \frac{d}{ds} \Big|_{s = 0} \; W\left(f (x + X_{s}) \exp \left(- \frac{1}{2} <q,l> - \frac{1}{2} \int_{0}^{s} q (x+ X_{r}) dr \right)\right)  \\
\lefteqn{\hbox{\big(from (1.3.18)\big)} } \\
&=& \exp \left(- \frac{1}{2} <q,l> \right) \cdot \left[\frac{1}{2} \, f'' (x) -\frac{1}{2} \, q(x) f(x)\right] \hspace*{4,4cm}(1.3.24) \\
&=& \frac{1}{2} \; \frac{\partial^{2} F}{\partial x^{2}} \, (x,l) - \frac{1}{2} \, q(x) F(x,l) \hspace*{7,2cm}{(1.3.25)}
\end{eqnarray*}

\noi Another way to prove Theorem 1.3.2 consists in showing that, for every $F$ of the form (1.3.17), we have :
$$
<\mathcal{G}F, \, 1>_{\widetilde{\text{\boldmath$\Lambda$\unboldmath}}^{a,b}} = 0 \eqno(1.3.26) $$

\noi {\bf Let us prove (1.3.26).}

\noi By symmetry, it suffices to prove (1.3.26) by replacing \text{\boldmath$\widetilde{\Lambda}$\unboldmath}$^{a,b}$ by \text{\boldmath$\widetilde{\Lambda}$\unboldmath}$^{+}$. Now, we obtain, for $F$ of the form (1.3.17) with $f$ of class $C^{2}$, with compact support :
\begin{eqnarray*}
<\mathcal{G}F, \, 1>_{\text{\boldmath$\widetilde{\Lambda}$\unboldmath}^{+}}
&=& \int_{\mathbb{R}} dx \int_{\widetilde{\Omega}} \text{\boldmath$\widetilde{\Lambda}$\unboldmath}_{x}^{+} (dl) \, e^{- \frac{1}{2} \, <q,l>} \left(\frac{1}{2} f'' (x)- \frac{1}{2} \, q(x) f(x) \right) \nonumber \\
\lefteqn {\hbox{\big(from (1.3.24)\big)}} \\
&=& \int_{\mathbb{R}} \varphi_{q}^{+} (x) dx \left(\frac{1}{2} \, f'' (x) - \frac{1}{2} \, q(x) f(x) \right) \\
\lefteqn {\hbox{ (from Lemma 1.3.3)} }\\
&=& \int_{\mathbb{R}} \frac{1}{2} \, f(x) \left[(\varphi_{q}^{+})''(x) -q (x) \varphi_{q}^{+} (x) \right] \, dx \\
\lefteqn {\hbox{(after integrating by parts)}} \\
&=& 0 \qquad \hbox{(from Lemma 1.3.3.)}
\end{eqnarray*}

\noi {\bf 2)} Theorem 1.3.2 invites to ask the following question : is the process $\big((X_{t}, \, L_{t}^{\bullet}), \; t \ge 0\big)$ reversible with respect to the measure \boldmath$\widetilde{\Lambda}$\unboldmath$^{a,b}$, i.e. : does the following hold :
$$
<Q_{s} F, \, G>_{\text{\boldmath$\widetilde{\Lambda}$\unboldmath}^{a,b}} = < F, \, Q_{s} G >_{{\text{\boldmath$\widetilde{\Lambda}$\unboldmath}}^{a,b}} \eqno(1.3.27) $$

\noi for every $F, G : \mathbb{R} \times \Omega \to \mathbb{R}_{+}$ measurable and positive ? The answer to this question is negative. In particular, the operator $\mathcal{G}$ is not symmetric, i.e., in general :
$$
<\mathcal{G}F, \, G>_{\text{\boldmath$\widetilde{\Lambda}$\unboldmath}^{a,b} } \neq <F, \, \mathcal{G} G>_{\text{\boldmath$\widetilde{\Lambda}$\unboldmath}^{a,b}} \eqno(1.3.28) $$

\noi {\bf We now show (1.3.28)}, with $F(x,l)=f(x)\; \exp \dis \left(-\frac{1}{2} <q,l>\right)$, $G (x,l)=g(x)$, $\;$ \boldmath$\widetilde{\Lambda}$\unboldmath$^{a,b} = $\;$ $\boldmath$\widetilde{\Lambda}$\unboldmath $\;$ := \boldmath$\widetilde{\Lambda}$\unboldmath$^{1,1}$. Assuming that the equality would hold in (1.3.28), we would obtain, after an elementary computation :
\begin{eqnarray*}
<\mathcal{G}F, \, G>_{\widetilde{\Lambda}}
&=& \int_{\mathbb{R}} \varphi_{q}(x) g(x) \left(\frac{1}{2} \, f'' (x) -\frac{1}{2} \,q(x) f(x)\right) dx \\
&=& \int_{\mathbb{R}} \varphi_{q} (x) f(x) \frac{1}{2} \, g'' (x) dx= <F, \, \mathcal{G} G >_{\text{\boldmath$\widetilde{\Lambda}$\unboldmath}}
\end{eqnarray*}

\noi Thus, the preceding equality would imply, after integrating by parts and use of the relation $\varphi''_{q} = q \, \varphi_{q}$ :
\begin{equation*}
-2 q(x) \varphi_{q} (x) f(x) = 2 \varphi'_{q} (x) f'(x)
\end{equation*}

\noi for every $f$ in class $C^{2}$, with compact support, which is absurd.

\smallskip

\noi {\bf 3)} Of course, the preceding point implies that the measure $\widetilde{\bf W}^{a,b}$ which is invariant for the process $(\mathcal{X}_{t}, \; t \ge 0)$ is not reversible.

\smallskip

\noi {\bf 4)} The following relation, which has been obtained from Lemma 1.3.3 and the definition of \boldmath${\Lambda}$\unboldmath$_{x}^{\pm}$ :
$$
W_{x} \left[ \varphi_{q}^{\pm} (X_{t}) \, \exp \left(- \frac{1}{2} \, A_{t}^{(q)}\right)\right] = \int_{\widetilde{\Omega}} \text{\boldmath$\widetilde{\Lambda}$\unboldmath}_{x}^{\pm} (dl) \, \exp \left(- \frac{1}{2} <q,l>\right) \eqno(1.3.29) $$

\noi is a  particular case of the following result, which is found in \big([RY, M], Chap. 2\big) : 

\noi Let $F : \widetilde{\Omega} \to \mathbb{R}_{+}$ measurable, and "sub-exponential at infinity", \Big(i.e. : there exists $q \in \mathcal{I}$ and $C>0$ such that, for every $l \in \widetilde{\Omega}$, $F(l) \le C$ exp$\big(-<q,l>\big)\Big)$, then :
$$
\left(\int_{\widetilde{\Omega}} \text{\boldmath$\widetilde{\Lambda}$\unboldmath}_{X_{t}}^{\pm} (dl) \, F (l+ L_{t}^{\bullet}), \; t \ge 0\right) \eqno(1.3.30) $$ 

\noi is a $\big((\mathcal{F}_{t}, \; t \ge 0), \; W\big)$ martingale ; hence :
\begin{eqnarray*}
W_{x} \left(\int_{\widetilde{\Omega}} \text{\boldmath$\widetilde{\Lambda}$\unboldmath}_{X_{t}}^{\pm} (dl) F (l+ L_{t}^{\bullet})\right)
&=& W_{x} \left(\int_{\widetilde{\Omega}} \text{\boldmath$\widetilde{\Lambda}$\unboldmath}_{X_{0}}^{\pm} (dl)  F(l)\right)  \\
&=& \int_{\widetilde{\Omega}} \text{\boldmath$\widetilde{\Lambda}$\unboldmath}_{x}^{\pm} (dl)  F(l) \hspace*{6cm} (1.3.31)
\end{eqnarray*}

\noi If $F(l)= \exp \dis \left(-\frac{1}{2} < q,l>\right)$, we have :
\begin{eqnarray*}
 \int_{\widetilde{\Omega}} \text{\boldmath$\widetilde{\Lambda}$\unboldmath}_{X_{t}}^{\pm} (dl)  F(l+ L_{t}^{\bullet})
 &=& \int_{\widetilde{\Omega}} \text{\boldmath$\widetilde{\Lambda}$\unboldmath}_{X_{t}}^{\pm} (dl) \exp \left(-\frac{1}{2} < q,l> - \frac{1}{2} \int_{\mathbb{R}} q(x) L_{t}^{x} dx \right) \\
 &=& \varphi_{q}^{\pm} (X_{t}) \exp \left(- \frac{1}{2} \, A_{t}^{(q)}\right)
 \end{eqnarray*}

\noi Thus, when : $F(l) = \exp \dis \left(-\frac{1}{2} < q,l>\right)$, (1.3.31) is nothing else but (1.3.29) since :
$$ \dis W_{x} \left(\varphi_{q}^{\pm} (X_{t}) \exp \left(- \frac{1}{2} \, A_{t}^{(q)}\right)\right)
 = \varphi_{q}^{\pm} (x). $$

\noi {\bf 5)} Theorem 1.3.2 also invites to ask the question : are the measures \boldmath$(\widetilde{\Lambda}$\unboldmath$^{a,b}, \; a,b \ge 0)$ the only invariant measures of the process $\big((X_{t}, \; L_{t}^{\bullet}), \; t \ge 0\big)$. Here is a partial answer to this question. Let \boldmath$\widehat{\Lambda}$\unboldmath $\;$ be an invariant measure for this process. 

\smallskip

\noi {\it i)} Since the first component of $\big((X_{t}, \; L_{t}^{\bullet}), \; t \ge 0\big)$ is a Brownian motion, and that process admits as its only invariant measure (up to a multiplicative factor) the Lebesgue measure on $\mathbb{R}$, the measure $\text{\boldmath$\widehat{\Lambda}$\unboldmath}$ admits a disintegration of the form : 
$$ \text{\boldmath$\widehat{\Lambda}$\unboldmath} (dx, dl) = dx \,�\text{\boldmath$\widehat{\Lambda}$\unboldmath}_{x} (dl) \eqno(1.3.32) $$

\noi and, denoting by $\widehat{\varphi}_{q}$ the function defined by :
$$  \widehat{\varphi}_{q} (x) = \text{\boldmath$\widehat{\Lambda}$\unboldmath}_{x} \left(\exp - \frac{1}{2} <q,l>\right)  \eqno(1.3.33)$$

\noi the computations which lead to (1.3.21) and to (1.3.22) imply, if \boldmath$\widehat{\Lambda}$\unboldmath $\;$ is invariant :
\begin{equation*}
\widehat{\varphi}_{q} (x) = W_{x} \left(\widehat{\varphi}_{q} (X_{t}) \exp \left(- \frac{1}{2} \int_{0}^{t} q (X_{s}) ds \right) \right)
\end{equation*}

\noi It follows from this formula, using It\^{o}'s Lemma, that :
$$    \widehat{\varphi}''_{q} = q \, \widehat{\varphi}_{q} \eqno(1.3.34) $$

\noi The vector space of the solutions of the Sturm-Liouville equation has dimension 2 ; \\ hence, there exist two constants $C_{\pm} (q)$ such that : 
$$
\widehat{\varphi}_{q} (x) = C_{+} (q) \varphi_{q}^{+} (x) + C_{-} (q) \varphi_{q}^{-} (x) \eqno(1.3.35) $$

\noi {\it ii)} The invariant measure \text{\boldmath$\widetilde{\Lambda}$\unboldmath}$^{a,b}$ which we described in Theorem 1.3.2, and which writes :
\hspace*{3cm}\begin{eqnarray*}
\text{\boldmath$\widetilde{\Lambda}$\unboldmath}^{a,b}
&=& \frac{1}{2} \int_{\mathbb{R}} dx (a \text{\boldmath$\Lambda$\unboldmath}_{x}^{+} + b \text{\boldmath$\Lambda$\unboldmath}_{x}^{-})  
= \int_{\mathbb{R}} dx \text{\boldmath$\Lambda$\unboldmath}_{x}^{a,b} \hspace*{5cm}(1.3.36) \\
{\rm with} \; \text{\boldmath$\Lambda$\unboldmath}_{x}^{a,b} 
&:=& \frac{1}{2} (a \text{\boldmath$\Lambda$\unboldmath}_{x}^{+} + b \text{\boldmath$\Lambda$\unboldmath}_{x}^{-}) \hspace*{8,1cm}(1.3.37)
\end{eqnarray*}

\noi enjoys the following property : both limits
$$
\mathop{\rm lim}_{x \to + \infty}^{} \frac{1}{x} \, \text{\boldmath$\Lambda$\unboldmath}_{x}^{a,b} \left(\exp - \frac{1}{2} <q,l>\right) \quad {\rm and} \quad \mathop{\rm lim}_{x \to - \infty}^{} \frac{1}{|x|} \, \text{\boldmath$\Lambda$\unboldmath}_{x}^{a,b}\left(\exp - \frac{1}{2} <q,l>\right) \eqno(1.3.38) $$

\noi do not depend on $q \in \mathcal{I}$. Indeed,
\begin{equation*}
\frac{1}{x} \, \text{\boldmath$\Lambda$\unboldmath}_{x}^{a,b} \left(\exp - \frac{1}{2} <q,l>\right) = \frac{1}{2x} \, \big(a \varphi_{q}^{+} (x) + b \varphi_{q}^{-} (x)\big) \mathop{\longrightarrow}_{x \to \infty} \frac{a}{2}
\end{equation*}

\noi from Lemma 1.3.3 and $\dis \frac{1}{|x|} \, \text{\boldmath$\Lambda$\unboldmath}_{x}^{a,b} \left(\exp - \frac{1}{2} <q,l>\right) \mathop{\longrightarrow}_{x \to - \infty} \frac{b}{2}\cdot$

\smallskip

\noi {\it iii)} We now assume that the invariant measure \text{\boldmath$\widehat{\Lambda}$\unboldmath}, which equals : $\text{\boldmath$\widehat{\Lambda}$\unboldmath} (dx,dl)=dx \text{\boldmath$\widehat{\Lambda}$\unboldmath}_{x} (dl)$ also satisfies that both limits :
\begin{equation*}
\mathop{\rm lim}_{x \to \infty}^{} \frac{1}{x} \, \text{\boldmath$\widehat{\Lambda}$\unboldmath}_{x}\left(\exp - \frac{1}{2} <q,l>\right) \quad {\rm and} \quad \mathop{\rm lim}_{x \to -\infty}^{} \frac{1}{|x|} \, \text{\boldmath$\widehat{\Lambda}$\unboldmath}_{x}\left(\exp - \frac{1}{2} <q,l>\right)
\end{equation*}

\noi exist and do not depend on $q \in \mathcal{I}$. Then, there exist $a$ and $b$ positive, such that :
$\text{\boldmath$\widehat{\Lambda}$\unboldmath} = \text{\boldmath$\widehat{\Lambda}$\unboldmath}^{a,b}$. Indeed, from (1.3.35), together with Lemma 1.3.3 and (1.3.33), we have :
\begin{gather*}
\mathop{\rm lim}_{x \to \infty}^{} \frac{1}{x} \, \text{\boldmath$\widehat{\Lambda}$\unboldmath}_{x}\left(\exp \left(- \frac{1}{2} <q,l>\right) \right) = \mathop{\rm lim}_{x \to \infty}^{}  \frac{\widehat{\varphi}_{q} (x)}{x} \nonumber \\
\qquad = \mathop{\rm lim}_{x \to \infty}^{} \frac{ C_{+} (q) \varphi_{q}^{+} (x)+ C_{-} (q) \varphi_{q}^{-} (x)} {x} = C_{+} (q) \hspace*{5cm}(1.3.39)
\end{gather*}

\noi Thus, $C_{+} (q)$ \big(and $C_{-} (q)$, by symmetry\big) are constants, which we shall denote respectively as $\dis \frac{a}{2}$ and $\dis \frac{b}{2}$. Thus, we have :
\hspace*{1cm}\begin{eqnarray*}
\text{\boldmath$\widehat{\Lambda}$\unboldmath}_{x} \left(\exp - \frac{1}{2} <q,l>\right) 
&=& \frac{a}{2} \, \varphi_{q}^{+} (x) + \frac{b}{2} \, \varphi_{q}^{-} (x) \\
&=& \text{\boldmath$\Lambda$\unboldmath}_{x}^{a,b} (e^{- \frac{1}{2} <q,l>})
\end{eqnarray*}

\noi Hence : $\text{\boldmath$\widehat{\Lambda}$\unboldmath}_{x} = \text{\boldmath$\Lambda$\unboldmath}_{x}^{a,b}$ and $\text{\boldmath$\widehat{\Lambda}$\unboldmath} = \text{\boldmath$\widetilde{\Lambda}$\unboldmath}^{a, b}$.

\bigskip

\noi {\bf 1.3.4}  \underline{Invariant measures of the process $(L_{t}^{X_t- \bullet}, \; t \ge 0)$.}

\smallskip

\noi {\bf 1.3.4.1}  \underline{For $t \ge 0$, we define the random measure $\mu_{t}$ via} :
$$  \mu_{t} (f) = \int_{0}^{t} f (X_{t} - X_{s}) ds \eqno(1.3.40) $$

\noi with $f$ positive, continuous and bounded. It is proven in [DMY] that $(\mu_{t}, \; t \ge 0)$ is a Markov process taking values in the space of positive measures on $\big(\mathbb{R}, \, \mathcal{B} (\mathbb{R})\big)$. Due to the density of occupation formula, we may write (1.3.40) in the form :
\begin{eqnarray*}
\hspace*{5cm}\mu_{t} (f)
&=& \int_{\mathbb{R}} f(X_{t} -y) L_{t}^{y} dy \\
&=& \int_{\mathbb{R}} f (z) L_{t}^{X_{t}- z} \; dz \hspace*{4,5cm}(1.3.41)
\end{eqnarray*}

\noi We deduce that :
$$    \mu_{t} (dz) = L_{t}^{X_{t}-z} \; dz \eqno(1.3.42) $$

\noi Hence, rather than working in the space of measures on $\mathbb{R}$, we shall consider the Markov process $(L_{t}^{X_{t}- \bullet}, \; t \ge 0)$ which takes values in $ \widetilde{\Omega} = \mathcal{C} (\mathbb{R} \to \mathbb{R}_+)$.

\noi {\bf 1.3.4.2} Of course, this Markov process is the image of the process $\big((X_{t}, L_{t}^{\bullet}), \; t \ge 0\big)$ by the application :
$$  \theta : \mathbb{R} \times \widetilde{\Omega} \to \widetilde{\Omega} $$

\noi defined by :
$$   \theta (x,l) (y) = l(x-y) \qquad x,y \in \mathbb{R}, \; l \in \widetilde{\Omega} \eqno(1.3.43) $$

\noi This application $\theta$ is not bijective since : 
$$  \theta (x,l) = \theta (x',l') $$

\noi as soon as : 
$$    l(x-x'+z)=l'(z) \eqno(1.3.44)   $$

\noi for every $z \in \mathbb{R}$ i.e. : as soon as $l'$ is an adequate translate of $l$. 

\smallskip

\noi {\it i)} \underline{We begin by verifying directly}, i.e. : without using the result of Donati-Martin-Yor recalled above - that the process $(L_{t}^{X_{t}-\bullet}, \; t \ge 0)$, which takes values in $\mathcal{C} (\mathbb{R} \to \mathbb{R}_{+})$ is Markovian, in the natural filtration of the process $\big((X_{t}, L_{t}^{\bullet}), \; t \ge 0\big)$. For this purpose, using Dynkin's criterion (see [D]), and denoting by $(Q_{t}, \; t \ge 0)$ the semi-group associated to the process $\big((X_{t}, L_{t}^{\bullet}), \; t \ge 0\big)$, one needs to verify that :
$$   Q_{t} (F \circ \theta) (x,l) = Q_{t} (F \circ \theta) (x', l')  \eqno(1.3.45) $$

\noi for every $t \ge 0$ and $F : \widetilde{\Omega} \to \mathbb{R}_{+}$ measurable, as soon as :
$$   \theta (x,l) = \theta (x',l')  $$

\noi Of course, from the monotone class theorem, it suffices to prove (1.3.45) for $F$ of the form $F_{q}, \; q \in \mathcal{I}$, with :
$$  F_{q} (l) := \exp \left(-\frac{1}{2} <q,l>\right) \qquad (l \in \widetilde{\Omega}) \eqno(1.3.46) $$

\noi We have, from (1.3.43) :
\begin{eqnarray*}
\hspace*{1cm}F_{q} \circ \theta (x,l)
&=& \exp \left(- \frac{1}{2} <q,l (x-\cdot)>\right)  \\
&=& \exp \left(- \frac{1}{2} \int_{\mathbb{R}} q(y) l(x-y) dy\right) = \exp \left(- \frac{1}{2} < \mathop{q_{x}}_{}^{\vee}, l>\right) \hspace*{1,5cm} (1.3.47) \\
{\rm with}
&& \mathop{q_{x}}_{}^{\vee} (y) = q(x-y) \hspace*{7,5cm}(1.3.48)
\end{eqnarray*}

\noi Thus, from (1.3.18) :
$$
Q_{t} (F_{q} \circ \theta) (x,l) = W \left(\exp \left(-\frac{1}{2} < {\mathop{ q}^{\vee}}_{x+X_{t}},l > - \frac{1}{2} \int_{0}^{t} q(x+X_{t} -(x + X_{r})) dr \right)\right) \eqno(1.3.49) $$

\noi However :
\begin{eqnarray*}
\hspace*{1cm}<{\mathop{ q}^{\vee}}_{x+X_{t}}, l>
&=& \int_{\mathbb{R}} q (x +X_{t}-y) l(y) dy  \\
&=& \int_{\mathbb{R}} q (X_{t}+z) l(x-z) dz \hspace*{6cm}(1.3.50)
\end{eqnarray*}

\noi Thus, from (1.3.44), if $\theta (x,l) = \theta(x',l')$, we have :
\begin{equation*}
l(x-z)=l'(x'-z) \quad {\rm hence} \quad < {\mathop{ q}^{\vee}}_{x+X_{t}},l> = <{\mathop{ q}^{\vee}}_{x'+X_{t}},l'>
\end{equation*}

\noi It now follows from (1.3.49) that :
\begin{equation*}
Q_{t} (F_{q} \circ \theta) (x,l) = Q_{t} (F_{q} \circ \theta) (x',l')
\end{equation*}

\noi{\bf 1.3.4.3} �\underline{Invariant measures for the process $(L_{t}^{X_{t}-\bullet}, \; t \ge 0)$.}

\noi Of course, from Theorem 1.3.2, the image of \boldmath$\widetilde{\Lambda}$\unboldmath$^{a,b}$ by $\theta$ \big(defined by (1.3.43)\big) is an invariant measure for the process $(L_{t}^{X_{t}-\bullet}, \; t \ge 0)$. Unfortunately, an elementary computation shows that this measure is identically infinite. Thus, we need to find directly - without refering to \boldmath$\widetilde{\Lambda}$\unboldmath$^{a, b}$ - invariant measures for $(L_{t}^{X_{t}-\bullet}, \; t \ge 0)$.

\smallskip

\noi {\bf Theorem 1.3.6.} {\it Let $a, b \ge 0$, and :
$$ \text{\boldmath$ \Lambda$\unboldmath}^{a,b} := a \text{\boldmath$\Lambda$\unboldmath}_{0}^{+} + b \text{\boldmath$\Lambda$\unboldmath}_{0}^{-} \eqno(1.3.51) $$

\noi Then, \boldmath$\Lambda$\unboldmath$^{a,b}$ is an invariant measure for $(L_{t}^{X_{t}-\bullet}, \; t \ge 0)$.} 

\noi We recall that \boldmath$\Lambda$\unboldmath$_{0}^{\pm}$ is the image of ${\bf W}_{0}^{\pm} = {\bf W}^{\pm}$ by the application $\mathcal{L}$. In particular :
$$\text{\boldmath$\Lambda$\unboldmath}_{0}^{\pm} \big(\exp -\frac{1}{2} < q,l>\big) = \varphi_{q}^{\pm} (0) \qquad (q \in \mathcal{I}) \eqno(1.3.52) $$

\noi {\bf We now show Theorem 1.3.6.}

\noi We denote by $(\overline{Q}_{s}, \; s \ge 0)$ the semi-group associated to the Markov process $(L_{t}^{X_{t}-\bullet}, \; t \ge 0)$. Thus, we have, for (1.3.49) :
$$
\overline{Q}_{s} (F_{q}) (l) = W \left(\exp \left(-\frac{1}{2} < q (X_{s}+ \cdot), l> - \frac{1}{2} \int_{0}^{s} q (X_{s} - X_{r}) dr \right)\right) \eqno(1.3.53) $$

\noi with : $F_{q} (l) = \exp \dis \left(- \frac{1}{2} <q,l>\right)$.

\noi On the other hand, by symmetry, it suffices to show that the measure \boldmath$\Lambda$\unboldmath$^{+}$ := \boldmath$\Lambda$\unboldmath$_{0}^{+}$ is invariant for $(L_{t}^{X_{t}-\bullet}, \; t \ge 0)$. We compute :
\begin{eqnarray*}
\lefteqn{\int_{\widetilde{\Omega}} \text{\boldmath$\Lambda$\unboldmath}^{+} (dl) \big(\overline{Q}_{s} (F_{q})\big) (l) } \\
&=& \int_{\widetilde{\Omega}} \text{\boldmath$\Lambda$\unboldmath}^{+} (dl)W \left(\exp\left(- \frac{1}{2} <q (X_{s}+ \cdot), l > - \frac{1}{2} \int_{0}^{s} q (X_{s}-X_{r})dr \right)\right) \\
&=& W\left\{\left(\exp - \frac{1}{2} \int_{0}^{s} q (X_{s}-X_{r}) dr\right) \cdot \int_{\widetilde{\Omega}} \text{\boldmath$\Lambda$\unboldmath}^{+} (dl) \exp \left(-\frac{1}{2} <q(X_{s}+ \cdot),l>\right)\right\} \\
\lefteqn{\hbox{(from Fubini)}} \\
&=& W \left\{\exp \left(- \frac{1}{2} \int_{0}^{s} q (X_{s}-X_{r})dr\right) \cdot \varphi_{q(X_{s}+\cdot)}^{+} (0) \right\}
\end{eqnarray*}

\noi \big(from (1.3.52)\big). Now, it is easy to check that :
$$   \varphi_{q(X_{s}+\cdot)}^{+} (0) = \varphi_{q}^{+} (X_{s}) \eqno(1.3.54) $$

\noi Thus :
\begin{eqnarray*}
\int_{\widetilde{\Omega}} \text{\boldmath$\Lambda$\unboldmath}^{+} (dl) \big(\overline{Q}_{s} (F_{q})(l)\big)
&=& W \left(\varphi_{q}^{+} (X_{s}) \exp \left(- \frac{1}{2} \int_{0}^{s} q (X_{s}-X_{r}) dr \right)\right) \\
&=& \varphi_{q}^{+} (0)
\end{eqnarray*}

\noi from Lemma 1.3.4 (replacing $(X_{t}, \; t \ge 0)$ by $(-X_{t}, \; t \ge 0)$\big)
\begin{equation*}
\hspace*{3.8cm} = \int_{\widetilde{\Omega}} \text{\boldmath$\Lambda$\unboldmath}^{+} (dl) F_{q} (l) \qquad \big({\rm from} \; (1.3.52)\big)
\end{equation*}

\noi This is Theorem 1.3.6. $\hfill \blacksquare$

\noi {\bf Remark 1.3.7.}

\noi {\bf 1)} Arguing as in point 2 of Remark 1.3.5, it is easily shown that none of the measures $\text{\boldmath$\Lambda$\unboldmath}^{a,b}$ is reversible for the process $(L_{t}^{X_{t}-\bullet}, \; t \ge 0)$.

\smallskip

\noi {\bf 2)} Here is another way to prove that $\text{\boldmath$\Lambda$\unboldmath}^{a,b}$ is invariant. (We give the details for $\text{\boldmath$\Lambda$\unboldmath}^{+}$). We have, with $\dis F_{q} (l) = \exp - \frac{1}{2} <q,l>$, from (1.3.53) :
$$
\overline{Q}_{s} (F_{q}) (l) = W \left(\exp \left( - \frac{1}{2} <q (X_{s}+ \cdot),l>- \frac{1}{2} \int_{0}^{s} q (X_{r}) dr \right)\right) \eqno(1.3.55) $$

\noi We proceeded from (1.3.53) to (1.3.55) by making the change of variable $r=s-u$ and using the fact that, under $W$, $(X_{s}-X_{s-r}, \; r \le s)\dis \mathop{=}_{}^{\rm (law)} (X_{r}, \; r \le s)$. Thus, denoting by $\overline{\mathcal{G}}$ the infinitesimal generator of the semi-group $(\overline{Q}_{s}, \; s \ge 0)$, we obtain :
\begin{eqnarray*}
\overline{\mathcal{G}} \, F_{q} (l)
&=& \frac{d}{ds} \Big|_{s = 0} \; \overline{Q}_{s} (F_{q}) (l) \\
&=&  \frac{d}{ds} \Big|_{s = 0} \; W \left[g (X_{s}) \exp \left(-\frac{1}{2} \int_{0}^{s} q (X_{r}) dr \right)\right] \\
\biggl({\rm with} \; g(x)
&:=& \exp \left(- \frac{1}{2} < q(x+ \cdot), l > \right)\biggl) \\
&=& \frac{1}{2} \; g'' (0) - \frac{1}{2} \, q(0) \, g(0) \\
&=& \frac{1}{2} \left[ \frac{\partial^{2}}{\partial x^{2}} \Big|_{x=0} \left(\exp \left(-\frac{1}{2} <q (x+ \cdot), l>\right) -q (0) \exp \left(-\frac{1}{2} <q,l>\right) \right] \right.
\end{eqnarray*}

\noi Thus : 
\begin{eqnarray*}
<\overline{\mathcal{G}} \, F_{q}, 1>_{\Lambda^{+}} 
&= & \int_{\widetilde{\Omega}} \overline{\mathcal{G}} \, F_{q} (l) \Lambda^{+} (dl) \\
&=& \frac{1}{2} \int_{\widetilde{\Omega}} \text{\boldmath$\Lambda$\unboldmath}^{+} (dl) \left(\frac{\partial^{2}}{\partial x^{2}} \Big|_{x=0} \exp \left(- \frac{1}{2} <q(x+\cdot), l> \right) 
 -q(0) \exp \left(-\frac{1}{2} <q,l>\right)\right)  \\
&=& \frac{1}{2} \big(\varphi^{+}_{q}{''} (0) - q(0) \varphi_{q}^{+} (0) \big) =0 \hspace*{7cm}(1.3.56)
\end{eqnarray*}

\noi after interverting the second derivative and integration with respect to $\text{\boldmath$\Lambda$\unboldmath}^{+} (dl)$, using Lemma 1.3.3 and the fact that $\varphi_{q(x+\cdot)}^{+} (0) = \varphi_{q}^{+} (x)$. From relation (1.3.56), we deduce of course that~: $<\overline{Q}_{s} \, F_{q}, 1>_{\text{\boldmath$\Lambda$\unboldmath}^{+}} = <F_{q}, 1>_{\text{\boldmath$\Lambda$\unboldmath}^{+}}$, i.e. that $\text{\boldmath$\Lambda$\unboldmath}^{+}$ is invariant.


\newpage

\noi {\bf{\Large Chapter 2. Existence and properties of the measure ${\bf W}^{(2)}$.}}

\bigskip

\noi We shall now establish a number of results similar to those of Chapter 1, but this time $(X_{t}, \; t \ge 0)$ is a 2-dimensional Brownian motion.

\smallskip

{\bf 2.1. Existence of ${\bf W}^{(2)}$.}

\smallskip

\noi {\bf 2.1.1} \underline{Notations and Feynman-Kac penalisations in two dimensions.}

\noi $\big(\Omega = \mathcal{C} (\mathbb{R}_{+} \to \mathbb{C}), (X_{t}, \; \mathcal{F}_{t})_{t \ge 0}, \; W_{z}^{(2)} (z \in \mathbb{C}) \big)$ denotes the two dimensional canonical Brownian motion, which takes its values in $\mathbb{C}$. We write $W^{(2)}$ for $W_{0}^{(2)}$. $\mathcal{I}$ denotes here the set of positive Radon measures on $\mathbb{C}$ admitting a density $q$ with compact support and such that $\dis \int q(x) dx >0$. Define : 
$$   A_{t}^{(q)} := \int_{0}^{t} q (X_{s}) ds \eqno(2.1.1)  $$

\noi Here is the analogue in dimension 2 of Theorem 1.1.1. A proof of this Theorem (in dimension 2) is found in [RVY, VI].

\smallskip

\noi {\bf Theorem 2.1.1.} {\it Let $q \in \mathcal{I}$ and, for every $t \ge 0$ and $z \in \mathbb{C}$ :
$$
W_{z,t}^{(2, q)} := \frac{ \exp \big(- \frac{1}{2} \, A_{t}^{(q)}\big) }{ Z_{z,t}^{(2,q)} } \cdot W_{z}^{(2)} \eqno(2.1.2) $$

\noi with
$$   Z_{z,t}^{(2,q)} := W_{z}^{(2)} \left( \exp -\frac{1}{2} \, A_{t}^{(q)} \right) \eqno(2.1.3) $$

\noi {\bf 1)} \underline{For every $s \ge 0$ and $\Gamma_{s} \in b(\mathcal{F}_{s})$} :

\noi $W_{z,t}^{(2,q)} (\Gamma_{s})$ admits a limit $W_{z, \infty}^{(2,q)} (\Gamma_{s})$ as $t \to \infty$ :
$$
W_{z,t}^{(2,q)} (\Gamma_{s}) \mathop{\longrightarrow}_{t \to \infty}^{} W_{z, \infty}^{(2,q)} (\Gamma_{s})
\eqno(2.1.4) $$

\noi {\bf 2)} \underline{$W_{z, \infty}^{(2,q)}$ is a probability on $(\Omega, \mathcal{F}_{\infty})$ such that} :
$$  W_{z, \infty}^{(2,q)} |_{\mathcal{F}_{s}} = M_{s}^{(2,q)} \cdot W_{z}^{(2)} |_{\mathcal{F}_{s}}  $$

\noi where $(M_{s}^{(2,q)}, \; s \ge 0)$ is the $\big((\mathcal{F}_{s}, \; s \ge 0), \; W_{z}^{(2)}\big)$ martingale defined by :
$$ 
M_{s}^{(2,q)} = \frac{ \varphi_{q} (X_{s}) }{ \varphi_{q} (z) } \; \exp \left(- \frac{1}{2} \, A_{s}^{(q)}\right)
\eqno(2.1.5) $$

\noi {\bf 3)} \underline{The function $\;\varphi_{q} : \mathbb{C} \to \mathbb{R}_{+}\;$ featured in (2.1.5) is strictly positive, continuous and}

 \underline{satisfies} :
$$   \varphi_{q} (z) \mathop{\sim}_{|z| \to \infty}^{} \frac{1}{\pi} \; \log \big(|z|\big) \eqno(2.1.6)  $$

\noi It may be defined via one or the other of the following descriptions :  

\noi {\it i)} $\varphi_{q}$ is the unique solution of the Sturm-Liouville equation :
$$    \Delta \varphi = q \cdot \varphi \quad \hbox{(in the sense of Schwartz distributions)}   $$

\noi which satisfies the limiting condition : 
$$
|z| \frac{\partial \varphi}{\partial r} \; (z) \mathop{\longrightarrow}_{r \to \infty}^{} \frac{1}{\pi} \qquad \big(r = |z|\big) \eqno(2.1.7) $$

\noi {\it ii)} $\dis \qquad \frac{1}{2 \pi} (\log \; t) \; W_{z}^{(2)} \left( \exp \left(- \frac{1}{2} \, A_{t}^{(q)} \right) \right) \mathop{\longrightarrow}_{t \to \infty}^{} \varphi_{q} (z) \hfill(2.1.8) $

\noi {\bf 4)} Under the family of probabilities $(W_{z, \infty}^{(2,q)}, z \in \mathbb{C})$, the canonical process $(X_{t}, \; t \ge 0)$ is a transient diffusion. More precisely, there exists a $\big(\Omega, (\mathcal{F}_{t}, \; t \ge 0), \; W_{z,\infty}^{(2,q)}\big)$ Brownian motion $(B_{t}, \; t \ge 0)$ valued in $\mathbb{C}$ and starting from 0 such that : 
$$
X_{t} = z + B_{t} + \int_{0}^{t} \; \frac{ \nabla \varphi_{q}}{\varphi_{q}} (X_{s}) ds \eqno(2.1.9) $$  }

\noi {\bf 2.1.2} { \underline{ Existence of the measure ${\bf W}^{(2)}$.}}

\noi {\bf Theorem 2.1.2.} {\it There exists on $\big(\Omega = \mathcal{C} (\mathbb{R}_{+} \to \mathbb{C}), \; \mathcal{F}_{\infty} \big)$ a $\sigma$-finite and positive measure ${\bf W}^{(2)}$ (with infinite total mass) such that, for every $q \in�\mathcal{I}$ :
$$ 
{\bf W}^{(2)} = \varphi_{q} (0) \; \exp \left(+ \frac{1}{2} \, A_{\infty}^{(q)} \right) \cdot W^{(2,q)}_{\infty} \eqno(2.1.10) $$

\noi In other terms, the RHS of (2.1.10) does not depend on $q \in \mathcal{I}$. }

\noi In fact, just as in the case of dimension 1, we show for every $z \in \mathbb{C}$, the existence of a measure ${\bf W}_{z}^{(2)}$, this measure being defined by : 
$$
{\bf W}_{z}^{(2)} \big(F (X_{s}, \; s \ge 0)\big) = {\bf W}^{(2)} \big(F (z+X_{s}, \; s \ge 0)\big) \eqno(2.1.11) $$

\noi {\bf Proof of Theorem 2.1.2.}

\noi It consists in showing that $\varphi_{q} (0)$ exp$\dis \left(+ \frac{1}{2} \, A_{\infty}^{(q)} \right) \cdot W_{\infty}^{(2,q)}$ does not depend on $q$. The proof is quite similar to that of Theorem 1.1.2. It hinges upon : 
\begin{itemize}
\item[$\bullet$] $\varphi_{q} (z) > 0$ for every $q \in \mathcal{I}$ and $z \in \mathbb{C}$ ;
\item[$\bullet$] $\dis \frac{\varphi_{q_{1}} (z)}{\varphi_{q_{2}} (z)} \mathop{\longrightarrow}_{|z| \to \infty}^{} 1$ for every $q_{1}$ and $q_{2} \in \mathcal{I}$ ;
\item[$\bullet$] $\dis \varphi_{q} (z) \mathop{\longrightarrow}_{|z| \to \infty}^{} + \infty$ and the $(W_{z, \infty}^{(2,q)}, \; z \in \mathbb{C})$ process $(X_{t}, \; t \ge 0)$ is transient.
\end{itemize}

\noi 

\noi These properties follow from Theorem 2.1.1. We also note, just as we did in Lemma 1.1.3 : 
\begin{eqnarray*}
W_{z, \infty}^{(2,q)} \left(\exp + \frac{\lambda}{2} \, A_{\infty}^{(q)} \right) &<& \infty \quad {\rm if} \;\; \lambda < 1 \hspace*{6cm}(2.1.12) \\
W_{z, \infty}^{(2,q)} \left(\exp + \frac{\lambda}{2} \, A_{\infty}^{(q)} \right) &=& \infty \quad {\rm if} \;\; \lambda \ge 1\hspace*{6cm} (2.1.13)
\end{eqnarray*}

\noi These two properties show that ${\bf W}^{(2)}$ is well defined via (2.1.10) \big(since $A_{\infty}^{(q)} < \infty \;\; W_{\infty}^{(2,q)}$ a.s.\big) and that ${\bf W}^{(2)}$ has infinite total mass ; it is $\sigma$-finite on $(\Omega, \mathcal{F}_{\infty})$ and it is such that ${\bf W}^{(2)} (\Gamma_{t}) = 0$ or $+ \infty$ for any $\Gamma_{t} \in b^{+} (\mathcal{F}_{t})$ depending whether $W^{(2)} (\Gamma_{t})$ is equal to 0 or is strictly positive.

\bigskip

{\bf 2.2 Properties of ${\bf W}^{(2)}$.}

\smallskip

\noi�{\bf 2.2.1} \underline{Some notation.}

\noi We shall now prepare for Theorem 2.2.1 - which plays for ${\bf W}^{(2)}$ a similar role as Theorem 1.1.5 for ${\bf W}$. However, in order to prepare for Theorem 2.2.1, we need the following notation :

\noi {\it i)} \underline{Denote by $C$ the unit circle in $\mathbb{C}$} : 
$$    C = \{ z \in \mathbb{C} \;;\; |z| =1\}    \eqno(2.2.1) $$

\noi and $(L_{t}^{(C)}, \; t \ge 0)$ the (continuous) local time process on $C$, which may be defined as :
$$
L_{t}^{(C)} := \mathop{\rm lim}_{\varepsilon \downarrow 0}^{} \; \frac{1}{2 \pi \varepsilon} \int_{0}^{t} 1_{C_{\varepsilon}} (X_{s}) ds  \eqno(2.2.2) $$

\noi where
$$
C_{\varepsilon} = \{ z \in \mathbb{C} \;;\; 1- \varepsilon \le |z| \le 1 + \varepsilon\}  $$

\noi so that, a.s. if $\;\;q_{0}\;\;$ denotes the uniform probability on $C$ : 
$$
\int_{C} f(z) \, q_{0} (dz) = \frac{1}{2 \pi} \; \int_{0}^{2 \pi} f (e^{i \theta}) d \theta \eqno(2.2.3) $$

\noi we have :
$$     (L_{t}^{(C)}, \; t \ge 0) = (A_{t}^{(q_{0})}, \; t \ge 0)  \eqno(2.2.4) $$

\noi In other terms, $(L_{t}^{(C)}, \; t \ge 0)$ is the additive functional which admits $q_0$ as 
Revuz's measure (see [Rev]).
\noi We denote by $(\tau_{l}^{(C)}, \; l \ge 0)$ the right continuous inverse of $(L_{t}^{(C)}, \; t \ge 0)$ :
$$    \tau_{l}^{(C)} := \inf \{t \ge 0 \;;\; L_{t}^{(C)} > l \} ,\qquad l \ge 0 \eqno(2.2.5) $$

\noi and we denote by $W_0^{(2, \tau_{l}^{(C)})}$ the law of a 2-dimensional Brownian motion starting from 0, considered up to $\tau_{l}^{(C)}$. 

\smallskip

\noi {\it ii)} \underline{We denote by $\; P_{1}^{(2, \log)}\;$ the law of the process $\;(R_{t}, \; t \ge 0)\;$ which solves the stochastic} 

\underline{differential equation} : 
$$
R_{t} = 1 + \beta_{t} + \int_{0}^{t} \, \frac{ds}{R_{s}} \left(\frac{1}{2} + \frac{1}{\log \, R_{s}}\right) \eqno(2.2.6) $$

\noi where $(\beta_{t}, \; t \ge 0)$ is a 1-dimensional Brownian motion starting from 0. We note that the process $(R_{t}, \; t \ge 0)$ starts from 1 and that $P ( R_{t} >1$ for every $t > 0) =1$.

\noi We adopted the notation $P_{1}^{(2, \log)}$ to indicate : 

\noi {\bf a)} that this process $R$ starts from 1 ;

\noi {\bf b)} that it "differs at infinity from a 2-dimensional Bessel process" by the presence of the term $\dis \frac{1}{\log \, R_{s}}$, in the drift part of equation (2.2.6). 

\smallskip

\noi {\it iii)} \underline{Here is another description of the process} $(R_{t}, \; t \ge 0)$ defined by (2.2.6) :
$$      (\log \, R_{t}, \; t \ge 0) \mathop{=}_{}^{(\rm law)} (\rho_{H_{t}}, \; t \ge 0) \eqno(2.2.7) $$

\noi with :

$\bullet$ $(\rho_{u}, \; u \ge 0)$ a 3-dimensional Bessel process starting from 0 ;
 
$ \bullet$ $ \dis H_{t} := \int_{0}^{t} \frac{ds}{R_{s}^{2}} \hfill (2.2.8)$

\smallskip

\noi \underline{We prove (2.2.7)}.

\smallskip

\noi We apply It\^{o}'s formula to the process $(R_{t})$ solution of (2.2.6) and we obtain : 
$$
\log \, R_{t} = \int_{0}^{t} \frac{d \beta_{s}}{R_{s}} + \int_{0}^{t} \frac{ds}{R_{s}^{2} \cdot \log \, R_{s}}
\eqno(2.2.9) $$

\noi We denote by $(\nu_{h}, \; h \ge 0)$ the inverse of the process $(H_{t}, \; t \ge 0)$ and we replace $t$ by $\nu_{h}$ in (2.2.9). Thus : 
\begin{eqnarray*}
\hspace*{1,7cm}\log \, R_{\nu_{h}} 
&=& \int_{0}^{\nu_{h}} \frac{d \beta_{s}}{R_{s}} + \int_{0}^{\nu_{h}} \frac{ds}{R_{s}^{2} \, \log \, R_{s}} \hspace*{5,8cm}(2.2.10) \\
&=& \widetilde{\beta}_{h} + \int_{0}^{h} \frac{du}{\log \, R_{\nu_{u}}} \hspace*{7,3cm}(2.2.11)
\end{eqnarray*}

\noi after the change of variable $s=\nu_{u}$ and with $(\widetilde{\beta}_{h}, \; h \ge 0) \dis := \left(\int_{0}^{\nu_{h}} \frac{d\beta_{s}}{R_{s}^{2}}, \; h \ge 0\right)$, which is a 1-dimensional Brownian motion since this - local - martingale admits as bracket $\dis \Big(\int_{0}^{\nu_{h}} \frac{ds}{R_{s}^{2}} =$

\noi $H_{\nu_{h}} = h,$  $h \ge 0\Big).$ Hence, from (2.2.11) (log $R_{\nu_{h}}, \; h \ge 0)$ is a 3-dimensional Bessel process starting from 0. 

\smallskip

\noi {\it iv)} \underline{Let now $\; (\alpha_{_{t}}, \; t \ge 0)\; $ be another 1-dimensional Brownian motion}, independent from 

\noi $(\beta_{t}, \; t \ge 0)$ \big(hence independent from $(R_{t}, \; t \ge 0)\big)$. We define the law $W^{(2, \tau_{l}^{(C)})} \circ \widetilde{P}_{1}^{(2, \log)}$ as the law of the 2-dimensional process $(Y_{t}, \; t \ge 0)$ satisfying to : 

a) $(Y_{t}, \; t \le \tau_{l}^{(C)})$ is a 2-dimensional Brownian motion starting from 0 and stopped in $\tau_{l}^{(C)}$~; its law, from point {\it i)}, is $W^{(2, \tau_{l}^{(C)})}$. Here, $\tau_{l}^{(C)}$ is the right-continuous inverse of $(L_t^{(C)}, t \geq 0)$, the local time on $C$ of the process $(Y_t, t \geq 0)$. 

b) after $\tau_{l}^{(C)}$, the process $(Y_{\tau_{l}^{(C)} + t} \; t \ge 0)$ writes :
$$    Y_{\tau_{l}^{(C)} + t} := R_{ t}  \cdot e^{ i \alpha_{H_{t}} } \qquad (t \ge 0) \eqno(2.2.12) $$

\noi where :
\begin{itemize}
\item[$\bullet$] the law of the process $(R_{t}, \; t \ge 0)$ is $P_{1}^{(2, \, \log)}$ 
\item[$\bullet$] $(\alpha_{t}, \; t \ge 0)$ is a 1-dimensional Brownian motion starting from $\alpha_{0}$, with $e^{i \alpha_{0}} = Y_{\tau_{l}^{(C)}}$   (we note that $Y_{\tau_{l}^{(C)}} \in C)$ 
\item[$\bullet$] $H_{t} = \dis \int_{0}^{t} \frac{ds}{R_{s}^{2}}$
\end{itemize}

c) $(\alpha_{t}, \; t \ge 0)$ and $(\beta_{t}, \; t \ge 0)$, the driving Brownian motion of $(R_{t}, \; t \ge 0)$ (see (2.2.6)) are, conditionally on $\alpha_{0}$, independent from the process $(Y_{t}, \; t \le \tau_{l}^{(C)})$. 

\smallskip

\noi Formula (2.2.7) - the second description of $(R_{t}, \; t \ge 0)$ - permits to write (2.2.12) in another form : 
$$
Y_{\tau_{l}^{(C)} + t} = \exp (\rho_{u} + i\alpha_{u}) |_{u=H_{t}} \qquad (t \ge 0)  \eqno(2.2.13) $$

\noi where $(\rho_{u}, \; u \ge 0)$ is a 3-dimensional Bessel process starting from 0 and $H_{t} \dis = \int_{0}^{t} \frac{ds}{R_{s}^{2}} \cdot$

\noi {\bf 2.2.2} \underline{Description of the canonical process $(X_{t}, \; t \ge 0)$ under $W_{\infty}^{(2, q_{0})}$.}

\noi In order to describe the measure ${\bf W}^{(2)}$, we shall use the formula :
$$
{\bf W}^{(2)} = \varphi_{q_{0}} (0) (e^{\frac{1}{2} \, L_{\infty}^{(C)}}) \cdot W_{\infty}^{(2, q_{0})} \eqno(2.2.14) $$

\noi This is formula (2.1.10), with $q=q_{0}$ \big(in fact, we use here a slight extension of (2.1.10) since $q_{0}$ is not absolutely continuous with respect to Lebesgue measure on $\mathbb{C}$\big). We now need to study the probability $W_{\infty}^{(2, q_{0})}$. This is the aim of the following Theorem :

\noi {\bf Theorem 2.2.1.} {\it With the notation of Theorem 2.1.1 : }
\begin{eqnarray*} 
 \hspace*{-1cm} {\bf 1)} \qquad \varphi_{q_{0}} (z)
 &=& 2 + \frac{1}{\pi} \; \log \, |z| \qquad {\rm if} \; |z| \ge 1\nonumber \\
 &=& 2 \hspace*{2,75cm} {\rm if} \; |z| \le 1 \hspace*{6,8cm}(2.2.15)
 \end{eqnarray*}

\noi {\it and $(M_{s}^{(q_{0})},\; s \ge 0)$ is the martingale defined by :}
\begin{eqnarray*} 
\hspace*{1cm}M_{s}^{(q_{0})} 
&=& \frac{\varphi_{q_{0}} (X_{s})} {\varphi_{q_{0}} (0)} \; \exp \left(- \frac{1}{2} \, L_{s}^{(C)}\right) �\hspace*{7cm}  (2.2.16) \\
&=& 1 + \frac{1}{\varphi_{q} (0)} \int_{0}^{s} < \nabla \varphi_{q_{0}} (X_{u}), \; dX_{u} > e^{- \frac{1}{2} \, L_{u}^{(C)}} \hspace*{4cm}(2.2.17)
\end{eqnarray*}

\noi {\bf 2)} {\it Let $g_{C} := \sup \{ t \ge 0 \;;\; X_{t} \in C\}$. Then $g_{C}$ is $W_{\infty}^{(2, q_{0})}$ a.s. finite and the r.v. $L_{\infty}^{(C)} (=L_{g_{C}}^{(C)})$ admits as density $f_{ L_{\infty}^{(C)}}^{W_{\infty}^{ (2,q_{0})} }$ with :
$$
f_{ L_{\infty}^{(C)}}^{W_{\infty}^{ (2,q_{0})} } (l) = \frac{1}{2} \; e^{- \frac{l}{2}} \, 1_{[0, \infty[} (l) \eqno(2.2.18) $$

\noi {\bf 3)} Under the probability $W_{\infty}^{(2, q_{0})}$ :

i) Conditionally on $X_{g_C}$, $(X_{s}, \; s \le g_{C})$ and $(X_{g_{C} + s}, \; s \ge 0)$ are independent 

ii) The law of the process $(X_{g_{C} + s}, \; s \ge 0)$ is $\widetilde{P}_{1}^{(2, \log)}$ \big(defined in point 2.2.1, iv)\big)

iii) Conditionally on $L_{g_{C}}^{(C)} =l$, the process $(X_{s}, \; s \le g_{C})$ is a 2-dimensional Brownian process stopped at $\tau_{l}^{(C)}$, and its law, from point 2.2.1 i), is $W_{0}^{(2, \tau_{l}^{(C)})}$. 

\noi In other terms :

{\it iv)} $ \dis \qquad W_{\infty}^{(2, q_{0})} = \frac{1}{2} \int_{0}^{\infty} e^{- \frac{l}{2}} dl
 \big(W_{0}^{(2, \tau_{l}^{(C)})} \circ \widetilde{P}_{1}^{(2, \log)}\big) \hfill(2.2.19) $   }

\noi We note, in particular, that $X_{\tau_{l}^{(C)}}$ under $W_{\infty}^{(2, q_{0})}$ is uniformly distributed on $C$. 

\smallskip

\noi {\bf Proof of Theorem 2.2.1.}

\noi In dimension 1, this Theorem is, essentially, proven in \big([RVY, II]\big). The only item which really differs from those of Theorem 8 in [RVY, II] is point 3, {\it ii)}. We shall emphasize the corresponding arguments. 

\noi \underline{We prove point 3, {\it ii)}}.

\smallskip

\noi We first recall and adapt to dimension 2 the notation and results of [RVY, II].

\smallskip

\noi {\it i)} \underline{Let $(\mathcal{G}_{t}, \; t \ge 0)$ be the smallest filtration containing} $(\mathcal{F}_{t},\; t \ge 0)$
and such that $g_{C}$ is a $(\mathcal{G}_{t}, \; t \ge 0)$ stopping time. Then, there exists a $\big((\mathcal{G}_{t}, \; t \ge 0), \; W_{\infty}^{2, q_{0}}\big)$ 2-dimensional Brownian motion $(B_{t}, \; t \ge 0)$ such that :
$$
X_{t} = B_{t} + \int_{t \wedge g_{_{C}}}^{t} \frac{{n}_{u}}{M_{u}^{(q_{0})} - \underline{M}_{u}^{(q_{0})}} \;\; du         \eqno (2.2.20) $$

\noi with :

\noi $ \bullet $ $n_{u} := e^{- \frac{1}{2} \, L_{u}^{(C)}} \cdot \frac{\nabla \varphi_{q_{0}} (X_{u})}{\varphi_{q_{0}} (0)} \hfill(2.2.21) $

\noi $\bullet$ $\; M_{u}^{(q_{0})}$ is defined by (2.2.16) and :
$$      \underline{M}_{u}^{(q_{0})} := \mathop{\inf}_{s \le u}^{} M_{s}^{(q_{0})} \eqno(2.2.22) $$

\noi {\it ii)} \underline{The function} $\underline{\varphi_{q_{0}} (z)} = 2 + \dis \frac{1}{\pi} \; \log \, |z| \;\; {\rm (for} \; |z| \ge 1)$ \big(see (2.2.15)\big) is increasing in $|z|$. On the other hand, for $u \ge g_{C}$, $L_{u}^{(q_{0})} = L_{g_{C}}^{(C)}$. Thus :
\begin{eqnarray*} 
\hspace*{2cm}\underline{M}_{u}^{(q_{0})} 
&= &M_{g_{C}}^{(q_{0})} =  \frac{ \varphi_{q_{0}} (X_{g_{C}}) }{\varphi_{q_{0}} (0)} \; e^{- \frac{1}{2} \, L_{g_{C}}^{C}}  \\
&=& e^{- \frac{1}{2} \, L_{g_{C}}^{(C)}} \hspace*{8cm}(2.2.23)
\end{eqnarray*}

\noi \big(from (2.2.15) and since $X_{g_{C}} \in C\big)$.

\smallskip

\noi {\it iii)} \underline{Gathering (2.2.20), (2.2.21) and (2.2.23), we obtain} :
\begin{eqnarray*}
X_{t} 
&=& B_{t} + \int_{t \wedge g_{_{C}}}^{t} du \; \frac{�\nabla�\varphi_{q_{0}} (X_{u}) \, e^{- \frac{1}{2} \, L_{g_{C}}^{(C)} }}{ \varphi_{q_{0}} (X_{u}) \, e^{- \frac{1}{2} \, L_{g_{C}}^{(C)}} - 2 e^{- \frac{1}{2} \, L_{g_{C}}^{(C)} } },  \\
&=& B_{t} + \int_{t \wedge g_{_{C}}}^{t} \; \frac{ \nabla \big(\log |\cdot|\big) (X_{u})}{ \log |X_{u}|} \; du \; \text{(after simplification by} \; e^{-\frac{1}{2} \, L_{g_{C}}^{(C)}})\hspace*{0,8cm}(2.2.24)
\end{eqnarray*}

\noi \big(from (2.2.15), since $\varphi_{q_{0}} (X_{u}) - 2 = \dis \frac{1}{\pi} \; \log |X_{u}| \;$ and $\nabla \varphi_{q_{0}} (X_{u}) = \dis \frac{1}{\pi} \big(\nabla \log |\cdot |\big) (X_{u})\big)$.

\smallskip

\noi {\it iv)} \underline{We now use It\^{o}'s formula to express} $|X_{g_{C} +t}|:= \widetilde{R}_{t}$. We obtain, from (2.2.24) :
$$
\widetilde{R}_{t} = (\widetilde{B}_{g_{C} +t} - \widetilde{B}_{g_{C}}) + \int_{0}^{t} \frac{ds}{\widetilde{R}_{s}} \left(\frac{1}{2} + \frac{1}{\log \, \widetilde{R}_{s}} \right) \eqno(2.2.25) $$

\noi where $\big(\widetilde{B}_{g_{C}+t} - \widetilde{B}_{g_{C} }, \; t \ge 0\big)$ is a 1-dimensional Brownian motion started at 1. Thus, from (2.2.6), the law of $\big(|X_{g_{C} + t}|, \; t \ge 0\big)$ is $P_{1}^{(2, \log)}$.

\noi Now, operating in an analogous manner to calculate Arg $(X_{g_{C}+t})$, we obtain :
$$     (X_{g_{C}+t}, \; t \ge 0) = (R_{t} \, e^{i \alpha_{H_{t}}}, \; t \ge 0)   \eqno(2.2.26) $$

\noi with notation of points 2.2.1, {\it ii)}, {\it iii)} and {\it iv)}.

\smallskip

\noi {\bf 2.2.3} \underline{Another description of the measure ${\bf W}^{(2)}$.}

\noi We now present a description of ${\bf W}^{(2)}$ which is analogous, in dimension 2, to the description of {\bf W} given by Theorem 1.1.6.

\smallskip

\noi {\bf Theorem 2.2.2.} 

\smallskip

\noi {\bf 1)} $\dis {\bf W}^{(2)} = \int_{0}^{\infty} dl \big(W_{0}^{(2, \tau_{l}^{(C)})} \circ \widetilde{P}_{1}^{(2, \log)}\big) \hfill(2.2.27)  $

\noi {\bf 2)} {\it For every $t \ge 0$ and $\Gamma_{t} \in b (\mathcal{F}_{t})$ : 
$$
{\bf W}^{(2)} [\Gamma_{t} \, 1_{g_{C} \le t}] = \frac{1}{\pi} \, W^{(2)} \big[\Gamma_{t} \, \log^{+} (|X_{t}|)\big]
\eqno(2.2.28) $$

\noi \big(Recall that $g_{C} := \sup\{s \ge 0 \;;\; X_{s} \in C\}$\big)  }
\begin{equation*} 
 \hspace*{-1cm} {\bf 3)} \quad {\it i)} \qquad {\bf W}^{(2)} (g_{C} \in dt) = e^{- \frac{1}{2 t}} \; \frac{dt}{2 \pi t} \quad (t \ge 0) \hfill(2.2.29)
 \end{equation*}
{\it ii)} {\it Conditionally on $g_{C} = t$, the law of the process $(X_{u}, \; u \le g_{C})$, under ${\bf W}^{(2)}$ is $\Pi_{0}^{(2,t,U)}$,

\noi where : 
\begin{itemize}
\item[$\bullet$] $U$ is a r.v. uniformly distributed on $C$ ;
\item[$\bullet$] Conditionally on $U=u, \; \Pi_{0}^{(2,t,U)}$ is the law of a 2-dimensional Brownian bridge 

$(b_{s}^{(2,t,u)}, \; 0 \le s \le t)$ of length $t$ such that $b_{0}^{(2,t,u)}=0$ and $b_{t}^{(2,t,u)}=u$.
\end{itemize}
\begin{equation*}
 \hspace*{-1cm} {\bf } \quad {\it iii)} \qquad {\bf W}^{(2)}= \int_{0}^{\infty} \frac{dt}{2 \pi t} \; e^{- \frac{1}{2t}} \big(\Pi^{2,t,U} \circ \widetilde{P}_{1}^{(2, \log)}\big) \hfill(2.2.30)
  \end{equation*} }

\noi {\bf Proof of Theorem 2.2.2.}

\noi {\it i)} Point 1) is an easy consequence of (2.2.14), (2.2.19) and (2.2.18).

\smallskip

\noi{\it ii)} \underline{We now show (2.2.28)}

\noi For this purpose, we use the definition (2.1.10) of ${\bf W}^{(2)}$ with $q=\lambda \,q_{0}$ (where $q_{0}$ is defined by (2.2.3), and $\lambda >0)$. We have :
$$     \varphi_{\lambda \, q_{0}} (z) = \frac{2}{\lambda} + \frac{1}{\pi} \; \log^{+} (|z|)   \eqno(2.2.31) $$

\noi \big(see (2.2.15)\big). Thus, for every $t \ge 0$ and $\Gamma_{t} \in b(\mathcal{F}_{t})$ :
\begin{eqnarray*}
 \hspace*{1cm}W^{(2)} \left(\Gamma_{t} \left( \frac{2}{\lambda} + \frac{1}{\pi} \, \log^{+} (|X_{t}|)\right) \right)
&=& \varphi_{\lambda \,q_{0}} (0) W_{\infty}^{(2, \lambda \,  q_{0})} \big(\Gamma_{t} \, e^{ \frac{\lambda}{2} \, L_{t}^{(C)})}\big)  \\
&=& {\bf W}^{(2)} \big(\Gamma_{t} \, e^{- \frac{\lambda}{2} (L_{\infty}^{(C)} - L_{t}^{(C)})}\big)  \hspace*{2,5cm} (2.2.32)
\end{eqnarray*}

\noi We then let $\lambda \to \infty$ in (2.2.32) and note that $L_{\infty}^{(C)} - L_{t}^{(C)} >0$ on the set $(g_{C} > t)$ (and equals to 0 on $g_{C} \le t$). The monotone convergence Theorem implies :
$$
\frac{1}{\pi} \; W^{(2)} \big(\Gamma_{t} \, \log^{+} (|X_{t}|) \big) = {\bf W}^{(2)} (\Gamma_{t} \, 1_{g_{C} \le t})$$

\noi This is (2.2.28). Note that we may replace $t$ by a stopping time $T$ in (2.2.28). We obtain :
$$
{\bf W}^{(2)} \big(\Gamma_{T} \, 1_{g_{C} \le T < \infty}) = \frac{1}{\pi} \, W^{(2)} \big(\Gamma_{T} \, \log^{+} (|X_{T}|) 1_{T < \infty}\big)  \eqno(2.2.33) $$

\noi with $\Gamma_{T} \in b(\mathcal{F}_{T})$.

\smallskip

\noi {\bf Remark 2.2.3.}

\noi We deduce from (2.2.32) and (2.2.28) :
\begin{eqnarray*}
 \hspace*{1cm}\frac{2}{\lambda} \; W^{(2)} (\Gamma_{t})
&=& {\bf W}^{(2)} \left(\Gamma_{t} \, 1_{g_{C} > t} \; \exp \big(-\frac{\lambda}{2} \big(L_{\infty}^{(C)} - L_{t}^{(C)}\big)\big) \right)  \\
&=& W^{(2)} (\Gamma_{t}) \left( \int_{0}^{\infty} e^{- \frac{\lambda}{2} \, l} dl \right) \hspace*{6cm} (2.2.34) 
\end{eqnarray*}

\noi and
$$
\frac{1}{\pi} W^{(2)} \big(\log^{+} |X_{t}|\big) = {\bf W}^{(2)} (g_{C} \le t) = {\bf W}^{(2)} (L_{\infty}^{(C)} -L_{t}^{(C)} =0)
\eqno(2.2.35) $$

\noi Then, operating as in the proof of Theorem 1.1.6, point 3) {\it i)} \big(see (1.1.45) and (1.1.46)\big), we obtain :

 \noi {\it i)} $ {\bf W}^{(2)} (L_{\infty}^{(C)} -L_{t}^{(C)} \in dl) = 1_{[0, \infty[} (l) dl + \frac{1}{\pi}W^{(2)} \big(\log^{+} (|X_{t}|)\big) \delta_{0}  (dl) \hfill(2.2.36) $
  
\noi {\it ii)} Conditionally on $L_{\infty}^{C} - L_{t}^{C} = l \quad (l > 0), \; (X_{u}, \; u \le t)$ is, under ${\bf W}^{(2)}$, a 2-dimensional Brownian motion indexed by $[0,t]$.

\smallskip

\noi {\bf Remark 2.2.4.} We can obtain (2.2.28) in the same manner as for point 2) of Remark 1.1.9. For this purpose, we need a scale function for the $W^{ (2, q_{0}) }$ process. The function $\dis z \to \frac{1}{1+ \frac{1}{\pi} \; \log (|z|)} \; (|z| \ge 1)$ is an adequate choice.

\smallskip

\noi {\it iii)} \underline{We now prove point 3 {\it i)} of Theorem 2.2.2.}

\noi We write (2.2.28) with $\Gamma_{t} \equiv 1$ :
$$       {\bf W}^{(2)} (g_{C} \le t) = \frac{1}{\pi} \; W^{(2)} (\log^{+} |X_{t}|)  \eqno(2.2.37) $$

\noi and we differentiate (2.2.37) with respect to $t$. Thus : 
\begin{eqnarray*}
{\bf W}^{(2)} (g_{C} \in dt)
&=& \frac{1}{\pi} \; \left(\frac{d}{dt} \, W^{(2)} (\log^{+} |X_{t}|)\right) \cdot dt \\
&=& \frac{1}{\pi} \; \frac{d}{dt} \, W^{(2)} \left( 1_{|X_{1}| > \frac{1}{\sqrt{t}}} \left(\log \, \sqrt{t} - \log \frac{1}{| X_{1}|} \right)\right) \cdot dt 
\quad \hbox{(by scaling)} \\
&=& \frac{1}{2 \pi t} \; W^{(2)} \left( \frac{ |X_{1}|^{2}}{2} > \frac{1}{2 t} \right) dt \\
&=& \frac{1}{2 \pi t} \; e^{- \frac{1}{2t}} dt \qquad (t \ge 0)
\end{eqnarray*}

\noi since $\dis \frac{|X_{1}|^{2}}{2}$ is a standard exponential r.v.

\noi The end of the proof of Theorem 2.2.2 is obtained by using arguments similar to those used for Theorem 1.1.6. We note, in particular, that conditionally on $X_{g_{C}}$, $(X_{g_{C} +t}, \; t \ge 0)$ and $(X_{s}, \; s \le g_{C})$ are independent.

\noi {\bf Remark 2.2.5.} From (2.2.29), we deduce :
$$
{\bf W}^{(2)} (e^{- \frac{\lambda^{2}}{2} g_{C}}) = \int_{0}^{\infty} \frac{dt}{2 \pi t} \; 
e^{ -\frac{\lambda^{2}}{2} \; t - \frac{1}{2 t} } = K_{0} ( \lambda) \eqno(2.2.38) $$

\noi where $K_{0}$ denotes the Bessel-Mc Donald function with index 0 \big(see [L], formula 5.10.25\big).

\newpage

 {\bf 2.3.} {\bf Study of the winding process under W$^{(2)}$.}

\noi Formula (2.2.12) : 
$$           X_{g_{C} + t} = R_{t} \; e^{i \alpha_{H_{t}}}, \; (t \ge 0) $$

\noi which provides a representation of $X$ after $g_{C}$ under $W^{(2)}$ invites to establish for this process a theorem similar to the classical theorem of Spitzer, which we recall : 

\smallskip

\noi {\bf 2.3.1} \underline{Spitzer's Theorem}.

\smallskip

\noi {\bf Theorem.} \big(Spitzer [S]\big)

\noi {\it Let $(X_{t}, \; t \ge 0)$ a $\mathbb{C}$ valued Brownian motion, starting from $z \neq 0$. We have : 
$$         X_{t} = |X_{t}| \, e^{i \alpha_{H_{t}}}  \eqno(2.3.1) $$

\noi with :    }

\noi {\it i)} $(\alpha_{u}, \; u \ge 0)$ a 1-dimensional Brownian motion independent from the 2-dimensional Bessel process $(|X_{t}|, \; t \ge 0)$ (one can also find a precise study of the winding process of planar Brownian motion in [PY1]). 

 \noi {\it ii)} $\dis  \qquad H_{t} = \int_{0}^{t} \frac{ds}{|X_{s}|^{2}} \hfill (2.3.2) $
 
\noi {\it Let $\dis (\theta_{t}, \; t \ge 0) := (\alpha_{H_{t}}, \; t \ge 0) = \left(\theta_{0} + {\rm Im} \int_{0}^{t} \frac{d X_{s}}{X_{s}}, \; t \ge 0\right)$ be the winding process. Then :
$$
\frac{2 \theta_{t}}{\log \, t} \; \mathop{\longrightarrow}_{t \to \infty}^{\rm (law)} \Gamma \mathop{=}_{}^{\rm (law)} \alpha_{T_{1} (\gamma)} \eqno(2.3.3) $$

\noi In (2.3.3), $(\gamma_{t}, \; t \ge 0)$ is a 1-dimensional Brownian motion started from 0 and independent from $(\alpha_{u}, \; u \ge 0)$. and :
$$            T_{1} (\gamma) := \inf \{ s \ge 0 \;;\; \gamma_{s} =1\}   \eqno(2.3.4) $$

\noi {\it iii)} Consequently $\Gamma$ is a standard Cauchy r.v.  }

\smallskip

\noi {\bf 2.3.2.} \underline{An analogue of Spitzer's Theorem}.

\smallskip

\noi Now, here is the analogue of the above (Spitzer) Theorem for the process $(X_{g_{C} + t}, \; t \ge 0)$ :

\smallskip

\noi {\bf Theorem 2.3.1.} {\it Under $\widetilde{P}_{1}^{(2, \log)}$, the winding process $(\theta_{t}, \; t \ge 0) = (\alpha_{H_{t}}, \; t \ge 0)$ satisfies :

 \noi {\bf 1)} $\dis \quad  \frac{4}{(\log \, t)^{2}} \; H_{t} \; \mathop{\longrightarrow}_{t \to \infty}^{\rm (law)} T_{1}^{(3)} \hfill (2.3.5) $
 
 \smallskip
 
\noi where $\dis  \quad  T_{1}^{(3)} := \inf \{ u \;;\; \rho_{u} =1\} \hfill(2.3.6)$ 

\noi is the first hitting time of level 1 by a 3-dimensional Bessel process $(\rho_{u}, \; u \ge 0)$ started at 0.

\smallskip

\noi {\bf 2)} $\dis \quad  \frac{2}{\log \, t} \; \theta_{t} \; \mathop{\longrightarrow}_{t \to \infty}^{\rm (law)}  \alpha_{T_{1}^{(3)}} \hfill(2.3.7) $

\smallskip
 
\noi where $(\alpha_{u}, \; u \ge 0)$ is a 1-dimensional Brownian motion independent from $(\rho_{u}, \; u \ge 0)$.     }

\smallskip

\noi \underline{We now recall our notation} (see Section 2.2.1)
\begin{itemize}
\item[$\bullet$] $(R_{t}, \; t \ge 0)$ is the process defined in (2.2.6)
\item[$\bullet$]  $\dis H_{t} = \int_{0}^{t} \frac{ds}{R_{s}^{2}}  \hfill(2.3.8) $
\item[$\bullet$] $(\alpha_{u}, \; u \ge 0)$ is a 1-dimensional Brownian motion independent from $(R_{t}, \; t \ge 0)$
\item[$\bullet$] $(\log \, R_{t}, \; t \ge 0) = (\rho_{H_{t}}, \; t \ge 0)$ and $(\rho_{u}, \; u \ge 0)$ is a 3-dimensional Bessel process started at 0.
\end{itemize} 

\smallskip

\noi {\bf Remark 2.3.2.}

\noi {\bf 1.} Theorem 2.3.1 differs from Spitzer's Theorem in that $T_{1}$ has been replaced by $T_{1}^{(3)}$.

\noi {\bf 2.} Let, for every $z \in \mathbb{C}, \; {\bf W}_{z}^{(2)}$ be defined by :
\begin{equation*}
{\bf W}_{z}^{(2)} \big(F (X_{s}, \; s \ge 0)\big) := {\bf W}^{(2)} \big(F (z+X_{s}, \; s \ge 0)\big)
\end{equation*}

\noi Theorem 2.3.1 then implies that, for $z \neq 0$, under ${\bf W}_{z}^{(2)}$ and conditionally on $g_{C} \le a$, the winding process $(\theta_{t}, \; t \ge 0)$ satisfies : 
$$
\frac{2}{\log \, t} \; \theta_{t} \mathop{\longrightarrow}_{t \to \infty}^{} \alpha_{T_{1}^{(3)}} \eqno(2.3.9) $$

\noi for all $a > 0$. This easily results from (2.3.7) and from the representation formula (2.2.6).

\smallskip

\noi \underline{Proof of Theorem 2.3.1.}

\noi {\it i)} \underline{We use the notation} (2.3.8). We admit for a moment that :
$$            H_{t}- H_{T_{\sqrt{t}} (R)} \qquad \text{converges in law as} \; t \to \infty, \; \text{with :}  \eqno(2.3.10) $$

\noi 
$$    T_{\sqrt{t}}(R)  := \inf \{ s \ge 0 \;;\; R_{s} \ge \sqrt{t}\}  \eqno(2.3.11) $$

\noi and we show that (2.3.10) implies Theorem 2.3.1. Indeed, from (2.3.10), we have :
$$
\frac{4}{(\log\, t)^{2}} \; H_{t} \mathop{\sim}_{t \to \infty}^{} \frac{1}{(\log \, \sqrt{t})^{2}} \; H_{T_{\sqrt{t}}(R)} 
\eqno(2.3.12) $$

\noi But :
$$
\frac{1}{(\log \, a)^{2}} \; H_{T_{a} (R)} = \frac{1}{(\log \, a)^{2}} \; T_{\log \, a} (\rho) \mathop{=}_{}^{\rm (law)} T_{1} (\rho) \eqno(2.3.13) $$

\noi with
$$    T_{\log (a)} (\rho) := \inf \{ t \ge 0 \;;\; \rho_{t} \ge \log \, a\} \eqno(2.3.14) $$

\noi The first equality in (2.3.13) results from definitions \big(see point 4 of (2.3.8)\big) and the second from the scaling property. Thus, from (2.3.10), we deduce :
$$
\frac{4}{(\log \, t)^{2}} \; H_{t} \; \mathop{\longrightarrow}_{t \to \infty}^{\rm (law)} T_{1}^{(3)} \eqno(2.3.15) $$

\noi and
\begin{eqnarray*}
\frac{2}{\log \, t} \; \theta_{t} 
&=&\frac{2}{\log \,t} \; \alpha_{H_{t}} \mathop{=}_{}^{\rm (law)} \frac{2 \sqrt{H_{t}}} {\log \, t} \; \alpha_{1} \quad \hbox{(by scaling)}\\
&& \mathop{\longrightarrow}_{t \to \infty}^{\rm (law)} \, \sqrt{T_{1}^{(3)}} \cdot \alpha_{1} \\
&&  \mathop{=}_{}^{\rm (law)}  \alpha_{ T_{1}^{(3)} } \qquad \text{(by scaling)}
\end{eqnarray*}

\noi which proves Theorem 2.3.1.

\smallskip

\noi {\it ii)} \underline{It remains to prove (2.3.10).}

\noi For this purpose, we start with the following Lemma : 

\noi {\bf Lemma 2.3.3.} {\it Let $(R_{t}, \; t \ge 0)$ be defined by (2.2.6). Then : 
$\dis \left( \frac{1}{\sqrt{t}} \, R_{tv}, \; v \ge 0 \right)$ converges in law, as $t \to \infty$, to a 2-dimensional Bessel process starting from 0. }

\smallskip

\noi {\bf Proof of Lemma 2.3.3.}

\noi From (2.2.6) we have :
$$
R_{t} = 1 + \beta_{t} + \int_{0}^{t} \left( \frac{1}{2 R_{s}} + \frac{1}{R_{s} \log \, R_{s}} \right) ds  $$

\noi Thus :
$$
\frac{1}{\sqrt{t}} \; R_{tv} = \frac{1}{\sqrt{t}} + \frac{1}{\sqrt{t}} \; \beta_{tv} + \frac{1}{\sqrt{t}} \int_{0}^{tv} \left( 
\frac{1}{2 R_{s}} + \frac{1}{R_{s} \log \, R_{s}} \right) ds \eqno(2.3.16) $$

\noi Denoting by $(\widetilde{\beta}_{v}, \; v \ge 0)$ the Brownian motion $\dis \left(\frac{1}{\sqrt{t}} \; \beta_{tv}, \; v \ge 0 \right)$ and making the change of variable $s=tv$, we obtain, with $\dis \left( \widetilde{R}_{v}^{(t)} = \frac{1}{\sqrt{t}} \; R_{tv}, \; v \ge 0 \right)$ :
$$
\widetilde{R}_{v}^{(t)} = \frac{1}{\sqrt{t}} + \widetilde{\beta}_{v} + \int_{0}^{v} \left( \frac{1}{2 \widetilde{R}_{u}^{(t)}} + \frac{1}{\widetilde{R}_{u}^{(t)} \big(\log \, \sqrt{t} + \log \widetilde{R}_{u}^{(t)}\big)} \right) du   \eqno(2.3.17) $$

Hence, as $t \to \infty$, $(\widetilde{R}_{v}^{(t)}, \; v \ge 0)$ converges in law to the law of the solution of the SDE :
$$\widetilde{R}_{v} = \widetilde{\beta}_{v} + \int_{0}^{v} \frac{du}{2 \widetilde{R}_{u} } $$
\noi i.e. to (the law of) a 2-dimensional Bessel process started at 0.

\smallskip

\noi {\it iv)} \underline{We may now end up the proof of (2.3.10).}

\noi We have, from (2.3.8) :
$$
 H_{T_{\sqrt{t}}(R)} - H_{t}= \int_{t}^{T_{\sqrt{t}}(R)} \frac{du}{R_{u}^{2}} = \int_{1}^{\frac{1}{t} \, T_{\sqrt{t}} (R) } \frac{dv}{\big( \frac{1}{t} \, R_{vt}^{2} \big)}    $$

\noi after making the change of variable $u=tv$. But, from Lemma 2.3.3, $\dis \left(\frac{1}{\sqrt{t}} \, R_{vt}, \; v \ge 0\right)$ converges in law to a 2-dimensional Bessel process $\big(R_{0}^{(2)} (v), \; v \ge 0\big)$ starting from 0. Thus : $H_{t}-H_{\sqrt{t}} (R)$ converges in law, as $t \to \infty$, to 
$$        \int_{1}^{T_{1} (R_{0}^{(2)})} \; \frac{du}{(R_{0}^{(2)} (u))^{2}} \eqno(2.3.18) $$

\noi with $T_{1} (R_{0}^{(2)}) = \inf \{ s \ge 0 \;;\; R_{0}^{(2)} (s) =1\}$.

\smallskip

\noi {\bf Remark 2.3.4.} (An extension of Theorem 2.3.1.)

\noi Let $(\beta_{t}, \; t \ge 0)$ denote a 1-dimensional Brownian motion starting at 0, $\delta >0$ and $(R_{t}^{(\delta)},$

\noi $t \ge 0)$ the solution of :
$$
R_{t}^{(\delta)} = 1 + \beta_{t} + \int_{0}^{t} \left( \frac{1}{2 R_{s}^{(\delta)}} + \frac{\delta}{R_{s}^{(\delta)} \log \, R_{s}^{(\delta)}} \right) ds   \eqno(2.3.19) $$

\noi The case we have just studied is that of $\delta =1$. Let :
$$     H_{t}^{(\delta)} := \int_{0}^{t} \frac{ds}{(R_{s}^{(\delta)})^{2}}   \eqno(2.3.20) $$

\noi and
$$         \theta_{t}^{(\delta)} = \alpha_{H_{t}^{(\delta)}}        $$

\noi where $(\alpha_{u}, \; u \ge 0)$ is a 1-dimensional Brownian motion independent from $(\beta_{t}, \; t \ge 0)$.

\noi The technique we have just developed allows to obtain : 

 \noi {\it i)} $\dis \quad  (\log \, R_{t}^{(\delta)}, \; t \ge 0) = (\rho_{H_{t}^{(\delta)}}^{(2 \delta +1)}, \; t \ge 0)
\hfill(2.3.21) $

\noi where $(\rho_{u}^{(2 \delta +1)}, \; u \ge 0)$ is a $(2 \delta +1)$-dimensional Bessel process starting at 0.
\begin{equation*}
 \hspace*{-1cm} {\it ii)} \quad  \frac{4}{(\log \, t)^{2}} \; H_{t}^{(\delta)} \mathop{\longrightarrow}_{t \to \infty}^{\rm (law)} T_{1}^{(2 \delta+1)}
 \end{equation*} 

\noi where $T_{1}^{(2 \delta+1)} := \inf \{u \ge 0 \;;\; \rho_{u}^{(2 \delta +1)} =1\}$.

\noi {\it iii)} $\dis \quad  \frac{2 \theta_{t}^{(\delta)}}{\log \, t} = \frac{2 \alpha_{H_{t}}^{(\delta)}}{\log \, t} \mathop{\longrightarrow}_{t \to \infty}^{\rm (law)} \alpha_{T_{1}^{(2 \delta +1)}} \hfill (2.3.22) $

\noi where $T_{1}^{(2 \delta +1)}$ is independent from the 1-dimensional Brownian motion $(\alpha_{u}, \; u \ge 0)$.

\bigskip

{\bf 2.4 $W^{(2)}$ martingales associated to ${\bf W}^{(2)}$.}

\noi Just as in Chapter 1, we associated to any r.v. $F \in L^{1} (\mathcal{F}_{\infty}, {\bf W})$ the $\big((\mathcal{F}_{t}, \; t \ge 0), \; W\big)$ martingale $\big(M_{t} (F), \; t \ge 0)$, we now associate to every r.v. $F \in L^{1} (\mathcal{F}_{\infty}, {\bf W}^{(2)})$ a $\big((\mathcal{F}_{t}, \; t \ge~0), \; W^{(2)}\big)$ martingale $\big(M_{t}^{(2)} (F), \; t \ge 0\big)$.

\smallskip

\noi {\bf 2.4.1} \underline{Definition of $\big(M_{t}^{(2)} (F), \; t \ge 0\big)$.}

\noi {\bf Theorem 2.4.1.} {\it Let $F \in L^{1} \big(\Omega = \mathcal{C} (\mathbb{R}_+ \to \mathbb{C}), \; \mathcal{F}_{\infty}, {\bf W}^{(2)}\big)$. There exists a $\big((\mathcal{F}_{t}, \; t \ge 0), \; W^{(2)}\big)$ martingale (which is necessarily continuous) $\big(M_{t}^{(2)} (F), \; t \ge 0\big)$, positive if $F \ge 0$, such that : 

\noi {\bf 1)} For every $t \ge 0$ and $\Gamma_{t} \in b (\mathcal{F}_{t})$ :
$$   {\bf W}^{(2)} (F \cdot \Gamma_{t}) = W^{(2)} \big(M_{t}^{(2)} (F) \cdot \Gamma_{t}\big) \eqno(2.4.1) $$

\noi In particular, for every $t \ge 0$ :
$$     {\bf W}^{(2)} (F) = W^{(2)} \big(M_{t}^{(2)} (F)\big)  \eqno(2.4.2) $$

\noi and, if $F$ and $G$ belong to $L_{+}^{1} (\mathcal{F}_{\infty}, {\bf W}^{(2)})$ :
$$
W^{(2)} \big(M_{t}^{(2)} (F) \cdot M_{t}^{(2)} (G)\big) = {\bf W}^{(2)} \big(F \cdot M_{t}^{(2)} (G)\big) = {\bf W}^{(2)} \big(M_{t}^{(2)} (F) \cdot G\big) \eqno(2.4.3) $$
  
\noi  {\bf 2)} $\dis \quad M_{t}^{(2)} (F) = \widehat{W}_{X_{t}(\omega_{t})}^{(2)} (F (\omega_{t}, \widehat{\omega}^{t})) \hfill(2.4.4) $
 
\noi {\bf 3)} $\dis \quad M_{t}^{(2)} (F) \underset{t \rightarrow \infty}{\longrightarrow} 0 \qquad W^{(2)} \; {\rm a.s.}  \hfill(2.4.5) $
  
 \noi In particular, the martingale $(M_{t}^{(2)} (F), \; t \ge 0)$ is not uniformly integrable if $F \neq 0$. 
 
 \smallskip
 
 \noi {\bf 4)} For every $q \in \mathcal{I}$ :
$$
 M_{t}^{(2)} (F) = \varphi_{q} (0) \; M_{t}^{(q)} \, W_{\infty}^{(2,q)} (F \, e^{\frac{1}{2}\, A_{\infty}^{(q)}} | \mathcal{F}_{t})  \eqno(2.4.6) $$

\noi where $M_{t}^{(q)}, \, \varphi_{q}$ and $W_{\infty}^{(2, q)}$ are defined in Theorem 2.1.1. }

\smallskip

\noi The proof of Theorem 2.4.1 is, mutatis mutandis, the proof of Theorem 1.2.1.  Here are some examples of martingales $(M_{t}^{(2)} (F), \; t \ge 0)$.

\smallskip

\noi {\bf Example 2.1.} Let $q \in \mathcal{I}$ and $\dis F_{q} = \exp \left(- \frac{1}{2} \, A_{\infty}^{(q)}\right)$. We have, from (2.1.10) :
$$    {\bf W}^{(2)} (F_{q}) = \varphi_{q} (0)   \eqno(2.4.7) $$

\noi and
$$
\left(M_{t}^{(2)} (F_{q}) = \varphi_{q} (X_{t}) \, \exp \left(- \frac{1}{2} \, A_{t}^{(q)}\right), \; t \ge 0\right)  \eqno(2.4.8) $$

\noi In particular, for $q = \lambda q_{0}\;$ \big(see (2.2.3) and (2.2.31)\big) :
$$
M_{t}^{(2)} \left( \exp - \frac{\lambda}{2} \,A_{\infty}^{(q_{0})} \right) = \left(\frac{2}{\lambda} + \frac{1}{\pi} \, \log^{+} \big(|X_{t}|\big)\right) \, \exp \left(- \frac{\lambda}{2} \, L_{t}^{(C)}\right)   \eqno(2.4.9) $$

\noi {\bf Example 2.2.} \big(see [RVY, VI]\big).

\noi We write the skew-product representation of the canonical 2-dimensional Brownian motion $(X_{t}, \; t \ge 0)$ starting at $z \neq 0$ as :
$$      X_{t} = |X_{t}| \cdot \exp (i \, \alpha_{H_{t}})   \eqno(2.4.10)    $$

\noi where : 

\noi {\it i)} $\big(|X_{t}|, \; t \ge 0\big)$ is a 2-dimensional Bessel process starting at $|z|$.

\noi {\it ii)} $\; \dis H_{t} = \int_{0}^{t} \frac{ds}{|X_{s}|^{2}}$

\noi {\it iii)} $(\alpha_{u}, \; u \ge 0)$ is a 1-dimensional Brownian motion, independent from $\big(|X_{u}|, \; u \ge 0\big)$.

\noi Let $(\theta_{t} := \alpha_{H_{t}}, \; t \ge 0)$ denote the winding process and introduce : 
$$
S_{t}^{\theta} := \mathop{\sup}_{s \le t}^{} \theta_{s} = \mathop{\sup}_{u \le H_{t}}^{} \alpha_{u}
\eqno(2.4.11) $$

\noi Let $\varphi : \mathbb{R}_{+} \to \mathbb{R}_{+}$ Borel and integrable. Then :
$$
\big(M_{t}^{(2)} \big(\varphi (S_{\infty}^{\theta})\big), \; t \ge 0\big) = \left( \varphi (S_{t}^{\theta}) (S_{t}^{\theta} - \theta_{t}) + \int_{S_{t}^{\theta}}^{\infty} \varphi (y) dy, \quad t \ge 0 \right)  \eqno(2.4.12) $$

\noi {\bf 2.4.2} \underline{A decomposition Theorem of positive $W^{(2)}$ supermartingales.}

\noi Just as in Theorem 1.2.5, we have obtained a decomposition Theorem for every $\big((\mathcal{F}_{t}, \; t \ge 0),$

\noi $W \big)$ positive supermartingale, we now present a decomposition theorem for every $\big((\mathcal{F}_{t}, \; t \ge 0),$

\noi $W^{(2)} \big)$ positive supermartingale.

\smallskip

\noi {\bf Theorem 2.4.2.} {\it Let $(Z_{t}, \; t \ge 0)$ denote a positive $\big(\Omega = \mathcal{C} (\mathbb{R}_{+}  \to \mathbb{C}), \; (\mathcal{F}_{t}, \; t \ge 0), \; W^{(2)}\big)$ supermartingale. We denote $Z_{\infty} := \dis \mathop{\rm lim}_{t \to \infty}^{} Z_{t}, \; W^{(2)}$ a.s. Then :

\noi  {\bf 1)} $\dis \quad z_{\infty} := \mathop{\rm lim}_{t \to \infty}^{} \;\; \pi \frac{Z_{t}}{1+ \log^{+} (|X_{t}|)} $  exists ${\bf W}^{(2)}$   a.s. $\hfill (2.4.13)$

\noi and : $\quad {\bf W}^{(2)} (z_{\infty}) < \infty \hfill (2.4.14) $
 
 \noi {\bf 2)} $(Z_{t}, \; t \ge 0)$ decomposes in a unique manner in the form :
$$
Z_{t} = M_{t}^{(2)} (z_{\infty}) + W^{(2)} (Z_{\infty}| \mathcal{F}_{t}) + \xi_{t} \quad (t \ge 0) \eqno(2.4.15) $$

\noi where $\big(M_{t}^{(2)} (z_{\infty}), \; t \ge 0\big)$ and $\big(W^{(2)} (Z_{\infty} |\mathcal{F}_{t}), \; t \ge 0\big)$ denote two $\big((\mathcal{F}_{t},\; t \ge 0), \; W^{(2)}\big)$ martingales and :

\noi $(\xi_{t}, \; t \ge 0)$ is a $\big((\mathcal{F}_{t}, \; t \ge 0), \; W^{(2)}\big)$ positive supermartingale such that : 

\noi i) $Z_{\infty} \in L_{+}^{1} (\mathcal{F}_{\infty}, W^{(2)})$, hence $W^{(2)} (Z_{\infty} | \mathcal{F}_{t})$ converges $W^{(2)}$ a.s. and in $L^{1} (\mathcal{F}_{\infty}, W^{(2)})$ towards $Z_{\infty}$.

\noi ii) $\quad  \;\;\;\; \dis \frac{W(Z_{\infty} | \mathcal{F}_{t}) + \xi_{t}}{1+ \log^{+} (|X_{t}|)} \mathop{\longrightarrow}_{t \to \infty}^{} 0 \quad {\bf W}^{(2)}$ a.s.
\begin{equation*}
 \hspace*{-0,9cm} {\it iii)} \qquad M_{t}^{(2)} (z_{\infty}) + \xi_{t}  \mathop{\longrightarrow}_{t \to \infty}^{} 0 \qquad W^{(2)} \; {\rm a.s.}
 \end{equation*}

\noi In particular, if $F \in L^{1} (\mathcal{F}_{\infty}, {\bf W}^{(2)})$, then :
$$
\pi \cdot \frac{M_{t} (F)}{1+ \log^{+} (|X_{t}|)} \mathop{\longrightarrow}_{t \to \infty}^{} F \qquad {\bf W}^{(2)} \;\; {\rm a.s.}  \eqno(2.4.16) $$

\noi and the map : $\; F \to \big(M_{t}^{(2)} (F), \; t \ge 0\big)$ is injective. }

\smallskip

\noi {\bf Corollary 2.4.3.} {\it (A characterisation of martingales of the form $\big(M_{t}^{(2)} (F), \; t \ge 0\big)$. A $\big((\mathcal{F}_{t}, \; t \ge 0), \; W^{(2)}\big)$ positive martingale $(Z_{t}, \; t \ge 0)$ is equal to $\big(M_{t}^{(2)} (F), \; t \ge 0\big)$ for an $F \in L^{1} (\mathcal{F}_{\infty}, {\bf W}^{(2)})$ if and only if : }
$$
Z_{0} = {\bf W}^{(2)} \left( \mathop{\rm lim}_{t \to \infty}^{} \pi \cdot \frac{Z_{t}}{1+ \log^{+} (|X_{t}|)} \right)
\eqno(2.4.17) $$

\noi Note that $\dis \mathop{\rm lim}_{t \to \infty}^{} \; \frac{Z_{t}}{1+ \log^{+} (|X_{t}|)}$ exists ${\bf W}^{(2)}$ a.s. from (2.4.13).

\smallskip

\noi {\bf Sketches of Proofs of Theorem 2.4.2 and of Corollary 2.4.3.}

\noi This proof is essentially the same as those of Theorem 1.2.5 and of Corollary 1.2.6. Two arguments need to be modified :

\noi {\it i)} The role of the r.v. $g$ in the proof of Theorem 1.2.5 is played here by that of the r.v. $g_{C}$.

\noi {\it ii)} The relation (1.1.41) : $ \dis {\bf W} (\Gamma_{t} \, 1_{g \le t}) = W(\Gamma_{t} |X_{t}|) $

\noi and the limiting result :
$$
\frac{\varphi_{q} (X_{t}) \, \exp (- \frac{1}{2} \, A_{t}^{(q)})} {1+ |X_{t}|} \mathop{\longrightarrow}_{t \to \infty}^{} \exp \; \left(- \frac{1}{2} \, A_{\infty}^{(q)}\right) \eqno(2.4.18) $$

\noi which were used in the proof of Lemma 1.2.8 need to be replaced respectively by :
$$
{\bf W}^{(2)} \left(\Gamma_{t} \, 1_{(g_{_{C}} \le t)} \right)= \left.\frac{1}{\pi} \; W^{(2)} (\Gamma_{t} \, \log^{+} |X_{t}|)\right).   $$

\noi (This is relation (2.2.28) of Theorem 2.2.2) and by :
$$
\pi \cdot  \frac{\varphi_{q} (X_{t}) \, \exp (- \frac{1}{2} \, A_{t}^{(q)})}{1+ \log^{+} (|X_{t}|)} \mathop{\longrightarrow}_{t \to \infty}^{} \exp \; \left(- \frac{1}{2} \, A_{\infty}^{(q)}\right) \quad {\bf W}^{(2)} \; {\rm a.s.} \eqno(2.4.19) $$

\noi The latter (2.4.19) follows easily from :
\begin{equation*}
\pi \cdot \varphi_{q} (z) \mathop{\sim}_{|z| \to \infty}^{} \log (|z|), \quad {\rm from} \; (2.1.6)
\end{equation*}

\noi and from : $|X_{t}| \dis \mathop{\longrightarrow}_{t \to \infty}^{} \infty \quad {\bf W}^{(2)}$ a.s.

\noi since the canonical process under $W_{\infty}^{(2,q)}$ is transient.

\noi {\bf 2.4.3} \underline{A decomposition Theorem for the martingales $\big(M_{t}^{(2)} (F), \; t \ge 0\big)$.}

\noi A difference with the preceding subsection is that the r.v.'s $F$ which we now consider belong to $L^{1} (\mathcal{F}_{\infty}, {\bf W}^{(2)})$ but are not necessarily positive. Here is the analogue, in dimension 2, of Theorem 1.2.11.

\smallskip

\noi {\bf Theorem 2.4.4.} {\it $F \in L^{1} (\mathcal{F}_{\infty}, {\bf W}^{(2)})$ and let $\big(M_{t}^{(2)} (F), \; t \ge 0\big)$ the $\big((\mathcal{F}_{t}, \; t \ge 0), \; W^{(2)}\big)$ martingale associated to $F$ by Theorem 2.4.1. Let $C$, $(L_{t}^{(C)}, \; t \ge 0)$ and $g_{C}$ be as in Section 2.2.1, {\it i)} and Section 2.2.2. Then :

\noi {\bf 1)} {\it i)} There exists a previsible process $\big(k_{s}^{(C)} (F), \; s \ge 0\big)$ which is defined $dL_{s}^{(C)} \cdot W^{(2)} (d \omega)$ a.s., positive if $F \ge 0$, and such that : 
$$
 W^{(2)} \left( \int_{0}^{\infty} \big| k_{s}^{(C)} (F) \big| dL_{s}^{(C)} \right) = {\bf W}^{(2)} \big(\big|k_{g_{C}}^{(C)} (F)\big|\big) \le {\bf W}^{(2)} (|F|) < \infty \eqno(2.4.20) $$

\noi and for every bounded previsible process $(\Phi_{s}, \; s \ge 0)$ : }
\begin{eqnarray*}
\hspace*{1cm}{\bf W}^{(2)} (\Phi_{g_{_{C}}} \cdot F)
&=& W^{(2)} \left( \int_{0}^{\infty} \Phi_{s} \, k_{s}^{(C)} (F) \, dL_{s}^{(C)}\right) \hspace*{4,8cm} (2.4.21) \\
&=& {\bf W}^{(2)} \big(\Phi_{g_{C}} \, k_{g_{C}}^{(C)} (F)\big) \hspace*{6,6cm} (2.4.22)
\end{eqnarray*}

\noi {\it Thus :}
$$
{\bf W}^{(2)} \big(F | \mathcal{F}_{g_{C}}\big) = k_{g_{C}}^{(C)} (F) \eqno(2.4.23) $$

 \noi {\it ii)} $\dis \big(k_{s}^{(C)} (k_{g_{C}}^{(C)} (F)), \; s \ge 0\big) = \big(k_{s}^{(C)} (F), \; s \ge 0\big) \hfill (2.4.24)$

\noi {\it iii) If $(h_{s}, \; s \ge 0)$ is a previsible process such that : } ${\bf W}^{(2)} \big(|h_{g_{_{C}}}|\big) < \infty$,
$$
\big(k_{s}^{(C)} (h_{g_{_{C}}}), \; s \ge 0\big) = (h_{s}, \; s \ge 0) \quad dL_{s}^{(C)} \cdot W^{(2)} (d \omega) \; {\rm a.s.}  \eqno(2.4.25) $$

\noi {\bf 2)} {\it There exist two continuous quasimartingales $\dis\big(\Sigma_{t}^{(2,C)}, \; t \ge 0\big)$ and $(\Delta_{t}^{(2,C)}, \; t \ge 0)$ such that, for every $t \ge 0$ :}
$$     M_{t}^{(2)} (F) = \Sigma_{t}^{(2, C)} (F) + \Delta_{t}^{(2,C)} (F)  \eqno(2.4.26) $$

\noi {\it with :}

\noi {\it i) For every $t \ge 0$ and $\Gamma_{t} \in b (\mathcal{F}_{t})$ :
\begin{gather*}
\hspace*{1cm}{\bf W}^{(2)} (\Gamma_{t} \, 1_{g_{C} \le t} \cdot F) = W^{(2)} \big(\Gamma_{t} \;\Sigma_{t}^{(2,C)} (F)\big)     \hspace*{5,5cm} (2.4.27) \\
\hspace*{1cm}{\bf W}^{(2)} (\Gamma_{t} \, 1_{g_{C} > t} \cdot F) = W^{(2)} \big(\Gamma_{t}\; \Delta_{t}^{(2,C)} (F)\big)  \hspace*{5,5cm} (2.4.28)
\end{gather*}

\noi In particular, from (2.4.27) applied with $\widetilde{\Gamma}_{t} = \Gamma_{t} \, 1_{|X_{t}| \le 1}$ and since $1_{g_{C} \le t} \cdot 1_{|X_{t}| \le 1} =0$, the process $\big(\Sigma_{t}^{(2,C)} (F), t \geq 0 \big)$ vanishes on the set $\big(|X_{t}| \le 1\big)$.

\noi ii) The Doob-Meyer decompositions of $\Sigma_{t}^{(2,C)} (F)$ and $\Delta_{t}^{(2,C)} (F)$ write :
\begin{eqnarray*}
\hspace*{1cm} \dis \Sigma_{t}^{(2,C)} (F) &=& - M_{t}^{\Sigma^{(2,C)}} (F) + \int_{0}^t k_{s}^{(C)} (F) dL_{s}^{(C)} \hspace*{5,3cm} (2.4.29) \\
\hspace*{1cm} \dis \Delta_{t}^{(2,C)} (F) &=& M_{t}^{\Delta^{(2,C)}} (F) - \int_{0}^t k_{s}^{(C)} (F) dL_{s}^{(C)} \hspace*{5,5cm} (2.4.30)
\end{eqnarray*}

\noi where $\big(M_{t}^{\Sigma^{(2,C)}} (F), \; t \ge 0\big)$ and $\big(M_{t}^{\Delta^{(2,C)}} (F), \; t \ge 0\big)$ are the martingale parts of the corresponding left-hand sides. The first martingale is not uniformly integrable ; the second one is uniformly integrable. In fact, we have : 
$$
M_{t}^{\Delta^{(2,C)}} (F) = W^{(2)} \left(\int_{0}^{\infty} k_{s}^{(C)} (F) dL_{s}^{(C)}| \mathcal{F}_{t}\right)
\eqno(2.4.31) $$

\noi with, from (2.4.20), $\dis \int_{0}^{\infty} k_{s}^{(C)} (F) dL_{s}^{(C)} \in L^{1} (\mathcal{F}_{\infty}, {\bf W}^{(2)})$. }

\smallskip

\noi {\it iii) The "explicit formula" :
$$
\Sigma_{t}^{(2,C)} (F) = \frac{1}{\pi} \; \log^{+} (|X_{t}|) \cdot \widehat{ \widetilde{E}}^{(2, \log)}_{X_{t} (\omega_{t})} \big(F (\omega_{t}, \widehat{\omega}^{t})\big)  \eqno(2.4.32) $$

\noi holds, where in (2.4.32) the expectation is taken with respect to $\widehat{\omega}^{t}$, and the argument $\omega_{t}$ is frozen. $\widetilde{E}^{(2, \log)}$ denotes the expectation with respect to the law $\;\widetilde{P}^{(2, \log)}\;$ defined in Theorem 2.2.2. \hbox{In particular :} 

\smallskip

\noi $\bullet$ $\Sigma_{t}^{(2,C)}$ vanishes on $\{t \;;\; |X_{t}| \le 1\}$, as we already observed,

\noi $\dis \bullet \quad \pi \; \frac{\Sigma_{t}^{(2,C)} (F)}{1 +\log^{+} (|X_{t}|)} \mathop{\longrightarrow}_{t \to \infty}^{} F \quad {\bf W}^{(2)}$  a.s. $\hfill(2.4.33)$

\smallskip

\noi and, from (2.4.16)

\smallskip

 $ \dis \pi \; \frac{\Delta_{t}^{(2,C)} (F)}{1 +\log^{+} (|X_{t}|)} \mathop{\longrightarrow}_{t \to \infty}^{} 0 \quad {\bf W}^{(2)}$  a.s. $\hfill(2.4.34)$  }
 
 \smallskip

\noi {\bf Corollary 2.4.5.} {\it Let $ F \in L^{1} (\mathcal{F}_{\infty}, {\bf W}^{(2)})$.

\noi One has $ M_{t}^{(2)} (F) = 0$ for every $t \ge 0$ such that $|X_{t}| \le 1 $, if and only if :    }

$$k_{g_{_{C}}}^{(C)} (F) = 0$$

\newpage

\noi { \bf{\Large Chapter 3. The analogue of the measure {\bf W} for a class of linear diffusions.}}

\bigskip

\noi In Chapters 1 and 2, we have, starting from penalisation results, associated to Wiener measure in dimensions 1 and 2 a positive and $\sigma$-finite measure ${\bf W}$ (resp. : ${\bf W}^{(2)}$ in dimension 2) on the canonical space $(\Omega, \mathcal{F}_{\infty})$. In this 3$^{\rm rd}$ Chapter, we shall prove the existence of a measure which is analogous to ${\bf W}$, in the more general situation of a large class of linear diffusions. This class is described in Section 3.2. Our approach in this Chapter does not use any penalisation result. Then, in Section 3.3, we shall particularize these results about linear diffusions to the situation of Bessel processes with dimension $d=2(1-\alpha)\; (0<d<2$, or $0< \alpha<1)$. Thus, we shall obtain the existence of the measure ${\bf W}^{(-\alpha)} (0<\alpha <1)$ on $\big(\mathcal{C} (\mathbb{R} \to \mathbb{R}_{+}), \; \mathcal{F}_{\infty}\big)$ and we shall then indicate its relationship with penalisation problems. Section 3.1 is devoted to a presentation of our hypotheses and notations. 

\smallskip

{\large\bf 3.1 Main hypotheses and notations.}

\smallskip

\noi \underline{3.1.1 Our framework is that of Salminen-Vallois-Yor.} [SVY], that is : 

\noi $(X_{t}, \; t \ge 0)$ is a $\mathbb{R}_{+} = [0, \infty[$ valued diffusion, with 0 an instantaneously reflecting barrier. The infinitesimal generator $\mathcal{G}$ of $(X_{t}, \; t \ge 0)$ is given by :
$$      \mathcal{G} \, f(x) = \frac{d}{dm} \; \frac{d}{dS} \; f(x) \qquad (x \ge 0)  \eqno(3.1.1)   $$

\noi where the scale function $S$ is a continuous, strictly increasing function s.t. : 
$$         S(0)=0, \quad S(+�\infty)= +\infty  \eqno(3.1.2) $$

\noi and $m(dx)$ is the speed measure of $X$ ; we assume $m(\{0\})=0$.

\smallskip

\noi {\bf 3.1.2} \underline{The semi-group of $(X_{t}, \; t \ge 0)$} admits $p(t,x,y)$ as density with respect to $m$ :
$$         P_{x} (X_{t} \in dy) =p(t,x,y) m (dy)  \eqno(3.1.3) $$

\noi with $p$ continuous in the 3 variables, and $p(t,x,y)=p(t,y,x)$. $\widehat{X} $ denotes the process $X$, killed at $T_{0} = \inf \{t \;;\; X_{t}=0\}$. We denote by $\widehat{p}$ its density with respect to $m$ :
$$
\widehat{P}_{x} (\widehat{X}_{t} \in dy) = P_{x} (X_{t} \in dy \;;\; 1_{t < T_{0}}) := \widehat{p} (t,x,y) m (dy)
\eqno(3.1.4)   $$ 
\noi with $\widehat{p} (t,x,y) = p(t,x,y) P_x (T_0 > t | X_t = y)$. 

\smallskip

\noi {\bf 3.1.3} \underline{The local time process} \\
\noi  We denote by $\{L_{t}^{y} \;;\; t \ge 0, \, y \ge 0\}$ the jointly continuous family of local times of $X$, which satisfy the density of occupation formula :
$$          \int_{0}^{t} h(X_{s}) ds = \int_{0}^{\infty} h(y) L_{t}^{y} \, m (dy)  \eqno(3.1.5)   $$

\noi for any $h: \mathbb{R}_{+} \to \mathbb{R}_{+}$, Borel. It is easily deduced from (3.1.5) and (3.1.3) that :
$$        E_{x} (d_{t} \, L_{t}^{y}) = p(t,x,y)dt    \eqno(3.1.6)   $$

\noi We denote by $P_{0}^{\tau_{l}}$ the law, under $P_{0}$, of $(X_{t}, \; t \le \tau_{l})$ with $\tau_{l} := \inf \{ t \ge 0 \;;\; L_{t}^{0} > l\}$. We have also :
$$          \big(S(X_{t})-L_{t}, \; t \ge 0\big) \quad \hbox{is a martingale}   \eqno(3.1.7)  $$

\noi a property which results from (3.1.1) and (3.1.5) \big(see [DM, RVY] for such a property in the context of Bessel processes\big).

\smallskip

\noi {\bf 3.1.4} \underline{The process $X$, conditioned not to vanish}, is a Doob $h$-transform of $\widehat{X}$, with $h(x)=S(x)$. In other terms : if $P_{x}^{\uparrow}$ is the law of $X$ conditioned not to vanish :
$$
P_{x}^{\uparrow} (F_{t}) = \frac{1}{S(x)} \; E_{x} [F_{t} \, S(X_{t}) 1_{t<T_{0}}]    \eqno(3.1.8)  $$

\noi for any $F_{t} \in b(\mathcal{F}_{t})$. In particular, the semi-group of the conditioned process is given by :
$$
P_{x}^{\uparrow} (X_{t} \in dy) = \frac{ \widehat{p}(t,x,y)}{S(x)} \; S(y) \, m(dy) \qquad (x \ge 0)  \eqno(3.1.9) $$
 
\noi Later, it will be interesting to use the following :
$$          P_{0}^{\uparrow} (X_{t} \in dy) = f_{y,0} (t) \; S(y) \, m(dy)  \eqno(3.1.10)  $$

\noi where $f_{y,0}(t)$ admits the following description :
\begin{gather*} 
\hspace*{1cm} \dis f_{y,0} (t) = \mathop{\rm lim}_{x \downarrow 0}^{} \; \frac{\widehat{p}(t,x,y)}{S(x)} \hspace*{8,5cm}(3.1.11) \\
\hspace*{1cm}\dis f_{y,0} (t) dt \mathop{=}_{(a)}^{} P_{y} (T_{0} \in dt) \mathop{=}_{(b)}^{} P_{0}^{\uparrow} (g_{y} \in dt) \hspace*{6cm}(3.1.12)
\end{gather*}

\noi with 
\begin{equation*}
g_{y} := \sup \{t \;;\; X_{t} = y\}
\end{equation*}

\noi We indicate here that (3.1.12) is a partial expression of the time reversal result :
$$
P_{y} \big(\{X_{T_{0}-t}, \; t \le T_{0} \}\big)= P_{0}^{\uparrow} \big(\{X_{u}, \; u \le g_{y}\}\big) \eqno(3.1.13) $$

\noi Furthermore :
$$
P_{0}^{\uparrow} \big( \{ X_{u}, \; u\le g_{y}\} \big| g_{y} =t\big) = P_{0}^{\uparrow} \big(\{X_{u}, \; u \le t\} \big| X_{t} =y\big)  \eqno(3.1.14) $$

\noi where in (3.1.13) and in (3.1.14) we have used the notation $P\big(\{X_{u}, \; u \le a\}\big)$ to denote the law of the process $(X_{u}, \; u \le a)$ under $P$. All these facts, as well as those presented in the following Proposition may be found in [SVY], [BS], [PY2], ... which all deal with properties of linear diffusions.

\smallskip

\noi {\bf 3.1.5} \underline{A useful Proposition} :

\noi We shall use the following :

\noi {\bf Proposition 3.1.1.} {\it Let $g^{(t)} := \sup \{s \le t, \, X_{s} =0\}$.

\noi {\bf 1)} Under $P_{0}$, conditionally on $g^{(t)}$, the processes $(X_{u}, \; u \le g^{(t)})$ and $(X_{g^{(t)} +u}, \; u \le t-g^{(t)})$ are independent.

\noi {\bf 2)} Conditionally on $g^{(t)} =s, \; (s \le t)$, the process $(X_{u}, \; u \le s)$ is distributed as $\Pi_{0}^{(s)}$, the law of the bridge of $X$ under $P_{0}$, with length $s$, ending at $x=0$ at time $s$.

\noi {\bf 3)} The law of the couple $\big(g^{(t)}, (X_{g^{(t)}+u}, \, u\le t-g^{(t)})\big)$ under $P_{0}$ may be described as follows :

{\it i)} $\dis  \qquad P_{0} (g^{(t)} \in ds, \; X_{t} \in dy) \, S(y) = p(s,0,0) \, 1_{s<t} P_{0}^{\uparrow} (X_{t-s} \in dy) ds  \hfill(3.1.15) $

\noi or equivalently, with the help of (3.1.10) :

{\it i')} $\dis  \qquad P_{0} (g^{(t)} \in ds, \; X_{t} \in dy) = p(s,0,0) \, f_{y,0} (t-s) \, 1_{s<t} ds \, m(dy) \hfill (3.1.16) $

\noi and, on the other hand :

{\it ii)} $\dis \qquad P_{0} \big(\{X_{g_{t}+u}, \; u \le t-g_{t} \} | X_{t}=y, \; g_{t} =s \big) $

$ \dis \hspace*{2cm} = P_{0}^{\uparrow} \big(\{X_{u}, \; u \le t-s\} | X_{t-s} =y \big) \hfill(3.1.17) $

}
\noi These different properties are established in [SVY], [Sa1] and [Sa2].     

\bigskip

  {\bf \large 3.2 The $\sigma$-finite measure ${\bf W}^{*}$.}
 
 \smallskip
 
 \noi {\bf 3.2.1} \underline{Definition of ${\bf W}^{*}$} :

\noi Here is the main result of this Section.

\noi {\bf Theorem 3.2.1.}

\noi {\bf 1)} {\it There exists a unique $\sigma$-finite measure, which we denote by ${\bf W}^{*}$, on $\big( \mathcal{C} \big(\mathbb{R}_{+} \to \mathbb{R}_{+}), \; \mathcal{F}_{\infty}\big)$ such that :
\begin{gather*}
\hspace*{1cm}\forall t \ge 0, \; \forall F_{t} \in b(\mathcal{F}_{t}) \, :  \\
\hspace*{2,4cm} E_{0} \big(F_{t} \, S(X_{t})\big) = {\bf W}^{*} (F_{t} \, 1_{g \le t}) \hspace*{6,3cm}(3.2.1)
\end{gather*}

\noi with $g:= \sup \{ t \ge 0 \;;\; X_{t}=0\}$

\noi {\bf 2)} $\dis \qquad {\bf W}^{*}= \int_{0}^{\infty} dl (P_{0}^{\tau_{l}} \circ P_{0}^{\uparrow}) \hfill(3.2.2)$

\smallskip

\noi {\bf 3)} $\dis \qquad {\bf W}^{*}= \int_{0}^{\infty} dt \, p(t,0,0) \, (\Pi_{0}^{(t)} \circ P_{0}^{\uparrow}) \hfill(3.2.3)$

 \noi In particular, if we denote ${\bf W}_{g}^{*}$ the restriction of ${\bf W}^{*}$ to $\mathcal{F}_{g}$, we have : 
$$
 {\bf W}_{g}^{*} = \int_{0}^{\infty} dl \, P_{0}^{\tau_{l}} = \int_{0}^{\infty} dt \, p(t,0,0) \, \Pi_{0}^{(t)} \eqno(3.2.4) $$      }
 
\noi Of course, this Theorem 3.2.1. has been guessed from the comparison with the Brownian situation described in Chapters 1 and 2. 

\smallskip   

\noi {\bf Proof of Theorem 3.2.1.}

\noi {\it i)} First of all, it is not difficult to show that, starting from equation (3.2.1), where ${\bf W}^{*}$ is the unknown, this problem admits at most one solution such that $g < \infty$,  ${\bf W}^{*}$ a.s.

\smallskip

\noi {\it ii)} Define 
$$
{\bf W}_{*} = \int_{0}^{\infty} dl \, (P_{0}^{\tau_{l}} \circ P_{0}^{\uparrow})  \eqno(3.2.5) $$

\noi We shall now prove that ${\bf W}_{*}$ satisfies (3.2.3) and (3.2.4). Since, under $P_{0}^{\uparrow}$, the process $(X_{t}, \; t \ge 0)$ remains in $\mathbb{R}_{+} \setminus \{0\}$, it follows immediately, from the definition (3.2.5) of ${\bf W}_{*}$ that
$$         {\bf W}_{*,g} = \int_{0}^{\infty} dl \, P_{0}^{\tau_{l}}    \eqno(3.2.6) $$
\noi where ${\bf W}_{*,g}$ denotes the restriction of ${\bf W}_{*}$ to $\mathcal{F}_g$. 

\noi On the other hand, a classical argument, which hinges on the fact that the random measure $(dL_{s})$ is carried by the zeros of $X$, allows to show easily that :
$$
\int_{0}^{\infty} dl \, P_{0}^{\tau_{l}} = \int_{0}^{\infty} dt \, p (t,0,0) \, \Pi_{0}^{(t)}  \eqno(3.2.7) $$

\noi Indeed, by integrating $F := \big(F_{t} := F (X_{u}, \; u \le t), \; t \ge 0\big)$ a positive measurable functional, we obtain on the LHS of (3.2.6)
$$
\int_{0}^{\infty} dl \, P_{0}^{\tau_{l}} (F) = \int_{0}^{\infty} dl \, P_{0} (F_{\tau_{l}}) = P_{0} \left( \int_{0}^{\infty} dL_{s} \cdot F_{s}\right)      $$

\noi (by time change $l=L_{s}$)
\begin{gather*}
= P_{0} \left(\int_{0}^{\infty} d L_{s} \, P_{0} (F_{s} | X_{s}=0) \right) \\
= \int_{0}^{\infty} P_{0} (dL_{s}) P_{0} (F_{s} | X_{s}=0) \\
= \int_{0}^{\infty} dt \, p(t,0,0) \, \Pi_{0}^{(t)} (F)
\end{gather*}

\noi by (3.1.6), with $x=y=0$.

\smallskip

\noi {\it ii)} We now prove that ${\bf W}_{*}$ satisfies (3.2.1), by showing this equality for the test functionals :
$$
F_{t} = \Phi (X_{u}, \; u \le g^{(t)}) \, \varphi (g^{(t)}) \, \psi (X_{g^{(t)}+u}, \; u \le t-g^{(t)}) \eqno(3.2.8) $$

\noi From (3.2.3), the RHS of (3.2.1) is equal to (with ${\bf W}_{*}$ instead of ${\bf W}^{*}$) :
\begin{eqnarray*}
\hspace*{1cm}R^{F} 
&:=& {\bf W}_{*} (F_{t} \, 1_{g \le t}) \\
&=& \int_{0}^{t} ds \, p(s,0,0) \, \Pi_{0}^{(s)} \big(\Phi (X_{u}, \; u \le s)\big) \varphi (s) \, P_{0}^{\uparrow} \big(\psi (X_{u}, \; u \le t-s)\big) \hspace*{1cm} (3.2.9) 
\end{eqnarray*}

\noi On the other hand, the LHS of (3.2.1) is equal to :
\begin{eqnarray*}
\hspace*{0,5cm}L^{F}
&:=& E_{0} \big[F_{t} \, S (X_{t})\big]  \\
&=& E_{0} \big[\Phi (X_{u}, \; u \le g^{(t)}) \, \varphi (g^{(t)}) \, \psi (X_{g^{(t)}+u}, \; u \le t-g^{(t)}) \, S(X_{t}) \big]  \\
&=& \int_{0}^{t} P_{0} (g^{(t)} \in ds) \varphi (s) \, E_{0} \big[\Phi (X_{u}, \; u \le s) |X_{s} =0\big] \, E_{0} \big[\psi (X_{s+u}, \; u \le t-s) \\
&& \hspace*{8cm} S(X_{t})�|g^{(t)}=s\big] \hspace*{1cm}(3.2.10)
\end{eqnarray*}

\noi where we have used a part of the results presented in the Proposition 3.1.1. Comparing (3.2.9) and (3.2.10), we now see that showing equality $R^{F}=L^{F}$ (i.e. the proof of (3.2.1)) has now been reduced to showing :
\begin{eqnarray*}
\hspace*{1cm}&&P_{0}^{\uparrow} \big(\psi (X_{u}, \; u \le t-s)\big) \, p(s,0,0) \, 1_{s < t}  ds\\
&=& P_{0} (g^{(t)} \in ds) \, E_{0} \big(\psi (X_{s+u} \;;\; u \le t-s) \cdot S(X_{t}) | g^{(t)}=s\big) \hspace*{3cm}(3.2.11)
\end{eqnarray*}

\noi But (3.2.11) is an easy consequence of point 3 of Proposition 3.1.1.

\bigskip

\noi {\bf 3.2.2} \underline{Some properties of ${\bf W}^{*}$.}

\noi The end of this subsection 3.2.2 is devoted to the statement of some results related to the measure ${\bf W}^{*}$. These results are presented without proofs since those are close to the ones found in Chapter 1. These theorems (below) are due to Christophe Profeta ([Pr], thesis in preparation).

\smallskip

\noi {\bf 3.2.2.1} \underline{The probabilities $P_{x, \infty}^{(\lambda)}$.}

\noi {\bf Theorem 3.2.2.}

\noi {\bf 1)} {\it Let, for $\lambda \ge 0$ and $x \ge 0$ :
$$
M_{t}^{(\lambda, x)} := \frac{ 1+ \frac{\lambda}{2} \, S (X_{t})}{1 + \frac{\lambda}{2} \,S(x)} \; e^{- \frac{\lambda}{2} \, L_{t}} = 1 + \frac{\lambda}{2 + \lambda S(x)} \int_{0}^{t} e^{- \frac{\lambda}{2} \, L_{s}} dN_{s} \eqno (3.2.12)$$

\noi where $\big(N_{s} := S (X_{s}) - L_{s}, \; s \ge 0\big)$ is the martingale defined by (3.1.7). Then, $(M_{t}^{(\lambda, x)}, \; t \ge 0)$ is a $\big((\mathcal{F}_{t}, \; t \ge 0), \; P_{x}\big)$ positive martingale such that : $\dis M_{t}^{(\lambda, x)} \mathop{\longrightarrow}_{t \to \infty} 0$, a.s.

\smallskip

\noi {\bf 2)} Let us define the probability $P_{x, \infty}^{(\lambda)}$ by :
$$
P_{x, \infty}^{(\lambda)} \big|_{\mathcal{F}_{t}} = M_{t}^{(\lambda, x)} \cdot P_{x}\big|_{\mathcal{F}_{t}} \eqno (3.2.13) $$

\noi Then, under $P_{x, \infty}^{(\lambda)}$ :
\begin{itemize}
\item[$\bullet$] The canonical process $(X_{t}, \; t \ge 0)$ is a transient diffusion with infinitesimal generator $\mathcal{G}_{\infty}^{(\lambda)}$ :
\begin{eqnarray*} 
\mathcal{G}_{\infty}^{(\lambda)} f(x)
&=& \frac{2}{2+ \lambda S(x)} \; \left( \mathcal{G}f (x) + \frac{\lambda}{2} \, \mathcal{G}(Sf) (x)\right) \\
&=& \mathcal{G}f(x) + \frac{2 \lambda}{2 + \lambda S(x)} \; \frac{df}{dm} \, (x) \hspace*{5,4cm} (3.2.14)
\end{eqnarray*}

and scale function $S_{\infty}^{(\lambda)}$ :
$$     S_{\infty}^{(\lambda)} := - \frac{2}{2+ \lambda S} \eqno(3.2.15) $$
\item[$\bullet$] If $\alpha < \lambda$ :
$$ E_{x, \infty}^{(\lambda)} \left(e^{\frac{\alpha}{2} \, L_{\infty}}\right) < \infty \eqno(3.2.16)$$

\noi and if $\alpha \ge \lambda$ :
$$ E_{x, \infty}^{(\lambda)} \left(e^{\frac{\alpha}{2} \, L_{\infty}}\right) = \infty \eqno(3.2.17)$$
\item[$\bullet$] The law of $L_{\infty}$ is given by :
$$  P_{x, \infty}^{(\lambda)} (L_{\infty} \in dl) = \frac{\lambda}{2 + \lambda S(x)} \; e^{- \frac{\lambda}{2} \, l} dl + \frac{\lambda S(x)}{2 + \lambda S(x)}  \delta_{0} (dl) \eqno (3.2.18) $$
\item[$\bullet$] $P_{x, \infty}^{(\lambda)}$ admits the following decomposition :
\begin{eqnarray*}
P_{x, \infty}^{(\lambda)}
&=& \frac{\lambda}{2 + \lambda S(x)} \int_{0}^{\infty} du \, p(u,x,0) e^{- \frac{\lambda}{2} \, L_{u}} \, . \, ( \Pi_{x,0}^{(u)} \circ P_{0}^{\uparrow} ) + \frac{\lambda S (x)}{2 + \lambda S(x)} \, P_{x}^{\uparrow} \hspace*{1,3cm} (3.2.19)\\
&=& \frac{\lambda}{2 + \lambda S(x)} \int_{0}^{\infty} e^{- \frac{\lambda l}{2}} dl \, (P_{x}^{\tau_{l}} \circ P_{0}^{\uparrow}) + \frac{ \lambda S(x)}{2+ \lambda S(x)} \, P_{x}^{\uparrow} \hspace*{3,4cm} (3.2.20)
\end{eqnarray*}
\end{itemize}  }

\noi {\bf 3.2.2.2} \underline{The measures $({\bf W}_{x}^{*}, \; x \in \mathbb{R}_{+})$.}

\noi {\bf Theorem 3.2.3.}

\noi {\bf 1)} {\it For any $\lambda >0$, the $\sigma$-finite measure $\dis \left(\frac{2}{\lambda} + S(x)\right) \cdot e^{\frac{\lambda}{2} \, L_{\infty}} \cdot P_{x, \infty}^{(\lambda)}$ does not depend on $\lambda$. We define :
$$ {\bf W}_{x}^{*} := \left(\frac{2}{\lambda} + S(x)\right) e^{\frac{\lambda}{2}\, L_{\infty}} \cdot P_{x, \infty}^{(\lambda)} \eqno (3.2.21)  $$

\noi We have the decompositions :
\begin{eqnarray*}
{\bf W}_{x}^{*} &=& \int_{0}^{\infty} du \, p(u,x,0) \, (\Pi_{x,0}^{(u)} \circ P_{0}^{\uparrow}) + S(x) P_{x}^{\uparrow} \hspace*{5,2cm} (3.2.22)  \\
{\rm and} \quad {\bf W}_{x}^{*} &=& \int_{0}^{\infty} dl \, (P_{x}^{\tau_{l}} \circ P_{0}^{\uparrow}) + S(x) P_{x}^{\uparrow} \hspace*{7cm} (3.2.23)
\end{eqnarray*}

\noi In particular, ${\bf W}^{*}_{0} = {\bf W}^{*}$, where ${\bf W}^{*}$ is defined by (3.2.2) or (3.2.3).

\smallskip

\noi {\bf 2)} i) For every $(\mathcal{F}_{t}, \; t \ge 0)$ stopping time $T$ and $\Gamma_{T} \in b(\mathcal{F}_{T})$ :
$$
E_{x} \big(\Gamma_{T} S (X_{T}) 1_{T < \infty}\big) = {\bf W}_{x}^{*} (\Gamma_{T} 1_{g \le T < \infty}) \eqno(3.2.24) $$

 where $g:= \sup \{ s \ge 0 \;;\; X_{s} = 0\}$
 
 ii) The law of $g$ under ${\bf W}_{x}^{*}$ is given by :
 $$
 {\bf W}_{x}^{*}  (g \in dt) = p(t,x,0) dt + S(x) \delta_{0} (dt) \quad (t�\ge 0) \eqno (3.2.25) $$
 
 and for every $(\mathcal{F}_{t}, \; t \ge 0)$ stopping time $T$, we have :
 $$
 {\bf W}_{x}^{*} (1_{T < \infty}, L_{\infty} -L_{T} \in dl) = P_{x} (T < \infty) \, 1_{[0, \infty[}(l) dl + E_{x} \big[S (X_{T}) 1_{T < \infty} \big] \delta_{0} (dl) \eqno (3.2.26) $$
 
 \noi {\bf 3)} For every previsible and positive process $(\Phi_{s}, \; s \ge 0)$, we have :
 $$
 {\bf W}_{x}^{*} (\Phi_{g}) = S(x) \Phi_{0} + E_{x} \left(\int_{0}^{\infty} \Phi_{s} d L_{s} \right) \eqno(3.2.27)$$  }
 
 \noi We note that, from (3.2.19), (3.2.20), (3.2.22) and (3.2.23), we have : 
$$ \underset{\lambda \rightarrow 0}{\lim} \, \frac{2}{\lambda} \, P_{x, \infty}^{(\lambda)} = {\bf W}_x^{*} \eqno(3.2.27')$$

 \noi {\bf 3.2.2.3} \underline{Martingales associated with $({\bf W}_{x}^{*}, \; x \in \mathbb{R}_{+}$).}
 
 \smallskip
 
 \noi {\bf Theorem 3.2.4.}
 
 \noi {\it Let $F \in L_{+}^{1} (\Omega, \mathcal{F}_{\infty}, {\bf W}_{x}^{*})$. There exists a positive $\big( (\mathcal{F}_{t}, \; t \ge 0), \, P_{x}\big)$ martingale $\big(M_{t}^{*} (F), \; t \ge 0\big)$ such that : 
 
 \smallskip
 
 \noi {\bf 1)} For every $t \ge 0$ and $\Gamma_{t} \in b (\mathcal{F}_{t})$ :
 $$
 {\bf W}_{x}^{*}  (F \cdot \Gamma_{t}) = E_{x} \big(M_{t}^{*} (F) \Gamma_{t}\big) \eqno(3.2.28) $$
 
 In particular, ${\bf W}_{x}^{*} (F) = E_{x} \big(M_{t}^{*} (F)\big)$
 
 \noi {\bf 2)} For every $\lambda > 0$ : 
 \begin{eqnarray*}
 M_{t}^{*} (F)
 &=& \left( \frac{2}{\lambda} + S (X_{t})\right) e^{- \frac{\lambda}{2}\, L_{t}} E_{x, \infty}^{(\lambda)} \big(F e^{\frac{\lambda}{2} \, L_{\infty}} | \mathcal{F}_{t}\big) \hspace*{5cm} (3.2.29)  \\
 &=& \widehat{\bf W}_{X_{t}}^{*} \big( F(\omega_{t}, \widehat{\omega}^{t})\big)
 \end{eqnarray*}
 
 \noi {\bf 3)} $M_{t}^{*} (F) \dis \mathop{\longrightarrow}_{t \to \infty} 0 \qquad P_{x} \qquad {\rm a.s.} \hfill (3.2.30) $  }
 
 \smallskip
 
 \noi {\bf Examples} : 
 
 \smallskip
 
 \noi $\bullet$ Let $h : \mathbb{R}_{+} \to \mathbb{R}_{+} $ such that $\dis \int_{0}^{\infty} h(u) du < \infty$. Then :
 $$
 M_{t}^{*} \big(h(L_{\infty})\big) = h(L_{t}) S(X_{t}) + \int_{L_{t}}^{\infty} h(l) dl \eqno (3.2.31) $$
 
In particular, if $h(y) = e^{- \frac{\lambda}{2} \,y} \; (y \ge 0)$ : 
 $$
 M_{t}^{*} (e^{- \frac{\lambda}{2} \, L_{\infty}}) = \left( \frac{2}{\lambda} + S(X_{t})\right) e^{- \frac{\lambda}{2} \, L_{t}} = \frac{2}{\lambda} \; M_{t}^{(\lambda, 0)} \qquad (x=0) $$
 
 \noi $\bullet$ Let $\Phi : \mathbb{R}_{+} \to \mathbb{R}_{+}$ Borel such that $\dis \int_{0}^{\infty} \Phi (u) p(u,x,0) du < \infty$. Then
 $$
 M_{t}^{*} \big(\Phi (g)\big) = \Phi (g^{(t)}) S(X_{t}) + \int_{0}^{\infty} \Phi (t+u) p(u,X_{t},0) du \eqno(3.2.32) $$
 
 \noi {\bf 3.2.2.4} \underline{A decomposition Theorem of $\big((\mathcal{F}_{t}, \; t \ge 0), \, P_{x}\big)$ positive supermartingales.}
 
 \smallskip
 
 \noi {\bf Theorem 3.2.5.}
 
 \noi {\it Let $(Z_{t} , \; t \ge 0)$ a positive $\big((\mathcal{F}_{t}, \; t \ge 0), \, P_{x}\big)$ supermartingale. We denote
 $$
 Z_{\infty} := \mathop{\lim}_{t \to \infty} Z_{t} \qquad P_{x} \qquad {\rm a.s.}  $$
 
 \noi Then :

 \noi {\bf 1)} $z_{\infty} := \dis \mathop{\lim}_{t \to \infty} \, \frac{Z_{t}}{1+S(X_{t})}$ exists ${\bf W}_{x}^{*}$ a.s. and ${\bf W}_{x}^{*} (z_{\infty}) < \infty$
 
 \smallskip
 
 \noi {\bf 2)} $(Z_{t}, \; t \ge 0)$ admits the following decomposition :
 $$
 Z_{t} = M_{t}^{*} (z_{\infty}) + E_{x} \big(Z_{\infty} | \mathcal{F}_{t}\big) + \xi_{t} \eqno(3.2.33) $$
 
 \noi where $\big(M_{t}^{*} (z_{\infty}), \; t \ge 0\big)$ and $\big( E_{x} \big(Z_{\infty} | \mathcal{F}_{t}\big), \, t \geq 0 \big)$ denote two positive $\big((\mathcal{F}_{t}, \; t \ge 0), \; P_{x}\big)$ martingales and $(\xi_{t}, \; t \ge 0)$ is a positive supermartingale such that :
 
 \noi $\bullet$ $Z_{\infty} \in L^{1}_+ (\mathcal{F}_{\infty}, P_{x})$, hence $\big(E_{x} (Z_{\infty} | \mathcal{F}_{t}), \; t \ge 0\big)$ is a uniformly integrable martingale converging towards $Z_{\infty}$.
 
 \smallskip
 
 \noi $\bullet$ $\dis \frac{E_{x} \big(Z_{\infty} | \mathcal{F}_{t}\big) + \xi_{t}}{1 + S(X_{t})} \mathop{\longrightarrow}_{t \to \infty} 0 \quad {\bf W}_{x}^{*}$ a.s.
 
 \smallskip
 
 \noi $\bullet$ $\dis M_{t}^{*} (z_{\infty}) + \xi_{t} \mathop{\longrightarrow}_{t \to \infty} 0 \quad P_{x}$ a.s.

\noi This decomposition (3.2.33) is unique. }

\smallskip

\noi {\bf Corollary 3.2.6.}

\noi {\it A positive martingale $(Z_{t}, \; t \ge 0)$ is equal to $\big(M_{t}^{*} (F), \; t \ge 0\big)$ for some $F \in L_{+}^{1} (\mathcal{F}_{\infty}, {\bf W}_{x}^{*})$ if and only if :
$$
Z_{0} = {\bf W}_{x}^{*} \left( \mathop{\lim}_{t \to \infty} \frac{Z_{t}}{1+S(X_{t})}\right) \eqno(3.2.34)  $$   }

\noi In the present framework of linear diffusions, it is possible to state a decomposition theorem for the martingales $\big(M_{t}^{*} (F), \; t \ge 0\big)$ $\big(F \in L^{1} (\mathcal{F}_{\infty}, {\bf W}_{x}^{*})\big)$ which is similar to the result stated in Theorem 1.2.11. We leave this task to the interested reader.

\smallskip

\noi {\bf 3.2.3} \underline{Relation between the measure ${\bf W}^{*}$ and penalisations.}

\noi In a submitted article, P. Salminen and P. Vallois (see [SV]) obtain the following result involving, as
weight functional, the local time of a diffusion, under a certain subexponentiality hypothesis. We now 
summarize their results.

\noi Let $(\tau_{l}, \; l \ge 0)$ denote the right continuous inverse of the local time process $(L_{t}, \; t \ge 0)$ at level 0 associated to $(X_{t}, \; t \ge 0)$ : 
$$ \tau_{l} := \inf \{ t \ge 0 \;;\; L_{t} > l\}     $$

\noi This subordinator $(\tau_{l}, \; l \ge 0)$ admits as its Levy measure a mesure $\nu$ with density, which we denote here by $\dis \mathop{\nu}^{\bullet}$ \big(see [KS]\big) :
$$  E(e^{- \lambda \tau_{l}}) = \exp \left\{ -l \int_{0}^{\infty} (1-e^{- \lambda x}) \mathop{\nu}^{\bullet} (x) dx\right\} \qquad (\lambda, l \ge 0) $$

\noi P. Salminen and P. Vallois then make the following hypothesis : the function $F : [1, \infty[ \to [0,1]$ defined by :
$$ F(x) := \frac{ \nu \big(]1,x[\big)}{\nu \big(]1, \infty[\big)} = \frac{\dis \int_{1}^{x} \mathop{\nu}^{\bullet} (y)dy} {\dis \int_{1}^{\infty} \mathop{\nu}^{\bullet} (y)dy} \eqno (3.2.35)  $$

\noi is sub-exponential \footnote{This notion has little to do with the sub-exponential functions, i.e. functions $f : \mathbb{R}_+ \to \mathbb{R}_+$ which satisfy : $f(x) \leq c_1 e^{-c_2 x}$ for some constants $c_1, c_2 > 0$, and which are considered in [RY, IX].}, i.e. :
$$  \mathop{\rm lim}_{x \to \infty} \; \frac{\overline{F * F}(x)}{\overline{F}(x)} = 2 \eqno (3.2.36)$$

\noi where $\overline{F} (x) := 1-F(x), \; x \ge 1$ and where $*$ indicates the convolution operation.

\noi One of the main consequences of the subexponentiality of $F$ is :
$$ \frac{\overline{F} (x+y)}{\overline{F} (x)} \mathop{\longrightarrow}_{x \to \infty} 1 \qquad \hbox{uniformly on compacts (in} \; y) $$

\noi Thus, here
$$ \frac{\nu \big(]x+y,  \infty[\big)}{\nu \big(]x, \infty[\big)} \mathop{\longrightarrow}_{x \to \infty} 1 \qquad 
\hbox{uniformly on compacts (in} \; y) \eqno (3.2.37)$$

\noi Under this subexponentiality hypothesis, P. Salminen and P. Vallois [SV] then prove the following Theorem.

\smallskip

\noi {\bf Theorem 3.2.7.} {\it (Penalisation by $(1_{(L_{t} < l)}, \; t \ge 0)$

\noi Let $l>0$ be fixed. Then, for every $s \ge 0$ and $\Gamma_{s} \in b (\mathcal{F}_{s})$ :
$$ \mathop{\lim}_{t \to \infty} \frac{E_{x} (\Gamma_{s} \, 1_{(L_{t} <l)})}{P_{x} (L_{t} <l)}  = E_{x} (\Gamma_{s} \cdot M_{s}^{(l)}) := P_{x, \infty}^{(l)} (\Gamma_{s}) \eqno (3.2.38)$$

\noi where $(M_{s}^{(l)}, \; s \ge 0)$ is the positive martingale defined by : 
$$ M_{s}^{(l)} := \frac{S(X_{s})-L_{s}+l}{S(x)+l} \cdot 1_{L_{s} < l} $$   }

\noi Let us remark that for $f : \mathbb{R}_+ \rightarrow \mathbb{R}_+$ such that 
$$ \int_0^{\infty} \left(1 + \frac{1}{l} \right) \, f(l) \, dl < \infty.$$
\noi Then, we have :

\begin{align*}
 M_s^{(f)} & := \int_{0}^{\infty} M_s^{(l)} \, f(l)\, dl \\ & = (S(X_s) -L_s) \int_{L_s}^{\infty} 
\frac{f(l) dl}{S(x)+l} \, + \int_{L_s}^{\infty}  \frac{f(l) l \, dl}{S(x) + l}
\end{align*}

\noi and, for $x=0$, 
\begin{align*}
M_s^{(f)} & = (S(X_s) - L_s) \int_{L_s}^{\infty} \frac{f(l)}{l} \, dl + \int_{L_s}^{\infty} f(l) \, dl \\
& = S(X_s) h(L_s) + \int_{L_s}^{\infty} h(y) dy
\end{align*}
with $$ h(y) := \int_{y}^{\infty} \frac{f(l)}{l} \, dl.$$

\noi Thus, $(M_s^{(f)}, s \geq 0)$ is the Az\'ema-Yor martingale associated to $h$ (see [AY1]). 

\noi The key point of the proof of Theorem 3.2.7 is the following

\smallskip

\noi {\bf Lemma 3.2.8.} \big([SV]\big)
$$ P_{x} (L_{t} < l) \mathop{\sim}_{t \to \infty} \big(S(x)+l\big) \, \nu \big(]t, \infty[\big) \eqno(3.2.39) $$

\noi Theorem 3.2.7 now follows easily from Lemma 3.2.8 and from relation (3.2.37). 

\noi From this Theorem 3.2.7, we deduce the following relation between the probability $P_{0, \infty}^{(l)}$ defined by (3.2.38) and the $\sigma$-finite measure ${\bf W}^{*}$ defined by (3.2.2) or (3.2.3) :
$$ 1_{L_{\infty} <l} \cdot {\bf W}^{*} = {\bf W}^{*} (L_{\infty} < l) \cdot P_{0, \infty}^{(l)} \eqno(3.2.39) $$

\noi \big(We note that $P_{0, \infty}^{(l)} (L_{\infty} < l)=1$\big). The reader may compare relation (3.2.39) with relation (1.1.107) of Theorem 1.1.11' and with relation (3.2.21) of Theorem 3.2.3. From (3.2.39), we also deduce, with the notation of Theorem 3.2.4, that :
$$ M_{t}^{*} (1_{(L_{\infty} <l)}) = {\bf W}^{*} (L_{\infty} <l) \cdot \left( \frac{S(X_{t}) + l-L_{t}}{l}\right) 1_{L_{t} <l} \qquad (x=0) \eqno(3.2.40) $$

\noi Finally, we indicate that, in further works in progress, C. Profeta (see [Pr]) studies the penalisation of a linear diffusion
reflected in $0$ and $1$ (the subexponentiality hypothesis is not satisfied) with the functional
$(e^{\alpha L_t}, t \geq 0)$ ($\alpha \in \mathbb{R}$). He proves that the penalised process is still a linear
diffusion reflected in $0$ and $1$ and computes the scale function and the speed measure of this new process. 
\bigskip

\noi {\bf \large 3.3 The example of Bessel processes with dimension $d$ ($0<d<2$)}

\smallskip

\noi {\bf 3.3.1} \underline{Transcription of our notation in the context of Bessel processes.}

\noi Let $d=2(1-\alpha)$ with $0<d<2$ (or $0<\alpha<1$). We now study the particular case of the process $(X_{t}, \; t \ge 0)$ described in Section 3.1 with :
\begin{eqnarray*}
\hspace*{1cm}m(dx) &=& \frac{x^{1-2 \alpha}}{\alpha} \; 1_{[0, \infty[} (x) dx \hspace*{7,5cm}(3.3.1) \\
S(x) &=& x^{2 \alpha} \quad (x \ge 0) \hspace*{8,3cm}(3.3.2)
\end{eqnarray*}

\noi Then, the process $(X_{t}, \; t \ge 0)$ described in Section 3.1 is a Bessel process with dimension $d$, and index $\dis \frac{d}{2} -1= - \alpha$. We denote by $(P_{x}^{(-\alpha)}, \; x \in \mathbb{R}_{+})$ the family of its laws. We note $\big(\Omega = \mathcal{C} (\mathbb{R}_{+} \to \mathbb{R}_{+}), \; (R_{t}, \mathcal{F}_{t}), \; t \ge 0, \; \mathcal{F}_{\infty}, \; P_{x}^{(-\alpha)} \; (x \in \mathbb{R}_{+})\big)$ the canonical realisation of the Bessel process with index $(-\alpha)$. Here, the probability $P_{x}^{\uparrow}$ defined in 3.1.4 is the law of Bessel process with dimension $4-d=2(1+\alpha)$, i.e. : index $\alpha$. We shall denote this law by $P_{x}^{(\alpha)}$. The formulae of subsection 3.1 now become :
\hspace*{1cm}\begin{gather*}
( R_{t}^{2 \alpha} - L_{t},  \; t \ge 0) \; \; \hbox{is a martingale} \hspace*{7,9cm}(3.3.3)\\
P_{x}^{\uparrow} = P_{x}^{(\alpha)} \hspace*{11,9cm}(3.3.4) \\
\int_{0}^{t} h (R_{s})ds = \frac{1}{\alpha} \int_{0}^{\infty} h(x) \, L_{t}^{x} \, x^{1-2 \alpha} dx \hspace*{7,3cm} (3.3.5)\\
E_{0}^{(-\alpha)} (L_{t}^{0}) = t^{\alpha} \, E_{0}^{(- \alpha)} (L_{1}) = \frac{2^{\alpha} \, t^{\alpha}}{\Gamma (1-\alpha)} \hspace*{7,2cm}(3.3.6)  \\
L^{(\alpha)} f (r) = \frac{1}{2} \, f'' (r) + \frac{1+2 \alpha}{2 r} \; f'(r) \hspace*{7,8cm}(3.3.7)
\end{gather*}

\noi The reader may refer to [DMRVY] for these formulae. 

\bigskip

\noi {\bf 3.3.2} \underline{The measure ${\bf W}^{(- \alpha)}$.}

\smallskip

\noi In this framework, Theorem 3.2.1 becomes :

\noi {\bf Theorem 3.3.1.} {\it For every $\alpha \in ]0,1[$ :

\noi {\bf 1)} There exists a unique positive and $\sigma$-finite measure ${\bf W}^{(-\alpha)}$ on $\big(\Omega = \mathcal{C} (\mathbb{R}_{+} \to \mathbb{R}_{+}), \; \mathcal{F}_{\infty}\big)$ such that, for
every $F_t \in b(\mathcal{F}_t)$ :
$$
{\bf W}^{(-\alpha)} (F_{t} \, 1_{g \le t}) = P_{0}^{(-\alpha)} (F_{t} \cdot R_{t}^{2 \alpha}) \eqno(3.3.8) $$

\noi {\bf 2)} $\dis \quad \; \;{\bf W}^{(-\alpha)}= \int_{0}^{\infty}  (P_{0}^{(-\alpha, \,\tau_{l})} \circ P_{0}^{(\alpha)}) dl \hfill(3.3.9)$

\noi {\bf 3)}   i) $\dis \qquad {\bf W}^{(-\alpha)} (g \in dt) = \frac{\alpha \, 2^{\alpha}}{\Gamma (1-\alpha)} \; t^{\alpha -1} dt \qquad (t \ge 0)$

 {\it ii)} Conditionally on $g=t$, under ${\bf W}^{(-\alpha)}, \; (R_{u}, \; u \le g)$ is a Bessel bridge with index $(-\alpha)$ and of length $t$ 

 {\it iii)} $\dis \qquad {\bf W}^{(-\alpha)}= \int_{0}^{\infty} \frac{\alpha \, 2^{\alpha} \, t^{\alpha -1}}{\Gamma (1-\alpha)} \; dt (\Pi_{0}^{(- \alpha, t)} \circ P_{0}^{(+ \alpha)} ) \hfill(3.3.10) $   }

\noi \underline{In this Theorem} :

\noi $\Pi_{0}^{(-\alpha, t)}$ denotes the law of a Bessel bridge with index $(-\alpha)$ and of length $t$.

\smallskip

\noi $P_{0}^{(-\alpha, \tau_{l})}$ denotes the law of a Bessel process with index $(-\alpha)$ starting at 0 and stopped at $\tau_{l}$, with :
$$     \tau_{l} = \inf \{ t \ge 0 \;;\; L_{t}^{0} > l\}   \eqno(3.3.11)   $$

\noi {\bf 3.3.3} \underline{Relations between ${\bf W}^{(-\alpha)} \big(d = 2(1-\alpha)\big)$ and Feynman-Kac penalisations.}

\noi {\bf Remark 3.3.2.} The measure ${\bf W}^{(-\alpha)}$ which we just described is also related to a penalisation problem. More precisely, one can prove \big(see [RVY, I or V]\big) : 

\noi {\it i)} Let $q$ be a positive Radon measure on $\mathbb{R}_{+}$, with compact support. Then :
$$
2^{\alpha} \Gamma (1+ \alpha) t^{\alpha} \, P_{r}^{(- \alpha)} \left(\exp \left(-\frac{1}{2} \, A_{t}^{(q)}\right) \right) \mathop{\longrightarrow}_{t \to \infty}^{} \varphi_{q}^{(-\alpha)} (r)  \eqno(3.3.12) $$

\noi with
$$
A_{t}^{(q)} := \int_{0}^{\infty} q (R_{s}) ds = \frac{1}{\alpha} \int_{0}^{\infty} L_{t}^{x} \, x^{1-2 \alpha} q(dx)
\eqno(3.3.13)  $$

\noi {\it ii)} The function $\varphi_{q}^{(-\alpha)}$ defined by (3.3.12) is characterised as the unique solution of :
\begin{equation*} 
\hspace*{2cm} \left\{ \begin{array} {ll}
\dis \frac{1}{2} \, f''(r) + \frac{1-2 \alpha}{2r} \, f'(r) = \frac{1}{2} \, f(r) q(r)\\
\hbox{(in the sense of Schwartz distributions)} \\
f(r) \dis \mathop{\sim}_{r \to \infty}^{} r^{2 \alpha} \hspace*{9,5cm}(3.3.14) 
\end{array} \right.
\end{equation*}

\noi {\it iii)} For every $s \ge 0$ and $\Gamma_{s} \in b (\mathcal{F}_{s})$ :
$$
E_{r}^{(-\alpha)} \left( \Gamma_{s} \; \frac{\exp - \frac{1}{2} \, A_{t}^{(q)}}{E_{r}^{(-\alpha)} \big(\exp - \frac{1}{2} \, A_{t}^{(q)}\big)} \right) \mathop{\longrightarrow}_{t \to \infty}^{} P_{r, \infty}^{(- \alpha,q)} (\Gamma_{s})  \eqno(3.3.15)  $$

\noi where the probability $P_{r, \infty}^{(- \alpha,q)}$ satisfies :
$$
P_{r, \infty}^{(- \alpha,q)}  |_{\mathcal{F}_{s}} = M_{s}^{(-\alpha, q)} P_{r}^{(-\alpha)} |_{\mathcal{F}_{s}}
\eqno(3.3.16)   $$

\noi with
$$
M_{s}^{(- \alpha, q)} = \frac{\varphi_{q}^{(- \alpha)} (R_{s})}{\varphi_{q}^{(-\alpha)} (r)} \; \exp \left(-\frac{1}{2} \, A_{s}^{(q)}\right)   \eqno(3.3.17)   $$

\noi and $(M_{s}^{(-\alpha, q)}, \; s \ge 0)$ is a $\big((\mathcal{F}_{s}, \; s \ge 0), \; P^{(-\alpha)}\big)$ martingale.

\smallskip

\noi {\it iv)} Under $P_{r, \infty}^{(- \alpha,q)} \; (r \ge 0)$, the canonical process $(R_{t}, \; t \ge 0)$ is a transient diffusion with infinitesimal generator $\mathcal{G}^{(-\alpha, q)}$ given by :
$$
\mathcal{G}^{(-\alpha, q)} f(r) = \frac{1}{2} \; f'' (r) + \left( \frac{1-2 \alpha}{2r} + \frac{(\varphi_{q}^{(- \alpha)})'}{ \varphi_{q}^{(-\alpha)} }\, (r) \right) f'(r)   \eqno(3.3.18)   $$

\noi {\bf Remark 3.3.3.}

\noi With the notation of Remark 3.3.2, in the particular case where $q$ is the measure $q_{0}$ such that $\dis \frac{1}{\alpha} \; x^{1-2 \alpha} q_{0} (dx)$ is Dirac mass in 0 (of course, this is somewhat informal : we need to choose a sequence $q_{0}^{(n)}$ such that $\dis \frac{1}{\alpha} \, x^{1-2 \alpha} q_{0}^{(n)} (dx)$ converges towards $\delta_{0}$ as $n \to \infty$), we obtain :
\begin{eqnarray*}
\hspace*{2cm}&&\varphi_{q_{0}}^{(-\alpha)} (r) = 2 + r^{2 \alpha}, \; \varphi_{q_{0}}^{(-\alpha)} (0) = 2 \hspace*{6,1cm}(3.3.19) \\
{\rm and}
&& M_{t}^{(-\alpha, q_{0})} = \left(1+ \frac{R_{t}^{2 \alpha}}{2}\right) \, e^{- \frac{1}{2} \, L_{t}} \hspace*{6,6cm} (3.3.20)
\end{eqnarray*}

\noi \underline{Now, the analogue of Theorem 1.1.5 is} : 

\noi {\bf Theorem 3.3.4.}

\noi {\it Under $P_{\infty}^{(-\alpha, q_{0})}$, the canonical process $(R_{t}, \; t \ge 0)$ satisfies :

\noi {\it i)} Let $g= \sup \{s \ge 0, \; R_{s}=0\}$. Then :
$$
g < \infty \qquad P_{\infty}^{(-\alpha, q_{0})} \qquad a.s. \qquad and   \eqno(3.3.21) $$

\noi {\it ii)} $L_{\infty}(=L_{g})$ admits as density :
$$
f^{P_{\infty}^{(-\alpha, q_{0})}}_{L_{\infty}} (l) = \frac{1}{2} \, e^{- \frac{l}{2}} \, 1_{[0, \infty[} (l)dl
\eqno(3.3.22) $$

\noi {\it iii)} Conditionally on $g$, $(R_{s}, \; s \le g)$ and $(R_{g+s}, \; s \ge 0)$ are independent.

\noi {\it iv)} $(R_{g+s}, \; s \ge 0)$ is a $(4-d)$ dimensional Bessel process starting at 0 \big(i.e. admits $P_{0}^{(+\alpha)}$ as its law\big).

\noi {\it v)} Conditionally on $L_{\infty} (=L_{g})=l, \; (R_{s}, \; s \le g)$ is a $d$-dimensional Bessel process stopped at $\tau_{l}$. Its law is $P_{0}^{(- \alpha, \tau_{l})}$.  }

\bigskip

\noi {\bf Remark 3.3.5.}

\noi {\bf 1)} Since, for $\dis \alpha = \frac{1}{2}, \; (R_{t}, \; t \ge 0)$ under $P^{(-\alpha)}$ is a reflected Brownian motion, one has :
$$
{\bf W} \big(F(|X_{s}|, \; s \ge 0)\big) = {\bf W}^{\big(-\frac{1}{2}\big)} \big(F(R_{s}, \; s \ge 0)\big)
\eqno(3.3.23)   $$

\noi (where ${\bf W}$ is defined by Theorem 1.1.2).

\noi {\bf 2)} In the same spirit, since the modulus of a 2-dimensional Brownian motion is a 2-dimensional Bessel process, hence has index 0, we conjecture that, in a sense to be made precise :
$$
{\bf W}^{(2)} \big(F (|X_{s}|, \; s \ge 0)\big) = \mathop{\lim}_{\alpha \downarrow 0}^{} \quad {\bf W}^{(- \alpha)} \big(F (R_{s}, \; s \ge 0)\big)   \eqno(3.3.24) $$
\noi (where ${\bf W}^{(2)}$ is defined by Theorem 2.1.2).

\noi {\bf Remark 3.3.6.}

\noi We have given, in Subsection 1.1.6, a proof of Theorem 1.1.6 (this is precisely Theorem 1.1.10) which hinges upon the disintegration of Wiener measure restricted to $\mathcal{F}_{t}$, with respect to the law of $g^{(t)}$ \big(see (1.1.82)\big). Formula (3.3.10) may be proven in a quite similar way by using the following :

\noi {\it i)} For fixed time $t$, the three following random elements are independent :
\begin{itemize}
\item[$\bullet$] $\dis \left(r_{u} := \frac{1}{\sqrt{g^{(t)}}} \, R_{u\,g^{(t)}}, \; u \le 1\right)$ which is a Bessel bridge with dimension $d=2(1-\alpha)$
\item[$\bullet$] $g^{(t)}:= \sup \{u<t \;;\; R_{u}=0\}$ which is distributed as :
$$ 
P_{0}^{(-\alpha)} (g^{(t)} \in du)= \frac{du}{b_{\alpha} u^{1-\alpha}(t-u)^{\alpha}} \qquad (0 \le u \le t)
\eqno(3.3.25) $$

with : $b_{\alpha} = B(\alpha, 1-\alpha)=\Gamma (\alpha) \Gamma (1-\alpha) = \frac{\pi}{\sin( \pi \alpha)}$.
\item[$\bullet$] $\dis \left(m_{u} := \frac{1}{\sqrt{t-g^{(t)}}} \, R_{g^{(t)} + u(1-g^{(t)})}, \; u \le 1\right)$, which is a Bessel meander (with dimension $d$).
\end{itemize}

\noi {\it ii)} Imhof's absolute continuity relationship between the laws of the Bessel meander $(m_{u}, \; u \le~1)$ and the Bessel process with dimension $2(1+\alpha)$ (i.e. : with index $\alpha$) is : 
$$
E^{(- \alpha)} \big(F(m_{u}, \; u \le 1)\big) = E_{0}^{(\alpha)} \left(F(R_{u}, \; u \le 1) \, \frac{2^{\alpha} \Gamma (1+\alpha)}{R_{1}^{2 \alpha}}\right)   \eqno(3.3.26)   $$

\bigskip

\noi {\large \bf 3.4 Another description of ${\bf W}^{(-\alpha)}$ (and of ${\bf W}_{g}^{*})$.}

\smallskip

\noi {\bf 3.4.1} \underline{We recall that \big(see (3.2.4), (3.3.9) and (3.3.10)\big)} :
$$
{\bf W}_{g}^{*} = \int_{0}^{\infty} dl \, P_{0}^{\tau_{l}} = \int_{0}^{\infty} dt \, p(t,0,0) \, \Pi_{0}^{(t)}
\eqno(3.4.1)   $$

\noi in the context of general linear diffusions and :
$$
{\bf W}_{g}^{(-\alpha)} = \int_{0}^{\infty} dl \, P_{0}^{(-\alpha, \tau_{l})} = \int_{0}^{\infty} \frac{\alpha 2^{\alpha} t^{\alpha -1}}{\Gamma (1-\alpha)} \; \Pi_{0}^{(-\alpha, t)} \, dt    \eqno(3.4.2)  $$

\noi in the context of the Bessel processes with index $(- \alpha) \;\; (0<\alpha<1)$

\noi We shall now give a new description of ${\bf W}_{g}^{(-\alpha)}$ (resp. ${\bf W}_{g}^{*})$ which is the restriction of ${\bf W}^{(-\alpha)}$ (resp. ${\bf W}^{*})$ to $\mathcal{F}_{g}$. This new description is simply the transcript in the Bessel framework of results found in Pitman-Yor \big(see [PY2]\big).

\smallskip

\noi {\bf 3.4.2} \underline{We begin by recalling in the framework of Bessel processes some of the results from [PY2].}

\noi We denote by $\widehat{\Omega}$ the space of continuous functions from $\mathbb{R}_{+}$ to $\mathbb{R}_{+}$ with finite lifetime $\xi$ :
\begin{gather*}
\hspace*{1,5cm}\widehat{\Omega} = \big\{ \omega \; : \mathbb{R}_{+} \to \mathbb{R}_{+} \;;\; \exists \xi (\omega) < \infty \; {\rm s.t.} \; \omega (0) = 0 = \omega (\xi),  \\
\hspace*{2cm}\qquad {\rm and} \qquad \omega (u) = 0 \quad {\rm for} \; {\rm every} \; u \ge \xi (\omega)\big\} 
\end{gather*}

\noi We denote by $(R_{t}, \; t \ge 0)$ the set of coordinates on this space :
\begin{equation*}
R_{t} (\omega) = \omega (t), \; \omega \in \widehat{\Omega}
\end{equation*}

\noi The result of Pitman-Yor which we use \big(Theorem 1.1 of [PY2]\big) asserts the existence, for every $\delta >0$, of a positive and $\sigma$-finite measure on $(\widehat{\Omega}, \mathcal{F}_{\infty})$, denoted as ${\bf \Lambda}_{0,0}^{(\delta)}$ and which may be described in either of the following manners : 

\smallskip

\noi \underline{First description}
$$
\text{\boldmath$\Lambda$\unboldmath}_{0,0}^{(\delta)} = \int_{0}^{\infty} \frac{2^{- \frac{\delta}{2}}}{\Gamma (\delta/2)} \; t^{- \frac{\delta}{2}} \quad \Pi_{0}^{\big (\frac{\delta}{2}-1, \, t\big)} dt \eqno(3.4.3) $$

\noi where $\Pi_{0}^{\big (\frac{\delta}{2}-1, \, t\big)}$ denotes the law of the Bessel bridge with index $\dis \frac{\delta}{2} - 1$, i.e. with dimension $\delta$, and length $t$.

\smallskip

\noi \underline{Second description}

\smallskip

\noi Let $m>0$ fixed and let $P_{0} ^{\big(\frac{\delta}{2} - 1, m, \nearrow \nwarrow\big)}$ denote the law of the process obtained by putting two Bessel processes with index $\dis \left(\frac{\delta}{2} - 1\right)$ (i.e. : with dimension $\delta$), back to back starting from 0, and stopped when they first reach level $m$. These two processes $R$ and $\widetilde{R}$ are assumed to be independent. In other terms, $P_{0} ^{\big(\frac{\delta}{2} - 1, m, \nearrow \nwarrow \big)}$ is the law of the process $(Y_{t}, \; t \ge 0)$ defined by :
$$ Y_{t} =
\left \lbrace
\begin{array}{lcl}
R_{t}  & {\rm if} & t \le T_{m}  \\
\widetilde{R}_{T_{m} + \widetilde{T}_{m} - t} &{\rm if} & T_{m} \le t \le T_{m} + \widetilde{T}_{m} \\
0 & {\rm if} & t \ge T_{m} + \widetilde{T}_{m} 
\end{array} \right. \eqno(3.4.4)$$

\noi where $T_{m}$ (resp $\widetilde{T}_{m})$ is the first hitting time of $m$ by $(R_{t}, \, t \ge 0)$ \big(resp. by $(\widetilde{R}_{t}, \; t \ge 0)$\big).

\noi Then : 
$$
\text{\boldmath$\Lambda$\unboldmath}_{0,0}^{(\delta)} = \int_{0}^{\infty} m^{1- \delta} dm \, P^{\big(\frac{\delta}{2} - 1, m, \nearrow \nwarrow\big)} \eqno(3.4.5) $$

\noi The measure \boldmath$\Lambda$\unboldmath$_{0,0}^{(\delta)}$ is called the "generalized excursion measure" in Pitman-Yor. When $\delta = 3$, ${\bf \Lambda}_{0,0}^{(3)}$ is the It\^{o} measure of (positive) Brownian excursions. Formula (3.4.3) is It\^{o}'s description of It\^{o}'s measure (see [ReY], Chap. XII), whereas formula (3.4.5) is Williams' description of that measure (see [Wi]).
 
\noi {\bf 3.4.3} \underline{Here is now, in the framework of the Bessel processes, the announced transcription} :

\noi {\bf Theorem 3.4.1} {\it For every $\alpha \in ]0,1[$ :
$$
{\bf W}^{(-\alpha)}\big|_{\mathcal{F}_{g}} = {\bf W}_{g}^{(-\alpha)} = 2 \alpha {\bf \Lambda}_{0,0}^{(2(1-\alpha))} \eqno(3.4.6) $$

\noi In particular :
\begin{gather*}
\hspace*{1cm}{\bf W}^{(-\alpha)} \big|_{\mathcal{F}_{g}} = \frac{\alpha 2^{\alpha}}{\Gamma (1-\alpha)} \, \int_{0}^{\infty} t^{\alpha -1} dt \, \Pi_{0}^{(-\alpha, t)} \hspace*{5,5cm}(3.4.7) \\
\hspace*{1cm}{\bf W}^{(-\alpha)} \big|_{\mathcal{F}_{g}}  = 2 \alpha \int_{0}^{\infty} m^{2 \alpha -1} dm \, P^{(- \alpha, m, \ \nearrow \nwarrow)} \hspace*{5,5cm}(3.4.8)
\end{gather*}

\noi Thus, formula (3.4.8) provides us with a new description of the measure ${\bf W}_{g}^{(-\alpha)}$. }

\smallskip

\noi {\bf Proof of Theorem 3.4.1} Of course, from (3.4.3), and (3.4.5), it suffices to show (3.4.6). Note that, from (3.3.8), for $\Gamma_{t} \in b(\mathcal{F}_{t})$, one has :
$$ {\bf W}^{(-\alpha)} (\Gamma_{t} \, 1_{g \le t}) = P_{0}^{(- \alpha)} (\Gamma_{t} \, R_{t}^{2 \alpha}) \eqno(3.4.9) $$

\noi Thus, for every $s \le t$ and $\Gamma_{s} \in b (\mathcal{F}_{s})$, since $(R_{t}^{2 \alpha} - L_{t}, \; t \ge 0)$ is a martingale \big(see (3.3.3)\big), we have :
\begin{eqnarray*}
 \hspace*{2cm} {\bf W}^{(- \alpha)} (\Gamma_{s} \, 1_{s \le g \le t}) 
  &= & P_{0}^{(-\alpha)} \big(\Gamma_{s} (R_{t}^{2 \alpha} - R_{s}^{2 \alpha})\big) \\
  &=& P_{0}^{(-\alpha)} \big(\Gamma_{s} (L_{t} - L_{s})\big) \hspace*{5cm} (3.4.10)
\end{eqnarray*}

\noi We deduce from the monotone class theorem and from (3.4.10) that, for every positive previsible process $(\Phi_{u}, \; u \ge 0)$, one has :
\begin{eqnarray*}
{\bf W}^{(-\alpha)} (\Phi_{g})
&=& P_{0}^{(-\alpha)} \left( \int_{0}^{\infty} \Phi_{u} \, dL_{u} \right) \\
&=& \int_{0}^{\infty} P_{0}^{(-\alpha)} (\Phi_{u} | R_{u} =0) P_{0}^{(-\alpha)} (dL_{u}) \\
&=& \int_{0}^{\infty} \Pi^{(- \alpha, u)} (\Phi_{u}) \, \frac{\alpha 2^{\alpha} u^{\alpha -1}}{\Gamma (1-\alpha)} \; du
\end{eqnarray*}

\noi from (3.3.6). Hence :
\begin{eqnarray*}
{\bf W}^{(-\alpha)} (\Phi_{g})
&=& \left( \int_{0}^{\infty} \Pi_{0}^{(- \alpha, u)} \frac{\alpha 2^{\alpha} u^{\alpha -1}}{\Gamma (1-\alpha)} \; du \right) (\Phi_{g}) \\
&=& \left(2 \alpha \int_{0}^{\infty} du \, \frac{2^{- \frac{\delta}{2}}}{\Gamma (\frac{\delta}{2})} \, u^{- \frac{\delta}{2}} \, \Pi_{0}^{( \frac{\delta}{2} -1,u)}\right) (\Phi_{g}) \\
&& \big({\rm since} \; \delta = 2 (1-\alpha)\big) \\
&=& 2 \alpha \text{\boldmath$\Lambda$\unboldmath}_{0,0}^{(2(1- \alpha))} (\Phi_{g}) \quad \text{(from (3.4.3))}
\end{eqnarray*}

\noi {\bf 3.4.4} \underline{In the general framework of linear diffusions, formulae (3.4.7) and (3.4.8) become} :
$$
{\bf W}^{*}_{g} = \int_{0}^{\infty} dt \, p (t, 0,0) \, \Pi_{0}^{(t)}  $$

\noi \big(this is formula (3.2.4)\big) and : 
$$ {\bf W}_{g}^{*} = \int_{0}^{\infty} P_{0}^{(m, \nearrow \nwarrow)} d S(m)   \eqno(3.4.11) $$

\noi The reader may refer to \big([PY2], 2.2, Corollary 2.1, p. 298\big) where the probability $P_{0}^{(m, \nearrow \nwarrow)}$ is defined in terms of the law $P_{0}^{\uparrow}$ \big(of the process $(X_{t}, \; t \ge 0)$ conditioned to remain $> 0$\big) just as $P_{0}^{(-\alpha, m, \nearrow \nwarrow)}$ is, in terms of the law $P_{0}^{(\alpha)}$. \big(see (3.4.4) with $\delta = 2 (1+ \alpha)$\big).

\bigskip

\noi {\large \bf 3.5 Penalisations of $\alpha$-stable symmetric L\'evy process ($1 < \alpha \leq 2$)}

\noi In this subsection, we summarize the results by K. Yano, Y. Yano and M. Yor [YYY] which bears upon the penalisation of the $\alpha$-stable symmetric L\'evy process, with $1< \alpha \le 2$. This summary is not exhaustive ; rather, it is an invitation to read [YYY].

\smallskip

\noi {\bf 3.5.1} \underline{Notation and classical results.} \big(see, e.g., [Be], [C], [SY]\big)

\smallskip

\noi {\bf 3.5.1.1} $\big(\Omega, (X_{t}, \mathcal{F}_{t})_{t \ge 0}, \mathcal{F}_{\infty}, P_{x}, \; x \in \mathbb{R}\big)$ denotes the canonical realization of the $\alpha$-stable symmetric L\'evy process, with $1<\alpha\le 2$. The notations are the same as in 1.0.1, with the difference that $\Omega$ now denotes the space of c\`{a}dl\`{a}g functions from $\mathbb{R}_{+}$ to $\mathbb{R}$. $\alpha$ being fixed once and for all, the dependency in $\alpha$ will be mostly omitted in our notation. This L\'evy process $(X_{t}, \; t \ge 0)$ is characterised via :
$$  E_{0} (e^{i \lambda X_{t}}) = \exp \big(-t |\lambda|^{\alpha}\big) \qquad (t \ge 0,\; \lambda \in \mathbb{R})  \eqno (3.5.1)   $$

\noi The case $\alpha =2$ corresponds to $(X_{t}, \; t \ge 0) \equiv (B_{2t}, \; t \ge 0)$ where $(B_{t}, \; t \ge 0)$ is a standard 1-dimensional Brownian motion. 

\smallskip

\noi {\bf 3.5.1.2} $p_{t} (x)$ denotes the density (with respect to Lebesgue measure on $\mathbb{R}$) of the law of the r.v. $X_{t}$ and $u_{\lambda} \; (\lambda > 0)$ the resolvent kernel :
$$   P_{x} (X_{t} \in dy) = p_{t} (x-y) dy = p_{t} (y-x) dy \eqno(3.5.2)   $$
$$   p_{t} (0) = \frac{1}{\alpha \pi} \; \Gamma \left(\frac{1}{\alpha}\right) \, t^{- \frac{1}{\alpha}} \qquad (t > 0) \eqno (3.5.3)    $$
$$   u_{\lambda} (x) := \int_{0}^{\infty} e^{- \lambda t} p_{t} (x) dt = \frac{1}{\pi} \int_{0}^{\infty} \frac{\cos (xy)}{\lambda + y^{\alpha}} \; dy \eqno (3.5.4) $$
$$ u_{\lambda} (0) = \frac{1}{\pi} \; B \left(1- \frac{1}{\alpha}, \frac{1}{\alpha}\right) \, \lambda^{\frac{1} {\alpha} -1} \eqno (3.5.5)   $$

\noi Let, for any $a \in \mathbb{R}$, $T_{a} := \inf \{t \ge 0 \;;\; X_{t} =a\}$. Then : 
$$  E_{x} [e^{- \lambda \, T_{0}}] = \frac{u_{\lambda} (x)}{u_{\lambda} (0)} \eqno (3.5.6) $$

\noi {\bf 3.5.1.3} We denote by $(L_{t}^{x}, \; t \ge 0, \; x \in \mathbb{R})$ the jointly continuous process of local times of $(X_{t}, \; t \ge 0)$, $(L_{t}, \; t \ge 0)$ stands for $(L_{t}^{0}, \; t \ge 0)$, the local time process at 0, and $(\tau_{l}, \; l \ge 0)$ its right continuous inverse. We have :
$$    E_{0} (e^{- \lambda \, \tau_{l}}) = \exp \left(- \frac{l}{u_{\lambda} (0)} \right) \eqno (3.5.7)  $$

\noi so that, from (3.5.5), $(\tau_{l}, \; l \ge0)$ is a stable subordinator with index $\dis 1- \frac{1}{\alpha}�\cdot$ On the other hand~:
\begin{eqnarray*}
\hspace*{1cm}E_{0} \left(\int_{0}^{\infty} e^{-\lambda t} dL_{t} \right)
&=& \int_{0}^{\infty} E_{0} (e^{- \lambda \, \tau_{l}}) dl = u_{\lambda} (0) \hspace*{5cm} (3.5.8) \\
{\rm and} \hspace*{1cm} E_{0} (dL_{t})
&=& p_{t} (0) dt = \frac{1}{\alpha \pi} \; \Gamma \left(\frac{1}{\alpha}\right) \, t^{- \frac{1}{\alpha}} dt \hspace*{4,5cm} (3.5.9) \\
\hbox{More generally :} \; E_{x} (dL_{t})
&=& E_{0} (d_{t} \, L_{t}^{x}) = p_{t} (x) dt \hspace*{5,8cm} (3.5.10)
\end{eqnarray*}

\noi {\bf 3.5.1.4} We denote by $h$ the function defined by :
$$     h(x) := \frac{1}{2 \, \Gamma(\alpha) \sin \left[ \frac{(\alpha - 1) \pi}{2} \right] }\; |x|^{\alpha -1} \qquad (x \in \mathbb{R}) \eqno (3.5.11)           $$

\noi This function is harmonic for the process $(X_{t}, \; t \ge 0)$ killed when it reaches 0, i.e. : for every $x \in \mathbb{R}$, and $t \ge 0$ :
$$    E_{x} \big[h(X_{t}) 1_{T_{0} \ge t}\big] = h(x) \eqno (3.5.12) $$

\noi Moreover, there exists a constant $c > 0$ such that, for every $x \in \mathbb{R}$ :
$$    \big(N_{t}^{x} := h (X_{t}) - h(x) - cL_{t}^{x}, \; t \ge 0 \big)  \eqno (3.5.13)    $$

\noi is a square integrable $P_{x}$-martingale \big(this formula may be compared with (3.1.7)\big). 

\smallskip

\noi {\bf 3.5.1.5} Since 0 is a regular and recurrent point for $(X_{t},\; t \ge 0)$, It\^{o}'s excursion theory may be applied. We denote by $\widetilde{\Omega}$ the excursions space, where $(Y_{t}, \; t \ge 0)$ is the process of coordinates, $\xi$ the lifetime of the generic excursion and ${\bf n}$ It\^{o}'s excursion measure. The master formula  from excursion theory implies :
$$    E_{0} \left[ \int_{0}^{\infty} e^{- \lambda t} f (X_{t}) dt \right] = E_{0} \left(\int_{0}^{\infty} e^{- \lambda \, \tau_{l}} dl \right) \cdot \int_{0}^{\infty} e^{- \lambda t} {\bf n} \big(f (Y_{t})\big) dt \eqno (3.5.14)   $$

\noi for any $f : \mathbb{R} \to \mathbb{R}_{+}$ Borel, such that $f(0)=0$. In particular :
$$   {\bf n} (\xi >t) = \frac{\alpha \pi}{B\left(1 - \frac{1}{\alpha}, \frac{1}{\alpha}\right) \Gamma \left(\frac{1}{\alpha}\right)} \; t^{\frac{1}{\alpha} -1} \eqno (3.5.15)    $$

\noi There exists a function $\rho (t,x)$ which is positive and jointly measurable such that : 
\begin{eqnarray*}
\hspace*{4cm}  {\bf n} (Y_{t} \in dx) &=& \rho (t,x) dx \hspace*{5,8cm} (3.5.16)  \\
{\rm and} \qquad 
P_{x} (T_{0} \in dt) &=& \rho (t,x) dt \hspace*{5,9cm} (3.5.17)
\end{eqnarray*}

\noi {\bf 3.5.2} \underline{Definition 
of the $\sigma$-finite measure ${\bf P}$} \footnote{We take up the notation from [YYY].}

\noi The measure ${\bf P}$ is defined on $(\Omega, \mathcal{F}_{\infty})$ by :
\begin{eqnarray*}
\hspace*{2cm} {\bf P}
&:=& \int_{0}^{\infty} P_{0} (dL_{u}) \, ( Q^{(u)} \circ P_{0}^{\uparrow} ) \hspace*{7cm}  (3.5.18) \\
&=& \frac{1}{\alpha \pi} \; \Gamma \left(\frac{1}{\alpha}\right) \int_{0}^{\infty} du \, u^{- \frac{1}{\alpha}} \, ( Q^{(u)} \circ P_{0}^{\uparrow}) \hspace*{5,6cm} (3.5.19)
\end{eqnarray*}

\noi \big(from (3.5.3)\big). We now explain the notations in (3.5.18) :
\begin{itemize}
\item $Q^{(u)}$ denotes the law of the $\alpha$-stable symmetric bridge with length $u$ :
$$ Q^{(u)} (\Gamma_{u}) = P_{0} \big(\Gamma_{u} | X_{u} =0\big) \qquad (\Gamma_{u} \in \mathcal{F}_{u}) \eqno(3.5.20) $$
\item We denote by $(P_{x}^{0}, \; x \neq 0)$ the law of the process $(X_{t}, \; t \ge 0)$ starting from $x$ and killed in $T_{0}$ :
$$   P_{x}^{0} (\Gamma_{t}) = E_{x} (\Gamma_{t} \, 1_{T_{0} > t}) \qquad \Gamma_{t} \in b(\mathcal{F}_{t})  $$
and by $P_{x}^{\uparrow}$ the law obtained from that of $P_{x}^{0}$ by Doob's $h$-transform \big(recall that $h$ is defined by (3.5.11) and that it is harmonic for the process $(X_{t}, \; t \ge 0)$ killed in $T_{0}$\big) :
$$     P^{\uparrow}_{x \; | \mathcal{F}_{t}} := \frac{h(X_{t})}{h(x)} \cdot P^{0}_{x \; | \mathcal{F}_{t}} \qquad x \neq 0 \eqno (3.5.21) $$
\noi Letting $x$ tend to 0 in (3.5.21), we obtain : 
$$   P^{\uparrow}_{0 \; |\mathcal{F}_{t}} := \mathop{\lim}_{x \to 0} \, \frac{h(X_{t})}{h(x)} \cdot P^{0}_{ x \; | \mathcal{F}_{t}} = h(X_{t}) \, {\bf n}_{\; | \mathcal{F}_{t}} \eqno (3.5.22)    $$
\item Another manner to define $P_{0}^{\uparrow}$ consists in first defining the law $M^{(t)}$ of the stable meander (with duration $t$) :
$$  M^{(t)} (\Gamma_{t}) : = {\bf n} \big(\Gamma_{t} | \xi >t\big) = \frac{{\bf n} \big(\Gamma_{t} \cap (\xi > t)\big)} {{\bf n} (\xi > t)} \qquad \big(\Gamma_{t} \in b(\mathcal{F}_{t})\big) \eqno (3.5.23)    $$
then to show that :
$$    M^{(t)} \mathop{\longrightarrow}_{t \to \infty} P_{0}^{\uparrow}  \eqno (3.5.24)   $$
with the preceding convergence taking place along $(\mathcal{F}_{s})$, i.e. : for every $s \ge 0$ and $\Gamma_{s} \in b(\mathcal{F}_{s})$ :
$$     M^{(t)} (\Gamma_{s}) \mathop{\longrightarrow}_{t \to \infty} P_{0}^{\uparrow} (\Gamma_{s}) \eqno (3.5.25)     $$
\item The measure ${\bf P}$ defined by (3.5.18) plays for the symmetric $\alpha$-stable L\'evy process the same role as the measure ${\bf W}$ for standard Brownian motion. Indeed, for $\alpha =2$, (3.5.18) becomes 
$$   {\bf P} = \frac{1}{2\sqrt{\pi}} \; \int_{0}^{\infty} \frac{du}{\sqrt{u}} \;( Q^{(u)} \circ P_{0}^{\uparrow}) = \frac{1}{\sqrt{2}} \; {\bf W}     $$
where ${\bf W}$ is defined by (1.1.16), or (1.1.43). The multiplication factor $\dis \frac{1}{\sqrt{2}}$ arises from the fact that, for $\alpha =2$, the 2-stable symmetric L\'evy process $(X_{t}, \; t \ge 0)$ is the process $(B_{2t}, \; t \ge 0)$ and not $(B_{t}, \; t \ge 0)$ \big(see (3.5.1)\big).
\end{itemize}

\noi {\bf 3.5.3} \underline{The martingales $\big(M_{t} (F), \; t \ge 0\big)$ associated with ${\bf P}$}

\smallskip

\noi {\bf 3.5.3.1} In the same manner that we have associated to the $\sigma$-finite measures ${\bf W}$, ${\bf W}^{(2)}$ and ${\bf W}^{*}$ introduced in Section 1.2, and in (3.2.2) and (3.2.3), a family of martingales, we associate here to every r.v. $F \in L_{+}^{1} (\Omega, \mathcal{F}_{\infty}, {\bf P})$ the $\big((\mathcal{F}_{t})_{t \ge 0}, \; P_{0}�\big)$ martingale $\big(M_{t} (F), \; t \ge 0\big)$ characterized by~:  for any $t \ge 0$ and $\Gamma_{t} \in b(\mathcal{F}_{t})$ :
$$    E_{\bf P} [F \cdot \Gamma_{t}] = E_{0} \big(M_{t} (F) \cdot \Gamma_{t}\big) \eqno (3.5.26)   $$

\noi In particular, for every $t \ge 0$ :
$$ E_{0} \big[M_{t} (F)\big] = E_{\bf P} (F) \eqno (3.5.27)    $$

\noi {\bf 3.5.3.2} \underline{Example 1.} Let $f: \mathbb{R}_{+} \to \mathbb{R}_{+}$ Borel such that $\dis \int_{0}^{\infty} f(y)dy < \infty$. Then :
$$   M_{t} \big(f (L_{\infty})\big) = f (L_{t}) h (X_{t}) + \int_{L_{t}}^{\infty} f(x) dx \qquad (t \ge 0) \eqno (3.5.28) $$

\noi where, in (3.5.28), the function $h$ is defined by (3.5.11). It is not difficult to see, thanks to (3.5.13), that $\big(M_{t} \big(f(L_{\infty})\big), \; t \ge 0\big)$ defined by (3.5.28) is indeed a martingale. We also note the analogy between (3.5.28) and formula (3.2.31) obtained in the framework of linear diffusions :
$$   M_{t}^{*} \big(f (L_{\infty})\big) = f(L_{t}) S(X_{t}) + \int_{L_{t}}^{\infty} f(y)dy \eqno (3.5.29)  $$

\noi Thus, we shift from (3.5.29) to (3.5.28) by replacing simply the scale function $S$ by the function $h$ \big(these two functions are such that, in both cases, $\big(S(X_{t}) 1_{t<T_{0}}, \; t \ge 0\big)$ and $\big(h (X_{t} ) 1_{t <T_{0}}\big)$ are martingales\big).

\smallskip

\noi {\bf 3.5.3.3} \underline{Example 2.} (Feynman-Kac martingales)

\noi Let $q$ denote a Radon measure on $\mathbb{R}$ such that :
$$  0 < \int_{\mathbb{R}} \big(1+ h(x)\big) \, q(dx) < \infty \quad \hbox{(with $h$ defined by (3.5.11))} \eqno (3.5.30)   $$

\noi Let 
$$ A_{t}^{(q)} := \int_{\mathbb{R}} L_{t}^{x} q(dx) \eqno (3.5.31)   $$

\noi and $\dis A_{\infty}^{(q)} := \mathop{\lim}_{t \to \infty} A_{t}^{(q)}$. Then
$$   M_{t} \big(\exp (-A_{\infty}^{(q)})\big) = \varphi_{q} (X_{t}) \cdot \exp (-A_{t}^{(q)}) \eqno (3.5.32)   $$

\noi with
$$ \varphi_{q} (x) := \mathop{\lim}_{t \to \infty} \, \frac{E_{x} (\exp - A_{t}^{(q)})}{{\bf n} (\xi > t)} \qquad (x \in \mathbb{R}) \eqno (3.5.33)   $$

\noi We note that : $E_{\bf P} \big(\exp (-A_{\infty}^{(q)})\big) = \varphi_{q} (0) $.

\noi Other descriptions of the function $\varphi_{q}$ are found in [YYY]. The reader will have noticed the complete analogy between the definition of $M_{t} (\exp - A_{\infty}^{(q)})$ given by (3.5.32) and that, in the Brownian case, of $M_{t} (\exp - A_{\infty}^{(q)})$ which is given by (1.2.19) :
$$   M_{t} \left( \exp - \frac{1}{2} \, A_{\infty}^{(q)} \right) = \varphi_{q} (X_{t}) \exp \left(- \frac{1}{2} \, A_{t}^{(q)}\right)    $$

\noi {\bf 3.5.4} \underline{Relations between ${\bf P}$ and penalisations}

\noi The following penalisation theorems, which we now present, are found in [YYY] : 

\noi {\bf Theorem 3.5.1} {\it Let $f : \mathbb{R}_{+} \to \mathbb{R}_{+}$ Borel such that $\dis \int_{0}^{\infty} f(y) dy < \infty$. Then :

\noi {\bf 1)} For every $s \ge 0, \; \Gamma_{s} \in b (\mathcal{F}_{s})$ :
$$ \mathop{\lim}_{t \to \infty} \, \frac{E_{0} \big[\Gamma_{s} f(L_{t})\big]}{{\bf n} (\xi >t)} = E_{0} \big[ \Gamma_{s} \, M_{s} \big(f(L_{\infty})\big)\big] \eqno(3.5.34)   $$

\noi where $\big(M_{t} \big(f (L_{\infty}), \; t \ge 0\big)$ is the positive martingale defined by (3.5.28).

\smallskip

\noi {\bf 2)} Let $P_{0, \infty}^{f (L)}$ the probability induced on $(\Omega, \mathcal{F}_{\infty})$ by : 
$$   P^{f(L)}_{0, \infty \; | \mathcal{F}_{t}} := \frac{M_{t} \big(f (L_{\infty})\big)}{ M_{0} \big(f (L_{\infty})\big)} \cdot P_{0 \; |\mathcal{F}_{t}} \eqno (3.5.35) $$

\noi Then, the absolute continuity formula :
$$   f(L_{\infty}) \cdot {\bf P} = E_{\bf P} \big(f (L_{\infty})\big) \cdot P_{0, \infty}^{f(L)} \quad holds \eqno (3.5.36) $$

\noi \Big(Note that : $\dis E_{\bf P} \big(f (L_{\infty})\big) = \int_{0}^{\infty} f(y)dy = E_{0} \big(M_{t} \big(f (L_{\infty})\big)$.\Big)       }

\noi Clearly, this formula (3.5.36) is formally identical to formula (1.1.107) obtained in the Brownian set-up (with $h^{+} = h^{-} =f$).   

\smallskip

\noi {\bf Theorem 3.5.2} {\it Let $q$ denote a Radon measure on $\mathbb{R}$ such that $\dis 0 < \int_{\mathbb{R}} \; \big(1+h (x)\big) q(dx) < \infty$ \big(with $h$ defined by (3.5.11)\big) and let $\dis A_{t}^{(q)} := \int_{\mathbb{R}} L_{t}^{x} q (dx)$. Then :

\noi {\bf 1)} For every $s \ge 0$ and $\Gamma_{s} \in b (\mathcal{F}_{s})$ : 
$$ \mathop{\lim}_{t \to \infty} \, \frac{E_{0} \big(\Gamma_{s} \exp \big(- A_{t}^{(q)})\big)}{ {\bf n} (\xi >t)} = E_{0} \big[\Gamma_{s} \, M_{s} \big(\exp (-A_{\infty}^{(q)})\big)\big] \eqno (3.5.37)   $$

\noi where $\big(M_{t} \big(\exp (-A_{\infty}^{(q)})\big), \; t \ge 0\big)$ is the positive martingale defined by (3.5.32).

\smallskip

\noi {\bf 2)} Let $P_{0, \infty}^{(q)}$ denote the probability induced on $(\Omega, \mathcal{F}_{\infty})$ by :
$$   P^{(q)}_{0, \infty \; | \mathcal{F}_{t}} = \frac{M_{t} \big(\exp (-A_{\infty}^{(q)})\big)}{M_{0} \big(\exp -(A_{\infty}^{(q)})\big)} \cdot P_{0 \; | \mathcal{F}_{t}} \eqno (3.5.38)   $$

\noi Then, the absolute continuity formula : 
$$  \exp (-A_{\infty}^{(q)}) \cdot {\bf P} = E_{\bf P} (\exp (-A_{\infty}^{(q)})\big) \cdot P_{0, \infty}^{(q)} \quad holds \eqno (3.5.39)   $$      }

\noi Of course, this formula is formally identical to formula (1.1.16') obtained in the Brownian framework (one should note that $E_{\bf P} (\exp -A_{\infty}^{(q)}) = \varphi_{q} (0) = E_{0} \big(M_{t} \big(\exp (- A_{\infty}^{(q)})\big)$ where $\varphi_{q}$ is defined by (3.5.33)).  

\smallskip

\noi Throughout the preceding discussion, a particular role was played by the point $x=0$. However, since any L\'evy process enjoys a property of invariance by translation, we may define, for every $x \in \mathbb{R}$ the $\sigma$-finite measure ${\bf P}_{x}$ by the formula : 
$$    E_{{\bf P}_{x}} \big[F (X_{t}, \; t \ge 0)\big] = E_{\bf P} \big[F (x+X_{t}, \; t \ge 0)\big]        $$

\noi for every positive measurable functional ; thus, the knowledge of ${\bf P}$ induces that of ${\bf P}_{x}$, for any $x \neq 0$.

\noi The reader will have noticed the quasi complete analogy between, on one hand, the results of [YYY] which we just described in the set-up of the $\alpha$-stable symmetric L\'evy process with $1 < \alpha \le 2$ and the results of Chapter 1 of this monograph, in the Brownian set-up. We refer the interested reader to [YYY] where the proofs of the results announced above are found, as well as many other informations.


\newpage

\noi {\bf{\Large Chapter 4. An analogue of ${\bf W}$ for discrete Markov chains.}}

\smallskip

\noi {\large \bf 4.0 Introduction.}

\smallskip

\noi In this chapter, we construct for Markov chains some $\sigma$-finite 
measures which enjoy similar properties as the measure $\bf{W}$ studied in Chapter 1. Very informally, these $\sigma$-finite measures are obtained by "conditioning a recurrent Markov
process to be transient". 

\smallskip

\noi Our construction applies to discrete versions of one- and two-dimensional
Brownian motion, i.e. simple random walk on $\mathbb{Z}$ and $\mathbb{Z}^2$,
but it can also be applied to a much larger class of Markov chains. 

\smallskip

\noi This chapter is divided into three sections; in Section 4.1, we give the construction of the $\sigma$-finite measures mentioned above ; in Section 4.2, we study the main properties of these measures, and in
Section 4.3,  we study some examples in more details.

\bigskip

\noi {\large \bf 4.1 Construction of the $\sigma$-finite measures $(\mathbb{Q}_x, x \in E)$}

\smallskip

\noi {\bf 4.1.1} \underline{Notation and hypothesis.}

\noi Let $E$ be a countable set, $(X_{n})_{n \geq 0}$ the canonical process on $E^
{\mathbb{N}}$, $(\mathcal{F}_n)_ {n \geq 0}$ its natural filtration, and $\mathcal{F}_{\infty}$ the $\sigma$-field generated by  $(X_n)_{n \geq 0}$. 

\noi Let us denote by $(\mathbb{P}_x)_{x \in E}$ the family of probability measures  on  $(E^{\mathbb{N}}, (\mathcal{F}_n)_{n \geq 0}, \mathcal{F}_{\infty})$ associated to a Markov chain ($\mathbb{E}_{x}$ below denotes the expectation with respect to $\mathbb{P}_{x})$ ; more precisely, we suppose there exist probability transitions $(p_{y,z})_{y, z \in E}$ such that :
$$
 \mathbb{P}_{x} (X_{0} = x_{0}, X_1=x_{1},..., X_k=x_k) = {\bf{1}}_{x_0 = x} p_{x_0,x_1}  p_{x_1,x_2} ... p_{x_{k-1},x_k} \eqno(4.1.1)$$

\noi for all $k \geq 0$, $x_0, x_1,..., x_k \in E$. 

\smallskip

\noi We assume three more hypotheses :
\begin{itemize}
\item For all $x \in E$, the set of $y \in E$ such that $p_{x,y} > 0$ is finite (i.e. the graph associated to the Markov chain is locally finite).
\item For all $x, y \in E$, there exists $n \in \mathbb{N}$ such that $\mathbb{P}_x (X_n = y) > 0$ (i.e. the graph of the Markov chain is connected). 
\item For all $x \in E$, the canonical process is recurrent under the probability $\mathbb{P}_x$.
\end{itemize}

\smallskip

\noi {\bf 4.1.2} \underline{A family of new measures.}

\noi From the family of probabilities $(\mathbb{P}_x)_{x \in E}$, we will construct families of $\sigma$-finite measures which should be informally considered to be the law of $(X_n)_{n \geq 0}$ under $\mathbb{P}_x$, after conditionning this process to be transient. 

\noi More precisely, let us fix a point $x_0 \in E$ and let us suppose there exists a function $\phi : E \rightarrow \mathbb{R}_{+}$ such that :
\begin{itemize}
\item $\phi(x) \geq 0$ for all $x \in E$, and $\phi (x_0) = 0$.
\item $\phi$ is harmonic with respect to $\mathbb{P}$, except at the point $x_0$, i.e. :

for all $x \neq x_0$, $\underset{y \in E}{\sum} p_{x,y} \phi(y) = \mathbb{E}_{x} [\phi (X_1)] = \phi(x)$.
\item $\phi$ is unbounded.
\end{itemize}

\noi As we will see in Section 4.2 (Lemma 4.2.9), if $\phi$ satisfies the two first conditions,
the third one is equivalent to the following (a priori weaker):
\begin{itemize}
\item $\phi$ is not identically zero. 
\end{itemize}

\noi In Section 4.3 (Proposition 4.3.1), we give some sufficient conditions for the existence of $\phi$. We
also study some examples. Generally, $\phi$ is not unique, but it will be fixed in this section. For any $r \in ]0,1[$, let us define: 
$$
\psi_r (x) = \frac{r}{1-r} \mathbb{E}_{x_0} [\phi(X_1)] + \phi(x). \eqno(4.1.2)$$
 
\noi From this definition, the following properties hold : 
\begin{itemize}
\item For all $x \neq x_0$, $\psi_r (x) =\mathbb{E}_{x} [\psi_r (X_1)]$. $\hfill(4.1.3)$
\item $\psi_r (x_0) = r \mathbb{E}_{x_0} [\psi_r (X_1)]$ $\hfill(4.1.4)$
\end{itemize}

\noi Now, for $y \in E$ and $k \geq -1$, let us denote by $L_k^{y}$ the local time of $X$ at point $y$ and time $k$, i.e. : 
$$        L_k^{y} =\sum_{m=0}^{k} {\bf 1}_{X_m = y} \eqno(4.1.5)$$

\noi (in particular, $L_{-1}^{y} = 0$ and $L_{0}^y = {\bf 1}_{X_0 = y}$). The properties of $\psi_r$ imply the following result :

\smallskip

\noi{\bf Proposition 4.1.1} 
{\it For every $x \in E$, $(\psi_r (X_n) r^{L_{n-1}^{x_0}},  n \geq 0)$ is a martingale under $\mathbb{P}_x$.}

\smallskip

\noi {\bf Proof of Proposition 4.1.1} For every $n \geq 0$, by Markov property : 
\begin{eqnarray*}
 \mathbb{E}_x \left[ \psi_r (X_{n+1}) r^{L_n^{x_0}} | \mathcal{F}_n \right] 
 & & = r^{L_n^{x_0}} \mathbb{E}_x [\psi_r (X_{n+1}) | \mathcal{F}_n ]  \\
& & = r^{L_n^{x_0}} \psi_r (X_{n}) \left( {\bf{1}}_{X_n \neq x_0} +
\frac{1}{r} {\bf{1}}_{X_n = x_0} \right) = r^{L_{n-1}^{x_0}} \psi_r (X_{n}). \hspace*{0,5cm}(4.1.6)
\end{eqnarray*}
\noi \big(from (4.1.3) and (4.1.4)\big). 

\noi {\bf Corollary 4.1.2}

\noi {\it There exists a finite measure $\mu_x^{(r)}$ on $(E^{\mathbb{N}}, \mathcal{F}_{\infty})$ such that :
$$
\mu^{(r)}_{x \, |\mathcal{F}_n} =  \psi_r (X_{n}) r^{L_{n-1}^{x_0}} \, .\mathbb{P}_{x \, |\mathcal{F}_n} \eqno(4.1.7)$$ }

\noi At this point, we remark that, for all $\sigma$, $0 < \sigma < 1/r$ :
\begin{itemize}
\item $ \psi_r (x) \leq \sup \left( \frac{1-\sigma r}{\sigma (1-r)}, 1\right) . \psi_{\sigma r} (x)$ for all $x \in E$.
\item Consequently, for $n \geq 1$ : 
\begin{eqnarray*} 
& & \mu_x^{(r)} (\sigma^{L_{n-1}^{x_0}} ) = \mathbb{P}_x [ \psi_r (X_n) (r\sigma)^{L_{n-1}^{x_0}} ] \quad \text{\big(from (4.1.7)\big)} \\ 
& & \leq  \sup \left( \frac{1-\sigma r}{\sigma (1-r)}, 1 \right) \mathbb{P}_x [ \psi_{\sigma r} (X_n) (r \sigma)^{L_{n-1}^{x_0}} ]  \\ 
& & \leq \sup \left( \frac{1-\sigma r}{\sigma (1-r)}, 1 \right) \mu_x^{(\sigma r)} (1) = C \hspace*{6,4cm}(4.1.8)
 \end{eqnarray*}
 where $C < \infty$ does not depend on $n$. 
 \end{itemize}

\noi Therefore, $\mu_x^{(r)} (\sigma^{L_{\infty}^{x_0}} ) <\infty $, with
$$L_{\infty}^{x_0} := \sum_{m=0}^{\infty} {\bf 1}_{X_m = x_0} = \underset{k \rightarrow \infty}{\lim}
 L_k^{x_0}.$$
\noi In particular, $L_{\infty}^{x_0} < \infty$,
  $\mu_x^{(r)}$-a.s. It is now possible to define a measure $\mathbb{Q}_x^{(r)}$, by :
$\mathbb{Q}_x^{(r)} = \left( \frac{1}{r} \right)^{L_{\infty}^{x_0}} . \mu_x^{(r)}$; this measure is
$\sigma$-finite since the sets $ \{ L_{\infty}^{x_0} \leq m \}$ increase to $\{ L_{\infty}^{x_0} < \infty\}$; moreover $\{L_{\infty}^{x_0} = \infty\}$ is $ \mathbb{Q}_x^{(r)}$-negligible, and
$$
 \mathbb{Q}_x^{(r)}  (  L_{\infty}^{x_0} \leq m) \leq \left(\frac{1}{r} \right)^m \mu_x^{(r)} (1) < \infty \eqno(4.1.9) $$

\noi {\bf 4.1.3} \underline{Definition of the measures $(\mathbb{Q}_{x}, \; x \in E)$.}
  
\noi Here is a remarkable result, which explains the interest of this construction : 
 
 \smallskip
 
 \noi {\bf Theorem 4.1.3} {\it The two following properties hold :

\noi {\it i)} For all $x \in E$, $\mathbb{Q}_x^{(r)}$ does not depend on $r \in ]0,1[$. 
 
\noi {\it ii)} Let $\mathbb{Q}_x$ denote the measure equal to $\mathbb{Q}_x^{(r)}$ for all $r \in ]0,1[$, and $F_n \geq 0$ a $\mathcal{F}_n$-measurable functional. If $q$ is a function from $E$ to $[0,1]$, such that $\{ q < 1 \}$ is a finite set, then :
 $$ \mathbb{Q}_x \left[ F_n \prod_{k=0}^{\infty} q(X_k) \right]  = \mathbb{E}_x  \left[ F_n  \psi_q (X_n)\prod_{k=0}^{n-1} q(X_k)  \right]  \eqno(4.1.10) $$

\noi  where for $y \in E$, $\dis \psi_q (y) := \mathbb{Q}_y \left[\prod_{k=0}^{\infty} q(X_k) \right]$. $\hfill (4.1.11)$ }

\smallskip

 \noi {\bf Remark 4.1.4} If we denote by $\mu_x^{(q)}$ the measure defined by : 
$$    \mu_x^{(q)} = \left( \prod_{k=0}^{\infty} q(X_k) \right) . \mathbb{Q}_x    \eqno(4.1.12)$$

\noi we obtain : 
$$
\mu^{(q)}_{x \, | \mathcal{F}_n} = \psi_q (X_n)  \left(\prod_{k=0}^{n-1} q(X_k) \right) .
\mathbb{P}_{x \, | \mathcal{F}_n}.  \eqno(4.1.13)$$

\noi These relations are similar to relations between ${\bf W}$ and Feynman-Kac penalisations of Wiener measure $W$ \big(see Chap. 1, Th. 1.1.2, formulae (1.1.7), (1.1.8), (1.1.16)\big). 

\noi Moreover, $\psi_q$ satisfies the "Sturm-Liouville equation" :
$$    \psi_q (x) = q(x) \mathbb{E}_x [\psi_q (X_1)]    \eqno(4.1.14) $$

\noi  The analogy between this situation and the Brownian case described in Chapter 1
 can be represented by the following correspondance :

\smallskip 

\begin{tabular}{|c|c|}
\hline
Markov chain & Brownian motion \\
\hline
$\mathbb{P}_{x_0}$ & $W_0$ \\
$\mathbb{P}_x$ & $W_x$ \\
$\mu_x^{(q)}$ & $W_{x, \infty}^{(q)}$ \\
$M_n^{(q)} = \psi_q (X_n) \prod_{k=0}^{n-1} q(X_k)$ & $M_t^{(q)} = \frac{\varphi_q (X_t)}{\varphi_q(x)}
\exp \left( - \frac{1}{2} A_t^{(q)} \right)$  \\
$\psi_q (x) = q(x) \mathbb{E}_x (\psi_q(X_1))$ & $\varphi''_q (x) = q(x) \varphi_q (x)$ \\
$\mu_x^{(q)} \,_{|\mathcal{F}_n} = M_n^{(q)}. \mathbb{P}_x \,_{|\mathcal{F}_n} $ & $ W_{x, \infty}^{(q)} \,_
{| \mathcal{F}_t} = M_t^{(q)} . W_x \,_{| \mathcal{F}_t}  $ \\
$ \mathbb{Q}_x $ &  ${\bf W}_x$ \\
$\mu_x^{(q)} = \left( \prod_{k=0}^{\infty} q(X_k) \right). \mathbb{Q}_x$ & $W_{x, \infty}^{(q)} = 
\frac{1}{\varphi_q(x)} \, \exp \left( - \frac{1}{2} A_{\infty}^{(q)} \right) . {\bf W}_x $ \\
\hline
\end{tabular}

\smallskip

\noi {\bf Proof of Theorem 4.1.3} To begin with, let us prove the point {\it ii)} (with $\mathbb{Q}_x^{(r)}$ instead of $\mathbb{Q}_x$) for a function $q$ such that $q(x_0) < 1$. Under the hypotheses of Theorem 4.1.3, for all $n \geq 0$, $F_n \prod_{k=0}^{N-1} q(X_k) \left( \frac{1}{r} \right) ^ {L_{N-1}^{x_0}}$
tends to $F_n \prod_{k=0}^{\infty} q(X_k) \left( \frac{1}{r} \right) ^{L_{\infty}^{x_0}}$ as $N \to \infty$ and is dominated by $\left(\frac{q(x_0)}{r} \vee 1\right)^{L_{\infty}^{x_0}}$, which is $\mu_x^{(r)}$-integrable because $\frac{q(x_0)}{r} \vee 1 < \frac{1}{r} \cdot$ (from (4.1.8)).

\noi By dominated convergence, if for $y \in E$, $k \geq 0$, we define :
$$  \chi_q^{r,k} (y) := \mathbb{E}_y \left[ \psi_r ( X_k) \prod_{m=0}^{k-1} q(X_m) \right], \eqno(4.1.15) $$

\noi for all $x \in E$ :
 \begin{eqnarray*}
\hspace*{1cm}& & \mathbb{E}_x \left[ F_n \, \chi_q^{r,N-n} (X_n) \prod_{k=0}^{n-1} q(X_k) \right]  =  \mathbb{E}_x \left[
F_n \, \psi_r (X_N) \prod_{k=0}^{N-1} q (X_k) \right] \nonumber \\ 
& & = \mu_x^{(r)} \left[ F_n \, \prod_{k=0}^{N-1} q (X_k) \, \left( \frac{1}{r} \right)^{L_{N-1}^{x_0}}
\right]  \\ 
& &\underset{N \rightarrow \infty}{\rightarrow} \mu_x^{(r)} \left[ F_n \, \prod_{k=0}^{\infty} q (X_k) \, \left( \frac{1}{r} \right)^{L_{\infty}^{x_0}} \right] = \mathbb{Q}_x^{(r)} \left[ F_n \,\prod_{k=0}^{\infty} q (X_k) \right].
\hspace*{2cm} (4.1.16)
\end{eqnarray*}

\noi In particular, if we take $n=0$ and $F_0 = 1$ : 
$$
\chi_q^{r,N} (y) \underset{N \rightarrow \infty}{\rightarrow} \mathbb{Q}_y^{(r)} \left[\prod_{k=0}^{\infty} q (X_k)\right] \eqno(4.1.17) $$

\noi for all $y \in E$. 

\noi Moreover : 
\begin{eqnarray*}  \hspace*{2cm}
& & \chi_q^{r,N-n} (y) \leq \mathbb{E}_y \left[ (q (x_0))^{L_{N-n-1}^{x_0}} \psi_r (X_{N-n}) \right]  \\
& & \leq \sup \left( \frac{r}{q(x_0)} \left( \frac{1-q(x_0)}{1-r} \right) ,1 \right) \mathbb{E}_y \left[ (q(x_0))
^ {L_{N-n-1}^{x_0}} \, \psi_{q(x_0)} (X_{N-n}) \right]   \\ 
& & = \sup \left( \frac{r}{q(x_0)} \left( \frac{1-q(x_0)}{1-r} \right) , 1 \right)  \psi_{q(x_0)} (y) \hspace*{3,8cm} (4.1.18)
\end{eqnarray*}

\noi where 
\begin{eqnarray*}
\hspace*{2cm} \mathbb{E}_x \left[ \psi_{q(x_0)} (X_{n}) \, \prod_{k=0}^{n-1} q(X_k) \right] 
&\leq& \mathbb{E}_x \left[ \psi_{q(x_0)} (X_{n}) (q (x_0))^{L_{n-1}^{x_0}} \right]  \\
& =&  \psi_{q(x_0)} (x) < \infty.  \hspace*{4cm} (4.1.19)
\end{eqnarray*}

\noi By dominated convergence : 
$$
\mathbb{E}_x \left[ F_n \,  \chi_q^{r,N-n} (X_n) \prod_{k=0}^{n-1} q(X_k) \right]  \underset{N \rightarrow \infty}{\rightarrow} \mathbb{E}_x \left[ F_n \,\psi_q^{(r)} (X_n) \prod_{k=0}^{n-1} q(X_k) \right], \eqno(4.1.20) $$
                          
\noi where $\psi_q^{(r)} (y) = \mathbb{Q}_y^{(r)} \left[ \prod_{k=0}^{\infty} q(X_k) \right]$.

\noi The two previous limits are equal; therefore :
$$
\mathbb{Q}_x^{(r)} \left[F_n \prod_{k=0}^{\infty} q(X_k) \right] = \mathbb{E}_x \left[ F_n \, \psi_q^{(r)}
(X_n) \prod_{k=0}^{n-1} q(X_k) \right]  , \eqno (4.1.21) $$
                                           
\noi as written in point {\it ii)} of Theorem 4.1.3 (with $\mathbb{Q}_x^{(r)}$ instead of $\mathbb{Q}_x$). 

\noi Now we can prove point {\it i)}, by taking for any $s \in ]0,1[$,  $q(x) = {\bf 1}_{x \neq x_0} + s
{\bf 1}_{x = x_0}$. 

\noi Let us first observe that $\dis \frac{\psi_r (X_n)}{\psi_s (X_n)}$ is $\mu_y^{(s)}$-a.s. well-defined for all $n \geq 0$; therefore, $\dis \mu_y^{(s)} \left[ \frac{\psi_r (X_n)}{\psi_s (X_n)} \right]$ is well-defined and :
 \begin{eqnarray*} \hspace*{2cm}
  & & \mu_y^{(s)} \left[ \frac{\psi_r (X_n)}{\psi_s (X_n)} \right] = \mathbb{E}_y \left[ s^{L_{n-1}^{x_0}}
 \psi_r (X_n) \right] = \mu_y^{(r)} \left[ \left( \frac{s}{r} \right) ^{L_{n-1}^{x_0}} \right] \\
 & & \underset{n \rightarrow \infty}{\rightarrow} \mu_y^{(r)} \left[ \left(\frac{s}{r} \right) ^ {L_{\infty}
 ^{x_0} } \right] = \mathbb{Q}_y^{(r)} [s ^{L_{\infty}^{x_0}}] =\psi_q^{(r)} (y).   \hspace*{4cm}(4.1.22)
 \end{eqnarray*}

\noi Moreover, for all $A > 0$ :
$$
 \mu_y^{(s)} \left[ \frac{\psi_r (X_n)}{\psi_s (X_n)} \right] = \mu_y^{(s)} \left[ \frac{\psi_r  (X_n)}{\psi_s (X_n)} {\bf 1}_ {\psi_s (X_n) \geq A} \right] + K_A, \eqno(4.1.23) $$
 
\noi where : 
 $$
 K_A \leq \sup \left( \frac{\psi_r}{\psi_s} \right) . \mu_y^{(s)} [ \psi_s (X_n) \leq A] \leq A \, \sup \left( \frac{\psi_r}{\psi_s} \right) \mathbb{E}_y [s ^{L_{n-1}^{x_0}}] \underset{n \rightarrow \infty}{\rightarrow} 0, \eqno (4.1.24) $$
\noi (from the definition (4.1.7) of $\mu_y^{(s)}$ and the fact that $(X_n)_{n \geq 0}$ is recurrent under
 $\mathbb{P}_y$). Hence :  
\begin{eqnarray*} \hspace*{2cm}
& & \underset{n \rightarrow \infty}{\lim \inf} \left( \underset{\psi_s (x) \geq A}{\inf} \frac{\psi_r(x)}
{\psi_s(x)} \right) \mu_y^{(s)} [ \psi_s (X_n) \geq A]  \\
 & & \leq \underset{n \rightarrow \infty}{\lim \inf} \, \mu_y^{(s)} \left[ \frac{\psi_r (X_n)}{\psi_s (X_n)} \right] \leq \underset{n \rightarrow \infty}{\lim \sup} \, \mu_y^{(s)} \left[ \frac{\psi_r (X_n)}{\psi_s (X_n)} \right] \\
 & & \leq \underset{n \rightarrow \infty}{\lim \sup} \left( \underset{\psi_s (x) \geq A}{\sup} \frac{\psi_r(x)}
{\psi_s(x)} \right) \mu_y^{(s)} [ \psi_s (X_n) \geq A]. \hspace*{4cm} (4.1.25)
\end{eqnarray*}

\noi Now, since $\phi$ (and hence, $\psi_s$) is unbounded, $\dis \underset{\psi_s (x) \geq A}{\inf} \frac{\psi_r(x)}{\psi_s(x)} $ and $\dis \underset{\psi_s (x) \geq A}{\sup} \frac{\psi_r(x)}{\psi_s(x)} $ tend to 1 when $A$ goes to infinity and :
$$
\mu_y^{(s)} [\psi_s (X_n) \geq A] \rightarrow \mu_y^{(s)} (1) = \psi_s (y). \eqno(4.1.26) $$

\noi Hence, $\dis \mu_y^{(s)} \left[ \frac{\psi_r (X_n)}{\psi_s (X_n)} \right] \underset{n \rightarrow \infty}
{\rightarrow} \psi_s (y)$, which implies that $\psi_q^{(r)} (y) = \psi_s(y)$. 

\noi By (4.1.21) : 
\begin{eqnarray*}
\hspace*{2cm}  \mathbb{Q}_x^{(r)} [F_n s^{L_{\infty}^{x_0}} ] 
&& = \mathbb{E}_x \left[ F_n s^{L_{n-1}^{x_0}} \psi_q^{(r)} (X_n) \right] = \mathbb{E}_x \left[ F_n s^{L_{n-1}^{x_0}} \psi_s (X_n) \right] \\
& & = \mu_x^{(s)} (F_n) = \mathbb{Q}_x^{(s)} [F_n s^{L_{\infty}^{x_0}} ].\hspace*{4cm} (4.1.27)
\end{eqnarray*}

 \noi By monotone class theorem, if $F$ is $\mathcal{F}_{\infty}$-measurable and positive :
 $$
 \mathbb{Q}_x^{(r)} (F . s^{L_{\infty}^{x_0}} ) = \mathbb{Q}_x^{(s)} (F . s^{L_{\infty}^{x_0}} ) \eqno(4.1.28)
 $$
 
\noi  for all $r, s \in ]0,1[$. Now, for all $r, s, t < 1$ :
$$
 \mathbb{Q}_x^{(r)} (F . t^{L_{\infty}^{x_0}} ) = \mathbb{Q}_x^{(t)} (F . t^{L_{\infty}^{x_0}} ) =
 \mathbb{Q}_x^{(s)} (F . t^{L_{\infty}^{x_0}} ). \eqno(4.1.29) $$

\noi Recall that $L_{\infty}^{x_0} < \infty$, $\mathbb{Q}_x^{(r)}$ and $\mathbb{Q}_x^{(s)}$-a.s. Therefore, by monotone convergence, $\mathbb{Q}_x^{(r)} (F) =  \mathbb{Q}_x^{(s)} (F)$ ; point {\it i)} of Theorem 4.1.3 is proven, and $\mathbb{Q}_x$ is well-defined. By (4.1.21), point {\it ii)} is proven if $q(x_0) < 1$. It is easy to extend this formula to the case $q(x_0)= 1$, again by monotone convergence ; the proof of Theorem 4.1.3 is now complete. $\hfill \blacksquare$

\bigskip

\noi {\bf Remark 4.1.5} The family $(\mathbb{Q}_x)_{x \in E}$ of $\sigma$-finite measures depends
on $x_0$ and $\phi$, which were assumed to be fixed in this section. In the sequel of the chapter, these parameters may vary; if some confusion is possible, we will write $(\mathbb{Q}_x^{(\phi,x_0)})_{x \in E}$ instead of $(\mathbb{Q}_x)_{x \in E}$.

\bigskip

\noi {\large \bf 4.2 Some more properties of $(\mathbb{Q}_x, x \in E)$.}

\smallskip

\noi {\bf 4.2.1} \underline{Martingales associated with $(\mathbb{Q}_{x}, x \in E)$.}

\noi At the beginning of this section, we extend the second point of Theorem 4.1.3 to more general functionals than functionals of the form $ \dis F_n \prod_{k=0}^{\infty} q(X_k)$. More precisely, the following result holds :

\noi {\bf Theorem 4.2.1}

\noi {\it Let $F$ be a positive $\mathcal{F}_{\infty}$-measurable functional. For $n \geq 0$, $y_0, y_1,..., y_n \in E$, let us define the quantity : 
$$
M(F, y_0, y_1,..., y_n) := \mathbb{Q}_{y_n} \left[ F(y_0, y_1,..., y_n=X_0, X_1, X_2,...) \right]. \eqno(4.2.1) $$

\noi Then, for every $(\mathcal{F}_n)_{n \geq 0}$-stopping time $T$, one has : 
$$
 \mathbb{Q}_x (F . {\bf 1}_{T < \infty}) = \mathbb{E}_x \left[ M(F, X_0, X_1,..., X_T) {\bf 1}_{T
 < \infty} \right]. \eqno(4.2.2) $$   }

\noi {\bf Proof of Theorem 4.2.1: } To begin with, let us suppose that $T=n$ for $n \geq 0$, and $F=r^{L_{\infty}^{x_0}} f_0 (X_0) f_1 (X_1) ... f_N (X_N)$ for $N>n$, $0 \leq f_i \leq 1$, $0<r<1$.

\noi One has : 
 \begin{eqnarray*}
\hspace*{1cm} & \mathbb{Q}_x (F) & = \mu_x^{(r)} \left[ f_0 (X_0) ... f_N (X_N) \right] \\
 & & = \mathbb{E}_x \left[ f_0 (X_0) ... f_N (X_N) r^{L_{N-1}^{x_0}} \psi_r (X_N)  \right] \hspace*{5cm} (4.2.3) \\
 & & = \mathbb{E}_x \left[ f_0 (X_0) ... f_{n-1} (X_{n-1}) r^{L_{n-1}^{x_0}} K(X_n)  \right],
 \end{eqnarray*}
 
\noi where : 
\hspace*{1cm} \begin{eqnarray*}
 & K(y) & = \mathbb{E}_y \left[ f_n (X_0) ... f_N (X_{N-n}) r^{L_{N-n-1}^{x_0}} \psi_r (X_{N-n}) \right]   \\
 & & = \mu_y^{(r)} \left[ f_n (X_0) ... f_N (X_{N-n}) \right] \hspace*{8cm} (4.2.4) \\
 & & = \mathbb{Q}_y \left[ f_n (X_0) ... f_N (X_{N-n}) r^{L_{\infty}^{x_0}}\right]. 
 \end{eqnarray*}
 
\noi Hence, for all $y_0,...,y_n$ :
 \begin{eqnarray*}
 \hspace*{1cm}& & f_0 (y_0) ... f_{n-1} (y_{n-1}) r^{ \sum_{k=0}^{ n-1} {\bf 1}_{y_k = x_0} } K(y_n)  \\
 & & = \mathbb{Q}_{y_n} \left[ f_0 (y_0) ... f_{n-1} (y_{n-1}) f_n(X_0) ... f_N (X_{N-n}) r^{\sum_{k=0}^{
 n-1} {\bf 1}_{y_k = x_0} + L_{\infty}^{x_0} } \right]  \hspace*{1cm} (4.2.5)\\
 & & = \mathbb{Q}_{y_n} \left[ F(y_0,...,y_n = X_0, X_1, ...) \right] = M(F,y_0,y_1,...,y_n).
 \end{eqnarray*}

 \noi Therefore : 
$$
 \mathbb{Q}_x (F) = \mathbb{E}_x \left[M (F,  X_0,..., X_n) \right],  \eqno(4.2.6)$$
 
\noi which proves Theorem 4.2.1 for these particular functionals $F$ and for $T=n$. 

\noi By monotone class theorem, we can extend (4.2.6) to the functionals $F=r^{L_{\infty}^{x_0}}. G$,  where $G$ is any positive functional, and by monotone convergence ($r$ increasing to $1$), Theorem 4.2.1 is proven for all $F$ and $T=n$. 

\noi Now, let us suppose that $T$ is a stopping time. 

\noi For $n \geq 0$, $M(F {\bf 1}_{T=n}, X_0, X_1,..., X_n ) = {\bf 1}_{T=n} M(F, X_0,..., X_n)$, because $\{ T =n \}$ depends only on $X_0, X_1, ..., X_n$; hence,
$$
\mathbb{Q}_x (F {\bf 1}_{T=n}) = \mathbb{E}_x \left[ {\bf1}_{T=n} M(F, X_0,..., X_n) \right]. \eqno(4.2.7) $$

\noi Summing from $n=0$ to infinity, we obtain the general case of Theorem 4.2.1. $\hfill \blacksquare$

\bigskip

\noi{\bf Corollary 4.2.2} {\it For any functional $F \in L^1 (\mathbb{Q}_x)$, $\dis \left( M(F, X_0, X_1,..., X_n) \right)_{n \geq 0}$ is a $\mathcal{F}_n$-martingale (with expectation
$\mathbb{Q}_x (F)$).  }

\noi The correspondance with the Brownian case is the following :

\smallskip 

\begin{tabular}{|c|c|}
\hline
Markov chain & Brownian motion \\
\hline
$F \in L_+^1 (\mathbb{Q}_x, \mathcal{F}_{\infty})$ & $F \in L_+^1 ({\bf W}_x, \mathcal{F}_{\infty})$ \\
\hline
$(M( F,X_0,...,X_n), \, n \geq 0) $ & $ (M_t (F), \, t \geq 0)$ a $(\mathcal{F}_t, \, t \geq 0, W_x)$ \\
a $(\mathcal{F}_n, \, n \geq 0, \mathbb{P}_x)$ martingale such that & martingale such that \\
$(*) \quad \mathbb{Q}_x [\Gamma_n F] = \mathbb{P}_x [\Gamma_n M(F, X_0,...,X_n)]$ ($\Gamma_n \in \mathcal{F}_n$)
& $ {\bf W}_x [\Gamma_t F] = W_x [ \Gamma_t M_t (F)]$ ($\Gamma_t \in \mathcal{F}_t$) \\
$\mathbb{Q}_x (F) = \mathbb{P}_x [M (F, X_0,..., X_n)]$ & ${\bf W}_x (F) = W_x (M_t(F))$ \\ 
\hline
\end{tabular}

\smallskip

\noi Here, $(*)$ is a consequence of (4.2.2) with $T = n. {\bf 1}_{\Lambda_n} + (+ \infty). {\bf 1}_{\Lambda_n^c}$. 
 
\noi Now, we are able to describe the properties of the canonical process under $\mathbb{Q}_x$. 

\smallskip

\noi {\bf 4.2.2.} \underline{Properties of the canonical process under $(\mathbb{Q}_{x}, x \in E)$.}

\noi We have already proven that $L_{\infty}^{x_0}$ is almost surely finite under $\mathbb{Q}_x$. In fact, the following proposition gives a more general result :

\smallskip

\noi {\bf Proposition 4.2.3} {\it Under $\mathbb{Q}_x$, the canonical process is a.s. transient, i.e $L_{\infty}^{y_0} < \infty$ for all $y_0 \in E$.}

\smallskip

\noi {\bf Proof of Proposition 4.2.3: } Let $y_0$ be in $E$, and $r$ be in $]0,1[$. If, for $k \geq 1$, $\tau_k^{(y_0)}$ denotes the $k$-th hitting time of $y_0$ for the canonical process $X$, then
for all $n \geq 0$ :
\begin{eqnarray*}
 \hspace*{1cm}& \mu_x^{(r)}[ L_{n-1}^{y_0} \geq k] & = \mu_x^{(r)} [\tau_k^{(y_0)} < n] = \mathbb{E}_x \left[ {\bf 1}_{  \tau_k^{(y_0)} < n} \, r^{L_{n-1}^{x_0}} \, \psi_r (X_n) \right] \\
& & = \mathbb{E}_x \left[ {\bf 1}_{  \tau_k^{(y_0)} < n} r^{L_{\tau_k^{(y_0)} -1}^{x_0}} \psi_r (y_0) \right]  \hspace*{4cm} (4.2.8)
\end{eqnarray*}

\noi by strong Markov property (applied at time $ \tau_{k}^{(y_{0})} \wedge n$), and by the fact that $\mathbb{E}_{y_{0}} [r^{L_{m-1}^{x_0}} \psi_r (X_m) ] =$ 

\noi $ \psi_r (y_0)$ for all $m \ge 0$ (from Proposition 4.1.1). 

\noi Hence : 
$$
\mu_x^{(r)} [L_{n-1}^{y_0}  \geq k ] \leq \psi_r (y_0) \mathbb{E}_x \left[r^{L_{ \tau_k^{(y_0)} -1}^{x_0}}
\right];  \eqno(4.2.9)$$

\noi and by monotone convergence :
$$
\mu_x^{(r)} [L_{\infty}^{y_0}  \geq k ] \leq \psi_r (y_0) \mathbb{E}_x\left[ r^{L_{ \tau_k^{(y_0)} -1}^{x_0}}
\right] \underset{k \rightarrow \infty}{\rightarrow} 0  \eqno(4.2.10)$$

\noi (since  $L_{ \tau_k^{(y_0)} }^{x_0} \underset{ k \rightarrow \infty}{\rightarrow} \infty$, $\mathbb{P}_x$-a.s.); this implies Proposition 4.2.3. $\hfill \blacksquare$

\bigskip

\noi Now, we have the following decomposition result which gives a precise description of the canonical process under $\mathbb{Q}_{y} \; (y \in E)$ :

\smallskip

\noi {\bf Proposition 4.2.4} {\it For all $y, y_0 \in E$, one has :
$$
\mathbb{Q}_y = \mathbb{Q}_y^{[y_0]} + \sum_{k \geq 1} \mathbb{P}_y^{\tau_k^{(y_0)}} \circ \widetilde{\mathbb{Q}}_{y_0}, \eqno(4.2.11) $$        

\noi where $ \dis \mathbb{Q}_y^{[y_0]} = {\bf 1}_{\forall n \geq 0, X_n \neq y_0} \mathbb{Q}_y $ is the restriction of $\mathbb{Q}_y$ to trajectories which do not hit $y_0$, $\widetilde{\mathbb{Q}}_{y_0} = {\bf 1}_{\forall n \geq 1, X_n \neq y_0} \mathbb{Q}_{y_0}$ is the restriction of $\mathbb{Q}_{y_0}$ to trajectories which do not return to $y_0$, and $\mathbb{P}_y^{\tau_k^{(y_0)}} \circ \widetilde{\mathbb{Q}}_{y_0}$ denotes the concatenation of $\mathbb{P}_y$ stopped at time $\tau_k^{(y_0)}$ and $ \widetilde{\mathbb{Q}}_{y_0}$, i.e. the image of $\mathbb{P}_y \otimes \widetilde{\mathbb{Q}}_{y_0}$ by the functional $\Phi$ from $E^{\mathbb{N}} \times E^{\mathbb{N}}$ such that : 
$$
\Phi ((z_0,z_1,...,z_n,...), (z'_0,z'_1,...,z'_n,...)) = (z_0,z_1,...,z_{\tau_k^{(y_0)}}, z'_1,...,z'_n). \eqno(4.2.12) $$    }
\noi \noi This formula (4.2.11) can be compared to (3.2.20) or (1.1.40). 

\noi {\bf Proof of Proposition 4.2.4 : } We apply Theorem 4.2.1 to the stopping time $T = \tau_k^{(y_0)}$, and to the functional :
$$
F=G H(X_{\tau_k^{(y_0)}}, X_{\tau_k^{(y_0)} + 1},...) {\bf 1}_{\forall u \geq 1, X_{\tau_k^{(y_0)} + u}
\neq y_0},  \eqno(4.2.13)$$

\noi where $G$, $H$ are positive functionals such that $G \in \mathcal{F}_{\tau_k^{(y_0)}}$. 

\noi For $k \geq 1$, we obtain :
\begin{eqnarray*}
 \hspace*{1cm} & & \mathbb{Q}_y \left[ G H(X_{\tau_k^{(y_0)}}, X_{\tau_k^{(y_0)} + 1},...) {\bf 1}_{L_{\infty}^{y_0} = k} \right]  \\ 
& & = \mathbb{E}_y \left[ {\bf 1}_ {\tau_k^{(y_0)} < \infty} G(X_0,...,X_{ \tau_k^{(y_0)}} ) \right] \, \widetilde{\mathbb{Q}}_{y_0} [H],  \hspace*{5cm}(4.2.14)
\end{eqnarray*}

\noi which implies : 
$$ \mathbb{Q}_y \left[ G H(X_{\tau_k^{(y_0)}}, X_{\tau_k^{(y_0)} + 1},...) {\bf 1}_{L_{\infty}^{y_0} = k} \right] = \mathbb{E}_y[G] \, \widetilde{\mathbb{Q}}_{y_0} [H], \eqno(4.2.15)$$

\noi because $\tau_k^{(y_0)} < \infty$, $\mathbb{P}_y$-a.s. (the canonical process is recurrent under $\mathbb{P}_y$). Moreover : 
$$
 \mathbb{Q}_y [H {\bf 1}_{L_{\infty}^{y_0} = 0} ] = \mathbb{Q}_y^{[y_0]}
(H) \eqno(4.2.16) $$

\noi by definition.  Now, $L_{\infty}^{y_0} < \infty$, $\mathbb{Q}_y$-a.s. by Proposition 4.2.3, so there exists $k \geq 0$ such that $L_{\infty}^{y_0} = k$ : the equalities (4.2.15) and (4.2.16) imply the Proposition 4.2.4 by monotone class theorem. $\hfill \blacksquare$

\bigskip

\noi {\bf 4.2.3} \underline{Dependence of $\mathbb{Q}_{x}$ on $x_0$.}

\noi The next Theorem shows that in the construction of the family $(\mathbb{Q}_x)_{x \in E}$, the choice of the point $x_0$ in $E$ is in fact not so important. More precisely, the following result holds :

\smallskip

\noi {\bf Theorem 4.2.5.} {\it For all $y_0 \in E$, let us define the function $\phi^{[y_0]}$ by :
$$       \phi^{[y_0]} (y) = \mathbb{Q}_y^{[y_0]} (1) \eqno(4.2.17) $$

\noi Then the following holds :

\noi {\it i)} $\phi^{[x_0]}$ is equal to $\phi$ and for all $y_0 \in E$,
$\phi^{[y_0]} - \phi$ is a bounded function.

\noi {\it ii)} For all $y_0 \in E$ :
\begin{itemize}
\item $\phi^{[y_0]}$ is finite and harmonic outside of $y_0$, i.e. for all $y \neq y_0$ : 
\begin{equation*}
\mathbb{E}_y [ \phi^{[y_0]} (X_1) ] = \phi^{[y_0]} (y).
\end{equation*}
\item $\phi^{[y_0]}  (y_0) = 0$.
\item $\widetilde{\mathbb{Q}}_{y_0} (1) = \mathbb{E}_{y_0} [ \phi^{[y_0]}
(X_1) ]$.
\end{itemize}
\noi {\it iii)} By point {\it ii)}, $y_0$ and the function $\phi^{[y_0]}$ can be used to construct a family $(\mathbb{Q}_x^{(\phi^{[y_0]}, y_0)})_{x \in E}$ of $\sigma$-finite measures by the method given in Section 4.1. Moreover, this family is equal to the family $(\mathbb{Q}_x = \mathbb{Q}_x^{(\phi, x_0)})_{x \in E}$ constructed with $\phi$ and $x_0$. 

\noi {\it iv)} For all $y_0, y \in E$, the image of the measure $\mathbb{Q}_y$ by the total local time at $y_0$ is given by the following expressions :
\begin{itemize}
\item $\mathbb{Q}_y [L_{\infty}^{y_0} = 0] = \phi^{[y_0]} (y)$.
\item For all $k \geq 1$, $\mathbb{Q}_y [L_{\infty}^{y_0} = k] = \mathbb{E}_{y_0} [\phi^{[y_0]} (X_1) ]$.
\end{itemize}   }

\smallskip

\noi {\bf Proof of Theorem 4.2.5. } Let $y_0$ and $y$ be in $E$. For all $r \in ]0,1[$, $n \geq 1$ :
\begin{eqnarray*}
\hspace*{1cm}  & \mu_y^{(r)} [L_{n-1}^{y_0} \geq 1] & = \mu_y^{(r)} [ \tau_1^{(y_0)} < n] =
\mathbb{E}_y  \left[ r^{L_{n-1}^{x_0}} . {\bf 1}_{\tau_1^{(y_0)} < n} . \psi_r (X_n) \right]  \\
& & = \mathbb{E}_y \left[ r^{L^{x_0}_{\tau_1^{(y_0)} -1}} . {\bf 1}_{\tau_1^{(y_0)} < n} \right] \, \psi_r(y_0)              \hspace*{4,5cm} (4.2.19)
\end{eqnarray*}
from (4.1.7) and the martingale property.  Hence : 
$$
\mu_y^{(r)} [L_{\infty}^{y_0} \geq 1] = \psi_r(y_0) \,\mathbb{E}_y \left[ r^{L^{x_0}_{\tau_1^{(y_0)} -1}} \right]. \eqno(4.2.20) $$

\noi If $y_0 = x_0$, this implies :
$$          \mu_y^{(r)} [L_{\infty}^{x_0} \geq 1] = \psi_r(x_0)   \eqno(4.2.21)   $$

\noi Therefore : 
\begin{eqnarray*}
\hspace*{1cm}  & \phi^{[x_0]} (y) & = \mathbb{Q}_y [ L_{\infty}^{x_0} = 0] = \mu_y^{(r)}
[L_{\infty}^{x_0} = 0]  \\
& & = \mu_y^{(r)} (1) -  \psi_r(x_0) = \psi_r (y) - \psi_r(x_0) = \phi(y) \hspace*{4cm} (4.2.22)
\end{eqnarray*}

\noi as written in Theorem 4.2.5. If $y_0 \neq x_0$, let us define the quantities :
$$
 p_{y,y_0}^{(x_0)} = \mathbb{P}_y [\tau_1^{y_0} < \tau_1^{x_0}], \eqno(4.2.23) $$

\noi and
$$
q_{y_0}^{(x_0)} = \mathbb{P}_{x_0} [ \tau_1^{y_0} > \tau_2^{x_0}]. \eqno(4.2.24)   $$

\noi We have : 
$$       \mathbb{P}_y \left[L^{x_0}_{\tau_1^{(y_0)} -1} = 0 \right] =p_{y,y_0}^{(x_0)}    \eqno(4.2.25) $$

\noi and, for $k \geq 1$, by strong Markov property :
$$
\mathbb{P}_y \left[L^{x_0}_{\tau_1^{(y_0)} -1} = k \right] = (1-p_{y,y_0}^{(x_0)})(q_{y_0}^{(x_0)})^{k-1} (1-q_{y_0}^{(x_0)}) \eqno(4.2.26) $$

\noi Summing all these equalities, one obtains :
$$
\mathbb{E}_y \left[ r^{L^{x_0}_{\tau_1^{(y_0)} -1}} \right] = p_{y,y_0}^{(x_0)} + \frac{r (1-p_{y,y_0}^{(x_0)}) (1-q_{y_0}^{(x_0)})}{ 1 - r q_{y_0}^{(x_0)}} \eqno(4.2.27) $$

\noi and from (4.2.21) and (4.2.27) :
\begin{eqnarray*}
\hspace*{1cm} & \mu_y^{(r)} [L_{\infty}^{y_0} \geq 1] & = \left[ \frac{r}{1-r} \mathbb{E}_{x_0} [\phi(X_1)] + \phi(y_0) \right] \\ 
& & \times \left[ p_{y,y_0}^{(x_0)} + \frac{r(1-p_{y,y_0}^{(x_0)}) (1-q_{y_0}^{(x_0)})}{ 1 - r q_{y_0}^{(x_0)}} \right]. \hspace*{4cm} (4.2.28)
\end{eqnarray*}

\noi (from (4.2.20) and (4.1.2)). Moreover : 
$$
 \mu_y^{(r)} (1) = \psi_r (y) = \frac{r}{1-r} \mathbb{E}_{x_0} [\phi(X_1)] + \phi(y). \eqno(4.2.29) $$

\noi By hypothesis, there exists $n \geq 0$ such that $\mathbb{P}_{x_0} (X_n = y_0) > 0$; it is easy to check that it implies : $q_{y_0}^{(x_0)} < 1$. 

\noi Hence, by considering the difference between (4.2.28) and (4.2.29) and taking $r \rightarrow 1$, one obtains : 
$$
\phi^{[y_0]} (y) = \mathbb{E}_{x_0} [\phi(X_1)] \frac{ 1 - p_{y,y_0}^{(x_0)}}{ 1- q_{y_0}^{(x_0)} }
+ [\phi(y) - \phi(y_0)]. \eqno(4.2.30) $$

\noi Therefore :
$$ 
\phi(y) - \phi(y_0) \leq \phi^{[y_0]} (y) \leq \frac{\mathbb{E}_{x_0} [\phi(X_1)]}{1- q_{y_0}^{(x_0)}} +
[\phi(y) - \phi(y_0)] \eqno(4.2.31) $$

\noi which implies point {\it i)} of the Theorem, and in particular the finiteness of $\phi^{[y_0]}$.  By applying Theorem 4.2.1 to $T=1$ and $F= {\bf 1}_{L_{\infty}^{y_0} = 0}$, one can easily check that $\phi^{[y_0]}$ is harmonic everywhere except at point $y_0$ (where it is equal to zero). \\
By taking $T=1$ and $F= {\bf 1}_{L_{\infty}^{y_0} = 1}$, one obtains the formula :
$\widetilde{\mathbb{Q}}_{y_0} (1) = \mathbb{E}_{y_0} [ \phi^{[y_0]} (X_1) ]$. Hence, we obtain point {\it ii)} of the Theorem, and the point {\it iv)} by formula (4.2.11). Now, by taking the notation : $\mu_y^{(r),y_0} = r^{L_{\infty}^{y_0}} .\mathbb{Q}_y$, one has (for all positive and $\mathcal{F}_n$-measurable
functionals $F_n$), by applying Theorem 4.2.1 to $T=n$ and $F= F_n \, r^{L_{\infty}^{y_0}}$ :
$$
  \mu_y^{(r),y_0} (F_n)  = \mathbb{Q}_y [F_n \, r^{L_{\infty}^{y_0}} ] = \mathbb{E}_y  \left[ F_n \, r^{L_{n-1}^{y_0}} \alpha (X_n) \right], \eqno(4.2.32) $$

\noi where  $\alpha (z) =  \mathbb{Q}_z [r^{L_{\infty}^{y_0}} ]$. By point {\it iv)} of the Theorem (already proven), one has :
\begin{eqnarray*}
\hspace*{1cm}  & \alpha(z) & = \phi^{[y_0]} (z) + \left(\sum_{k=1}^{\infty} r^k \right) \mathbb{E}_{y_0} [\phi^{[y_0]} (X_1) ] \\
& & = \frac{r}{1-r}  \, \mathbb{E}_{y_0} [\phi^{[y_0]} (X_1) ] + \phi^{[y_0]} (z)  \hspace*{6cm} (4.2.33)
\end{eqnarray*}

\noi Hence : 
$$
\mu_y^{(r),y_0} (F_n) = \mathbb{E}_y \left[ F_n \, r^{L_{n-1}^{x_0}} \left( \frac{r}{1-r}  \, \mathbb{E}_{y_0}
[\phi^{[y_0]} (X_1) ] +  \phi^{[y_0]} (X_n) \right) \right]  \eqno(4.2.34) $$

\noi This formula implies that $\mu_y^{(r),y_0}$ is the measure defined in the same way as $\mu_y^{(r)}$, but from the point $y_0$ and the function $\phi^{[y_0]}$, instead of the point $x_0$ and the function $\phi$. By considering the new measure with density $r^{- L_{\infty}^{y_0}}$ with respect to $\mu_y^{(r), y_0}$, one obtains
the equality : 
$$     \mathbb{Q}_y = \mathbb{Q}_y^{(\phi^{[y_0]}, y_0)} \eqno(4.2.35) $$

\noi which completes the proof of Theorem 4.2.5. $\hfill \blacksquare$

\bigskip

\noi There is also an important formula, which is a direct consequence of (4.2.1), (4.2.5) and Theorem 4.2.5. : 

\smallskip

\noi {\bf Corollary 4.2.6} {\it Let $F_n$ be a positive $\mathcal{F}_n$-measurable functional, $y, y_0$ be
in $E$ and $g_{y_0}$ be the last hitting time of $y_0$ for the canonical process. Then the following formula holds :
$$
\mathbb{Q}_y \left[ F_n {\bf 1}_{g_{y_0} < n} \right] = \mathbb{E}_y [F_n \phi^{[y_0 ]} (X_n)]
\eqno(4.2.36) $$

\noi In particular, one has :
$$
\mathbb{Q}_y \left[ F_n {\bf 1}_{g_{x_0} < n} \right] = \mathbb{E}_y [F_n \phi(X_n)] \eqno(4.2.37) $$

\noi and $\big(\phi^{[y_{0}]} (X_{n}), \; n \ge 0\big)$, $\big(\phi (X_{n}), \; n \ge 0\big)$ are two $\mathbb{P}$ submartingales.  }

\noi The correspondance with the Brownian case is the following :

\smallskip

\begin{tabular}{|c|c|}
\hline
Markov chain & Brownian motion \\
\hline
$ \mathbb{Q}_y [ F_n {\bf 1}_{g_{x_0} < n}] = \mathbb{E}_y [F_n \phi (X_n)]$ & $ {\bf W}_x 
(F_t 1_{g < t} ) = W_x (F_t |X_t|) $ \\
$ \mathbb{Q}_y [ F_n {\bf 1}_{g_{y_0} < n}] = \mathbb{E}_y [F_n \phi^{[y_0]} (X_n)]$ & $ {\bf W}_x 
(F_t 1_{\sigma_a < t} ) = W_x (F_t (|X_t|-a)_+) $ \\
$F_n \in \mathcal{F}_n$ & $F_t \in \mathcal{F}_t$ \\
\hline
\end{tabular}

\smallskip

\noi By 
 Theorem 4.2.5, the construction of a given family $(\mathbb{Q}_x)_{x \in E}$ can be obtained by taking any point $y_0$ instead of $x_0$, if the corresponding harmonic function $\phi^{[y_0]}$ is well-chosen.

\smallskip

\noi {\bf 4.2.4} \underline{Dependence of $\mathbb{Q}_{x}$ on $\phi$.}

\noi In fact, this family of $\sigma$-finite measures depends only upon the equivalent class of the function $\phi$, for an equivalence relation which will be described below.  More precisely, if $\alpha$ and $\beta$ are two functions from $E$ to $\mathbb{R}_+$, let us write : $\alpha \simeq \beta$, iff
$\alpha$ is equivalent to $\beta$ when $\alpha + \beta$ tends to infinity ; i.e, for all $\epsilon \in ]0,1[$, there exists $A > 0$ such that for all $x \in E$, $\alpha(x) + \beta(x) \geq A$ implies $ \dis 1- \epsilon <
\frac{\alpha(x)}{\beta(x)} < 1 + \epsilon$. With this definition, one has the following result :

\noi {\bf Propostion 4.2.7} {\it The relation $\simeq$ is an equivalence relation.}

\smallskip

\noi {\bf Proof of Proposition 4.2.7 } The reflexivity and the symmetry are obvious, so let us prove the transitivity. 

\noi We suppose that there are three functions $\alpha, \beta, \gamma$ such that $\alpha \simeq \beta$ and $\beta \simeq \gamma$. 

\noi There exists $\epsilon : \mathbb{R}_+ \rightarrow \mathbb{R}_+ \cup \{\infty \} $, tending to zero at
infinity, such that $\alpha + \beta \geq A$ implies $ \left| \frac{\alpha}{\beta} - 1 \right| \leq \epsilon (A)$, and $\beta + \gamma \geq A$ implies $ \left| \frac{\beta}{\gamma} - 1 \right| \leq \epsilon(A)$.  For a given $x \in E$, let us suppose that $\alpha(x) + \gamma(x)  \geq A$ for $A >  4 \sup \{ z, \epsilon (z) \geq 1/2 \}$. There are two cases :

\begin{itemize}
\item $\alpha(x) \geq A/2$. In this case, $\alpha(x) + \beta(x) \geq A/2$; hence, $ \left|\frac{ \alpha(x)}
{\beta(x)} - 1 \right| \leq \epsilon(A/2) \leq 1/2$, which implies : $\beta (x) + \gamma (x) \geq \beta(x) 
\geq \alpha(x) /2 \geq A/4$. 

\noi Therefore :  $ \left|\frac{ \beta(x)} {\gamma(x)} - 1 \right| \leq \epsilon(A/4)$. Consequently, there exist $u$ and $v$, $|u| \leq \epsilon(A/2) \leq 1/2$, $|v| \leq \epsilon(A/4) \leq 1/2$, such that $\frac{\alpha(x)}{\gamma(x)} =  (1+u) (1+v)$, which implies : 
\begin{eqnarray*}
\hspace*{1cm} & \left| \frac{\alpha(x)}{\gamma(x)} - 1 \right| & \leq |u| + |v| + |uv| \leq \epsilon (A/2) +  \epsilon (A/4) + \epsilon (A/2)  \epsilon (A/4)  \\ 
& & \leq \frac{3}{2} \left(  \epsilon (A/2) +  \epsilon (A/4)   \right) \hspace*{6cm}(4.2.38)
\end{eqnarray*}
\item $\alpha(x) \leq A/2$. In this case, $\gamma(x) \geq A/2$, hence we are in the same situation as in the first case if we exchange $\alpha (x)$ and $\gamma(x)$
\end{itemize}

\noi The above inequality  implies : $\alpha \simeq \gamma$, since $ \epsilon (A/2) +  \epsilon (A/4)$ tends to zero when $A$ goes to infinity.  Hence, $\simeq$ is an equivalence relation. $\hfill \blacksquare $

\bigskip

\noi This equivalence relation satisfies the following property : 

\smallskip

\noi {\bf Lemma 4.2.8} {\it Let $\phi_1$ and $\phi_2$ be two functions from $E$ to $\mathbb{R}_+$ which 
are equal to zero at a point $y_0 \in E$ and which are harmonic at the other points i.e. for all $y \neq y_{0}$, $E_{y} [\phi_{i} (X_{1})] = \phi_{i} (y), \; i=1,2$. If $\phi_1 \simeq \phi_2$, then $\phi_1 = \phi_2$.   }

\smallskip

\noi {\bf Proof of Lemma 4.2.8 }  By the martingale property, for all $x \in E$, $A > 0$ : 
\begin{eqnarray*}
\hspace*{1cm} & \phi_1 (x) & = \mathbb{E}_x \left[ \phi_1 (X_{n \wedge \tau_1^{(y_0)}})  \right]  \\
& & = \mathbb{E}_x \left[ \phi_1 (X_{n \wedge \tau_1^{(y_0)}}) {\bf 1}_{\phi_1 (X_{n \wedge \tau_1^{(y_0)}})+ \phi_2 (X_{n \wedge \tau_1^{(y_0)}}) \geq A} \right] + K, \hspace*{3cm}(4.2.39)
\end{eqnarray*}

\noi where $|K| \leq A \, \mathbb{P}_x (\tau_1^{(y_0)} > n)$. Now, if $\phi_1(y) + \phi_2(y) \geq A$, one has : 
$$    (1- \epsilon(A)) \phi_1(y) \leq \phi_2 (y) \leq (1 + \epsilon(A)) \phi_1(y), \eqno(4.2.40)  $$

\noi where $\epsilon(A)$ tends to zero when $A$ tends to infinity. Therefore : 
$$ 
\phi_1 (x)= \alpha \, \mathbb{E}_x \left[ \phi_2 (X_{n \wedge  \tau_1^{(y_0)}}) {\bf 1}_{  \phi_1 (X_{n \wedge \tau_1^{(y_0)}})+ \phi_2 (X_{n \wedge \tau_1^{(y_0)}}) \geq A} \right] + K, \eqno(4.2.41)   $$

\noi where $1 - \epsilon(A) \leq \alpha \leq 1 + \epsilon(A)$. Moreover :  
$$
\phi_2 (x)= \mathbb{E}_x \left[ \phi_2 (X_{n \wedge \tau_1^{(y_0)}}) {\bf 1}_{  \phi_1 (X_{n \wedge \tau_1^{(y_0)}})+ \phi_2 (X_{n \wedge \tau_1^{(y_0)}}) \geq A} \right] + K', \eqno(4.2.42) $$

\noi where $|K'| \leq A \, \mathbb{P}_x (\tau_1^{(y_0)} > n)$. Hence : 
$$     \phi_1 (x) = \alpha \left( \phi_2 (x) - K' \right) + K. \eqno(4.2.43) $$

\noi Now, if $A$ is fixed, $|K| + |K'|$ tend to zero when $n$ goes to infinity. Therefore : 
$$      (1- \epsilon(A)) \phi_1(x) \leq \phi_2 (x) \leq (1 + \epsilon(A)) \phi_1(x). \eqno(4.2.44) $$

\noi This inequality is true for all $A \geq 0$; hence : $\phi_1 = \phi_2$, which proves Lemma 4.2.8. $\hfill \blacksquare $

\noi We now give another lemma, which is quite close to Lemma 4.2.8 : 

\smallskip

\noi {\bf Lemma 4.2.9} {\it Let $\phi$ be a function from $E$ to $\mathbb{R}$ which 
is equal to zero at a point $y_0 \in E$ and harmonic at the other points.
If $\phi$ is bounded, it is identically zero.}
\smallskip

\noi {\bf Proof of Lemma 4.2.9 } Since $\phi$ is bounded, there exists $A>0$ such that $|\phi(x)| < A$. 
The harmonicity of $\phi$ implies, for every $n \geq 0$ and $x \neq y_0$ : 
$$ \phi(x) = \mathbb{E}_{x} [ \phi(X_{n \wedge \tau_1^{y_0}})]$$
\noi Consequently, since $\phi(y_0) = 0$, we get :
$$|\phi (x)| \leq A \, \mathbb{P}_x (\tau_1^{y_0} > n) \underset{n \rightarrow \infty}{\longrightarrow}
 0 $$ 
\noi since $(X_n, \, n \geq 0)$ is recurrent. Hence, $\phi$ is identically zero. 

\noi If $\phi$ is bounded and positive, then $\phi$ is equivalent to zero (by definition of $\simeq $). Hence,
in this case, Lemma 4.2.9 may be considered as a particular case of Lemma 4.2.8.
 
\smallskip

\noi Now, let us state the following result, which explains why we have defined the previous equivalence relation :

\smallskip

\noi {\bf Proposition 4.2.10} {\it Let $x_0$, $y_0$ be in $E$, $\phi$ a positive function which is harmonic 
except at $x_0$ and equal to zero at $x_0$, $\psi$ a positive function which is harmonic except at $y_0$ and equal to zero at $y_0$. In these conditions, the family $(\mathbb{Q}_x^{(\phi, x_0)})_{x \in E}$ of 
$\sigma$-finite measures is identical to the family  $(\mathbb{Q}_x^{( \psi, y_0)})_{x \in E}$ if and only if $\phi \simeq \psi$. Therefore this family can also be denoted by $(\mathbb{Q}_x^{[\phi]})_{x \in E}$, where $[\phi]$ is the equivalence class of $\phi$.  }
\smallskip

\noi {\bf Proof of Proposition 4.2.10 } If the two families of measures are equal, for all $x \in E$, $\mathbb{Q}_x^{(\phi, x_0)} =  \mathbb{Q}_x^{(\psi, y_0)}$. Now, it has been proven that $\psi(x) = \mathbb{Q}_x^{(\psi, y_0)} (L_{\infty}^{y_0} = 0)$. Hence, if $\phi^{[y_0]} (x) = \mathbb{Q}_x^{(\phi, x_0)} (L_{\infty}^{y_0} = 0)$, one has $\psi = \phi^{[y_0]}$.
 
\noi Since $\phi - \phi^{[y_0]}$ is bounded (point {\it i)} of Theorem 4.2.5), $\phi - \psi$ is bounded, which implies that $\phi$ is equivalent to $\psi$. On the other hand, if $\phi$ is equivalent to $\psi$, and if $\phi^{[y_0]} = \mathbb{Q}_x^{(\phi, x_0)} (L_{\infty}^{y_0} = 0)$, $\psi$ and $\phi^{[y_0]}$ are two equivalent functions which are harmonic except at point $y_0$, and equal to zero at $y_0$. Hence, by Lemma 4.2.8, $\psi = \phi^{[y_0]}$, and by Theorem 4.2.5, for all $x \in E$, $\mathbb{Q}_x^{(\phi, x_0)} = \mathbb{Q}_x^{(\phi^{[y_0]}, y_0)}$. 
 
\noi Therefore, $\mathbb{Q}_x^{(\phi, x_0)} = \mathbb{Q}_x^{(\psi, y_0)}$, which proves Proposition 4.2.10.

\noi In the next Section, we give some examples of the above construction.

\bigskip

\noi {\large \bf 4.3 Some examples.}

\smallskip






\noi {\bf 4.3.1} \underline{The standard random walk.}

\noi In this case, $E = \mathbb{Z}$ and for all $x  \in E$, $\mathbb{P}_x$ is the law of the standard random walk. The functions $\phi_+ : x \rightarrow x_+$, $\phi_- : x \rightarrow x_-$ and their sum $\phi : x \rightarrow |x|$ satisfies the harmonicity conditions above at point $x_0 = 0$.

\noi Let $(\mathbb{Q}_x^{+})_{x \in \mathbb{Z}}$, $(\mathbb{Q}_x^{-})_{x \in \mathbb{Z}}$ and 
$(\mathbb{Q}_x)_{x \in \mathbb{Z}}$ be the associated $\sigma$-finite measures. For all $a \in \mathbb{Z}$, let us take the notations : $\phi_+^{[a]} (x) = \mathbb{Q}^+_x [L_{\infty}^{a} = 0]$,
 $\phi_-^{[a]} (x) = \mathbb{Q}^-_x [L_{\infty}^{a} = 0]$ and $\phi^{[a]} (x) = \mathbb{Q}_x [L_{\infty}^{a} = 0]$.

\noi The function $\phi_+^{[a]}$ satisfies the harmonicity conditions at point $a$ and is equivalent to $\phi_+$. Now, these two properties  are also satisfied by the function $x \rightarrow (x-a)_+$; hence, by Lemma 4.2.8, $\phi_+^{[a]} (x) = (x-a)_+$. By the same argument, $\phi_-^{[a]} (x) = (x-a)_-$ and $\phi^{[a]} (x) = |x-a|$. 

\noi Therefore, we have the equalities for every positive and $\mathcal{F}_n$-measurable functional $F_n$, and for every $x, a \in  \mathbb{Z}$ : 
\begin{eqnarray*} 
\hspace*{1cm} \mathbb{Q}_x^{+} [F_n \, {\bf 1}_{g_a < n} ] &=& \mathbb{E}_x [F_n (X_n - a)_+], \label{+} \hspace*{7cm}(4.3.1) \\
\mathbb{Q}_x^{-} [F_n \,{\bf 1}_{g_a < n} ] &=& \mathbb{E}_x [F_n (X_n - a)_-], \label{-}
\hspace*{7cm} (4.3.2) \\
\mathbb{Q}_x [F_n \, {\bf 1}_{g_a < n} ] &=& \mathbb{E}_x [F_n |X_n - a|]. \label{+-} \hspace*{7,4cm} (4.3.3) 
\end{eqnarray*}

\noi These equations and the fact that the canonical process is transient under $\mathbb{Q}_x^{+}$, $\mathbb{Q}_x^{-}$, $\mathbb{Q}_x$ characterize these measures. Moreover, by using equations (4.3.1), (4.3.2) and (4.3.3), it is not difficult to prove that for all $x \in \mathbb{Z}$, these measures are the 
images of $\mathbb{Q}_0^{+}$, $\mathbb{Q}_0^{-}$ and $\mathbb{Q}_0$ by the translation by $x$. 

\noi Now, for all $a, x \in \mathbb{Z}$, and for all positive and $\mathcal{F}_n$-measurable functional $F_n$ : 
$$
\mathbb{Q}_x^{+,[a]} [F_n] := \mathbb{Q}_x^{+} [F_n \, {\bf 1}_{L_{\infty}^{a} = 0} ] =  \mathbb{E}_x [F_n 
 (X_{n \wedge \tau_1^{(a)}} - a)_+ ]. \eqno (4.3.4)  $$

\noi Hence, if $x \leq a$, $\mathbb{Q}_x^{+,[a]} = 0$, and if $x > a$, $\mathbb{Q}_x^{+,[a]}$ is $(x-a)$ times the law of a Bessel random walk strictly above $a$, starting at point $x$ \big(cf [LG] for a definition of the Bessel random walk\big).

\noi By the same arguments, if $x \geq a$, $\mathbb{Q}_x^{-,[a]} = 0$, and if $x < a$, $\mathbb{Q}_x^{-,[a]}$ is $(a-x)$ times the law of a Bessel random walk strictly below $a$, starting at point $x$. Moreover, $\mathbb{Q}_x^{[a]}$ is the $|x-a|$ times the law of a Bessel random walk above or below $a$, depending on the sign of $x-a$. The same kind of arguments imply that (with obvious notations) : 
\begin{itemize}
\item $\widetilde{\mathbb{Q}}^+_{a}$ is $1/2$ times the law of a Bessel 
random walk strictly above $a$.
\item $\widetilde{\mathbb{Q}}^-_{a}$ is $1/2$ times the law of a Bessel 
random walk strictly below $a$.
\item $\widetilde{\mathbb{Q}}_{a}$ is the law of a symmetric Bessel random 
walk, strictly above or below $a$ with equal probability.
\end{itemize}

\noi The equalities given by Proposition 4.2.4 are the following :
$$
\mathbb{Q}^+_x = \mathbb{Q}_x^{+,[a]} + \sum_{k \geq 1} \mathbb{P}_x^{\tau_k^{(a)}} \circ \widetilde{\mathbb{Q}}^+_{a}, \eqno(4.3.5) $$
$$\mathbb{Q}^-_x = \mathbb{Q}_x^{-,[a]} + \sum_{k \geq 1} \mathbb{P}_x^{\tau_k^{(a)}} \circ \widetilde{\mathbb{Q}}^-_{a}, \eqno(4.3.6) $$
$$
\mathbb{Q}_x \;= \mathbb{Q}_x^{[a]} + \sum_{k \geq 1} \mathbb{P}_x^{\tau_k^{(a)}} \circ \widetilde{\mathbb{Q}}_{a}. \eqno(4.3.7) $$

\noi Moreover, one has :
\begin{itemize}
\item $\mathbb{Q}_x^{+} [L_{\infty}^{a} = 0] = (x-a)_+$ and $ \mathbb{Q}_x^{+} [L_{\infty}^{a} = k] = 1/2$ for all $k \geq 1$.
\item $\mathbb{Q}_x^{-} [L_{\infty}^{a} = 0] = (x-a)_-$ and $ \mathbb{Q}_x^{-} [L_{\infty}^{a} = k] = 1/2$ for all $k \geq 1$.
\item $\mathbb{Q}_x [L_{\infty}^{a} = 0] = |x-a|$ and $ \mathbb{Q}_x [L_{\infty}^{a} = k] = 1$ for all $k \geq 1$.
\end{itemize}

\noi Hence, by applying Theorem 4.2.1 and Corollary 4.2.2 to the functional $F=h(L_{\infty}^{a})$ for a positive function $h$ such that $\sum_{n \in \mathbb{N}} h(n) < \infty$, and for $a \in \mathbb{Z}$, one 
obtains that for all $x \in \mathbb{Z}$ :
$$
M_n^+ = (X_n - a)_+ \, h(L_{n-1}^{a}) + \frac{1}{2} \sum_{k = L_{n-1}^{a} + 1}^{\infty} h(k), \eqno(4.3.8)
$$
$$
M_n^- = (X_n - a)_- \, h(L_{n-1}^{a}) + \frac{1}{2} \sum_{k = L_{n-1}^{a} + 1}^{\infty} h(k), \eqno(4.3.9) $$

\noi and their sum
$$
M_n \;= |X_n - a| \, h(L_{n-1}^{a}) + \sum_{k = L_{n-1}^{a} + 1}^{\infty} h(k) \eqno(4.3.10) $$

\noi are martingales under the probability $\mathbb{P}_x$. Other  martingales can be obtained by taking other functionals $F$. 

\smallskip

\noi {\bf 4.3.2} \underline{The "bang-bang random walk".} 

\noi In this case, we suppose that $E = \mathbb{N}$ and that $(\mathbb{P}_x)_{x \in \mathbb{N}}$ is the family of measures associated to transition probabilities : $p_{0,1} = 1$, $p_{y,y+1} = 1/3$ and $p_{y,y-1}=2/3$ for all $y \geq 1$. Informally, under $\mathbb{P}_x$ (for any $x \in \mathbb{N}$), the canonical 
process is a Markov process which tends to decrease when it is strictly above zero, and which increases if it is equal to zero. 

\noi The family of measures $(\mathbb{Q}_x)_{x \in \mathbb{N}}$ can be constructed by taking $x_0  = 0$ and $\phi(x) = 2^x - 1$ for all $x \in \mathbb{N}$. 

\noi For $y \in \mathbb{N}$, the function $\phi^{[y]} : x \rightarrow \mathbb{Q}_x [L_{\infty}^{y} = 0]$ is harmonic except at $y$ where it is equal to zero, and it is equivalent to $\phi$. 

\noi Since the function : $x \rightarrow (2^x - 2^y). {\bf 1}_{x \geq y}$ satisfies the same properties, by Lemma 4.2.8, we get :  $\phi^{[y]} (x) =  (2^x - 2^y). {\bf 1}_{x \geq y}$. For all $x \in \mathbb{N}$, the measure $\mathbb{Q}_x$ is characterized by the transience of the canonical process, and by the formula : 
$$
\mathbb{Q}_x [F_n \, {\bf 1}_{g_a < n} ] = \mathbb{E}_x [F_n \, (2^{X_n} - 2^a)_+ ], \eqno(4.3.11) $$

\noi which holds for all $a$, $n \in \mathbb{N}$ and for every positive $\mathcal{F}_n$-measurable functional $F_n$. 

\noi Adopting obvious notations, it is not difficult to prove the formula : 
$$
\mathbb{Q}_x^{[a]} (F_n) = \mathbb{E}_x [ F_n \, (2^{X_{n \wedge\tau_1^{(a)}}} - 2^a)] \, {\bf 1}_{x \geq a}, 
\eqno(4.3.12) $$

\noi and for $n \geq 1$ : 
$$
\widetilde{\mathbb{Q}}_{a} (F_n) = \mathbb{E}_a \left[ F_n \, (2^{X_{n \wedge \tau_2^{(a)}}} - 2^a) \, {\bf 1}_{X_1 = a+1} \right]. \eqno(4.3.13) $$

\noi Moreover : 
\begin{itemize}
\item The total mass of $\mathbb{Q}_x^{[a]}$ is zero if $x \leq a$, and $2^x - 2^a$ if $x>a$. 
\item The total mass of $ \widetilde{\mathbb{Q}}_{a} $ is 1 if $a=0$, and $2^a/3$ if $x \geq 1$. 
\item For $x > a$ and under the probability $\bar{\mathbb{P}}_x^{[a]} = \mathbb{Q}_x^{[a]}/ (2^x-2^a)$, the canonical process is a Markov process with probability transitions : $\dis \bar{p}_{x, x+1} = \frac{2. 2^{x-a} - 1}{3. 2^{x-a}-3}$ and $\dis \bar{p}_{x, x-1} = \frac{ 2^{x-a} - 1}{3. 2^{x-a}-3}$. We remark that $\bar{p}_{x, 
x+1}$ tends to $2/3$ when $x$ goes to infinity, and $\bar{p}_{x, x-1}$ tends to $1/3$ (the opposite case as the initial transition probabilities).
\item Under the probability $\dis \frac{  \widetilde{\mathbb{Q}}_{a} }{(2^a/3) {\bf 1}_{a \geq 1} + {\bf 1}_{a=0}}$, the canonical process is a Markov process with the same transition probabilities as 
under $\bar{\mathbb{P}}_x^{[a]}$, with $X_1 = a+1$ almost surely. 
\end{itemize}

\noi For all $a, x \in \mathbb{N}$, the image of $\mathbb{Q}_x$ by the total local times is given by the equalities : 
$$
\mathbb{Q}_x [L_{\infty}^{a} = 0] = (2^x - 2^a) \, {\bf 1}_{x > a}, \eqno (4.3.14) $$

\noi and for all $k \geq 1$ : 
$$
\mathbb{Q}_x [L_{\infty}^{a} = k] = K(a), \eqno (4.3.15) $$

\noi where $K(0)=1$ and $K(a)= 2^a/3$ for $a \geq 1$.

\noi Moreover, if $h$ is an integrable function from $\mathbb{N}$ to $\mathbb{R}_+$, and if $a, x \in \mathbb{N}$, 
$$
M_n = h(L_{n-1}^a) \, (2^{X_n}-2^a)_+ \, + \, K(a) \sum_{k = L_{n-1}^{a} + 1}^{\infty} h(k) \eqno(4.3.16) $$

\noi is a martingale under the initial probability $\mathbb{P}_x$. 

\smallskip

\noi {\bf 4.3.3} \underline{The random  walk on a tree.}

\noi We consider a random walk on a binary tree, which can be represented by the set $E = \{ \varnothing, (0), (1), (0,0), (0,1), (1,0), (1,1), (0,0,0), ... \}$ of $k$-uples of elements in $\{0,1\}$ 
($k \in \mathbb{N}$).  

\noi Obviously, $k$ is the distance to the origin $\varnothing$ of the tree.

\noi The probability transitions of the Markov process associated to the starting 
family of probabilities $(\mathbb{P}_x)_{x \in E}$ are $p_{\varnothing, (0)} = p_{\varnothing, (1)} = 1/2$, and for $k \geq 1$ : $p_{(x_1,x_2,...,x_k), (x_1, x_2,..., x_{k-1})} = 1/2$, $p_{(x_1,...,x_k), (x_1,...,x_k,0)} = 
p_{(x_1,...,x_k), (x_1,...,x_k,1)} = 1/4$. 

\noi In particular, under $\mathbb{P}_x$ (for all $x \in E$), the distance to the origin is a standard reflected random walk. 

\noi If the reference point $x_0$ is $\varnothing$, we can take for $\phi$ the 
distance to the origin of the tree. But there are other possible functions $\phi$ for the same point $x_0$. 
For example, if $(a_0,a_1,a_2,...)$ is an infinite sequence of elements in $\{0,1\}$ it is possible to take for $\phi$ the function such that for all $(x_0,x_1,...,x_k) \in E$, one has $\phi (x_0,x_1,...,x_k) = 2^p-1$, where $p$ is the smallest element of $\mathbb{N}$ such that $p > k$ or $x_p \neq a_p$. In particular, if $a_p = 0$ for all $p$, one has $\phi(\varnothing) = 0$, $\phi((0)) = 1$, $\phi((1)) = 0$, $\phi((0,0)) = 3$, $\phi((0,1)) = 1$, $\phi((1,0)) = \phi((1,1)) = 0$, $\phi((0,0,0)) = 7$, etc. 

\noi Each choice of the sequence $(a_p)_{p \in \mathbb{N}}$ gives a different function $\phi$ and hence a different family $(\mathbb{Q}^{[\phi]}_x)_{x \in E}$ of $\sigma$-finite measures. 

\smallskip

\noi {\bf 4.3.4} \underline{Some more general conditions for existence of $\phi$.}

\noi The following proposition gives some sufficient conditions for the existence of a function $\phi$ which satisfies the hypothesis of Section 4.1.2 : 

\smallskip

\noi {\bf Proposition 4.3.1} {\it Let $(\mathbb{P}_x)_{x \in E}$ be the family of probabilities associated to 
a discrete time Markov process on a countable set $E$. We suppose that for all $x \in E$, the set of $y \in E$ such that the transition probability $p_{x,y}$ is strictly positive is finite.  Furthermore, let us consider a function $\phi$ which satisfies one of the two following conditions (for a given point $x_0 \in E$) : 
\begin{itemize}
\item There exists a function $f$ from $\mathbb{N}$ to $\mathbb{R}_+^*$ such that $f(n)/f(n+1)$ tends to $1$ when $n$ goes to infinity, and such that for all $x \in E$ :
$$
\mathbb{E}_x [\tau_1^{(x_0)} \geq n] \underset{n \rightarrow \infty}{\sim} f(n) \, \phi(x).  \eqno(4.3.17) $$
\noi where $\tau_1^{(x_0)}$ is the first hitting time of $x_0$, for the canonical process. 
\item For all $x \in E$, $\mathbb{P}_x (X_k = x_0)$ tends to zero when $k$ tends to infinity, and : 
$$
\sum_{k=0}^{N} \left[ \mathbb{P}_{x_0} (X_k = x_0) -  \mathbb{P}_{x} (X_k = x_0) \right] \underset{N \rightarrow \infty}{\rightarrow} \phi(x). \eqno(4.3.18) $$
\end{itemize}

\noi In these conditions, $\phi$ is harmonic, except at point $x_0$ where this function is equal to zero.  }

\noi {\bf Proof of Proposition 4.3.1 } Let us suppose that the first condition is satisfied. For all $x \neq x_0$ and for all $y \in E$ such that $p_{x,y}  > 0$ :
$$  
\mathbb{E}_y \left[ \tau_1^{(x_0)} \geq n \right] \underset{n \rightarrow \infty}{\sim} f(n) \, \phi(y). \eqno (4.3.19) $$

\noi By adding the equalities obtained for each point $y$ and multiplied by $p_{x,y}$, we obtain : 
$$
\sum_{y \in E} p_{x,y} \mathbb{E}_y \left[ \tau_1^{(x_0)} \geq n \right] \underset{n \rightarrow \infty}{\sim}
f(n) \, \sum_{y \in E} p_{x,y} \, \phi(y), \eqno (4.3.20) $$

\noi which implies : 
$$
\mathbb{E}_x \left[ \tau_1^{(x_0)} \geq n+1 \right] \underset{n \rightarrow \infty}{\sim} f(n) \, \mathbb{E}_x 
[\phi(X_1)]. \eqno(4.3.21) $$

\noi Moreover :
$$
\mathbb{E}_x \left[ \tau_1^{(x_0)} \geq n+1 \right] \underset{n \rightarrow \infty}{\sim} f(n+1) \, \phi(x). \eqno (4.3.22) $$

\noi By comparing these equivalences and by using the fact that $f(n)$ is equivalent to $f(n+1)$ and is strictly positive, one obtains : 
$$        \phi(x) = \mathbb{E}_x [\phi(X_1)].  \eqno (4.3.23) $$

\noi Since $\phi(x_0)$ is obviously equal to zero ($\mathbb{E}_{x_0} \left[ \tau_1^{(x_0)} \geq n \right] = 0$), Proposition 4.3.1 is proven if the first condition holds. 

\noi Now let us assume the second condition holds. 

\noi If $x \neq x_0$, for all $y$ such that $p_{x,y} > 0$ :
$$
\sum_{k=0}^{N} \left[ \mathbb{P}_{x_0} (X_k=x_0) - \mathbb{P}_y (X_k = x_0) \right] \underset{N \rightarrow \infty}{\rightarrow} \phi(y). \eqno (4.3.24) $$

\noi Therefore : 
$$
\sum_{y \in E} p_{x,y} \, \left[ \sum_{k=0}^{N} \left( \mathbb{P}_{x_0} (X_k=x_0) - \mathbb{P}_y (X_k = x_0) \right)\right] \underset{N \rightarrow \infty} {\rightarrow} \sum_{y \in E} p_{x,y} \, \phi(y). \eqno (4.3.25) $$

\noi This equality implies : 
$$
\sum_{k=0}^{N} \left[ \mathbb{P}_{x_0} (X_k=x_0) \right] -  \sum_{k=1}^{N+1} \left[ \mathbb{P}_{x} (X_k=x_0) \right]\underset{N \rightarrow \infty} {\rightarrow} \mathbb{E}_x [\phi(X_1)]. \eqno (4.3.26) $$

\noi Now, $\mathbb{P}_{x} (X_0 = x_0) = 0$ (since $x \neq x_0$) and when $N$ goes to infinity, $\mathbb{P}_x(X_{N+1} = x_0)$ tends to zero by hypothesis. Hence : 
$$
\sum_{k=0}^{N} \left[ \mathbb{P}_{x_0} (X_k=x_0) - \mathbb{P}_x (X_k = x_0) \right] \underset{N \rightarrow \infty} {\rightarrow} \mathbb{E}_x [\phi(X_1)], \eqno (4.3.27) $$

\noi which implies : 
$$     \phi(x) = \mathbb{E}_x [\phi(X_1)], \eqno (4.3.28) $$

\noi as written in Proposition 4.3.1. $\hfill \blacksquare$

\bigskip

\noi {\bf Remark 4.3.2} {\it If the condition : 
$$
\sum_{k=0}^{N} \left[ \mathbb{P}_{x_0} (X_k=x_0) - \mathbb{P}_x (X_k = x_0) \right] \underset{N \rightarrow \infty}{\rightarrow} \phi(x) \eqno (4.3.29) $$

\noi is satisfied for a function $\phi$, then $\phi$ is automatically positive. Indeed : 

$$ \sum_{k=0}^{N} \left[ \mathbb{P}_{x_0} (X_k=x_0) - \mathbb{P}_x (X_k = 
x_0) \right] 
 =  \mathbb{E}_{x_0} \left[ \sum_{k = 0}^{N} {\bf 1}_{X_k  = x_0} \right] 
\, - \,  \mathbb{E}_{x} \left[ \sum_{k = 0}^{N} {\bf 1}_{X_k  = x_0}  \right], \eqno (4.3.30)$$

\noi where, by the strong Markov property : 
\begin{eqnarray*}
\hspace*{1cm} & \mathbb{E}_{x_0} \left[ \sum_{k = 0}^{N} {\bf 1}_{X_k  = x_0}  \right] & \geq \mathbb{E}_x \left[ 
\sum_{k = 0}^{  \tau_1^{(x_0)} + N} {\bf 1}_{X_k  = x_0}  \right]  \\ 
& & \geq \mathbb{E}_{x} \left[  \sum_{k = 0}^{N} {\bf 1}_{X_k  = x_0}   \right]. \hspace*{4cm} (4.3.31)
\end{eqnarray*}   }

\noi {\bf 4.3.5} \underline{The standard random walk on $\mathbb{Z}^2$.}

\noi In this case, $E = \mathbb{Z}^2$ and $(\mathbb{P}_x)_{x \in \mathbb{Z}^2}$ is the family of probabilities associated to the standard random walk. If we take $x_0 = (0,0)$, the problem is to find a function $\phi$ which satisfies the hypothesis of Section 4.1.2 : it can be solved by using 
Proposition 4.3.1.

\noi More precisely, by doing some classical computations (see for example [Spi]), we can prove that for all $(x,y) \in \mathbb{Z}^2$, and for all $k \in \mathbb{N}$ : 
$$
\mathbb{P}_{(x,y)} \left[ X_k = (0,0) \right] = {\bf 1}_{k \equiv x+y \,  (mod. \, 2)} \, \frac{C}{k+1}  \, + \, \epsilon_{(x,y)} (k),  \eqno (4.3.32) $$

\noi where for all $(x,y)$, $k^2 \, \epsilon_{(x,y)} (k)$ is bounded and $C$ is a universal constant. 

\noi Therefore, for all $N$ : 
\begin{eqnarray*}
\hspace*{1cm} &  \sum_{k=0}^{N} \mathbb{P}_{(x,y)} \left[ X_k = (0,0) \right] & = C \, \sum_{k \leq N, \, k \equiv x+y \,  (mod. \, 2)} \, \frac{1}{k+1} \, + \, \sum_{k=0}^N \epsilon_{(x,y)} (k)  \\ 
& & = \frac{C}{2} \log (N) + \eta_{(x,y)} (N), \hspace*{4cm} (4.3.33) 
\end{eqnarray*}

\noi where for all $(x,y) \in \mathbb{Z}^2$, $\eta_{(x,y)} (N)$ converges to a limit $\eta_{(x,y)} (\infty)$ when $N$ goes to infinity. Therefore : 
$$
\sum_{k=0}^{N} \left[ \mathbb{P}_{(0,0)} \left( X_k = (0,0) \right) -  \mathbb{P}_{(x,y)} \left( X_k = (0,0) \right)\right] \underset{N \rightarrow \infty}{\rightarrow} \phi((x,y)) := \eta_{(0,0)} (\infty) - \eta_{(x,y)} (\infty).
\eqno (4.3.34) $$

\noi By Proposition 4.3.1, the function $\phi$ is harmonic except at $(0,0)$, and can be used to construct the family of probabilities $(\mathbb{Q}_{(x,y)})_{(x,y) \in \mathbb{Z}^2}$, as in dimension one. Moreover, it is not difficult to check that $\mathbb{Q}_{(x,y)}$ is the image of $\mathbb{Q}_{(0,0)}$ by the translation of $(x,y)$. 


\newpage

{\bf{\LARGE About our bibliography}}

\bigskip

It contains several references which are not directly cited in our text. However, they all concern topics which are 
closely related to those which are developed in the present monograph. They are : 
\begin{itemize}
\item articles whose subject is "close" to the penalisations described here : [Kn], [Ko], [LG2], [N1], 
[N2], [N3], [N4], [RVY, J], [W1], [W2], [W3];
\item articles concerning the Az\'ema-Yor martingales; we encountered those martingales in Chapter 1, Examples 2
and 3 : [AY1], [O]; 
\item articles discussing enlargements of filtrations, an indispensable tool in several penalisation problems : 
[J], [JY]; 
\item the article [PY2] by J. Pitman and M. Yor, which is a key to several multidimensional penalisations 
(see [RVY, VI]) and which helps to understand Chapter 2 of this monograph; 
\item the article [Pi] by J. Pitman, where he establishes his $2S-X$ representation of the BES(3) process; this
result has been adequately extended to processes obtained via a penalisation procedure in [RVY, IV]. 
\end{itemize}

\newpage

{\bf{\LARGE Index of main notations}}

\bigskip

\noi {\bf Chap. 1}

\noi $\Omega = \mathcal{C} ( \mathbb{R}_{+} \to \mathbb{R})$ : the space of continuous functions from $\mathbb{R}_{+}$ to $\mathbb{R}$

\noi $(X_{t}, \; t \ge 0)$ : the set of coordinates on this space

\noi $(\mathcal{F}_{t}, \; t \ge 0)$ : the natural filtration of $(X_t, \, t \geq 0)$

\noi $\mathcal{F}_{\infty} = \dis \mathop{\vee}_{t \ge 0} \mathcal{F}_{t}$

\noi $b(\mathcal{F}_{t})$ : the space of bounded real valued $\mathcal{F}_{t}$ measurable functions

\noi $(W_{x}, \; x \in \mathbb{R})$ : the set of Wiener measures on $(\Omega, \mathcal{F}_{\infty})$

\noi $W = W_{0}$

\noi $W_{x} (Y)$ : the expectation of the r.v. $Y$ with respect to $W_{x}$ 

\noi $ (L_{t}^{y}, \; y \in \mathbb{R}, \; t \ge 0)$ : the bicontinuous process of local times 

\noi $(L_{t} := L_{t}^{0}, \; t \ge 0)$ the local time at level $0$

\noi $\big(\tau_{l} := \inf \{t \ge 0 \;;\; L_{t} > l\}, \; l \ge 0\big)$ : the right continuous inverse of $(L_{t}, \; t \ge 0)$

\noi $q$ : a positive Radon measure on $\mathbb{R}$

\noi $\mathcal{I}$ : the set of positive Radon measures on $\mathbb{R}$ s.t. $\dis 0 < \int_{-\infty}^{\infty} \big( 1+|x|\big) q (dx) < \infty$

\noi $\delta_{a}$ : the Dirac measure at $a$

\noi $\dis \left( A_{t}^{(q)} := \int_{0}^{t} q(X_{s}) ds = \int_{\mathbb{R}} L_{t}^{y} q(dy), \; t \ge 0\right)$ : the additive functional associated with $q$

\noi $(W_{x, \infty}^{(q)}, \; x \in \mathbb{R})$ : the family of probabilities on $(\Omega, \mathcal{F}_{\infty})$ obtained by Feynman-Kac penalisation

\noi $(M_{x,s}^{(q)}, \; s \ge 0)$ : the martingale density of $W_{x, \infty}^{(q)}$ with respect to $W_{x}$

\noi $\gamma_{q}$ : a scale function

\noi $\varphi_{q}, \varphi_{q}^{\pm}$ : solutions of the Sturm-Liouville equation $\varphi'' =q \varphi$

\noi $({\bf W}_{x}, \; x \in \mathbb{R})$ : a family of positive $\sigma$-finite measures on $(\Omega, \mathcal{F}_{\infty})$

\noi $L^{1} (\Omega, \mathcal{F}_{\infty}, {\bf W})$ \big(resp. $L_{+}^{1} (\Omega, \mathcal{F}_{\infty}, {\bf W})\big)$ : the Banach space of 

\noi ${\bf W}$-integrable r.v.'s (resp. the cone of positive and ${\bf W}$-integrable r.v.'s)

\noi $\big(M_{t} (F), \; t \ge 0\big)$ : a martingale associated with $F \in L^{1} (\Omega, \mathcal{F}_{\infty}, {\bf W})$ 

\noi $g_{a} := \sup \{s \ge 0 \;;\; X_{s} = a\} \quad ; \quad  g_{0} = g$

\noi $g_{a}^{(t)} := \sup \{ s \le t, \; X_{s} = a\} \quad ; \quad g_{0}^{(t)} = g^{(t)}$

\noi $\sigma_{a} := \sup \big\{ s \ge 0 \;;\; X_{s} \in [-a,a]\big\} \;;\; \sigma_{a,b} := \sup \big\{ s \ge 0 \;;\; X_{s} \in [a,b]\big\}$

\noi $f_{Z}^{(P)}$ : density of the r.v. $Z$ under $P$

\noi $T$ : a $(\mathcal{F}_{t}, \; t \ge 0)$ stopping time

\noi $P_{0}^{(3)}$ (resp. $\widetilde{P}_{0}^{(3)}$) : the law of a 3-dimensional Bessel process (resp. of the opposite of a 3-dimensional Bessel process) started at 0

\noi $\dis P_{0}^{(3, {\rm sym})} = \frac{1}{2} (P_{0}^{(3)} + \widetilde{P}_{0}^{(3)})$

\noi $W_{0}^{\tau_{l}}$ : the law of a 1-dimensional Brownian motion stopped at $\tau_{l}$

\noi $\Pi_{0,0}^{(t)}$ : the law of the Brownian bridge $(b_{u}, \; 0 \le u \le t)$ of length $t$ and s.t. $b_{0} = b_{t} =0$

\noi $\omega \circ \widetilde{\omega}$ :  the concatenation of $\omega$ and $\widetilde{\omega}$ $(\omega, \widetilde{\omega} \in \Omega)$

\noi $\omega = (\omega_{t}, \omega^{t})$ : decomposition of $\omega$ before and after $t$

\noi $\dis \Gamma^{+} = \big\{\omega \in \Omega \;;\; X_{t} \mathop{\longrightarrow}_{t \to \infty} \infty \big\}$, $\dis \Gamma^{-} = \big\{\omega \in \Omega \;;\; X_{t} (\omega) \mathop{\longrightarrow}_{t \to \infty} -\infty \big\}$

\noi ${\bf W}^{+} = 1_{\Gamma^{+}} \cdot {\bf W}, \quad {\bf W}^{-} = 1_{\Gamma^{-}} \cdot {\bf W}$

\noi $W^{F} \big(F \in L_{+}^{1} (\Omega, \mathcal{F}_{\infty}, {\bf W})\big)$: the finite measure defined on $(\Omega, \mathcal{F}_{\infty})$ by : $W^{F} (G) = {\bf W}(F \cdot G)$

\noi $\mathcal{C}$ : the class of "good" weight processes for which Brownian penalisation holds

\noi $(\nu_{x}^{(q)}, \; x \in \mathbb{R})$ : a family of $\sigma$-finite measures associated with the additive functional $(A_{t}^{(q)}, \; t \ge~0)$

\noi $(Z_{t}, \; t \ge 0)$ : a positive Brownian supermartingale

\noi $\dis Z_{\infty} := \mathop{\lim}_{t \to \infty} Z_{t} \quad W$ a.s. ; $\dis z_{\infty} := \mathop{\lim}_{t \to \infty} \frac{Z_{t}}{1+ |X_{t}|} \quad {\bf W}$ a.s.

\noi $\big(\Delta_{t} (F), \; t \ge 0\big), \; \big(\Sigma_{t} (F), \; t \ge 0\big)$ : two quasimartingales associated with $F \in L^{1} (\Omega, \mathcal{F}_{\infty}, {\bf W})$

\noi $(\Phi_s, \; s \ge 0)$ : a predictable positive process

\noi $\big(k_{s} (F), \; s \ge 0 \big)$ a predictable process such that ${\bf W} \big(F | \mathcal{F}_{g}\big) = k_{g} (F)$ $\big(F \in L_{+}^{1} (\Omega, \mathcal{F}_{\infty}, {\bf W}) \big)$

\noi $(\chi_{t}, \; t \ge 0)$ : a $\mathcal{C} (\mathbb{R}_{+} \to \mathbb{R})$ valued Markov process 

\noi $(\mathbb{P}_{t}, \; t \ge 0)$ : the semigroup associated to $(\chi_{t}, \; t \ge 0)$

\noi ${\bf W}_{x}^{a,b} = a {\bf W}_{x}^{+} + b {\bf W}_{x}^{-}$

\noi $\dis \widetilde{\bf W}^{a,b} = \int \; dx {\bf W}_{x}^{a,b}$ : is an invariant measure for $(\chi_{t}, \; t \ge 0)$

\noi $\widetilde{\Omega} = \mathcal{C} (\mathbb{R} \to \mathbb{R}_{+})$ : the space of continuous functions from $\mathbb{R}$ to $\mathbb{R}_{+}$

\noi $<q,l> := \dis \int_{\mathbb{R}} l(x) q (dx), \; q \in \mathcal{I}, \; l \in \widetilde{\Omega}$

\noi $\mathcal{L} : \Omega \to \widetilde{\Omega}$ defined by $\mathcal{L} (X_{t}, \; t \ge 0) = (L_{\infty}^{y}, \; y \in \mathbb{R})$

\noi $(Q_{t}, \; t \ge 0)$ : the semigroup associated with the Markov process $\big((X_{t}, \; L_{t}^{\bullet}) , \; t \ge 0\big)$ which is $\mathbb{R} \times \widetilde{\Omega}$ valued

\noi $\mathcal{G}$ : the infinitesimal generator of $(Q_{t}, \; t \ge 0)$

\noi $(\widetilde{\bf \Lambda}^{a,b}, a,b \ge 0)$ : a family of invariant measures for $\big((X_{t}, \; L_{t}^{\bullet}) , \; t \ge 0\big)$

\noi $({\bf \Lambda}_{x}, \; x \in \mathbb{R})$ : a family of positive and $\sigma$-finite measures on $\widetilde{\Omega}$

\noi $\theta \;:\; \mathbb{R} \times \widetilde{\Omega} \to  \widetilde{\Omega}$ defined by $\theta (x,l) (y) = l(x-y)$ \quad $(x,y \in \mathbb{R}, \; l \in  \widetilde{\Omega})$

\noi $(L_{t}^{ X_{t^{- \bullet} }}, \;  t \ge 0)$ : a $ \widetilde{\Omega}$ valued Markov process

\noi $(\overline{Q}_{t}, \; t \ge 0)$ : the semigroup associated with $(L_{t}^{ X_t - \bullet} , \;  t \ge 0)$

\noi $\overline{\mathcal{G}}$ :  the infinitesimal generator of $(\overline{Q}_{t}, \; t \ge 0)$

\noi \boldmath$\Lambda$\unboldmath$^{a,b} = a$ \boldmath$\Lambda$\unboldmath$^{+} + b$ \boldmath$\Lambda$\unboldmath$^{-}$

\bigskip

\noi {\bf Chap. 2}

\noi $\Omega = \mathcal{C} (\mathbb{R}_{+} \to \mathbb{C})$ : the space of continuous functions from $\mathbb{R}_{+}$ to $\mathbb{C}$

\noi $(X_{t}, \; t \ge 0)$ : the coordinate process on $\Omega$

\noi $(W_{x}^{(2)}, \; x \in \mathbb{C})$ the set of Wiener measures ; $W_{0}^{(2)} = W^{(2)}$

\noi $\mathcal{J}$ : the set of positive Radon measures on $\mathbb{C}$ with compact support 

\noi $\dis (A_{t}^{(q)} := \int_{0}^{t} q (X_{s}) ds, \; t \ge 0)$ : the additive functional associated with $q \in \mathcal{J})$

\noi $(W_{z, \infty}^{(2,q)}, \; z \in \mathbb{C})$ : the set of probabilities obtained by Feynman-Kac penalisations associated with $q \in \mathcal{I}$ ; $W_{0, \infty}^{(2,q)}= W_{\infty}^{(2,q)}$

\noi $(M_{s}^{(2,q)}, \; s \ge 0)$ : the martingale density of $W_{z, \infty}^{(2,q)}$ with respect to $W_{z}^{(2)}$

\noi $\varphi_{q}$ : a solution of Sturm-Liouville equation $\Delta \varphi = q \varphi$

\noi $\Delta$ : the Laplace operator

\noi $({\bf W}_{z}^{(2)}, \; z \in \mathbb{C})$ : a family of positive and $\sigma$-finite measures on $(\Omega, \mathcal{F}_{\infty})$

\noi ${\bf W}_{0}^{(2)} = {\bf W}^{(2)}$

\noi $C$ : the unit circle in $\mathbb{C}$

\noi $(L_{t}^{(C)}, \; t \ge 0)$ : the continuous local time process on $C$

\noi $(\tau_{l}^{(C)}, \; l \ge 0)$ : the right continuous inverse of $(L_{t}^{(C)}, \; t \ge 0)$

\noi $(R_{t}, \; t \ge 0)$ : the process solution of (2.2.6)

\noi $P_{1}^{(2, \log)}$ : the law of process $(R_{t}, \; t \ge 0)$

\noi $(\rho_{u}, \; u \ge 0)$ : a 3-dimensional Bessel process starting from 0.

\noi $\dis \left(H_{t} := \int_{0}^{t} \frac{ds}{R_{s}^{2}}, \; t \ge 0\right)$

\smallskip

\noi $g_{C} := \sup \{s \ge 0 \;;\; X_{t} \in C\}$

\noi $W_{0}^{\big(2, \tau_{l}^{(C)}\big)}$ : the law of a $\mathbb{C}$-valued Brownian motion stopped at $\tau_{l}^{(C)}$

\noi $\widetilde{P}_{1}^{(2, \log)}$ : the law of $(X_{g_{C} + s}, \; s \ge 0)$

\noi $\nabla$ : the gradient operator

\noi $K_{0}$ : the Bessel Mc Donald function with index 0

\smallskip

\noi $T_{1}^{(3)} := \inf \{u \;;\; \rho_{u}=1\}$

\noi $(R_{t}^{(\delta)}, \; t \ge 0)$ : the process solution of (2.3.19)

\noi $\big(M_{t}^{(2)} (F), \; t \ge 0\big)$ : the Brownian martingale associated with $F \in L^{1} (\Omega, \mathcal{F}_{\infty}, {\bf W}^{(2)})$

\bigskip

\noi {\bf Chap. 3}

\noi $\Omega = \mathcal{C} (\mathbb{R}_{+} \to \mathbb{R}_{+})$ : the space of continuous functions from $\mathbb{R}_{+}$ to $\mathbb{R}_{+}$ 

\noi $S$ : the scale function

\noi $m$ : the speed measure

\noi $(X_{t}, \; t \ge 0, \; P_{x}, \; x \in \mathbb{R}_{+})$ : the canonical process associated with $S$ and $m$

\noi $(\mathcal{F}_{t}, \; t \ge 0)$ : the natural filtration of $(X_t, \, t \geq 0)$; $\dis \mathcal{F}_{\infty} = \mathop{\vee}_{t \ge 0} \mathcal{F}_{t}$

\noi $\dis L = \frac{d}{dm} \; \frac{d}{dS}$ : the infinitesimal generator of $(X_{t}, \; t \ge 0)$

\noi $p(t,x,\bullet)$ : the density of $X_{t}$ under $P_{x}$ with respect to $m$

\noi $(L_{t}^{y}, \; t \ge 0, \; y \ge 0)$ : the jointly continuous family of local times of $X$
     
\noi $(L_{t}, \; t \ge 0)$ : the local time process at level 0

\noi $(\tau_{l}, \; l \ge 0)$ : the right continuous inverse of $(L_{t}, \; t \ge 0)$

\noi $P_{x}^{\tau_{l}}$ : the law of the process $(X_{t}, \; t \ge 0)$ started at $x$ and stopped at $\tau_{l}$

\noi $g_{y} := \sup \{ t \ge 0 \;;\; X_{t}=y\} \quad ; \quad g := g_{0}$

\noi $g_{y}^{(t)} := \sup \{ s \le t \;;\; X_{s}=y\} \quad ; \quad g^{(t)} := g_{0}^{(t)}$

\noi $T_{0} := \inf \{ t \ge 0 \;;\; X_{t} = 0\}$

\noi $(\widehat{X}_{t}, \; t \ge 0)$ : the process $(X_{t}, \; t \ge 0)$ killed at $T_{0}$

\noi $\widehat{p} (t, x, \bullet)$ : the density of $\widehat{X}_{t}$ under $P_{x}$ with respect to $m$

\noi $(P_{x}^{\uparrow}, \; x \in \mathbb{R}_{+})$ : the laws of $X$ conditionned not to vanish ; $P^{\uparrow} := P_{0}^{\uparrow}$

\noi $f_{y,0} (t)$ defined by : $f_{y,0}(t)dt = P_{y} (T_{0} \in dt) = P_{0}^{\uparrow} (g_{y} \in dt)$

\noi ${\bf W}^{*}$ a $\sigma$-finite measure on $(\Omega, \mathcal{F}_{\infty})$

\noi $\Pi_{0}^{(t)}$ : the law of the bridge of length $t$

\noi ${\bf W}_{g}^{*}$ : the restriction of ${\bf W}^{*}$ to $\mathcal{F}_{g}$

\noi $\left(M_{t}^{(\lambda, x)} = \dis \frac{1 + \frac{\lambda}{2} \, S(X_{t})}{1 + \frac{\lambda}{2} \, S(x)} \cdot e^{- \frac{\lambda}{2} \, L_{t}}, \; t \ge 0\right)$ : the martingale density of $P_{x, \infty}^{(\lambda)}$ with respect to $P_{x}$

\noi $(M_{t}^{*} (F), \; t \ge 0)$ : the positive $\big((\mathcal{F}_{t}, \, t \ge 0), \; P_{0}\big)$ martingale associated with $F \in L^{1} (\Omega, \mathcal{F}_{\infty}, {\bf W}^{*})$

\noi $\big(P_{x}^{(-\alpha)}, \; x \ge 0\big)$ : the family of laws of a Bessel process with dimension $d=2(1-\alpha)$ ($0<d<2$, or equivalently $0<\alpha <1$) started at $x$

\noi ${\bf W}^{(-\alpha)}$ : the measure ${\bf W}^{*}$ in the particular case of a Bessel process with index $(- \alpha)$ 

$(0< \alpha <1)$

\noi $\Pi_{0}^{(-\alpha,t )}$ : the law of the Bessel bridge with index $(- \alpha)$ and length $t$

\noi $P_{x}^{(-\alpha, \tau_{l})}$ : the law of a Bessel process with index $(- \alpha)$ started at $x$ and stopped at $\tau_{l}$

\noi $\varphi_{q}$ : a particular solution of the Sturm-Liouville equation :
$$
\frac{1}{2} \, \varphi'' (r) + \frac{1-2 \alpha}{2 r} \, \varphi' (r) = \frac{1}{2} \, \varphi (r) \, q(r), \quad  r \ge 0$$

\noi with $q$ a positive Radon measure with compact support

\noi $(m_{u}, \; 0 \le u \le 1)$ : the Bessel meander with dimension $d$

\noi $P_{0}^{\big( \frac{\delta}{2}-1, \; m, \; \nearrow \nwarrow)}$ : the law of the process obtained by putting two Bessel processes with index $\dis \left(\frac{\delta}{2}-1\right)$ back to back; these processes start from 0 and are stopped when they first reach level $m$

\bigskip

\noi {\bf Chap. 4}

\noi $E$ : a countable set

\noi $(X_{n}, \; n \ge 0)$ : the canonical process on $E^{\mathbb{N}}$

\noi $(\mathcal{F}_{n}, \; n \ge 0)$ : the natural filtration, $\mathcal{F}_{\infty} = \dis \mathop{\vee}_{n \ge 0} \mathcal{F}_{n}$

\noi $(\mathbb{P}_{x}, \; x \in E)$ : the family of probabilities associated to Markov process $(X_{n}, \; n \ge 0)$ s.t. $\mathbb{P}(X_{n+1} = z | X_{n}=y)=p_{y,z}$ and $\mathbb{P}_{x}(X_{0}=x)=1$

\noi $\left(L_{k}^{y} = \dis \sum_{m=0}^{k} 1_{X_{m}=y}, \; k \ge 0 \right)$ : the local time of $(X_{n}, \; n \ge 0)$ at level $y$ (with $L_{-1}^{y}=0)$

\noi $\phi$ : a positive function from $E$ to $\mathbb{R}_{+}$, harmonic with respect to $\mathbb{P}$, except at the point $x_{0}$ and such that $\phi (x_{0})=0$

\noi $\psi_{r} (x) := \dis \frac{r}{1-r} \, \mathbb{E}_{x_{0}} \big(\phi(X_{1})\big) + \phi (x) \quad \big(r \in ]0,1[, \; x \in E)$

\noi $\big(\mu_{x}^{(r)}, \; x \in E, \; r \in ]0,1[\big)$ : a family of finite measures on $(E^{\mathbb{N}}, \mathcal{F}_{\infty})$

\noi $\mathbb{Q}_{x} = \dis \left(\frac{1}{r}\right)^{L_{\infty}^{x_{0}}}\, \mu_{x}^{(r)}$, independent of $r \in ]0,1[$

\noi $\mathbb{Q}_{x}^{(\psi, y_{0})}$ : the $\sigma$-finite measure $\mathbb{Q}_{x}$ constructed from the point $y_{0}$ and the function $\psi$

\noi $q$ : a function from $E$ to $[0,1]$ such that $\{q<1\}$ is a finite set

\noi $\big(M(F, X_{0}, X_{1}, \cdots, X_{n}), \; n \ge 0\big)$ : the  $\big((\mathcal{F}_{n}, \; n \ge 0), \; \mathbb{P}_{x}\big)$ martingale associated with $F \in L^{1}(\Omega, \mathcal{F}_{\infty}, \mathbb{Q}_{x})$

\noi $\tau_{k}^{(y)}$ : the $k$-th hitting time of $y$

\noi $(\tau_{k}^{(y)}, \; k \ge 0)$ : the inverse of $(L_{k}^{y}, \; k \ge 0)$ 

\noi $\mathbb{Q}_{y}^{[y_{0}]}$ : the restriction of $\mathbb{Q}_{y}$ to trajectories which do not hit $y_{0}$

\noi $\widetilde{\mathbb{Q}}_{y}$ : the restriction of $\mathbb{Q}_{y}$ to trajectories which do not return to $y$

\noi $\mathbb{P}_{x}^{\tau_{k}^{(y_{0})}}$ : the law of the Markov chain $(X_{n}, \; n \ge 0)$ starting from $x$ and stopped at $\tau_{k}^{(y_{0})}$

\noi $2 \widetilde{\mathbb{Q}}_{a}^{+}$ : the law of a Bessel random walk strictly above $a$

\noi $2 \widetilde{\mathbb{Q}}_{a}^{-}$ : the law of a Bessel random walk strictly below $a$

\noi $\widetilde{\mathbb{Q}}_{a} :=  \widetilde{\mathbb{Q}}_{a}^{+} + \widetilde{\mathbb{Q}}_{a}^{-}$

\noi $g_{a} := \sup \{n \ge 0 \;;\; X_{n} =a\}$

\noi $\phi^{[y_{0}]}$ defined by $\phi^{[y_{0}]} (y) = \mathbb{Q}_{y}^{[y_{0}]} (1)$

\noi $\simeq$ : the equivalence relation defined in Subsection 4.2.4

\noi $\mathbb{Q}_{x}^{[\psi]}$ : the measure $Q_{x}^{(\psi, y_{0})}$ where $[\psi]$ denotes the equivalence class of $\psi$


\newpage

{\bf{\LARGE Acknowledgment}}

\bigskip

The authors are grateful to the Mathematical Society of Japan for accepting to publish this monograph, and to 
Dr. K. Yano (Kobe University) who carefully read the manuscript. The third author (M. Yor) is grateful for the 
hospitality he received at RIMS during the month of October 2007 and for the series of lectures he could give
in the Mathematical Department of Kyoto University at that time.


\newpage

\begin{thebibliography}{99}

   \bibitem[AY1]{clefReference1}
      J. Az\'ema, M. Yor  {\em Une solution simple au probl\`eme de Skorokhod}. In S\'eminaire de Probabilit\'es XIII (Univ. Strasbourg, Strasbourg 1977/78). LNM {\bf 721}, Springer, Berlin,  p. 90-115, 1979
 
   \bibitem[AY2]{clefReference2}
      J. Az\'ema, M. Yor  {\em Sur les z\'eros des martingales continues}. In S\'eminaire de Probabilit\'es XXVI. LNM {\bf 1526}, Springer, Berlin, Heidelberg, New-York, p. 248-306, 1992
 
 \bibitem[Be]{clefReference3}
      J. Bertoin  {\em L\'evy processes}. Cambridge Tracts in Math. Cambridge, Vol. {\bf 121}, Univ. Press, Cambridge 1996

   \bibitem[Bi]{clefReference3}
      P. Biane  {\em Decomposition of Brownian trajectories and some applications}. In A. Badrikian, P.A. Meyer, J.A. Yan (Eds). Probability and Statistics. Rencontres franco-chinoises en Probabilit\'es et Statistiques, Proceedings of Wuhan meeting, World Scientific, Singapore, p. 51-76, 1993
      
    
 \bibitem[BY]{clefReference3}
      P. Biane, M. Yor {\em Valeurs principales associ\'ees aux temps locaux browniens}. Bul. Sci. Math. {\bf 111}  p. 23-101 (1987)
 
 \bibitem[BS]{clefReference3}
      A.N. Borodin, P. Salminen {\em Handbook of Brownian motion - facts and formulae - }. Prob. and its Applications. Birkh\"auser, Basel, 2d edition 2002
 
 \bibitem[C]{clefReference3}
      L. Chaumont {\em Excursion normalis\'ee, m\'eandre et pont pour les processus de L\'evy stables}. Bull. Sci. Math.  {\bf 121(5)}  p. 377-403 (1997)

\bibitem[D]{clefReference3} 
     E.B. Dynkin {\em Markov processes. Vols I. II.} Translated with the authorization and assistance of the author
              by J. Fabius, V. Greenberg, A. Maitra, G. Majone. Die
              Grundlehren der Mathematischen Wissenschaften {\bf 121}, {\bf 122}. Springer-Verlag, Berlin-G\"ottingen-Heidelberg (1965)
 
 \bibitem[DMY]{clefReference3}
      C. Donati-Martin, M. Yor {\em Some measure valued Markov processes attached to occupation time of Brownian motion}. Bernoulli {\bf 6(1)}  p. 63-72  (2000)
 
 
 \bibitem[DMRVY]{clefReference3}
      C. Donati-Martin, B. Roynette, P. Vallois, M. Yor {\em On constants related to the choice of the local time at $0$, and the corresponding  It\^o measure for Bessel processes with dimension $d = 2(1-\alpha),  0<\alpha<1$ }. Studia Sci. Math. Hungarica {\bf 44(2)}, p. 207-221 (2008)
 
 \bibitem[Fe]{clefReference3}
      W. Feller {\em An introduction to probability theory and its applications. Vol. II}. J. Wiley and Sons Inc., New-York 1966
 
 \bibitem[F]{clefReference3}
      H. F\"ollmer {\em The exit measure of a supermartingale}. Z. Wahrscheinlichkeits\-theorie und Verw. Gebiete  {\bf 21} p. 154-166 (1972)
 
 \bibitem[J]{clefReference3}
      T. Jeulin {\em Semi-martingales et grossissement d'une filtration}. LNM  {\bf 833}, Springer, Berlin, Heidelberg, New-York, 1980
 
  \bibitem[JY]{clefReference3}
    T. Jeulin, M. Yor {\em Grossissement de filtrations : exemples et applications}. LNM  {\bf 1118}, Springer, Berlin, Heidelberg, New-York. Papers from the seminar on stochastic calculus held at University Paris VI, Paris 1982/83

 
 \bibitem[Kn]{clefReference3}
    F.B. Knight {\em Brownian local time and taboo processes}. Trans. Amer. Math. Soc. {\bf 143} p. 173-185 (1969)
 
 \bibitem[Ko]{clefReference3}
      S. Kotani {\em Asymptotics for expectations of multiplicative functionals of one-dimensional Brownian motion}. preprint (nov. 2006)
 
\bibitem[KS]{clefReference3}
      U. K\"uchler, P. Salminen {\em On spectral measures of strings and excursions of quasi-diffusions}. In J. 
Az\'ema, P.A. Meyer and M. Yor, editors. S\'eminaire de Probabilit\'es XXIII. LNM {\bf 1372}, Springer, Berlin,  p. 490-502, 1989
    
       \bibitem[L]{clefReference3}
      N.N. Lebedev {\em Special functions and their applications}. Revised English Ed., Dover Publications Inc, New-York 1965

 \bibitem[LG]{clefReference3}
      J.F. Le Gall {\em Une approche \'el\'ementaire des th\'eor\`emes de d\'ecomposition de Williams}. In S\'eminaire de Probabilit\'es XX. LNM {\bf 1204}, Springer, Berlin,  p. 447-464, 1986
 
\bibitem[LG2]{clefReference3}
      J.F. Le Gall {\em Sur le temps local d'intersection du mouvement brownien plan et la m\'ethode de renormalisation de Varadhan}. In S\'eminaire de Probabilit\'es XIX. LNM {\bf 1123}, Springer, Berlin, Heidelberg, New-York,  p. 314-331, 1985

 
 \bibitem[M]{clefReference3}
      P.A. Meyer {\em Probabilit\'es et potentiel}. Publication de l'Institut de Math. de l'Univ. de Strasbourg XIV, Actualit\'es Scientifiques et Industrielles {\bf 1318}, Hermann, Paris 1966
 
 \bibitem[N1]{clefReference3}
        J. Najnudel {\em Temps locaux et p\'enalisations browniennes}. Th\`ese de l'Universit\'e Paris 6, juin 2007
      
      \bibitem[N2]{clefReference3}
        J. Najnudel {\em P\'enalisations de l'araign\'ee brownienne}. Ann. Inst. Fourier  {\bf 57(4)},  p. 1063-1093 (2007)

 
  \bibitem[N3]{clefReference3}
       J. Najnudel {\em On the construction of an Edwards' probability measure on $\mathcal{C} (\mathbb{R}_{+}, \mathbb{R})$}. Submitted (2008)

   \bibitem[N4]{clefReference3}
       J. Najnudel {\em Penalizations of the Brownian motion with a functional of its local times}. S. P. A. {\bf 118}, p. 1407-1433 (2008)

 \bibitem[NRY]{clefReference3}
      J. Najnudel, B. Roynette, M. Yor {\em A remarkable $\sigma$-finite measure on $\mathcal{C} (\mathbf{R}_+, \mathbf{R})$ related to many Brownian penalisations}. C. R. Acad. Sci. Paris {\bf 345}, p. 459-466 (2007)
  
   
   \bibitem[O]{clefReference3}
      J. Obl\'oj {\em A complete characterization of local martingales which are functions of brownian motion and its supremum}. Bernoulli {\bf 12(6)} p. 955-969 (2006)
   
    \bibitem[Pi]{clefReference3}
       J. Pitman {\em One dimensional Brownian motion and the three dimensional Bessel process}. Adv. in Appl. Prob.   {\bf 7(3)}, p. 511-526 (1975)
    
    \bibitem[PY1]{clefReference3}
      J. Pitman, M. Yor {\em Asymptotic laws of planar Brownian motion}. Ann. Prob. {\bf 14}, p. 733-779 (1986)
     
      \bibitem[PY2]{clefReference3}
      J. Pitman, M. Yor  {\em Decomposition at the maximum for excursions and bridges of one dimensional diffusions}. In N. Ikeda, S. Watanabe, M. Fukushima, H. Kunita (Eds). It\^o's stochastic calculus and probability theory. Springer, Berlin, Heidelberg, New-York, p. 293-310, 1996
      
       \bibitem[Pr]{clefReference3}
      C. Profeta, {\em }Thesis in preparation (2008-2010)
       
        \bibitem[R]{clefReference3}
      K.M. Rao {\em Quasi-martingales}. Math. Scand. {\bf 24}, p. 79-92 (1969)
        
        \bibitem[Rev]{clefReference3}
        D. Revuz {\em Mesures associ\'ees aux fonctionnelles additives de Markov I}. Trans. Amer. Math. Soc. {\bf 148}, p. 501-531 (1970) 

         \bibitem[ReY]{clefReference3}
      D. Revuz, M. Yor {\em Continuous martingales and Brownian motion}. Grundlehren der Math. Wissenschaften {\bf 293}, Springer, Berlin, Third Edition 1999
         
          \bibitem[RVY, I]{clefReference3}
    B. Roynette, P. Vallois, M. Yor {\em Limiting laws associated with Brownian motion perturbed by normalized exponential weights, I}. Studia Sci. Math. Hungarica  {\bf 43(2)}, p. 171-246 (2006)
          
           \bibitem[RVY, II]{clefReference3}
      B. Roynette, P. Vallois, M. Yor {\em Limiting laws associated with Brownian motion perturbed by its maximum, minimum and local time II}.  Studia Sci. Math. Hungarica  {\bf 43(3)}, p. 295-360 (2006)
           
            \bibitem[RVY, III]{clefReference3}
      B. Roynette, P. Vallois, M. Yor {\em Limiting laws for long Brownian bridges perturbed by their one-sided maximum III}.  Period. Math. Hungar. {\bf 50(1-2)}, p. 247-280 (2005)
            
             \bibitem[RVY IV]{clefReference3}
      B. Roynette, P. Vallois, M. Yor {\em Some extensions of Pitman and Ray-Knight theorems for penalized Brownian motion and their local time IV}. Studia Sci. Math. Hungarica  {\bf 44(4)}, p. 469-517 (2007)
             
              \bibitem[RVY, V]{clefReference3}
B. Roynette, P. Vallois, M. Yor {\em Penalizing a BES(d) process $(0< d <2$) with a function of its local time V}. Studia Sci. Math. Hungarica {\bf 45(1)}, p. 67-124 (2008)
              
               \bibitem[RVY, VI]{clefReference3}
      B. Roynette, P. Vallois, M. Yor {\em Penalisation of multidimensional Brownian motion VI}. ESAIM P.S. (2009)
               
                \bibitem[RVY, VII]{clefReference3}
      B. Roynette, P. Vallois, M. Yor {\em Brownian penalisations related to excursions lengths VII}. To appear in Annales de l'IHP (2009)
                
                 \bibitem[RY, VIII]{clefReference3}
      B. Roynette, M. Yor {\em Ten penalisations of Brownian motion involving its one-sided supremum until first and last passage times VIII}. Journ. Funct. Anal {\bf 255(9)}, p. 
2606-2640 (2008)   
                  
                  \bibitem[RY, IX]{clefReference3}
      B. Roynette, M. Yor {\em Local limit theorems for Brownian additive functionals and penalizations of Brownian paths IX}. To appear in ESAIM P.S. (2009)                 
                         
         \bibitem[RVY, X]{clefReference3}
      B. Roynette, P. Vallois, M. Yor {\em Penalisations of Brownian motions with its maximum and minimum processes as weak forms of Skorokhod embeddings}. To appear in Theory of Proba. and its App. (2009)
      
         \bibitem[RVY, J]{clefReference3}
         B. Roynette, P. Vallois, M. Yor  {\em Some penalisations of the Wiener measure}. Japan Jour. of Math. {\bf 1(1)}, p. 263-290 (2006)
      
         \bibitem[RY, M]{clefReference3}
      B. Roynette, M. Yor {\em Penalising Brownian paths}. LNM {\bf 1969}, Springer (2009)
       
        \bibitem[Sa1]{clefReference3}
      P. Salminen {\em One dimensional diffusions and their exit spaces}. Math. Scand.  {\bf 54}, p. 209-220 (1984)
      
         \bibitem[Sa2]{clefReference3}
      P. Salminen {\em On last exit decomposition of linear diffusions}. Studia Sci. Math. Hungarica  {\bf 33}, p. 251-262 (1997)
      
      \bibitem[SV]{clefReference3}
      P. Salminen, P. Vallois {\em On subexponentiality of the L\'evy measure of the diffusion inverse local time, with applications to penalisations. }. Pr\'epublication IEC Nancy (2008)
      
         \bibitem[SVY]{clefReference3}
      P. Salminen, P. Vallois, M. Yor  {\em On the excursion theory for linear diffusions}. Japan Jour. of Math. {\bf 2}, p. 91-127 (2007)
      
        \bibitem[SY]{clefReference3}
      P. Salminen, M. Yor {\em Tanaka formula for symmetric L\'evy processes}. In S\'eminaire de Probabilit\'es XL. LNM {\bf 1899}, Springer, Berlin,  p. 265-285, 2007   
      
      \bibitem[S]{clefReference3}
      F. Spitzer {\em Some theorems concerning $2$-dimensional Brownian motion}. Trans. Amer. Math. Soc. {\bf 87}, p. 187-197 (1958)

      \bibitem[Spi]{clefReference3}
      F. Spitzer {\em Principles of random walks}. Van Nostrand Princeton (1964)
        
         \bibitem[W1]{clefReference3}
        J. Westwater {\em On Edwards model for long polymer chains}. Comm. Math. Phys. {\bf 72}, p. 131-174 (1980)
      
      \bibitem[W2]{clefReference3}
        J. Westwater {\em On Edwards model for polymer chains. II. The self-consistent potential}. Comm. Math. Phys. {\bf 79}, p. 53-73 (1981)
      
      \bibitem[W3]{clefReference3}
        J. Westwater {\em On Edwards model for polymer chains. III. Borel summability}. Comm. Math. Phys. {\bf 84(4)}, p. 459-470 (1982)

      \bibitem[Wi]{clefReference3}
       D. Williams {\em Diffusions, Markov processes and Martingales, vol. 1}. Foundations, Wiley and Sons, New York (1979)
      \bibitem[YYY]{clefReference3}
      K. Yano, Y. Yano, M. Yor {\em Penalising symmetric stable L\'evy paths}. To appear in J. Math. Soc. Japan. (2009)
 
\end{thebibliography}
\end{document}